\documentclass[12pt,centertags]{amsart}
\usepackage{amsmath,amsfonts,amssymb,amsthm, mathtools}
\usepackage{graphicx,color, dsfont}
\usepackage{epsfig, tikz, amscd}
\usepackage{url}
\usepackage{fixme}
\fxsetup{status=draft} 
\fxsetup{theme=color}
\fxsetup{layout=inline}

\usepackage{enumerate}
\usepackage[toc,page]{appendix}

  \usepackage[margin=1.3in]{geometry}
  \usepackage{tikz}
  \usetikzlibrary{calc,through,backgrounds, arrows, cd}
  \usetikzlibrary{decorations.pathmorphing, decorations.pathreplacing,positioning,fit}

\usepackage{imakeidx}
\makeindex[title=Index Notation,columns=3]

\usepackage{hyperref}
\hypersetup{
  hidelinks,
  pdftitle={Thurston measure is stationary with exponential Teichm\"uller moment},
  pdfauthor={Alex Eskin and Kasra Rafi}
}
\DeclareMathRadical{\sqrtsign}{symbols}{"70}{largesymbols}{"70}

\newcounter{kcomments}

\newcommand{\bb}{\mathbb}

\newlength{\figboxwidth}             
  \newcommand{\param}{{\mathchoice{\mkern1mu\mbox{\raise2.2pt\hbox{$
  \centerdot$}}
  \mkern1mu}{\mkern1mu\mbox{\raise2.2pt\hbox{$\centerdot$}}\mkern1mu}{
  \mkern1.5mu\centerdot\mkern1.5mu}{\mkern1.5mu\centerdot\mkern1.5mu}}}

\newcommand{\integers}{{\bb Z}}
\newcommand{\natls}{{\bb N}}

\newcommand{\reals}{{\bb R}}

  \newcommand{\sB}{{\sf B}}   
  \newcommand{\sC}{{\sf C}}   
  \newcommand{\sD}{{\sf D}}

  \newcommand{\sL}{{\sf L}}     
  \newcommand{\sM}{{\sf M}}

  \newcommand{\sR}{{\sf R}}   
       
  \newcommand{\sT}{{\sf T}}   
     
  \newcommand{\sV}{{\sf V}}   
  \newcommand{\sW}{{\sf W}}

  \newcommand{\kk}{{\sf k}}   
       
  \newcommand{\mm}{{\sf m}}

  \newcommand{\rr}{{\sf r}}   
  \renewcommand{\ss}{{\sf s}}     
  \renewcommand{\tt}{{\sf t}}

\renewcommand{\Re}{\mbox{Re}\:}
\renewcommand{\Im}{\mbox{Im}\:}

\newcommand{\st}{\;\: : \;\:}         

\newcommand{\zed}{\integers}

\newcommand{\axis}{\operatorname{axis}}

\newcommand{\area}{\operatorname{area}}

\newcommand{\length}{\operatorname{length}}

\newcommand{\odd}{{\operatorname{odd}}}
\newcommand{\even}{{\operatorname{even}}}
\newcommand{\GM}{{\operatorname{GM}}}

\newcommand{\GL}{\operatorname{GL}}

\newcommand{\AGY}{{\scriptscriptstyle  \rm AGY}}

\newcommand{\DDM}{{\scriptscriptstyle  \rm DDM}}

\def\@ifundefined#1#2#3%
  {\expandafter\ifx\csname#1\endcsname\relax#2\else#3\fi}

\@ifundefined{theoremstyle}{
}{
\theoremstyle{plain} 
}
\newtheorem{theorem}{Theorem}[section]

\newtheorem{proposition}[theorem]{Proposition}
\newtheorem{lemma}[theorem]{Lemma}

\newtheorem{corollary}[theorem]{Corollary}
\newtheorem{claim}[theorem]{Claim}

\@ifundefined{theoremstyle}{
}{
\theoremstyle{definition} 
}
\newtheorem{definition}[theorem]{Definition}

\newtheorem{notation}[theorem]{Notation}

\theoremstyle{remark}
\newtheorem{remark}[theorem]{Remark}
\newtheorem*{remark*}{Remark}

\newcommand{\cA}{{\mathcal A}}
\newcommand{\cB}{{\mathcal B}}
\newcommand{\cC}{{\mathcal C}}

\newcommand{\cF}{{\mathcal F}}
\newcommand{\cG}{{\mathcal G}}

\newcommand{\cJ}{{\mathcal J}}

\newcommand{\cM}{{\mathcal M}}
\newcommand{\cN}{{\mathcal N}}

\newcommand{\cP}{{\mathcal P}}
\newcommand{\cQ}{{\mathcal Q}}

\newcommand{\cT}{{\mathcal T}}
\newcommand{\cU}{{\mathcal U}}
\newcommand{\cV}{{\mathcal V}}

\newcommand{\cZ}{{\mathcal Z}}

\mathchardef\GG="321D
\newcommand{\bfv}{{\mathbf v}}
\newcommand{\bfw}{{\mathbf w}}

\newcommand{\mc}[1]{{}}  
\newcommand{\mcc}[1]{{}}

\numberwithin{equation}{section}

\newcommand{\QT}{{\cQ\cT}}
\newcommand{\OT}{{\cQ_1\cT}}
\newcommand{\QM}{{\cQ\cM}}
\newcommand{\OM}{{\cQ_1\cM}}
\newcommand{\PQT}{{\cP\cQ\cT}}
\newcommand{\POT}{{\cP\cQ_1\cT}}
\newcommand{\PQM}{{\cP\cQ\cM}}
\newcommand{\POM}{{\cP\cQ_1\cM}}

\newcommand{\MF}{{\cM\cF}}
\newcommand{\NF}{{\cN\cF}}
\newcommand{\PMF}{{\cP\cM\cF}}

\newcommand{\from}{\colon\thinspace} 
\newcommand{\Ext}{\operatorname{Ext}}
\newcommand{\Hyp}{\operatorname{Hyp}}
\newcommand{\Flat}{\operatorname{Flat}}
\newcommand{\Mod}{\operatorname{Mod}}

\newcommand{\C}{\mathbb{C}}
\newcommand{\R}{\mathbb{R}}

\newcommand{\Z}{\mathbb{Z}}
\newcommand{\PP}{\mathbb{P}}

\newcommand{\supp}{\operatorname{supp}}

\newcommand{\norm}[1]{\left\lVert #1 \right\rVert }

\DeclareMathOperator{\sys}{sys}

\def\la{ \lambda}

\def\<{\langle}
\def\>{\rangle}

\def\ep{\epsilon}

\usepackage{xfrac}
  \newcommand{\oof}{{\tfrac{1}{50}}}
  \newcommand{\ooh}{{\tfrac{1}{100}}} 
  \newcommand{\ooth}{{\tfrac{1}{200}}}

  \newcommand{\EE}{\mathbb{E}}
  \newcommand{\tV}{{\widetilde V}}
  \newcommand{\I}{\mathrm{i}}
  \newcommand{\muMV}{{\mu_{\rm MV}}}
  \newcommand{\muG}{\mu_{\Gamma}}
  \newcommand{\muT}{{\mu_{\rm Th}}}

\newcommand{\Teich}{Teich\-m\"u\-ller\ } 
\newcommand{\Hor}{\operatorname{Hor}}

\begin{document}

  \title[Thurston measure is stationary]{Thurston measure is stationary\\
    with exponential Teichm\"uller moment}


 \author   {Alex Eskin}
 \address{Department of Mathematics, University of Chicago, Chicago, IL }
 \email{eskin@math.uchicago.edu}
 
  

 \author   {Kasra Rafi}
 \address{Department of Mathematics, University of Toronto, Toronto, ON }
 \email{rafi@math.toronto.edu}
 
 
  \date{November 18, 2026}

\begin{abstract}
In 1967, Furstenberg proved that every lattice in a connected semisimple
Lie group with finite center admits a probability measure whose stationary
measure on the associated Furstenberg boundary belongs to the Lebesgue
measure class. We prove an analogue of this theorem
for the action of the mapping class group on \Teich space.
Every regular probability measure in the Thurston measure class on the
space of projective measured foliations is the stationary measure of a
non-elementary random walk with a finite exponential moment in the \Teich
metric. 
\end{abstract}

\maketitle 

\tableofcontents 

\section{Introduction}\label{sec:introduction}

A $G$--space is a metric space with a continuous action of a group
$G$.  
A useful first example of a $G$--space is the Poincar\'e disk: the group
$G=\mathrm{PSL}_2(\mathbb R)$ acts on the hyperbolic disk by isometries, and the
circle at infinity records the directions in which points of the disk can escape. The
action extends naturally to this boundary circle. Thus the boundary circle is
itself a compact $G$--space. This picture has a broad generalization. For a semisimple Lie group, the
canonical boundary is the Furstenberg boundary $G/P$, where $P$ is a minimal parabolic
subgroup; it plays, in higher rank, the role played by the ideal circle for the
Poincar\'e disk.

\subsection*{Boundaries and random walks: Lie groups}

Furstenberg's theory explains how the boundary $G/P$ enters
probability. Let $\mu$ be a probability measure on $G$, or on a
lattice $\Gamma<G$, and let
\[
g_n=h_1\cdots h_n
\]
be the associated random walk. Under suitable nondegeneracy
hypotheses, a sample path determines an almost sure limiting flag
$\xi\in G/P$. If $\nu$ denotes the law of $\xi$, then $\nu$ is
$\mu$--stationary:
\[
\nu=\mu*\nu=\int g_*\nu\,d\mu(g).
\]
This identity reflects the Markov property. After the first step, the
remaining tail of the walk has the same distribution as a fresh copy
of the walk, while the endpoint at infinity is translated by the first
increment. Thus $\nu$ can be viewed as the hitting measure of the
random walk on the boundary: it records the distribution of the point
at infinity hit, in the limit, by a typical sample path.

\subsection*{Boundaries and random walks: mapping class groups}
The same picture has a mapping-class-group analogue. 
Let $S$ be an oriented surface of finite type, let $\Mod=\Mod(S)$ be its mapping
class group, and let $\cT=\cT(S)$ be its \Teich space. The Thurston compactification
of $\cT$ has boundary $\PMF$, the space of projective measured foliations. The action
of $\Mod$ on $\cT$ extends continuously to an action on $\PMF$.

Fix a basepoint $Y\in\cT$ and
let
\[
\gamma_n=h_1\cdots h_n
\]
be the random walk on $\Mod$ driven by $\muG$. 
Kaimanovich and Masur, in
\cite{kaimanovich:masur:rennes,kaimanovich:masur:PB},  
proved that, assuming that the semigroup generated by
$\supp(\muG)$ is non-elementary, the orbit points $\gamma_nY$ converge almost
surely to a projective measured foliation; moreover, the limiting foliation is
almost surely uniquely ergodic. If $\nu$ denotes the law
of this limit, then $\nu$ is $\muG$--stationary:
\begin{equation}
\label{eq:stationary}
\muG*\nu=\sum_{\gamma\in\Mod}\muG(\gamma)\,\gamma_*\nu=\nu.
\end{equation}
In this setting, $\nu$ is the unique non-atomic $\muG$--stationary measure on
$\PMF$.

Thus $\nu$ is the hitting measure of the random walk on the Thurston boundary.

\subsection*{Inverse Theorems.}
For lattices in semisimple Lie groups, Furstenberg also gives an
inverse statement. Let $\Gamma$ be a lattice in $G$. Then
there are probability measures $\mu$ on $\Gamma$
with finite first moment for displacement in the associated symmetric
space, for which the
resulting stationary measure is equivalent to the natural smooth measure on $G/P$
\cite{furstenberg:PB,furstenberg:RW}. This result, which goes back to
1967, played a key role in
later developments. In particular, it was used by Furstenberg to prove
an early (pre-Mostow) rigidity theorem for lattices, and later
it was used by Margulis in the first proof of
the Margulis superrigidity theorem. 

In this paper, we are considering the analogue of Furstenberg's
inverse theorem in the Teichm\"uller space setting. Instead of starting with
a random walk and asking what boundary measure it produces, we start with a measure
$\nu$ in the Thurston measure class on $\PMF$ and ask whether one can find a
probability measure $\muG$ on $\Mod$ such that
\[
\nu=\muG*\nu.
\]
In other words, we ask whether a prescribed Thurston-class measure can be realized as
a stationary measure for a random walk on the mapping class group.

\subsection*{Main theorem}

We exclude the once-punctured torus throughout. This is the classical
flat-torus case of Furstenberg's theory
\cite{furstenberg:PF,furstenberg:RW}.

We work with the following regular representatives of the Thurston
measure class.
Let $\rho:\MF\to\R_{\ge 0}$ be continuous, homogeneous of degree one,
positive away from $0$, and piecewise smooth in train-track coordinates.
Here and below, piecewise smoothness is understood on $\MF\setminus\{0\}$:
in each train-track coordinate domain there is a locally finite
polyhedral subdivision on whose pieces $\rho$ is smooth up to the
faces, meaning that each restriction extends smoothly to a neighborhood
of the closed piece. Together with continuity, this implies that
$\rho$ is locally Lipschitz away from $0$.
The Thurston measure on $\MF$ can be disintegrated using the
level sets of $\rho$, and this gives a probability measure $\nu_\rho$ on $\PMF$.

\begin{theorem}
\label{theorem:main}
Let $\nu=\nu_\rho$ be a probability measure on $\PMF$ in the Thurston measure class,
where $\rho$ is continuous, homogeneous of degree one, positive away from $0$,
and piecewise smooth in train-track coordinates. Then there exist $a > 0$
and a probability measure $\muG$ on $\Mod$ such that the semigroup generated by the support of $\muG$ is
non-elementary and
\[
\nu=\muG*\nu.
\]
Moreover, $\nu$ is the unique non-atomic
$\muG$--stationary probability measure on $\PMF$.
In addition $\muG$ has a small exponential moment:  for every
$P\in\cT$,
\[
\sum_{\gamma\in\Mod}
\muG(\gamma)\,e^{a \,  d_\cT(P,\gamma P)}<\infty.
\]
\end{theorem}

\begin{remark}[Regularity and limits]
The piecewise smoothness assumption on $\rho$ is used to obtain uniform spike and
Lipschitz estimates for the Radon--Nikodym derivatives of the mapping class group
action on $\PMF$. This assumption is not expected to be essential. Suppose
$\nu_n=\nu_{\rho_n}$ are measures of the above regular type, $\nu_n\to\nu$ weakly,
and the measures $\mu_n$ constructed from
\[
\nu_n=\mu_n*\nu_n
\]
satisfy a uniform exponential \Teich moment bound. Then the sequence $\mu_n$ is tight on
the countable group $\Mod$. Passing to a subsequence, we get $\mu_n\to\mu$, and by
testing against continuous functions on $\PMF$ and using continuity of the $\Mod$-action,
we obtain
\[
\nu=\mu*\nu.
\]
Thus the construction extends to limits of regular densities whenever the required
moment bounds are uniform.
\end{remark}

\begin{remark} 
In fact, our proof proves the following more general statement. Let 
\[
\nu_1=\nu_{\rho_1}
\qquad\text{and}\qquad
\nu_2=\nu_{\rho_2}
\]
be two probability measures on $\PMF$ in the Thurston measure class, where
$\rho_1$ and $\rho_2$ are continuous, homogeneous of degree one, positive away
from $0$, and piecewise smooth in train-track coordinates. Then, there exists
a probability measure $\muG$ on $\Mod$ such that
\[
\nu_1=\muG*\nu_2.
\]
Moreover, $\muG$ has finite exponential
moment with respect to the \Teich metric. 
\end{remark} 

\subsection*{Moment assumptions and singularity}

The moment condition in Theorem~\ref{theorem:main} is stated using the \Teich
metric, rather than a word metric on $\Mod$, and this distinction is essential.
For many natural random walks on $\Mod$, the resulting hitting measure on
$\PMF$ is known to be singular with respect to the Thurston measure class.
Gadre proved this for finitely supported measures \cite{gadre:singular}, and
Gadre--Maher--Tiozzo proved more generally that, for a non-elementary random
walk with finite first moment in a word metric, the hitting measure on $\PMF$ is
singular with respect to the Lebesgue, equivalently Thurston, measure class
\cite{gadre:maher:tiozzo}.

Thus the conclusion of Theorem~\ref{theorem:main} cannot hold under a finite first
word-metric moment assumption. The theorem is close to optimal in this respect.
Fix a basepoint $Y\in\cT$. Word length is at most exponentially distorted
relative to the orbit metric induced by the \Teich metric: there are constants
$A,B>0$ such that
\[
|\gamma|_{\mathrm{word}}
\le
A\exp\bigl(B\,d_{\cT}(Y,\gamma Y)\bigr).
\]
Consequently, any exponential moment in the \Teich metric strong enough to
dominate this distortion would imply finite first moment in the word metric, and
hence would force singularity of the hitting measure. The random walks
constructed here must therefore have heavy tails in every word metric. The point
of Theorem~\ref{theorem:main} is that this obstruction does not
prevent a small exponential moment in the \Teich metric.

\subsection*{Relation to Lyapunov exponents}

Another motivation comes from the study of Lyapunov exponents.  Eskin and Matheus used 
Furstenberg's inverse theorem for lattices in $SL(2, \R)$ to give a coding-free simplicity criterion 
for the Lyapunov exponents of the Kontsevich--Zorich cocycle over \Teich curves 
\cite{eskin:matheus:LC}. The random walk constructed by Furstenberg's inverse theorem 
has Lyapunov exponents that are identical (up to scaling) to those of the \Teich geodesic flow
(this was an observation of Alex Furman). 

%

Theorem~\ref{theorem:main} is designed to provide an analogous source of random walks
in the mapping class group setting. It can be plugged into the proof of
\cite{eskin:matheus:LC} to give another proof of the simplicity of the odd and
even parts of the Lyapunov spectrum for the principal stratum of quadratic
differentials. The corresponding
simplicity result for Abelian differentials is due to Avila and Viana
\cite{Avila-Viana}. For quadratic differentials, the plus and minus spectra are
known to be simple in all components
\cite{bell:delecroix:gadre:gutierrezromo:schleimer:DF}. See also the references
in \cite{bell:delecroix:gadre:gutierrezromo:schleimer:DF} for a more complete
history. A broader boundary-theoretic approach to
Lyapunov spectra developed in \cite{bader:furman:boundary} also recovers 
the result for the principal stratum of quadratic differentials. 

\subsection*{Historical remarks and prior work}
Furstenberg's inverse theorem dates to 1967. Since then,
many analogues of statements from homogeneous dynamics have been
proved in the \Teich setting, but Furstenberg's inverse theorem has
been a notable holdout. The original (highly non-trivial) proof is
based on a discretization of Brownian motion.  In the homogeneous
setting, it is clear that the hitting measure for Brownian motion is
in the Lebesgue measure class. Furstenberg then develops an iterative argument
for approximating Brownian motion using a measure supported on the
lattice $\Gamma$. This approach was considered by Kaimanovich and Masur
in \cite[Section 3]{kaimanovich:masur:rennes}. They study invariant
Markov chains on \Teich space and show that, under their regularity and
positive recurrence assumptions, these can be discretized to random walks
on the mapping class group without changing the hitting measure.
However, this does not yield Theorem~\ref{theorem:main} since it is not
clear whether the hitting measure of their process belongs to the
Thurston measure class.

\subsection*{Outline of the proof}

Our approach is strongly influenced by the work of Connell and Muchnik
\cite{connell:muchnik:QC,connell:muchnik:Gibbs}. They studied the corresponding
inverse problem for groups acting on hyperbolic spaces. In \cite{connell:muchnik:QC},
they show that, under their geometric hypotheses, a continuous
quasiconformal measure on the limit set can be realized as a stationary measure for
a random walk on the group, provided the radial limit set has full measure.
In the CAT$(-1)$ cocompact case, for a Lipschitz quasiconformal measure the resulting random
walk has finite first moment and the boundary with this measure is the Poisson
boundary. In \cite{connell:muchnik:Gibbs}, they extend the construction to measures
coming from generalized Gibbs states.

The Connell--Muchnik argument is based on a positive-basis construction using
boundary functions with controlled decay, called spikes. Hyperbolicity supplies sharp
estimates for these spikes: a spike is large on a small boundary ball and decays at a
controlled rate away from that ball.

In our setting, the spikes are the Radon--Nikodym derivatives $H_\rho$ on the Thurston
boundary defined in \eqref{eq:H}. Our task is to find nonnegative coefficients
$\muG(\gamma)$ satisfying
\[
\sum_{\gamma\in\Mod}\muG(\gamma)H_\rho(\gamma, \param)= \mathds{1}
\]
$\nu_\rho$--almost everywhere. By
\eqref{eq:main:goal:spikes}, this is precisely the stationarity equation,
so such a positive decomposition gives the desired random walk.
For the pseudo-Anosov elements $\gamma$ constructed below, the function
$H_\rho(\gamma,\param)$ has a peak near the attracting foliation.
The peak becomes narrower as the translation length grows. With enough
spike centers, the coarse estimates in Theorem~\ref{theorem:h:estimates}{\rm(a)}
already let us approximate $\mathds{1}$ by positive sums of spikes. We
subtract one such sum, then approximate the positive remainder using
narrower spikes.

However, to obtain even a finite first moment, one needs finer estimates than
hyperbolicity alone provides. Connell and Muchnik call these $Q$-spike
estimates; we formulate the version needed here in Section~\ref{sec:spikes}.
The main technical question is how narrow the next set of spikes must be
after each subtraction. Without a quantitative answer, the translation
lengths needed at successive steps could grow too quickly to obtain a
finite first moment.
The Lipschitz estimate in Theorem~\ref{theorem:h:estimates}{\rm(b)} provides this 
control. It bounds the variation of the remainder and hence the scale needed for 
the next approximation.

We use extremal length to control the size of $H_\rho$. To estimate its
variation, we use homogeneous train-track coordinates in which the mapping
class action on measured foliations is linear (Lemma~\ref{lem:V_of_xi}).

The difficulty is obtaining uniform control near measured foliations that
are endpoints of \Teich geodesics with long excursions into the thin part.
To address this, we work at times when geodesics have good statistics
(Definition~\ref{def:good-statistics}). At such a time the differential
lies in a fixed compact subset of the principal stratum in moduli space,
and the geodesic spends most of every backward interval ending at that
time in this thick part. We also need a boundary metric
adapted to this control. The usual separation time of two rays may occur
during a thin excursion. We instead use the last time at which the rays
are still close and both have good statistics, as in \eqref{eq:def-tau}.
From these times we construct the Gromov metric $d_G$ 
(see Section~\ref{sec:gromov-distance}). Choosing the metric this way is a key
part of the argument as it gives a notion of boundary distance for which
the needed size and Lipschitz estimates hold.

We also need enough pseudo-Anosov elements with these statistical properties.
We obtain them by combining the horospherical large-deviation and mixing
estimates of Appendix~\ref{app:large:deviation}
(Proposition~\ref{prop:appendix-final}) with the closing argument in
Proposition~\ref{prop:properties:gamma:n}. These elements have controlled
translation lengths. Proposition~\ref{prop:extra:properties:gamma:n} shows
that the small balls around their attracting foliations, where the
corresponding spikes are large, cover each boundary point at a large
proportion of the first $n$ approximation scales, except for a set of
$\nu_\rho$--measure exponentially small in $n$.

In Sections~\ref{sec:basis}--\ref{sec:first_moment}, we use these families
to approximate the constant function $1$ by positive sums of spikes.
The Lipschitz estimates allow the translation lengths to grow only linearly
with the step number, while the frequency of coverage just described makes
the integral of the positive remainder decrease exponentially. The resulting
coefficients satisfy the stationarity equation and define a probability
measure with a small exponential \Teich moment.

\subsection*{Order of choosing the constants}
To say that a geodesic has good statistics, we have to decide what part
of the principal stratum we will call thick and what proportion of each
backward interval the geodesic may spend outside it. The less time we
allow outside, the larger the thick part may have to be. The contraction
rate in Corollary~\ref{cor:AGY:half} depends on this choice of thick part,
and we use that rate to define the boundary metric and estimate how the
spikes vary. We therefore fix the allowed proportion first. The order of
choices below ensures that the later estimates work with that proportion.

The surface $S$, base point $P$
in \Teich space and prescribed measure $\nu_\rho$ on $\PMF$ are fixed throughout, and
$w=\dim(\MF)=6g-6+2n$.

\begin{enumerate}
\item \emph{Charts and closing accuracy.}
Fix the chart system and the $\Mod$--compact set $K_0$ containing all
backward chart rays, as in Theorem~\ref{theorem:charts}. This fixes
$\sD_{\rm chart}$, $\rr_{\rm chart}$, $\sD_{\rm temp}$, and the
change-of-charts constant $\kk$ of Lemma~\ref{lem:uniform-change-of-charts}.
The closing lemma (Proposition~\ref{prop:closing:lemma}) approximates a
sufficiently recurrent segment that nearly returns to its starting point
in moduli space by a periodic orbit. Its constants $\ep_{\rm clos},C_{\rm clos}$
depend only on $K_0$. Choose the allowed return error $\epsilon$ as in the proof of
Proposition~\ref{prop:properties:gamma:n}, small enough that closing preserves
the required inclusions between the nested chart domains. This choice
precedes the return-time bounds and the larger compact sets below.
The bound $\sD_0$ on the distance from $P$ to the resulting pseudo-Anosov
axes also depends only on these fixed geometric data.

\item \emph{Time outside compact sets.}
Fix once and for all
\[
\eta_2=\frac14,\qquad
0<\eta_1<\min\left\{\frac1{100},\frac1{512w}\right\},
\qquad \eta=256w\eta_1,
\]
as required in \eqref{eq:spike-error-choice}. Here $\eta_1$ bounds
the proportion of time outside $K_1$ on each backward segment, $\eta_2$
is the larger tolerance for time outside $K_2$ used in defining the Gromov metric, and $\eta$ is
the error in the spike exponents.
The notation $\eta=\eta_1/2$ used locally in the proof of
Proposition~\ref{prop:properties:gamma:n} denotes its input to the appendix,
not this spike error. These choices do not depend on $K_1,K_2$
or on the later contraction rate $\beta$.

\item \emph{Return data and $K_1$.}
Apply Proposition~\ref{prop:appendix-final} with tolerance $\eta_1/2$
to obtain a $\Mod$--compact set $K$, a sampling interval $\tau$, an exponential
decay rate $\kappa$ for the exceptional measure, and a tolerance $\theta_0$
for the measure omitted from the return set. Set $\theta=\theta_0/2$.
With $\epsilon$ fixed, Claim~\ref{claim:ergodic:small:set} gives a set
$\Omega_\theta$ of normalized Masur--Veech measure at least $1-\theta$
whose points flow into the required neighborhoods of the initial chart points within a bounded
time $\rho_\theta=\rho_1$. Choose $K_1$ as in the
proof of Proposition~\ref{prop:properties:gamma:n}, containing the
specified AGY neighborhoods of $K_0$, the initial chart sets, and
$\{g_tq:q\in K,\ 0\le t\le\rho_1\}$.
Next choose the approximation error $\delta_\theta$ small enough to
remain within the closing accuracy after flowing for up to $\rho_1$,
and obtain the lower bound $\tau'_0>\rho_1$ on the sampling interval
used later. The closing lemma now gives the required total time $R_1$
in $K_1$; this threshold affects only how many initial indices are
omitted. None of these choices changes $\epsilon$.

\item \emph{$K_2$ and the metric.}
Choose $K_2$ by Definition~\ref{def:K_2}, large enough to accommodate
the time shifts bounded by $\sD_0$, $\sD_{\rm chart}$, and
$\sD_{\rm temp}$ and the AGY perturbations of radius $1/100$
that arise when comparing rays through $K_1$.
With $K_2$ fixed, Corollary~\ref{cor:AGY:half} gives the contraction
rate $\beta=\beta_{K_2}$ and multiplicative constant $C_\AGY=C_{K_2}$.
We then choose the radius $\rr_1$ and minimum time $\sR_1$ used to
compare the AGY and Gromov distances. In particular, \eqref{eq:R1} uses
the convergence threshold from Proposition~\ref{cor:rays-converge-new}
and the positive margin $\eta_2-\eta_1$. This $\sR_1$ is distinct from the closing
threshold $R_1$.

\item \emph{Spike estimates.}
The multiplicative bounds and additive time delays in the geometric and
Gauss--Manin estimates of \S\ref{sec:Lipschitz} may now depend on the
fixed compact sets. In the Lipschitz proof, first choose the accuracies
and delays that keep the projections inside their chart domains
(\eqref{eq:domain-buffer-choice}--\eqref{eq:domain-delay-choice}),
and only then the step size in \eqref{eq:epsilon-step-choice}.
The coefficients $8$ and $16$ multiplying the time errors in
Proposition~\ref{prop:thin} remain numerical, so the compact-dependent
constants do not force a smaller choice of $\eta_1$. The approximation
argument uses the volume exponent $W$ and the error exponent $\Xi$ given by
\[
W=\frac w\beta,\qquad \Xi=\frac{2\eta}\beta.
\]
The required inequality $\Xi<W$ follows from
$\Xi/W=512\eta_1<1$. No condition comparing $\eta$ with $\beta$ is needed.

\item \emph{Iteration and moment exponent.}
Fix the spike and volume bounds, the bound $\sB$ on the number of balls
covering any point at one stage, and the bound
$\sW_0=e^{\beta(\rho_1+2)}$ on their radius ratios.
Proposition~\ref{prop:subfunction} then gives the required separation
factor $\mm$ between successive scales and a lower bound $\epsilon_0$
on the fraction subtracted from the current remainder on the covered set.
Only now choose $\tau'\in\tau\natls$ above the recurrence, metric,
spike, and scale-separation thresholds listed in
\S\ref{sec:first_moment}, and set $R_j=j\tau'$.
The frequency of coverage and this subtraction make the integrals of
the positive remainders decay at a rate $0<\kappa_1<1$.
Finally choose $0<a<-\log(\kappa_1)/\tau'$, so that this decay
outweighs the factor $e^{aj\tau'}$ in the exponential moment.
These final choices may depend on $\rho$, but require no further change
to $\eta_1$ or to $K_0,K_1,K_2$.
\end{enumerate}

\subsection*{Use of AI tools}

This project began as joint work with Maryam Mirzakhani around 2008. The overall strategy 
and main ideas of the proof have been in place for many years, but completing 
the manuscript took much longer.

We used OpenAI's Codex to help finish the manuscript. It assisted with
routine editing, reference searches, and checking and completing technical
details, including repairing gaps in proofs and checking how the constants
depend on one another.
It also helped simplify some existing proofs, for example in Section~\ref{sec:GM-new}. 
We take full responsibility for the mathematical content of the paper. 

\section{Preliminaries}
\label{sec:prelims}

In this section, we recall the notation and some facts about \Teich space and measured foliations
that we will use. We also rewrite the stationarity equation
\eqref{eq:stationary} as \eqref{eq:main:goal:spikes}, which asks us to express
the constant function $1$ as a nonnegative sum of the functions defined in
\eqref{eq:H}. The rest of the paper constructs these coefficients with the
required exponential moment bound.

Let $S$ be a surface of finite hyperbolic type, i.e., a surface of genus $g$
with $n$ punctures such that
\[
3g - 3 + n > 0.
\]
We also assume $(g,n)\ne(1,1)$. Let $\cT(S)=\cT_{g,n}$ be the \Teich space of $S$. 
When the surface $S$ is
fixed, we simply write $\cT$ for \Teich space. Points in $\cT$ correspond to
marked Riemann surfaces $(\phi_X \from S \to X)$, where $\phi_X$ is defined up
to homotopy. However, we often simply write $X\in\cT$, keeping the marking map
implicit. For two points $(\phi_X \from S \to X)$ and
$(\phi_Y \from S \to Y)$ in $\cT$, the map
$\phi_X\circ \phi_Y^{-1}\from Y\to X$ is referred to as the difference of
markings map, again defined up to isotopy (see \cite{hubbard:TT, IT:Teich}
for background on \Teich space).

Let $\QT=\QT(S)$ be the space of nonzero meromorphic quadratic differentials
on $S$, where we always assume that quadratic differentials have at most
simple poles. We denote a quadratic differential by $q$, or, when emphasizing
the underlying Riemann surface $X$, by $(X,q)$. We denote the space of
unit-area quadratic differentials by $\OT$.

Let $\PQT(S)\subset\QT(S)$ denote the principal stratum: the stratum of
quadratic differentials whose zeros are simple and which have a simple pole at
every marked point. Let $\POT$ be the principal stratum in $\OT$ (see
\cite{strebel:QD} for background on quadratic differentials).

Let $\Mod=\Mod(S)$ be the mapping class group of $S$ (see \cite{farb:MCG} for
background on the mapping class group). The group $\Mod$ naturally acts on all
of the above spaces by changing the marking map. For example, for
$(\phi_X \from S \to X)\in\cT$ and $\gamma\in\Mod$,
\[
\gamma(\phi_X \from S \to X)
=
(\phi_X\circ\gamma^{-1}\from S\to X).
\]
We denote the quotients of $\cT$, $\QT$, and $\PQT$ by $\Mod$, respectively, as
$\cM$, $\QM$, and $\PQM$. Similarly, $\OM$ and $\POM$ denote the quotients of
$\OT$ and $\POT$ by $\Mod$. Abusing notation, we denote all such quotient maps
by $\pi$. For example, if $(X,q)\in\POT$, then
\[
\pi(X)\in\cM
\qquad\text{and}\qquad
\pi(X,q)=\pi(q)\in\POM.
\]

\begin{definition}
\label{def:Mod-compact}
A set $K\subset\POT$ is called \emph{$\Mod$--compact} if it is
$\Mod$--invariant and $\pi(K)$ is compact.
\end{definition}

\subsection*{Measured foliations}

Let $\MF=\MF(S)$ denote the space of measured foliations on $S$, defined up to
homotopy and Whitehead moves. For $F,F'\in\MF$, let
\[
\I(F,F')
\]
denote the intersection number between the two measured foliations (see
\cite{flp:TTS} for background on measured foliations).

For $c>0$, let $c\cdot F$ denote the measured foliation obtained by scaling the
transverse measure of $F$ by a factor of $c$. The space of projectivized
measured foliations is
\[
\PMF=(\MF\setminus\{0\})/\R_{>0}.
\]
For $F\in\MF\setminus\{0\}$, the notation $[F]\in\PMF$ represents the
projective class of $F$.

Given a quadratic differential $(X,q)$, one obtains two transverse measured
foliations, called the horizontal and vertical foliations, by looking at the
local square-root coordinate of $q$. Concretely, away from the zeros of $q$
there is a local coordinate chart from $X$ to $\C$ such that $q=dz^2$. In this
coordinate, the horizontal leaves are given by
\[
\Im(z)=\text{constant}
\]
and have transverse measure $\bigl|\mathrm d(\Im z)\bigr|$, while the vertical
leaves are given by
\[
\Re(z)=\text{constant}
\]
and have transverse measure $\bigl|\mathrm d(\Re z)\bigr|$. Define the maps
\[
F_\pm\from \QT\to\MF
\]
by assigning to a quadratic differential $q$ its vertical measured foliation
$F_+(q)$ and its horizontal measured foliation $F_-(q)$.

This process can often be reversed. We say measured foliations $F$ and $F'$
fill $S$ if for every nonzero measured foliation $G\in\MF-\{0\}$,
\[
\I(F,G)+\I(F',G)>0.
\]
For $F\in\MF$, let $\PMF(F)$ be the set of points $x\in\PMF$ such that $F$
and any, equivalently every, representative of $x$ fill $S$. Let
\[
\NF
=
\Bigl\{
(F,F')\in\MF\times\MF :
\exists\, G\in\MF-\{0\}
\text{ such that }
\I(F,G)+\I(F',G)=0
\Bigr\}
\]
be the set of pairs of measured foliations that do not fill $S$.

\begin{theorem}[Gardiner and Masur \cite{masur:EG}]
\label{theorem:GM}
The map
\[
F_+\times F_-\from \QT\longrightarrow \MF\times\MF-\NF
\]
is a homeomorphism. That is, for $F,F'\in\MF$, if $F$ and $F'$ fill the
surface $S$, then there exists a unique quadratic differential $q$ such that
\[
F_+(q)=F
\qquad\text{and}\qquad
F_-(q)=F'.
\]
We denote this quadratic differential by $Q(F,F')$.
\end{theorem}

The area of the differential determined by a filling pair is
\[
\area\bigl(Q(G,F)\bigr)=\I(G,F).
\]
Thus, for $x=[G]\in\PMF(F)$, the unique unit-area differential with
horizontal foliation $F$ and projective vertical foliation $x$ is
\[
Q_1(x,F)=Q\left(\frac{G}{\I(G,F)},F\right).
\]
This does not depend on the representative $G$ of $x$. In this notation,
the first foliation is vertical, the second is horizontal, and the
subscript $1$ indicates unit area.

Hubbard and Masur \cite{masur:QDF} showed that a point $X\in\cT$ and a
nonzero measured foliation $F$ determine a unique quadratic differential
on $X$ with vertical foliation $F$, and a unique one with horizontal
foliation $F$. We denote them by $Q(F,X)$ and $Q(X,F)$, respectively.
Their unit-area rescalings depend only on $x=[F]$ and are denoted by
$Q_1(x,X)$ and $Q_1(X,x)$.

\subsection*{Extremal length and \Teich distance}

By a curve, we always mean a free homotopy class of non-trivial,
non-peripheral, simple closed curves. Given a quadratic differential $(X,q)$,
there are several notions of length one can associate to a curve $\alpha$. The
extremal length $\Ext_X(\alpha)$ of a curve $\alpha$ is defined to be
\[
\Ext_X(\alpha)
=
\sup_\sigma
\frac{\length_\sigma(\alpha)^2}{\area(\sigma)},
\]
where the supremum is over all metrics $\sigma$ in the conformal class of $X$.

The uniformization theorem implies that there is a unique hyperbolic metric on
$X$. Let $\Hyp_X(\alpha)$ denote the hyperbolic length of $\alpha$ in $X$,
representing the length of the geodesic representative of $\alpha$ in this
metric. Also, $|q|$ defines a singular Euclidean metric on $X$. In this
metric, every curve has a unique geodesic representative except when there is a family
of geodesic representatives foliating a flat cylinder. We denote the length of
the geodesic representative of $\alpha$ in this metric by $\Flat_q(\alpha)$.

For weighted curves, flat and hyperbolic length scale linearly with the
weight, while extremal length scales quadratically. Weighted curves are
dense in $\MF$, and these functions extend continuously to measured
foliations.

\begin{theorem}[Corollary 28 in \cite{rafi:LD}]
\label{theorem:Length-Continuous}
The length functions
\[
\Flat\from \QT\times\MF\to\R_{\ge 0},
\qquad
\Hyp,\ \Ext\from \cT\times\MF\to\R_{\ge 0}
\]
are continuous with respect to all parameters.
\end{theorem}

For a fixed $\ep>0$, let $\cM_{\ge \ep}\subset\cM$ denote the space of Riemann
surfaces where, for every $X\in\cM_{\ge \ep}$ and every curve $\alpha$, we
have $\Hyp_X(\alpha)\ge\ep$. We refer to the set
\[
\cT_{\ge \ep}=\pi^{-1}(\cM_{\ge \ep})
\]
as the \emph{thick part} of \Teich space. Also, let $\OM_{\ge \ep}$ be the
restriction of $\OM$ to $\cM_{\ge \ep}$ and let $\OT_{\ge \ep}$ be the
restriction of $\OT$ to $\cT_{\ge \ep}$. Note that $\cM_{\ge \ep}$ and
$\OM_{\ge \ep}$ are compact. Compactness, continuity, homogeneity, and
$\Mod$--equivariance give the following comparison:

\begin{proposition}
\label{prop:Thick-Length}
For every $(X,q)\in\OT_{\ge \ep}$ and every curve $\alpha$, we have
\[
\sqrt{\Ext_X(\alpha)}
\asymp \Hyp_X(\alpha)
\asymp \Flat_q(\alpha),
\]
where the implied constants depend only on $\ep$.
\end{proposition}

\Teich space can be equipped with several natural $\Mod$--invariant metrics. We
work with the usual \Teich metric. We will use Kerckhoff's formula
\cite{kerckhoff:AG}:
\begin{equation}
\label{eq:kerckhoff}
d_\cT(X,Y)
=
\frac12\sup_\alpha
\log\frac{\Ext_Y(\alpha)}{\Ext_X(\alpha)},
\end{equation}
where the supremum is taken over all curves $\alpha$.

\subsection*{\Teich geodesics}

The \Teich geodesic flow
\[
g_t\from \QT\to\QT,
\qquad t\in\R,
\]
is defined by
\[
g_t(q)
=
\begin{bmatrix}
e^t & 0\\
0 & e^{-t}
\end{bmatrix}
\cdot q
=
Q\bigl(e^t\cdot F_+(q),\,e^{-t}\cdot F_-(q)\bigr).
\]
A \Teich geodesic is the trajectory of a quadratic differential under this
flow. To simplify notation, for $q\in\QT$, we write $g_t(q)$ as $q(t)$. In
particular, $q=q(0)$.

We call $[F_+(q)]$ and $[F_-(q)]$ the forward and backward endpoints of
$q(t)$. Note that these names do not assert that the underlying \Teich geodesic
converges in the Thurston compactification of \Teich space.
For a quadratic differential $(X,q)$, its area equals the extremal length
of either of its measured foliations. In particular, if $(X,q)\in\OT$, then
\[
\Ext_X(F_+(q))=\Ext_X(F_-(q))=1.
\]
Since extremal length is quadratic under scaling,
\[
\Ext_X(cF)=c^2\Ext_X(F),
\]
for unit-area quadratic differentials $q$ and $(X_t,q_t)=g_t(q)$ we have
\begin{equation}
\label{eq:F_-}
\Ext_{X_t}(e^{-t}\cdot F_-(q))=1
\qquad\Longrightarrow\qquad
\Ext_{X_t}(F_-(q))=e^{2t},
\end{equation}
and
\begin{equation}
\label{eq:F_+}
\Ext_{X_t}(e^t\cdot F_+(q))=1
\qquad\Longrightarrow\qquad
\Ext_{X_t}(F_+(q))=e^{-2t}.
\end{equation}
Thus, along a \Teich geodesic, the extremal length of $F_-(q)$, the
horizontal foliation, increases exponentially, while the extremal length of
$F_+(q)$, the vertical foliation, decreases exponentially.

Sometimes, we describe a \Teich geodesic using its horizontal and vertical
foliations. For any $F\in\MF$ and $x\in\PMF(F)$, define
\[
x_F\from\R\to\OT
\]
as the \Teich geodesic given by
\begin{equation}
\label{eq:def:x(t)}
x_F(t)
=
Q_1\bigl(x,e^{-t}\cdot F\bigr)
=
g_t\bigl(Q_1(x,F)\bigr).
\end{equation}

\subsection*{Horospheres}

In $\OT$, a horosphere refers to the set of quadratic differentials whose
vertical or horizontal foliation is a given measured foliation. Hence, there
are two types of horospheres. For $F\in\MF$, define $\Hor_+(F)$ as the set of
area-one quadratic differentials $p\in\OT$ such that
\[
F_+(p)=F,
\]
and define $\Hor_-(F)$ as the set of $p\in\OT$ such that
\[
F_-(p)=F.
\]
When no sign is indicated, we write
\[
\Hor(F)=\Hor_-(F).
\]
Under the geodesic flow,
\[
x_F(t)\in\Hor_-(e^{-t}F).
\]

We may also define a horosphere using a quadratic differential. For $q\in\OT$,
set
\[
\Hor_+(q)=\Hor_+(F_+(q))
\qquad\text{and}\qquad
\Hor_-(q)=\Hor_-(F_-(q)).
\]
We call $\Hor_+(q)$ the stable horosphere and $\Hor_-(q)$ the
unstable horosphere. Fixing only the projective vertical foliation
gives the corresponding center-stable leaf.
Alternatively, for $x\in\PMF$ and $X\in\cT$, define
\[
\Hor_+(x,X)=\Hor_+\bigl(Q_1(x,X)\bigr)
\qquad\text{and}\qquad
\Hor_-(X,x)=\Hor_-\bigl(Q_1(X,x)\bigr).
\]

\subsection*{The $\GL^+(2,\R)$-action}

The space $\QT$ carries a natural $\GL^+(2,\R)$--action. Its positive
diagonal subgroup
\begin{equation}
\label{eq:def:D}
D=\left\{\begin{bmatrix}a&0\\0&b\end{bmatrix}:a,b>0\right\}
\subset\GL^+(2,\R)
\end{equation}
acts by scaling the vertical and horizontal foliations separately:
\[
\begin{bmatrix}a&0\\0&b\end{bmatrix}\cdot Q(G,F)=Q(aG,bF).
\]
We denote the $D$--orbit of a set $E\subset\QT$ by $DE$. For $q\in\OT$,
\begin{equation}
\label{eq:D-orbit}
D\Hor_\pm(q)
=
\{p\in\QT:[F_\pm(p)]=[F_\pm(q)]\}.
\end{equation}
Thus the area and the scale of the fixed foliation may vary, while its
projective class stays fixed.

\subsection*{The Thurston measure class}

The Thurston measure $\muT$ on $\MF$ is invariant under the mapping class
group and homogeneous of degree
\begin{equation} \label{eq:def:w}
w=\dim(\MF)=6g-6+2n.
\end{equation} 
Thus scaling a measurable set by $c>0$ multiplies its measure by $c^w$.
We refer to \cite{ph:TT,bonahon:WT} for the construction of $\muT$.

Let $\rho\from\MF\to\R_{\ge0}$ be a continuous function that is positive away from $0$ 
and is homogeneous of degree one:
\[
\rho(cF)=c\,\rho(F)\qquad(c>0).
\]
Write
\[
S_\rho=\{F\in\MF:\rho(F)=1\},
\qquad
B_\rho=\{F\in\MF:\rho(F)\le1\}.
\]
The set $B_\rho$ is compact and contains a neighborhood of the origin.

The function $\rho$ defines a probability measure on $\PMF$ as follows. 
For a measurable set $E\subset\PMF$, let
\[
\R_+E=\{F\in\MF\setminus\{0\}:[F]\in E\}
\]
be the cone over $E$ and define 
\[
\widetilde\nu_\rho(E)=\muT(\R_+E\cap B_\rho).
\]
Then set 
\[
m_\rho=\widetilde\nu_\rho(\PMF)=\muT(B_\rho),
\qquad \text{and} \qquad
\nu_\rho=\frac{1}{m_\rho}\widetilde\nu_\rho.
\]
Since $0<m_\rho<\infty$, this defines a probability measure on $\PMF$.
Homogeneity gives the polar-coordinate formula (recall the definition of $w$ from 
\eqref{eq:def:w})
\begin{equation} \label{eq:polar}
d\muT=m_\rho\,w\,r^{w-1}\,dr\,d\nu_\rho,
\qquad r=\rho(F).
\end{equation} 

Any two such functions $\rho_1,\rho_2$ have a ratio bounded above and
below by positive constants: the ratio is well defined and continuous on 
the compact space $\PMF$. The associated measures are
therefore equivalent. We call their common measure class the
\emph{Thurston measure class on $\PMF$}.

A natural example for such a function is $\rho(F)=\sqrt{\Ext_P(F)}$ for a fixed $P\in\cT$.
To control the variation of our spike functions, we fix $\rho$ to be
piecewise smooth in train-track coordinates, with the regularity
convention stated before Theorem~\ref{theorem:main}: the subdivision is
locally finite and each restriction extends smoothly across its faces.
In particular, $\rho$ is locally Lipschitz on $\MF\setminus\{0\}$,
and $S_\rho$ is the piecewise smooth boundary of $B_\rho$.

\subsection*{The stationarity equation.}
Recall that the pushforward of $\nu_\rho$ by a mapping class $\gamma$ is defined as: 
\[
(\gamma_*\nu_\rho)(E)=\nu_\rho(\gamma^{-1}E).
\]
By $\Mod$--invariance of $\muT$, this equals
\[
\frac{1}{m_\rho}\,
\muT\bigl(\{F\in\R_+E:\rho(\gamma^{-1}F)\le1\}\bigr).
\]
On the ray through $F\in S_\rho$, the radial cutoff is now
$r\le\rho(\gamma^{-1}F)^{-1}$. Integrating the polar density \eqref{eq:polar} up to this
cutoff we get that, for $x=[F]$ and any nonzero representative $F$,
\begin{equation}
\label{eq:H}
H_\rho(\gamma,x)
=
\frac{d\,(\gamma_*\nu_\rho)}{d\nu_\rho}(x)
=
\left(\frac{\rho(F)}{\rho(\gamma^{-1}F)}\right)^w.
\end{equation}
The ratio is unchanged when $F$ is rescaled. These functions satisfy
the cocycle identity
\[
H_\rho(\gamma_1\gamma_2,x)
=
H_\rho(\gamma_1,x)\,
H_\rho(\gamma_2,\gamma_1^{-1}x).
\]
For example, if $\gamma F=\lambda F$ with $\lambda>0$, then
\[
H_\rho(\gamma,[F])=\lambda^w.
\]
Section~\ref{sec:Lipschitz} gives the estimates that allow us to use
suitable functions $H_\rho(\gamma,\param)$ as spikes. We often suppress
$\rho$ and write $H(\gamma,x)$.

Since $\gamma_*\nu_\rho$ is a probability measure, each of these densities
has integral one:
\[
\int_{\PMF}H_\rho(\gamma,x)\,d\nu_\rho(x)=1.
\]
For a probability measure $\muG$ on $\Mod$, stationarity
\eqref{eq:stationary} is therefore equivalent to
\begin{equation}
\label{eq:main:goal:spikes}
\sum_{\gamma\in\Mod}\muG(\gamma)H_\rho(\gamma, \param)= \mathds{1}
\end{equation}
in $L^1(\nu_\rho)$, or equivalently $\nu_\rho$--almost everywhere.
Conversely, if nonnegative coefficients $\muG(\gamma)$ satisfy this
identity, integrating and applying Tonelli's theorem gives
\[
\sum_{\gamma\in\Mod}\muG(\gamma)=1.
\]
Thus our task is to write $\mathds{1}$ as a nonnegative sum of these densities,
with coefficients satisfying the required exponential moment bound.

\section{Period Coordinates and the Gauss--Manin Connection}
\label{sec:GM}

Period coordinates describe nearby differentials by complex
vectors. Their real and imaginary parts give local coordinates for the
vertical and horizontal measured foliations. We explain this description,
including when these foliations have saddle connections, and show how
Gauss--Manin transport compares the resulting coordinates along a path.
We begin with Abelian differentials, where the coordinates are obtained
by integration.

\subsection*{Period coordinates for Abelian differentials}

Let $\Sigma$ be a closed oriented surface. An Abelian differential
$(Y,\omega)$ consists of a Riemann surface $Y$ homeomorphic to $\Sigma$
and a nonzero holomorphic $1$--form $\omega$. Let
\[
Z=\{z_1,\ldots,z_m\}
\]
consist of its zeros together with any regular marked points we
wish to retain. The stratum $\Omega\mathcal M(\kappa)$ fixes the
orders $\kappa=(k_1,\ldots,k_m)$ at these points, with $k_j=0$
for a regular marked point. If $\omega$ has no zeros, we retain
one regular marked point for the triangulation below.

Away from $Z$, integrating $\omega$ gives flat coordinates whose
transition maps are translations. For an oriented saddle connection
$\delta$, its relative period is
\[
\operatorname{Per}(\omega,\delta)=\int_\delta\omega.
\]
Its absolute value
is the length of $\delta$. The transverse measures of the vertical
and horizontal foliations across this straight edge are, respectively,
the absolute values of its real and imaginary parts.

Choose a triangulation by saddle connections with vertices in $Z$.
Each oriented edge has a relative period, and the three
periods around an oriented triangle sum to zero. Nearby edge
vectors satisfying these relations and the edge identifications give
new triangles. For sufficiently small changes, all triangles remain
nondegenerate with the same orientation, and gluing them in the same
pattern gives a nearby differential in the same stratum. This is the
geometric reason that periods provide local coordinates.

The full list of edge vectors is redundant. To remove this
redundancy, we use relative homology:
\[
H_1(Y,Z;\mathbb Z).
\]
Its elements are represented by relative $1$--cycles, that is, chains whose
boundary lies in $Z$. A basis consists of absolute cycles around the
handles, together with arcs joining one point of $Z$ to the others.
Fix a small enough marked neighborhood $\cV(\omega)$ of $(Y, \omega)$, so that these cycles and arcs
can be followed on nearby surfaces $(Y', \omega') \in \cV(\omega)$. For a basis
$\{\delta_1,\ldots,\delta_N\}$ of $H_1(Y,Z;\mathbb Z)$, the period coordinates are
\[
(Y',\omega')
\longmapsto
\left(
\int_{\delta_1}\omega',\ldots,\int_{\delta_N}\omega'
\right)
\in \mathbb C^N .
\]

Equivalently, the periods define a linear functional on relative
homology, hence a relative cohomology class. Using the local marking
to identify the homology groups, we write the period map as
\[
\cV(\omega)\longrightarrow H^1(Y,Z;\mathbb C),
\qquad
\omega'\longmapsto
\left[\delta\mapsto \int_\delta \omega'\right].
\]
When $\cV(\omega)$ is sufficiently small, this gives
complex coordinates onto an open subset of $H^1(Y,Z;\mathbb C)$.
The triangle relations are the cocycle equations for this cohomology
class, and the gluing construction gives the local inverse. The same
triangulation persists throughout a smaller neighborhood, where both
the periods and their inverse vary continuously. These charts give 
$\Omega\mathcal M(\kappa)$ the structure of an orbifold (see \cite[Section~3.4]{zorich:FS})
(when $(Y, \omega)$ has symmetries, we can pass to the finite quotient by automorphisms). 

\subsection*{Period coordinates for quadratic differentials}

The flat coordinates of a quadratic differential $(X,q)$ have transition maps
$z\mapsto\pm z+c$. Thus integrating the local forms $dz=\sqrt q$ on $X$ 
does not in general give a cohomology class.

To remedy this, we consider the \emph{orienting double cover}
\[
\pi:\widetilde X\to X
\]
of $(X,q)$. Here we fill in the punctures of $X$ and include the
points over them in $\widetilde X$. Let $\cZ(q)$ denote the set of
singularities of $q$, that is, its
zeros and poles. Continuing a square root around an odd-order singularity
reverses its sign. In the principal stratum, the zeros and poles are
all simple, so the orienting double cover is branched exactly over
$\cZ(q)$. By abuse of notation, we also write
$\cZ(q)\subset\widetilde X$ for the ramification locus, which is naturally
identified with the set of singularities of $q$.

The pullback $\tilde q=\pi^*q$ has orientable horizontal and vertical
foliations. Therefore $\tilde q$ is the square of a holomorphic Abelian
differential, well-defined up to sign. We write
\[
\tilde q=\omega_q^2.
\]
The points over simple zeros of $q$ are double zeros of $\omega_q$,
while the points over simple poles are regular marked points.
Indeed, writing $z=w^2$ at a branch point transforms the local models
$z\,dz^2$ and $dz^2/z$ into $4w^4\,dw^2$ and $4\,dw^2$,
respectively. We keep all points of $\cZ(q)$ as marked points in
the relative cohomology, so that saddle connections ending at poles
are included. Since the deck involution $i: \widetilde{X} \to \widetilde{X}$
sends $\omega_q$ to $-\omega_q$, the corresponding period classes 
are odd for the deck involution.

After choosing the sign of $\omega_q$ and choosing a local marking of the
orienting double cover, period coordinates give a local map
\[
\cV(q)\hookrightarrow H^1(\widetilde X,\cZ(q);\C),
\qquad
q'\longmapsto
\left[\gamma\mapsto \int_\gamma \omega_{q'}\right],
\]
for some neighborhood $\cV(q)$ of $q$ in the marked principal stratum
$\PQT(S)$. Changing the sign of the square root multiplies all coordinates by
$-1$.

\begin{proposition}
\label{prop:odd-cohom-dim}
Let $(X,q) \in \PQT(S)$, where $S$ is a surface of genus $g$ with $n$ marked
points, equivalently punctures. Denote by $\widetilde X$ the orienting double
cover of $(X,q)$. Then $\widetilde X$ is a closed surface of genus
\[
   \widetilde{g} \;=\; 4g \;-\; 3 \;+\; n .
\]
Let $i: \widetilde{X} \to \widetilde{X}$ be the deck involution associated to
this double cover. Then the induced map on cohomology
\[
   i^* : H^1(\widetilde{X};\C) \;\longrightarrow\; H^1(\widetilde{X};\C)
\]
splits $H^1(\widetilde{X};\C)$ into a direct sum of $+1$ and $-1$ eigenspaces,
\[
   H^1(\widetilde{X};\C)
   =
   H^1_\even(\widetilde{X};\C)
   \oplus
   H^1_\odd(\widetilde{X};\C).
\]
The \emph{odd} part $H^1_\odd(\widetilde{X};\C)$, the $-1$ eigenspace, has
complex dimension
\[
   \dim_{\C}\!\big( H^1_\odd(\widetilde{X};\C) \big)
   =
   6g \;-\; 6 \;+\; 2n .
\]
\end{proposition}

\begin{proof}
There are $4g-4+n$ simple zeros and $n$ simple poles, hence
$b=4g-4+2n$ branch points. Riemann--Hurwitz gives
\[
2\widetilde g-2=2(2g-2)+b=8g-8+2n,
\]
which proves the genus formula.

The invariant subspace is
$H^1_\even(\widetilde X;\C)=\pi^*H^1(X;\C)$: pullback classes are
invariant, and invariant classes descend to the quotient
$\widetilde X/i=X$. Its dimension is therefore $2g$. Since the
total dimension is $2\widetilde g$, the odd part has dimension
\[
2\widetilde g-2g=6g-6+2n. \qedhere 
\]
\end{proof}

Similarly, the relative cohomology splits into even and odd parts:
\[
H^1(\widetilde X,\cZ(q);\C)
=
H^1_\even(\widetilde X,\cZ(q);\C)
\oplus
H^1_\odd(\widetilde X,\cZ(q);\C).
\]
The long exact sequence of the pair $(\widetilde X,\cZ(q))$ gives
\begin{equation}
\label{eq:Exact}
0
\longrightarrow
H^0(\widetilde X;\C)
\longrightarrow
H^0(\cZ(q);\C)
\xrightarrow{\;\delta\;}
H^1(\widetilde X,\cZ(q);\C)
\xrightarrow{\;r\;}
H^1(\widetilde X;\C)
\longrightarrow
0.
\end{equation}
This sequence is equivariant with respect to the deck involution. Since the
deck involution fixes every point of $\cZ(q)$, it acts trivially on
$H^0(\cZ(q);\C)$. Hence the image of $\delta$ lies in the even part of
$H^1(\widetilde X,\cZ(q);\C)$. Taking the odd part of \eqref{eq:Exact}, we
obtain an isomorphism
\[
r:
H^1_\odd(\widetilde X,\cZ(q);\C)
\longrightarrow
H^1_\odd(\widetilde X;\C).
\]
Thus, for the principal stratum, the odd relative period coordinates may be
viewed canonically as absolute odd cohomology classes.
This uses that every marked singularity is a branch point; in an
arbitrary stratum the relative and absolute odd spaces can differ.
We will also use the real vector space
\begin{equation}
\label{eq:def-Lq}
L(q):=H^1_\odd(\widetilde X;\R)
\cong H^1_\odd(\widetilde X,\cZ(q);\R).
\end{equation}
The relative description lets us evaluate these classes on saddle
connections. The absolute description lets us use the Hodge norm
and cup product on the closed orienting cover.

\begin{theorem}[Masur \cite{masur:IE}]
\label{theorem:Masur}
For any $(X,q)\in\PQT(S)$, there exists a neighborhood $\cV(q)$ on
which a local marking of the orienting covers gives a period coordinate map
\[
\Phi_q\from
\cV(q)
\to
H^1_\odd(\widetilde X,\cZ(q);\C)
\cong
H^1_\odd(\widetilde X;\C),
\]
well-defined up to multiplication by $-1$, that is a homeomorphism onto
an open subset of the indicated cohomology space.
\end{theorem}

Here the remaining sign is the choice of square root. As in the
Abelian case, a triangulation persists on a small neighborhood and
provides a common local marking for the period construction. This
gives continuous coordinates and inverse throughout that neighborhood.
The marked preimages of poles are retained in this construction.
See also \cite[Theorem~5.4]{dumas:skinning} for the period-coordinate
description in the closed-surface case.

\subsection*{Gauss--Manin transport}

To compare periods at different surfaces, we must identify their
cohomology groups. Gauss--Manin transport does this by following
cycles with the surface and keeping the value of a transported
class on each such cycle constant.

First consider marked Riemann surfaces $Y_1,Y_2\in\cT(\Sigma)$,
with markings $\phi_{Y_j}:\Sigma\to Y_j$. Their difference-of-markings
map is
\[
f_{12}=\phi_{Y_2}\circ\phi_{Y_1}^{-1}:Y_1\to Y_2,
\]
defined up to isotopy. The induced transport is
\[
\GM_{Y_1,Y_2}:H^1(Y_1;\C)\to H^1(Y_2;\C),
\qquad
\GM_{Y_1,Y_2}(\alpha)=(f_{12}^{-1})^*\alpha.
\]
Equivalently, it is characterized by
\[
\int_{(f_{12})_*\delta}\GM_{Y_1,Y_2}(\alpha)=\int_\delta\alpha
\qquad\text{for every }\delta\in H_1(Y_1;\Z).
\]
If a finite set $\cZ\subset\Sigma$ is marked as well and if the marking maps 
are taken up to homotopy rel $\cZ$, the same construction gives transport between
$H^1(Y_j,\phi_{Y_j}(\cZ);\C)$. These maps are also defined over $\R$.

The markings make transport independent of the path in \Teich
space. In moduli space however, a loop can return a cycle to a different
homology class. This is called \emph{monodromy}. 

For the orienting covers, a marking of $X$ by $S$ does not determine
a marking of $\widetilde X$ relative to $\cZ(q)$. The zeros of $q$
are not part of the marking of $X$, and choosing a lift also leaves
a deck-involution ambiguity. Thus transport on odd cohomology must
keep track of a path in $\PQT(S)$.

For a path $\xi(t)=(X_t,q_t)$, $0\le t\le1$, in $\PQT(S)$ from
$(X,q)$ to $(Y,p)$, write $i_t$ for the deck involution of
$\widetilde X_t$. On a sufficiently short subinterval, a
triangulation used for period coordinates persists. Matching the
corresponding lifted triangles identifies nearby covers, carries
their ramification points along, and commutes with the deck
involutions. Subdivide $[0,1]$ into finitely many such subintervals
and compose these local identifications, matching them at the
subdivision points and starting with the identity. This gives a
continuous family of homeomorphisms
\[
\widetilde f_t:\widetilde X\longrightarrow\widetilde X_t,
\qquad
\widetilde f_0=\mathrm{Id},
\qquad
\widetilde f_t\circ i_0=i_t\circ\widetilde f_t.
\]
This is what it means to trivialize the family along $\xi$.
The endpoint map $\widetilde f_\xi=\widetilde f_1$ respects the
deck involutions, so its inverse pullback preserves odd cohomology:
\[
\GM_\xi=(\widetilde f_\xi^{-1})^*:
H^1_\odd(\widetilde X;\C)\to H^1_\odd(\widetilde Y;\C).
\]

\begin{lemma}[Gauss--Manin map is well-defined]
\label{lem:GM}
For a path $\xi$ connecting $(X,q)$ to $(Y,p)$ in $\PQT(S)$, the map
\[
\GM_\xi\from
H^1_\odd(\widetilde X;\C)
\to
H^1_\odd(\widetilde Y;\C)
\]
is well-defined by Gauss--Manin transport along $\xi$. It depends only on the
homotopy class of $\xi$ in $\PQT(S)$ relative to the endpoints. Furthermore,
$\GM_\xi$ is defined over $\R$ and commutes with complex conjugation.
\end{lemma}

\begin{proof}
Two choices of the family $\widetilde f_t$ above differ by a
continuous family of self-homeomorphisms of $\widetilde X$ starting
at the identity. Their endpoint maps are therefore isotopic and induce the same map
on cohomology. The same argument over a homotopy of paths proves
homotopy invariance. Pullback by a homeomorphism is defined over
$\R$, which also gives compatibility with complex conjugation.
\end{proof}

For example, consider the loop $q_t=e^{2\pi it}q$, $0\le t\le1$.
The underlying surface $X$ stays fixed. Away from the singularities,
the orienting cover consists of the two choices of a square root
$\zeta$ with $\zeta^2=q_t$. Continuing these choices along the loop
sends $\zeta$ to $e^{\pi it}\zeta$. At $t=1$ this is the deck
involution, so transport around the loop is $-\mathrm{Id}$ on odd
cohomology. Thus even a loop over one fixed base surface can have
nontrivial monodromy.

Within a contractible neighborhood, transport identifies the period
vectors of the same differential in different charts. This is the
content of the next lemma. The period vector of the differential
itself can change as the differential moves; it need not be a
parallel class.

\begin{lemma}
\label{lem:GM-identity}
Let $\xi$ be a path in $\PQT(S)$ connecting $(X,q)$ to $(Y,p)$, and suppose
there is a contractible set $\cU\subset\PQT(S)$ containing $\xi([0,1])$ as
well as the period-coordinate neighborhoods $\cV(q)$ and $\cV(p)$. Consider
the maps
\[
\Phi_q\from \cV(q)\to H^1_\odd(\widetilde X;\C),
\qquad
\Phi_p\from \cV(p)\to H^1_\odd(\widetilde Y;\C),
\]
and the Gauss--Manin transport
\[
\GM_\xi\from H^1_\odd(\widetilde X;\C)\to H^1_\odd(\widetilde Y;\C).
\]
After choosing the signs of $\Phi_q$ and $\Phi_p$ consistently over $\cU$, we
have, on the overlap $\cV(q)\cap\cV(p)$,
\[
\GM_\xi\circ\Phi_q=\Phi_p.
\]
\end{lemma}

\begin{proof}
Let $\widetilde\cU$ be the family of orienting double covers of points in
$\cU$. Since $\cU$ is contractible, we may choose a lift of this family to the
marked Abelian stratum. This choice identifies all odd cohomology groups over
$\cU$ with a single fixed vector space. Under this identification, the two
period-coordinate maps $\Phi_q$ and $\Phi_p$ are given by the same period
vector, and the Gauss--Manin map $\GM_\xi$ is the identity on the fixed vector
space. The claim follows.
\end{proof}

\subsection*{Measured foliations and the local product structure}

To control the variation of $H_\rho$ in \eqref{eq:H}, we will need to
compare nearby measured foliations and their images under a mapping class
in linear coordinates. Train-track charts provide such coordinates on
$\MF$, identified with the space of measured laminations. A foliation in
such a chart is represented by nonnegative branch weights satisfying the
switch equations: at each switch, the total weight on one side equals that
on the other. These charts define the standard piecewise-linear structure;
see Penner--Harer \cite[Chapters~1 and~3]{ph:TT}. In suitable source and
target charts, the mapping class action is linear on these weights. We
now relate these coordinates to period coordinates.

Near a principal differential $q$, the real and imaginary edge periods
record the transverse measures of the vertical and horizontal foliations.
We choose a saddle-connection triangulation with no horizontal or vertical
edge, so that the corresponding train-track weights are linear functions
of the real and imaginary periods throughout a small neighborhood. The
two foliations can then be varied independently, giving the local product
structure in Proposition~\ref{prop:Local-Product}. Fixing the horizontal
foliation, for example, gives real period coordinates for the vertical
foliation. This identification will let us use Gauss--Manin transport to
compare foliation coordinates along paths. The size of the neighborhood
$\cV(q)$ may depend on $q$.

\begin{proposition}
\label{prop:Local-Product}
At every $q_0\in\PQT(S)$ there is a saddle-connection
triangulation with no horizontal or vertical edge. It persists in
a period neighborhood
\[
\cV(q_0)=\Phi^{-1}(U+iV),
\]
where $U,V$ are open convex subsets of a fixed real odd cohomology
space $L$. There are fixed train-track coordinate maps $f,h$ such
that
\[
F_+\bigl(\Phi^{-1}(u+iv)\bigr)=f(u),\qquad
F_-\bigl(\Phi^{-1}(u+iv)\bigr)=h(v).
\]
Their train-track weights are linear functions of $u$ and $v$,
respectively. In particular, with $\cV_+(q_0)=f(U)$ and
$\cV_-(q_0)=h(V)$, the foliation map is a homeomorphism
\[
\cV(q_0)\longrightarrow\cV_+(q_0)\times\cV_-(q_0).
\]
\end{proposition}

\begin{proof}
Suppose first that $\Delta$ is a saddle-connection triangulation
with no horizontal or vertical edge. We explain how it gives the
product coordinates, then construct such a triangulation.

The triangles of $\Delta$ persist in a small period neighborhood of
$q_0$. Shrink it until the real and imaginary parts of every edge
period remain nonzero with fixed signs. In each triangle, the largest
absolute real edge period is the sum of the other two. Place a
trivalent switch in the triangle, with its large branch dual to the
edge with the largest weight and its two small branches dual to the
other edges, as in Figure~\ref{fig:period-switch}.
Gluing dual branches across edges gives a train track
$\tau_+$ carrying the vertical foliation; the branch crossing an edge
has weight equal to the absolute value of that edge's real period.
The horizontal construction gives $\tau_-$. A simple zero gives a
complementary trigon, and a simple pole gives a punctured monogon.
Thus these are maximal train tracks, and the weights of the two
foliations are strictly positive.

\begin{figure}[ht]
\centering
\begin{tikzpicture}[>=stealth, line width=.7pt, font=\small]
  \coordinate (A) at (0,0);
  \coordinate (B) at (2.1,.65);
  \coordinate (C) at (3.1,2.4);
  \draw[->] (A) -- node[below right, pos=.45] {$z_1$} (B);
  \draw[->] (B) -- node[right, pos=.5] {$z_2$} (C);
  \draw[->] (C) -- node[above left, pos=.5] {$z_3$} (A);
  \node at (1.55,-.65) {$z_1+z_2+z_3=0$};

  \coordinate (S) at (6.2,1.2);
  \draw (4.7,1.2) -- node[above] {$x_3$} (S);
  \draw (S) .. controls (6.85,1.2) and (7.15,2.1) ..
    (7.8,2.1) node[right] {$x_1$};
  \draw (S) .. controls (6.85,1.2) and (7.15,.3) ..
    (7.8,.3) node[right] {$x_2$};
  \node at (6.4,-.65) {$x_3=x_1+x_2$};
\end{tikzpicture}
\caption{From triangle periods to switch weights:
$x_j=|\Re z_j|$, with $\Re z_1,\Re z_2>0$ and $\Re z_3<0$.}
\label{fig:period-switch}
\end{figure}
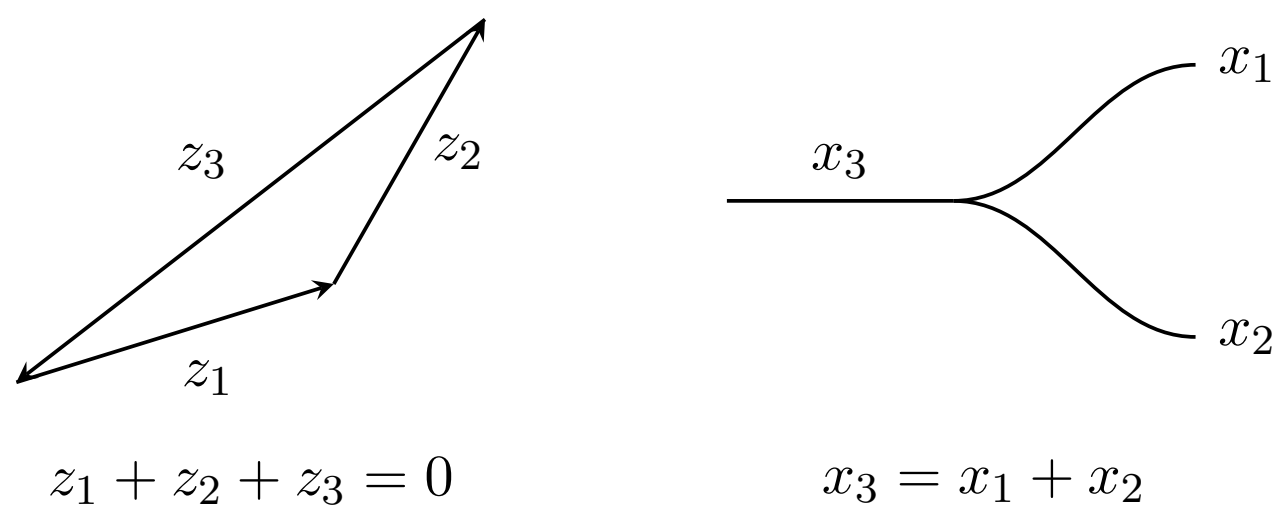

Lift $\Delta$ to the orienting cover and orient one lift of each edge
so that its real period is positive. An odd relative cohomology class
is represented by an odd cellular $1$-cocycle on this triangulation:
there are no relative $0$-cochains, since all vertices belong to the
relative set. Its cocycle equation on each triangle is exactly the
switch equation for $\tau_+$. Conversely, switch weights assign values
to all lifted edges, with the required signs, and give such a cocycle.
This defines a fixed linear isomorphism
\[
E_+:L\longrightarrow W(\tau_+),
\]
where $W(\tau_+)$ is the real vector space of switch weights. We used
the identification of odd relative and absolute cohomology from
\eqref{eq:Exact}; every zero and pole in the principal stratum is a
branch point. Repeated edges in a triangle are counted with their
multiplicities. The imaginary periods give $E_-$ and $\tau_-$ in
the same way.

Positive weights in these maximal track charts give local coordinates
on $\MF$. Choose a small convex product $U+iV$ inside the period
neighborhood. The two fixed linear maps $E_+$ and $E_-$ give the
asserted maps $f,h$. Their images are open, and both maps are
homeomorphisms onto their images. This gives the claimed coordinates
from $\Delta$.

It remains to choose $\Delta$. We follow \cite[Lemma~5.7]{dumas:skinning}.
We avoid both coordinate directions and include surfaces with punctures.
First rescale $q_0$ to area one; rescaling does not change its saddle
connections or their directions. There are only finitely many horizontal
and vertical saddle connections on $q_0$, since each starts along one of
finitely many prongs. By \cite[Theorem~1.2(iii)]{kleinbock:weiss:bounded},
choose a rotation $r_\theta q_0$ whose forward orbit in moduli space has
compact closure inside its stratum, and such that neither original
coordinate direction is contracted by $g_t r_\theta$. The theorem
applies to punctured surfaces, including one-prong singularities at
simple poles. Each of the finitely many original horizontal and vertical
saddle connections has length tending to infinity under this flow.

Choose local marked lifts of the compact orbit closure. At each lift,
choose a triangulation by saddle connections, with vertices at the zeros
and poles, that persists nearby in the stratum. A finite cover by smaller
such neighborhoods provides a number $M$ such that every surface in the
orbit closure has one of these triangulations with all edge lengths at
most $M$. For large $t$, none of the original horizontal or vertical
saddle connections can be an edge of that triangulation on
$g_t r_\theta q_0$. Pull it back by the affine map, and undo the initial
area normalization, to obtain the required triangulation $\Delta$ of
$q_0$. The parameters $\theta$, $t$, and $M$ are chosen for this $q_0$;
the size of the resulting product neighborhood may depend on $q_0$.
\end{proof}

It follows in particular that, for $p_1,p_2\in\cV(q)$,
\begin{equation}
\label{eq:same-real}
F_+(p_1)=F_+(p_2)
\quad\Longleftrightarrow\quad
\Re\Phi_q(p_1)=\Re\Phi_q(p_2).
\end{equation}
The corresponding equivalence holds for the horizontal foliation and
the imaginary periods. These are local assertions. Along a path they
apply by continuation through successive product boxes; they do not
require one product box to contain the entire path.

\subsection*{Homogeneous foliation coordinates and their transport}

We now assign a real vector to a measured foliation, keeping track
of its scale. To do this near $F_+(q)$, fix the horizontal foliation
$F_-(q)$ and take the real period of the resulting differential.
For $G\in\cV_+(q)$, define
\[
\Psi_q(G)=\Re\Phi_q\bigl(Q(G,F_-(q))\bigr).
\]
In the preceding coordinates this is $f^{-1}(G)$. We extend it
homogeneously by setting
\[
\Psi_q:\R_+\cV_+(q)\longrightarrow L(q),\qquad
\Psi_q(cG)=c\,f^{-1}(G).
\]
The extension is well-defined because the train-track weights and
their inverse linear map both respect positive scaling. Its domain
is the cone over the local foliation neighborhood.

\begin{lemma}
\label{lem:Homogeneous}
For $c>0$ and $F\in\R_+\cV_+(q)$,
\[
\Psi_q(cF)=c\Psi_q(F).
\]
For $s>0$, identify the orienting covers of $q$ and $s^2q$
canonically and choose compatible period signs. Then the maps
$\Psi_q$ and $\Psi_{s^2q}$ agree as germs
near their common projective vertical foliation: they agree on
a common conical neighborhood of that ray.
\end{lemma}

\begin{proof}
The first assertion follows from the homogeneous extension. Scalar
multiplication by $s^2$ multiplies periods and both measured foliations
by $s$, so $\Psi_{s^2q}(sG)=s\Psi_q(G)$ near $F_+(q)$.
The first assertion then gives $\Psi_{s^2q}(G)=\Psi_q(G)$ on
a common neighborhood of that ray.
\end{proof}

To describe these coordinates on a unit-area horosphere, first
write the period of $q$ as $a+ib$ and set
\[
A(u,v)=\frac12\int_{\widetilde X}u\wedge v
\qquad (u,v\in L(q)).
\]
The area form on the orienting cover is
$\Re\omega_q\wedge\Im\omega_q$. Since the cover has twice the
area of $q$, we have
\[
\area(q)=A(a,b).
\]
In coordinates continued along the \Teich flow, the period becomes
$e^ta+i e^{-t}b$. This describes the stretching and contraction
while also showing that area is preserved.

Now suppose $q$ has area one. Fixing its horizontal foliation
fixes $b$, so a small neighborhood in its unit-area horosphere
has affine coordinates $a+v+ib$ with $A(v,b)=0$. For a nonzero
measured foliation $F$ such that $p=Q_1([F],F_-(q))$ lies in
this neighborhood, the normalization in Section~\ref{sec:prelims}
gives $F_+(p)=F/\I(F,F_-(q))$. Homogeneity therefore gives
\[
\Psi_q(F)=\I(F,F_-(q))\,\Re\Phi_q(p).
\]
The point $p$ records the projective class $[F]$; the intersection
factor restores the scale of $F$. The same formula defines an extension
whenever these affine leaf coordinates continue injectively to a
larger domain. Lemma~\ref{lem:AGY:vs:cV} will supply such a domain
of uniform intrinsic size. We shrink the original local period box
so that its homogeneous coordinate is the restriction of that
extension.

We next compare the vectors representing a fixed measured foliation
as the base differential moves. The straight paths in the leaf
charts specify how these coordinates are continued. Requiring those
paths to vary continuously makes their continuation agree with
Gauss--Manin transport along the base path.

\begin{lemma}
\label{lem:Psi-is-identity}
Let $\xi:[0,1]\to\POT$ be a continuous path from $q$ to $p$,
and let $F\ne0$. Suppose that
\[
p_t=Q_1([F],F_-(\xi(t)))
\]
belongs to the affine leaf-coordinate domain based at $\xi(t)$
used to define $\Psi_{\xi(t)}(F)$. Suppose also that the straight
leaf paths from $\xi(t)$ to $p_t$ vary continuously with $t$.
With signs chosen consistently along these paths,
\[
\GM_\xi\bigl(\Psi_q(F)\bigr)=\Psi_p(F).
\]
\end{lemma}

\begin{proof}
Parametrize the straight leaf path from $\xi(t)$ to $p_t$ by
$H(t,s)$, $0\le s\le1$. Thus $H(t,0)=\xi(t)$ and
$H(t,1)=p_t$. By assumption this is a continuous homotopy in
the principal stratum; its sides $H(0,s)$ and $H(1,s)$ are the
leaf paths based at $q$ and $p$.
Write $c_t=\I(F,F_-(\xi(t)))$. By the definition of the
homogeneous coordinate, transport along the leaf path identifies
$\Psi_{\xi(t)}(F)$ with $c_t\Re\omega_{p_t}$.

The vertical foliation of $p_t$ is $F/c_t$. Consequently, the
path $c_t^2p_t$ has the fixed vertical foliation $F$. In each local
product box its real period is constant, by
\eqref{eq:same-real}. Continuing this identity through finitely
many boxes shows that $c_t\Re\omega_{p_t}$ is parallel along
$p_t$. Scalar multiplication gives the canonical identification of
the covers used here. Finally, homotopy invariance of Gauss--Manin
transport along the boundary of this homotopy identifies the
endpoint coordinates as claimed.
\end{proof}

The construction following Lemma~\ref{lem:AGY:vs:cV} supplies
uniform leaf domains with the continuous paths required here.
In applications it will therefore suffice to place the points $p_t$
in those domains.

Finally, recall that $\Mod(S)$ acts by changing markings, with the
convention fixed in Section~\ref{sec:prelims}. For
$\gamma \in \Mod(S)$, the induced linear map
$\gamma^*: L(q) \to L(\gamma q)$ is defined as follows: for each
$v \in L(q) = H^1_\odd(\widetilde X;\R)$ and each homology class
$h \in H_1(\widetilde{\gamma X};\R)$, set
\[
\big(\gamma^* v\big)(h) = v\!\big(\gamma^{-1}_* h\big).
\]
The map $\gamma^*$ is compatible with the Gauss--Manin map in the following
sense:

\begin{lemma}
\label{lem:V_of_xi}
For $q\in\PQT(S)$, $\gamma\in\Mod(S)$, and
$F\in\R_+\cV_+(q)$, the homogeneous local coordinates satisfy
\[
\gamma^*\bigl(\Psi_q(F)\bigr)
=\Psi_{\gamma q}(\gamma F).
\]
Here the local charts are taken to be translates, and the
period signs are chosen compatibly. The same identity holds on the
extended unit-area domains constructed after Lemma~\ref{lem:AGY:vs:cV}.
If a path $\xi$ from $q$ to $\gamma q$ satisfies the hypotheses
of Lemma~\ref{lem:Psi-is-identity} for $\gamma F$, then also
\[
\gamma^*\bigl(\Psi_q(F)\bigr)
=\Psi_{\gamma q}(\gamma F)
=\GM_\xi\bigl(\Psi_q(\gamma F)\bigr).
\]
\end{lemma}

\begin{proof}
Changing the marking sends $Q(G,F_-(q))$ to
$Q(\gamma G,F_-(\gamma q))$ and sends its period vector to
$\gamma^*$ of that vector. Taking real parts proves the local
identity, and homogeneity extends it to the cone
$\R_+\cV_+(q)$.
The affine leaf continuation, its norm domain, and intersection
number are likewise equivariant, so the same proof applies on the
extended domain. The final equality is
Lemma~\ref{lem:Psi-is-identity} applied to $\gamma F$.
\end{proof}

\section{The Avila-Gou\"ezel-Yoccoz metric}
\label{sec:AGY}

We use the Avila--Gou\"ezel--Yoccoz norm to measure changes in period
coordinates. It gives uniform charts on small horospherical neighborhoods
and controls projection to a common horosphere. We also record how
horospherical distances change under the \Teich flow, with strict
contraction when the orbit spends enough time in a compact set.

The norm was introduced in \cite{AGY} and its further properties were
developed in \cite{AG}. We work on the marked principal stratum; the
construction is $\Mod(S)$--invariant and descends to moduli space.
By Theorem~\ref{theorem:Masur}, for $(X,q)\in\PQT$ we have
\[
T_q\PQT\cong H^1_\odd(\widetilde X,\cZ(q);\C).
\]
For a tangent vector $\omega$ in this space, define
\begin{equation}
\label{eq:AGY-norm}
\|\omega\|_q
=
\sup_{\substack{\gamma\subset(\widetilde X,\tilde q)\\
\gamma\text{ a saddle connection}}}
\frac{\bigl|\omega([\gamma])\bigr|}{\Flat_{\tilde q}(\gamma)}.
\end{equation}
Here the supremum runs over all saddle connections $\gamma$ of the flat
surface $(\widetilde X,\tilde q)$; equivalently, the denominator is
$\big|\int_\gamma \omega_q\big|$. The points of $\cZ(q)$ are retained as
marked endpoints of saddle connections. In particular, this includes the
preimages of simple poles, which are regular points of $\omega_q$ on the
orienting cover. The norm is continuous
\cite[Propositions~2.10--2.11]{AGY}.

The Abelian strata in \cite[Section~2.2.1]{AGY} allow regular marked
points, so they include these marked preimages of poles. For $q\in\POT$,
the cover has area $2$. To apply their unit-area results, rescale its
Abelian differential by $1/\sqrt2$. Scaling the tangent vector by the
same factor leaves the ratio in \eqref{eq:AGY-norm} unchanged. We use their estimates
on the locus carrying the orienting involution, restricted to odd
cohomology.

For a piecewise $C^1$ path $\xi\colon[0,1]\to\POT$, its length is
\begin{equation}
\label{eq:def:length}
\ell(\xi)
=
\int_0^1 \|\dot\xi(t)\|_{\xi(t)}\,dt.
\end{equation}
The \emph{ambient AGY-metric} on $\POT$ is
\[
D_{\AGY}(q,q')
=
\inf \big\{\ell(\xi):\xi(0)=q,\; \xi(1)=q',\;
                    \xi([0,1])\subset\POT\big\}.
\]
Here ambient paths lie in the unit-area principal stratum, while
horospherical paths will be constrained to a single horosphere.
Completeness follows from \cite[Corollary~2.13]{AGY} on the marked
orienting cover. The locus on which the fixed covering involution
acts by $-1$ is closed; local period charts give convergence in its
intrinsic length metric as well. Its moduli quotient is a complete,
locally compact length space, hence is proper. In particular, closed neighborhoods
of $\Mod$--compact sets are again $\Mod$--compact.

\begin{notation}
For $q \in \POT$ and $r>0$, let
\[
B_{\AGY}(q,r)
=
\bigl\{p\in\POT:\;D_{\AGY}(p,q)< r\bigr\}
\]
be the ambient $AGY$--ball of radius $r$ around $q$. Also, for a set
$E\subset\POT$, let $N_\AGY(E,r)$ be the $r$--neighborhood of $E$ in the
ambient $AGY$--metric.
\end{notation}

\subsection*{AGY norm for $L(q)$}

The same norm restricts to the real space $L(q)$ from
\eqref{eq:def-Lq}. Recall the area pairing
\[
A(u,v)=\frac12\int_{\widetilde X}u\wedge v
\qquad (u,v\in L(q)).
\]
We will need to control the change of area when periods are varied.
The following comparison gives a bound independent of the differential.
Here $\|u\|_{H,q}$ is the $L^2$ norm of the harmonic representative of the
absolute class $u$ on the compact orienting cover $\widetilde X$.

\begin{lemma}
\label{lem:AGY-area}
For $q\in\POT$ and $u\in L(q)$, we have
\begin{equation}
\label{eq:AGY-Hodge-upper}
\|u\|_{H,q}\le\sqrt2\,\|u\|_q.
\end{equation}
Consequently,
\begin{equation}
\label{eq:AGY-area-bound}
|A(u,v)|\le\|u\|_q\,\|v\|_q
\qquad (u,v\in L(q)).
\end{equation}
\end{lemma}

\begin{proof}
The periods of $u$ define a function on the lifts of the marked points.
We extend it to a Lipschitz function on the universal cover. Its
differential will represent $u$ with controlled energy.

Let $Y$ be the universal cover of $\widetilde X$ with its lifted flat
metric, let $\widehat\Sigma$ be the set of lifts of $\cZ(q)$, and fix
$z_0\in\widehat\Sigma$. For $z\in\widehat\Sigma$, define $f(z)$ by
evaluating $u$ on the projection of a path from $z_0$ to $z$. This is
independent of the path because $Y$ is simply connected. A minimizing
flat geodesic between two points of $\widehat\Sigma$ splits into
finitely many saddle connections: the lifted flat metric is proper
and its marked set is locally finite. The definition of the AGY norm,
with $N=\|u\|_q$, therefore gives
\[
|f(z)-f(z')|\le N\,d_Y(z,z').
\]
For a deck transformation $h$ of $Y$, we also have
$f(hz)=f(z)+u(h)$, where $u(h)$ is the absolute period of $u$.

The extension
\[
\widetilde f(y)
=\inf_{z\in\widehat\Sigma}\bigl(f(z)+N\,d_Y(y,z)\bigr)
\]
is finite, is $N$--Lipschitz, and agrees with $f$ on
$\widehat\Sigma$. Reindexing the infimum gives
$\widetilde f(hy)=\widetilde f(y)+u(h)$. Hence $d\widetilde f$
descends to a weakly closed real one-form on $\widetilde X$ representing
the absolute image of $u$, with pointwise flat norm at most $N$ almost
everywhere. Since the orienting cover has area $2$, its energy is at
most $2N^2$. Energy of one-forms on a surface is conformally invariant,
and the harmonic representative
minimizes energy, also among weakly closed $L^2$ representatives by
smooth approximation. This proves \eqref{eq:AGY-Hodge-upper}.

The Hodge Cauchy--Schwarz inequality now gives
\[
|A(u,v)|\le\frac12\|u\|_{H,q}\|v\|_{H,q}
\le\|u\|_q\|v\|_q. \qedhere
\]
\end{proof}

If the period of $q\in\POT$ is $a+ib$, then
$\|a\|_q,\|b\|_q\le1$ directly from \eqref{eq:AGY-norm}. Since
$A(a,b)=1$, the lemma gives
\[
1=A(a,b)\le\|a\|_q\|b\|_q\le1.
\]
Thus $\|a\|_q=\|b\|_q=1$. We will use this normalization when
comparing foliation coordinates along the flow.

We also use the following norm comparison along a path. If $\xi$ is
a piecewise $C^1$ path in $\POT$, of length $L$, and $v$ is an odd
complex cohomology class at its initial point, then
\begin{equation}
\label{eq:AGY-GM-path}
e^{-L}\|v\|_{\xi(0)}
\le\|\GM_\xi(v)\|_{\xi(1)}
\le e^L\|v\|_{\xi(0)}.
\end{equation}
This is \cite[Proposition~5.5]{AG}, restricted to the odd subspace
as above.

\subsection*{The metric on horospheres}

Suppose $x\in\POT$ and $y\in\Hor(x)\cap\POT$. The
\emph{horospherical AGY-metric} is
\[
d_\AGY(x,y)
=
\inf \big\{\ell(\xi):\xi(0)=x,\; \xi(1)=y, \; \xi(t) \in \Hor(x)\cap\POT, 0\le
t \le 1\big\}.
\]
The infimum is over piecewise $C^1$ paths, with length defined by
\eqref{eq:def:length}. This is the metric denoted by $d_{W^u}$ in
\cite[Section~5.1]{AG}.

We use the same definition and notation on $\Hor_+(x)\cap\POT$.

We now continue affine coordinates a uniform distance along a single
horosphere. Write the period of $q\in\POT$ as $a+ib$, and put
\[
E_q=\{v\in L(q):A(v,b)=0,\ \|v\|_q<1/2\}.
\]
The condition $A(v,b)=0$ makes the segment $a+tv+ib$ stay in
the area-one hyperplane.

\begin{lemma}
\label{lem:AGY:vs:cV}
For every $q\in\POT$, continuation of the periods $a+tv+ib$,
$0\le t\le1$, defines an injective coordinate map
\[
\operatorname{Exp}_q:E_q\longrightarrow\Hor(q)\cap\POT.
\]
Its image contains every $p\in\Hor(q)\cap\POT$ with
$d_\AGY(q,p)\le1/50$.
\end{lemma}

\begin{proof}
The radial continuation is the one in
\cite[Proposition~5.3]{AG}, restricted to the odd real subspace
and the area-one hyperplane. Here is the length estimate that
ensures it exists. If $\ell(t)$ is the length of the radial lift,
then \eqref{eq:AGY-GM-path} gives
\[
\ell'(t)\le e^{\ell(t)}\|v\|_q,
\qquad
\ell(t)\le-\log(1-t\|v\|_q).
\]
Completeness continues the lift to $t=1$, and the local product
description in Proposition~\ref{prop:Local-Product} keeps its
horizontal foliation equal to $F_-(q)$. Covering a lifted segment
by finitely many period charts also gives continuous dependence
on $q,v,t$, with compatible local markings.

We will use local inverses for this family before proving global
injectivity. Fix a radial segment, and transport its endpoint
period chart back to $L(q)$. In these coordinates the endpoint
period is $a(q)+v+ib(q)$. The horizontal-coordinate maps at the
two ends agree near $b(q)$: nearby radial lifts preserve the
horizontal foliation, and Proposition~\ref{prop:Local-Product}
identifies that foliation with the imaginary period. Thus, for
a nearby pair $(q',p')$ with the same horizontal foliation and
area one, the imaginary period of $p'$ is $b(q')$, and its inverse
vector is its real period minus $a(q')$. This vector satisfies
$A(v',b(q'))=0$ and remains in the open norm ball of radius $1/2$.
The inverse depends continuously on the pair. This argument
uses local period charts, without a smoothness assumption on
foliation coordinates.

First suppose that $H=F_-(q)$ has the principal singularity type
and no saddle connections. Every differential with horizontal
foliation a positive multiple of $H$ is then principal.
Hubbard--Masur gives $D\Hor_-(H)\cong\cT(S)\times\R_{>0}$.
If $r$ in this set has area $B$ and horizontal foliation $cH$,
the normalization
\[
r\longmapsto
\begin{bmatrix}c/B&0\\0&1/c\end{bmatrix}r
\]
together with the parameters $B,c$, identifies this set with
$\Hor_-(H)\times\R_{>0}^2$.
Hence $\Hor_-(H)$ is contractible, and every loop in it has
trivial Gauss--Manin transport with the square-root sign
continued consistently. Two radial paths with the same endpoint
therefore have the same endpoint period in $L(q)$, so their
displacement vectors agree.

For an arbitrary $q$, suppose that distinct $v,w\in E_q$ have
the same image $p$. Choose disjoint local inverse branches
about $(q,v)$ and $(q,w)$ as above. Take $q_n\to q$ with
horizontal foliations having no saddle connections; small
rotations avoiding the countably many saddle-connection
directions suffice. Set
\[
p_n=Q_1\bigl([F_+(p)],F_-(q_n)\bigr).
\]
One local product neighborhood at $p$ shows that $p_n$ is
principal for large $n$ and tends to $p$. The pair $(q_n,p_n)$
then lies in the domains of both inverse branches. This gives
distinct vectors in $E_{q_n}$ with the same image $p_n$,
contradicting the generic case. Thus $\operatorname{Exp}_q$
is injective at every center.

Finally, let $d_\AGY(q,p)\le1/50$. Choose an intrinsic path
from $q$ to $p$ of length less than $1/40$ and continue its
periods from $q$. Its imaginary period is $b$, and its real
displacement $u(s)$ satisfies
\[
\|u(s)\|_q\le e^{\ell(s)}-1<e^{1/40}-1<1/2
\]
by \eqref{eq:AGY-GM-path}. Local uniqueness of period lifts
identifies the path with $\operatorname{Exp}_q(u(s))$.
In particular, $p$ belongs to the image.
\end{proof}

These charts give the larger domains for $\Psi_q$ promised in
Section~\ref{sec:GM}. For $F\ne0$ such that
$p=Q_1([F],F_-(q))=\operatorname{Exp}_q(v)$, define
\begin{equation}
\label{eq:Psi-leaf-domain}
\Psi_q(F)=\I(F,F_-(q))\,(a+v).
\end{equation}
The factor restores the scale of $F$, since
$F_+(p)=F/\I(F,F_-(q))$. This agrees with the local homogeneous
coordinate after its original neighborhood is shrunk.

The joint local inverses in the proof show that the image of the
family of charts is open in the space of pairs with the same
horizontal foliation and area one. Consequently the set of pairs
$(q,F)$ on which \eqref{eq:Psi-leaf-domain} is defined is open.
On this set, $\Psi_q(F)$ and the radial paths
$t\mapsto\operatorname{Exp}_q(tv)$ depend continuously on $(q,F)$
in local markings. In particular, these paths satisfy the
continuity requirement in Lemma~\ref{lem:Psi-is-identity}.
All the constructions commute with changing markings. Reversing
the square-root sign reverses every period vector, so they also
respect the sign convention in Lemma~\ref{lem:V_of_xi}.

We will also use the stable version of the construction. Rotating
periods by $a+ib\mapsto-b+ia$ preserves the AGY norm and exchanges
the two foliations. Thus, if $A(a,u)=0$ and $\|u\|_q<1/2$,
the segment $a+i(b+tu)$ continues in $\Hor_+(q)\cap\POT$,
has length at most $-\log(1-\|u\|_q)$, and depends continuously
on its parameters with the same local inverse properties.

\begin{proposition}
\label{prop:AGY-L(q)}
Let $q\in\POT$, and suppose $p_1,p_2\in\Hor(q)\cap\POT$ satisfy
$d_\AGY(q,p_i)\le1/200$ for $i=1,2$. Put $F_i=F_+(p_i)$. Then
\[
\frac12 d_\AGY(p_1,p_2)
\le\|\Psi_q(F_1)-\Psi_q(F_2)\|_q
\le2d_\AGY(p_1,p_2).
\]
\end{proposition}

\begin{proof}
Use Lemma~\ref{lem:AGY:vs:cV} to write
$p_i=\operatorname{Exp}_q(v_i)$. Since the two points have area one
and the same horizontal foliation as $q$, their foliation vectors
are $a+v_i$ in these coordinates. Norms throughout the affine domain
are comparable to the norm at $q$ by a factor of two: the radial
paths have length less than $\log2$, so this follows from
\eqref{eq:AGY-GM-path}.

The straight segment from $v_1$ to $v_2$ lies in the convex set
$E_q$. Its image has length at most $2\|v_1-v_2\|_q$, which gives
the first inequality. Conversely, a nearly minimizing path from
$p_1$ to $p_2$ has length less than $1/100+\varepsilon$. Every
point on it is within $3/200+\varepsilon<1/50$ of $q$, for small
$\varepsilon$, so the entire path lies in this one leaf chart.
Integrating its period velocity and letting the length error tend
to zero gives $\|v_1-v_2\|_q\le2d_\AGY(p_1,p_2)$.
\end{proof}

\subsection*{The AGY metric and the \Teich geodesic flow}

\Teich geodesics are $1$--Lipschitz paths in the AGY metric
\cite[Lemma~5.1]{AG}: for $q\in\POT$ and $t\in\reals$,
\begin{equation}
\label{eq:gt-Lip}
D_\AGY(g_t(q),q)\le |t|.
\end{equation}
On an unstable horosphere, backward flow does not increase distances.
The forward flow has the following upper bound \cite[Lemma~5.2]{AG}.

\begin{theorem}
\label{theorem:properties:AGY}
For $F\in\MF$ and $x,y\in\PMF(F)$ with
$x_F(0),y_F(0)\in\POT$, and for $s\le t$, we have
\begin{equation}
\label{eq:Monotone}
d_\AGY\big(x_F(s),y_F(s)\big)
\le
d_\AGY\big(x_F(t),y_F(t)\big)
\le
e^{2(t-s)}d_\AGY\big(x_F(s),y_F(s)\big).
\end{equation}
\end{theorem}

\begin{proof}
For $h\ge0$, the flow multiplies an unstable tangent vector by $e^h$,
while the length of each saddle connection changes by a factor between
$e^{-h}$ and $e^h$. The definition of the AGY norm therefore gives
\[
\|v\|_q\le\|Dg_h(v)\|_{g_hq}\le e^{2h}\|v\|_q
\qquad(v\in T_q\Hor(q)).
\]
Integrating along paths in the horosphere, and using that the flow
maps these paths bijectively to paths in the image horosphere, proves
the two distance inequalities.
\end{proof}

Gauss--Manin transport keeps a cohomology vector fixed while the norm
changes. For the vector representing the vertical foliation, homogeneity
gives an exact formula.

\begin{lemma}
\label{lem:GM:near:geodesic}
Let $q\in\POT$, $\tau\ge0$, and let $\xi$ be the segment
$q(t)=g_t(q)$, $0\le t\le\tau$. For $\bfv\in L(q)$,
\begin{equation}
\label{eq:GM-exp}
e^{-\tau}\|\bfv\|_q
\le
\|\GM_\xi(\bfv)\|_{q(\tau)}
\le
e^\tau\|\bfv\|_q.
\end{equation}
If $[F]=[F_+(q)]$, then
\begin{equation}
\label{eq:F-decay}
\|\Psi_q(F)\|_q
=e^\tau\|\Psi_{q(\tau)}(F)\|_{q(\tau)}.
\end{equation}
\end{lemma}

\begin{proof}
Under the flow, saddle connections correspond bijectively and their
lengths change by factors between $e^{-\tau}$ and $e^\tau$. Their
evaluations on the transported class are unchanged, proving
\eqref{eq:GM-exp}. For \eqref{eq:F-decay}, write
$F=cF_+(q)=ce^{-\tau}F_+(q(\tau))$. The real period has AGY norm
one throughout $\POT$, as observed after
Lemma~\ref{lem:AGY-area}. Homogeneity of $\Psi$ therefore gives
$\|\Psi_q(F)\|_q=c$ and
$\|\Psi_{q(\tau)}(F)\|_{q(\tau)}=ce^{-\tau}$.
\end{proof}

Time spent in a compact set gives strict contraction. We first record
the tangent-vector estimate from \cite[Proposition~4.3]{AG}.

\begin{proposition}
\label{prop:AG43}
For every $\Mod$--compact set $K'\subset\POT$, there is $T=T(K')>0$
such that, if $x\in\POT$ and $t\ge0$ satisfy
\begin{equation}
\label{eq:stays:in:Kprime}
\bigl|\{s\in[0,t]:g_{-s}x\in K'\}\bigr|\ge T,
\end{equation}
then
\[
\|Dg_{-t}v\|_{g_{-t}x}\le\frac12\|v\|_x
\qquad(v\in T_x\Hor(x)).
\]
Here $T_x\Hor(x)$ is the real period subspace
$\{u\in L(x):A(u,\Im\Phi_x(x))=0\}$.
\end{proposition}

\begin{proof}
The cited statement assumes that both endpoints lie in the compact
set. Trim the interval to its first and last visits to $K'$. These
visits are attained because $K'$ is closed
(Definition~\ref{def:Mod-compact}), and trimming removes no time spent
in $K'$. Apply the cited contraction estimate on this interval.
Backward flow does not increase the unstable tangent norm on either
omitted piece, by the calculation above.
\end{proof}

To pass from tangent vectors to distances, we allow the points of a
short connecting path to visit a slightly larger compact set.

\begin{corollary}
\label{cor:AG43}
For every $\Mod$--compact set $K\subset\POT$, there is $T=T(K)>0$
such that the following holds. Suppose $x\in K$,
$y\in\Hor(x)\cap\POT$, and $d_\AGY(x,y)\le\ooh$. If $t>0$ and
\begin{equation}
\label{eq:stays:in:K}
\bigl|\{s\in[0,t]:g_{-s}x\in K\}\bigr|\ge T,
\end{equation}
then
\[
d_\AGY(g_{-t}x,g_{-t}y)\le\frac12d_\AGY(x,y).
\]
\end{corollary}

\begin{proof}
Let $K'$ be the closed $1/100$--neighborhood of $K$ in the ambient
AGY metric. This is $\Mod$--compact by properness of the quotient metric.
Choose $T=T(K')$ from
Proposition~\ref{prop:AG43}.

First suppose $d_\AGY(x,y)<1/100$. Choose a piecewise $C^1$ path
$\gamma$ in $\Hor(x)\cap\POT$ from $x$ to $y$, of length less than
$d_\AGY(x,y)+\varepsilon<1/100$. Backward nonexpansion implies
that every point of $g_{-s}\gamma$ is within $1/100$ of $g_{-s}x$.
Whenever $g_{-s}x\in K$, all these points therefore lie in $K'$.
Apply Proposition~\ref{prop:AG43} to the tangent vectors along
$\gamma$ and integrate:
\[
d_\AGY(g_{-t}x,g_{-t}y)
\le\ell(g_{-t}\gamma)
\le\frac12\ell(\gamma)
<\frac12\bigl(d_\AGY(x,y)+\varepsilon\bigr).
\]
Let $\varepsilon$ tend to zero. If $d_\AGY(x,y)=1/100$, truncating
nearly minimizing paths gives points $y_j$ with
$d_\AGY(x,y_j)<1/100$ and $d_\AGY(y_j,y)\to0$.
Backward nonexpansion also gives
$d_\AGY(g_{-t}y_j,g_{-t}y)\to0$, so the same estimate follows
by passage to the limit.
\end{proof}

Repeated halving gives exponential contraction when the orbit spends
a fixed proportion of its time in $K$.

\begin{corollary}
\label{cor:AGY:half}
For every $\Mod$--compact set $K\subset\POT$, there exist constants
$0<\beta_K<1$ and $C_K\ge1$, depending only on $K$, with the following
property. Let $F\in\MF$ and $x,y\in\PMF(F)$ with
$x_F(0),y_F(0)\in\POT$. If $s\le t$,
\[
d_\AGY(x_F(t),y_F(t))\le\ooh,
\qquad
\bigl|\{u\in[s,t]:x_F(u)\in K\}\bigr|\ge\frac12(t-s),
\]
then
\begin{equation}
\label{eq:Decay}
d_\AGY\big(x_F(s),y_F(s)\big)
\le C_K e^{-\beta_K(t-s)}
d_\AGY\big(x_F(t),y_F(t)\big).
\end{equation}
\end{corollary}

\begin{proof}
Let $T=T(K)$ be given by Corollary~\ref{cor:AG43}, and put
\[
m(u)=\bigl|\{v\in[s,u]:x_F(v)\in K\}\bigr|,
\qquad M=m(t),\qquad N=\lfloor M/T\rfloor.
\]
For $1\le j\le N$, let $s_j$ be the first time at which
$m(s_j)=jT$, and put $s_0=s$. Continuity of $m$ gives these
times. Since $K$ is closed, minimality of $s_j$ implies
$x_F(s_j)\in K$: otherwise $m$ would be constant on a neighborhood
of $s_j$. Each interval $[s_{j-1},s_j]$ contains $T$ units of
time in $K$.

Starting at time $t$, apply backward nonexpansion down to $s_N$
and Corollary~\ref{cor:AG43} on each of the $N$ full intervals.
All distances remain at most $1/100$. With no full interval, use
nonexpansion on $[s,t]$. Since $N\ge M/T-1$, in either case
\begin{align*}
d_\AGY(x_F(s),y_F(s))
&\le 2e^{-(\log2)M/T}d_\AGY(x_F(t),y_F(t))\\
&\le 2e^{-(\log2)(t-s)/(2T)}d_\AGY(x_F(t),y_F(t)).
\end{align*}
Thus, for example, $C_K=2$ and
$\beta_K=\min\{1/2,(\log2)/(2T)\}$ satisfy the assertion.
\end{proof}

For points farther apart, each full block of compact time removes a
fixed amount of distance, as long as the final distance is at least
$1/100$.

\begin{corollary}
\label{cor:AGY:endpoint:contraction}
Let $K\subset\POT$ be $\Mod$--compact, and let $T=T(K)$ be as in
Corollary~\ref{cor:AG43}. Suppose $x\in\POT$,
$y\in\Hor(x)\cap\POT$, and $t\ge0$ satisfy
\[
d_\AGY(g_{-t}x,g_{-t}y)\ge\ooh,
\qquad
\bigl|\{s\in[0,t]:g_{-s}x\in K\}\bigr|\ge2T.
\]
Then
\[
d_\AGY(g_{-t}x,g_{-t}y)
\le d_\AGY(x,y)
-\frac{1}{400T}\bigl|\{s\in[0,t]:g_{-s}x\in K\}\bigr|.
\]
\end{corollary}

\begin{proof}
We may assume $d_\AGY(x,y)<\infty$, since otherwise the conclusion
is immediate. Put
\[
m(u)=\bigl|\{s\in[0,u]:g_{-s}x\in K\}\bigr|,
\qquad M=m(t),\qquad N=\lfloor M/T\rfloor.
\]
Let $t_0$ be the first visit to $K$, and for $1\le j\le N$
let $t_j$ be the first time with $m(t_j)=jT$. These times exist,
and $g_{-t_j}x\in K$ by closedness of $K$ and the same minimality
argument as above. Each interval $[t_j,t_{j+1}]$ contains $T$
units of compact time. Write $x_j=g_{-t_j}x$ and
$y_j=g_{-t_j}y$.

Backward nonexpansion and the hypothesis at time $t$ give
$d_\AGY(x_j,y_j)\ge1/100$. For $0\le j<N$, take a path in
$\Hor(x_j)\cap\POT$ from $x_j$ to $y_j$ of length less than
$d_\AGY(x_j,y_j)+\varepsilon$, and let $z_j$ be its point at
length $1/100$. Since $d_\AGY(x_j,z_j)\le1/100$,
Corollary~\ref{cor:AG43} sends these two endpoints to distance
at most $1/200$ over the next block. The rest of the path does
not increase in length. Therefore
\[
d_\AGY(x_{j+1},y_{j+1})
\le\frac1{200}
+d_\AGY(x_j,y_j)+\varepsilon-\frac1{100}.
\]
Letting $\varepsilon$ tend to zero gives a decrease of $1/200$
on each block. The initial and final omitted pieces are
nonexpanding, so
\[
d_\AGY(g_{-t}x,g_{-t}y)
\le d_\AGY(x,y)-\frac{N}{200}
\le d_\AGY(x,y)-\frac{M}{400T},
\]
where $N\ge M/T-1\ge M/(2T)$ because $M\ge2T$.
\end{proof}

\subsection*{Lipschitz projection}

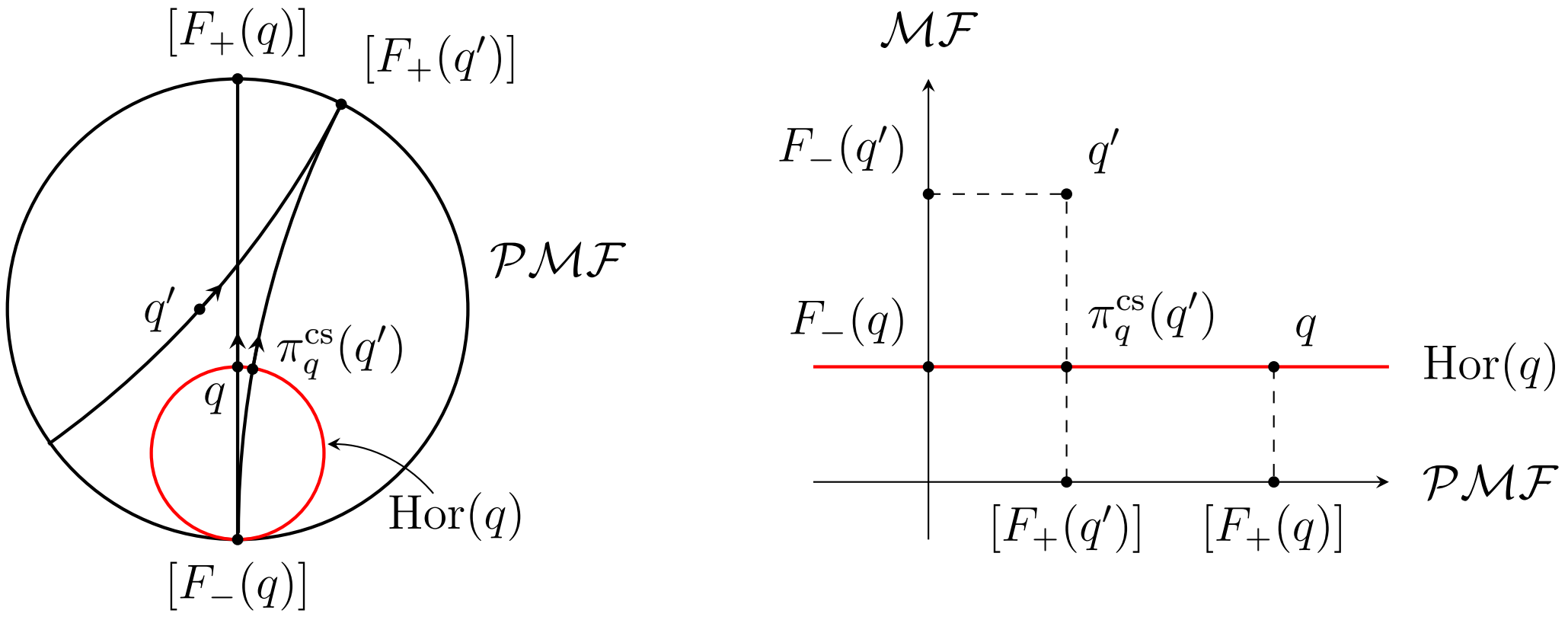
\begin{figure}[ht]
\setlength{\unitlength}{0.01\linewidth}
\begin{picture}(100, 32)

\put(5,-3){
\begin{tikzpicture} [>=stealth]

\tikzstyle{vertex} =[circle,draw,fill=black,thick, inner sep=0pt,minimum size= .7 mm]

\draw[thick] (0,0) circle (2);
\draw[thick] (0,-2) -- (0,2);

\draw[thick] (.9,1.78) arc (153:180:8.2);
\draw[thick] (.9,1.78)  arc(-27: -54.5: 8.2);
\draw[red, thick] (0,-1.25) circle (.75);

\node [vertex] (q) at (0,-.5)  {};
\node at (-.2 ,-1.2)    [label=$q$] {};
\draw [thick, ->] (q) -- (0, -.2);
\node [vertex] (p) at (.13,-.52)  {};
\draw [thick, ->] (p) -- (.19, -.22);
\node at (.9 ,-.9)    [label=$\pi_q^{\rm cs}(q')$] {};

\node [vertex] (qq) at (-.33,0)  [label=left:$q'$] {};
\draw [thick, ->] (qq) -- (-.13, .22);

\node [vertex] at (0,2)  [label=90:${[F_+(q)]}$] {};
\node [vertex] at (0,-2)  [label=-90:${[F_-(q)]}$] {};
\node [vertex] at (.9,1.78)  [label=45:${[F_+(q')]}$] {};
\node at (1.9,0) [label=45:$\PMF$] {};
\node at (1.9 ,-1.3)   [label=-90:$\Hor(q)$] {};

\draw[->]  (1.7, -1.6) arc (40: 90: 1.2);


\draw[->] (5,-1.5) -- (10,-1.5);
\node at (10 ,-1.5)   [label=right:$\PMF$] {};
\draw[->] (6,-2) -- (6,2);
\node at (6 ,2)   [label=90:$\MF$] {};

\draw[red, thick] (5,-.5) -- (10,-.5);
\node at (10 ,-.5)   [label=right:$\Hor(q)$] {};

\node[vertex] at (6 ,-.5)  [label=135:$F_-(q)$] {};
\node[vertex] at (6 ,1)  [label=135:$F_-(q')$] {};

\draw[dashed] (6, 1) -- (7.2,1) -- (7.2,-1.5);
\node[vertex] at (7.2 ,1)  [label=45:$q'$] {};
\node[vertex] at (7.2 ,-.5)   [label=45:$\pi_q^{\rm cs}(q')$] {};
\node[vertex] at (7.2 ,-1.5)  [label=-90:{$[F_+(q')]$}] {};

\draw[dashed] (9,-.5) -- (9,-1.5);
\node[vertex] at (9 ,-.5)  [label=45:$q$] {};
\node[vertex] at (9 ,-1.5)   [label=-90:{$[F_+(q)]$}] {};


\end{tikzpicture}}
\end{picture}
\caption{The map $\pi_q^{\rm cs}$.}
\label{fig:CS}
\end{figure}

We now estimate projection along the center-stable foliation. Given
a quadratic differential $q \in \POT$, let $\cQ_q \subset \POT$ be the
set of all quadratic differentials $q'$ such that $F_+(q')$ and $F_-(q)$ fill
$S$. There is a map
\[
\pi_q^{\rm cs} \from \cQ_q \to \Hor(q),
\]
defined by
\[
\pi_q^{\rm cs}(q')
=
Q_1\bigl([F_+(q')],F_-(q)\bigr).
\]
Recall that $\Hor(q)=\Hor_-(q)$. Thus $\pi_q^{\rm cs}(q')$ is the unique
point in $\Hor(q)$ with the same projective vertical foliation as $q'$. In
other words, $\pi_q^{\rm cs}(q')$ is the unique point in the intersection of
the center-stable foliation of $q'$ with the unstable foliation of $q$; see
Figure~\ref{fig:CS}. Recall that $B_{\AGY}(q',r)\subset\POT$ is the ball of
radius $r$ in the $D_\AGY$--metric centered at $q'$. The following
estimate also ensures that this projection is defined on the stated
neighborhood and remains in the principal stratum.

\begin{theorem}
\label{theorem:cs-projection}
For $q,q'\in\POT$, assume $D_\AGY(q,q')\le\oof$. Then the
projection $\pi_q^{\rm cs}$ is defined on $B_\AGY(q',\oof)$, its
image lies in $\Hor(q)\cap\POT$, and
\[
d_\AGY\bigl(\pi_q^{\rm cs}(q_1),\pi_q^{\rm cs}(q_2)\bigr)
\le4D_\AGY(q_1,q_2)
\qquad\text{for }q_1,q_2\in B_\AGY(q',\oof).
\]
\end{theorem}

\begin{proof}
In local period coordinates, the projection keeps the real period up
to the factor needed to restore area one, and replaces the imaginary
period by that of $q$. We first show that this operation can be
continued along a short path. This proves that the projected
differential exists in the principal stratum, including when one of
its foliations has saddle connections. We then estimate its derivative.

Let $\kappa$ be a piecewise $C^1$ path in $\POT$ starting at $q$,
of length less than $\log(11/10)$. Transport the odd cohomology groups
along $\kappa$, and write the period at $r_s=\kappa(s)$ as
$a_s+ib_s$. Put $b=b_0$. If $L_s$ is the length of the initial
segment of $\kappa$ ending at $r_s$, then
\eqref{eq:AGY-GM-path}, applied to the imaginary part of its velocity,
gives
\[
\delta_s:=\|b-b_s\|_{r_s}
\le\int_0^s e^{L_s-L_u}\|\dot\kappa(u)\|_{r_u}\,du
=e^{L_s}-1<\frac1{10}.
\]
By Lemma~\ref{lem:AGY-area},
\[
h_s:=A(a_s,b),\qquad |h_s-1|\le\delta_s.
\]
The candidate for the projected period is therefore $a_s/h_s+ib$.

We construct this point in two steps: move along a stable leaf, then
flow to the required normalization. Set
\[
u_s=\frac{b}{h_s}-b_s.
\]
Then $A(a_s,u_s)=0$, and
\[
\|u_s\|_{r_s}
\le\frac{\|b-b_s\|_{r_s}+|1-h_s|\|b_s\|_{r_s}}{h_s}
\le\frac{2\delta_s}{1-\delta_s}<\frac12.
\]
The stable construction following Lemma~\ref{lem:AGY:vs:cV}
therefore gives a principal differential with period $a_s+ib/h_s$.
Flow for time $-\log h_s$ to obtain the period $a_s/h_s+ib$.
Parametrize this concatenation by $H_s:[0,1]\to\POT$. Its length satisfies
\[
\ell(H_s)
\le-\log(1-\|u_s\|_{r_s})+|\log h_s|
\le-\log(1-3\delta_s),
\]
where we used $|\log h_s|\le-\log(1-\delta_s)$.

These connecting paths depend continuously on $s$ and form a
homotopy from $\kappa$ to the path of their endpoints. Transport
over the parameter square gives compatible period coordinates
throughout this family. Each connecting path preserves the
projective vertical foliation $[F_+(r_s)]$. The endpoint path
has constant imaginary period $b$ and starts at $q$, so its
horizontal foliation is $F_-(q)$, by
Proposition~\ref{prop:Local-Product}. These local foliation
identities hold also at saddle connections. Uniqueness of the
differential determined by its two foliations now gives
\[
H_s(1)=Q_1\bigl([F_+(r_s)],F_-(q)\bigr)
=\pi_q^{\rm cs}(r_s).
\]
In particular, the foliations fill, the projected differential is
principal, and it is independent of the path chosen from $q$.
The same uniqueness allows us to forget the auxiliary markings
of zeros and orienting covers.

We next estimate the derivative, using the local period formula
$a/h+ib$. At $r=r_s$, write $a=a_s$, $h=h_s$, and $\delta=\delta_s$.
For a tangent vector $\xi=u+iv$ to the unit-area locus, the
derivative of the projection has real period
\[
w=\frac{u}{h}-\frac{a\,A(u,b)}{h^2}.
\]
Since $\|a\|_r=1$, $\|b\|_r\le1+\delta$, and
$\|u\|_r\le\|\xi\|_r$, Lemma~\ref{lem:AGY-area} gives
\[
\|w\|_r\le\frac{2}{(1-\delta)^2}\|\xi\|_r.
\]
Transport along $H_s$ and \eqref{eq:AGY-GM-path} then give
\[
\|D\pi_q^{\rm cs}(\xi)\|_{\pi_q^{\rm cs}(r)}
\le\frac{2}{(1-\delta)^2(1-3\delta)}\|\xi\|_r
\le\frac{2000}{567}\|\xi\|_r
<4\|\xi\|_r.
\]
This is an intrinsic horospherical tangent norm, since the image
has fixed horizontal foliation and area one.

Finally, let $q_1,q_2\in B_\AGY(q',\oof)$. Join $q$ to $q_1$
by a path of length arbitrarily close to $D_\AGY(q,q_1)$, and then
join $q_1$ to $q_2$ by a path of length arbitrarily close to
$D_\AGY(q_1,q_2)$. Their total length can be chosen less than
$4/50+\varepsilon<\log(11/10)$. The preceding construction and
derivative bound apply along this whole concatenation, even if the
second path leaves the stated ball. Integrating along the second
path, then letting its length error tend to zero, proves the required
inequality.
\end{proof}

Projection immediately compares the ambient and intrinsic metrics
for nearby points on the same horosphere.

\begin{lemma}
\label{lemma:DAGY:vs:dAGY}
For $x\in\POT$ and $y\in\Hor(x)\cap\POT$ with
$D_\AGY(x,y)\le\ooh$, we have
\[
D_\AGY(x,y)\le d_\AGY(x,y)\le4D_\AGY(x,y).
\]
\end{lemma}

\begin{proof}
The first inequality follows from the definitions. For the second,
apply Theorem~\ref{theorem:cs-projection} with $q=q'=x$.
Both points lie in $B_\AGY(x,\oof)$ and are fixed by the projection.
\end{proof}

The same comparison holds on stable horospheres, by the rotation
$q\mapsto-q$.

\section{A system of charts for $\PMF$} \label{sec:Charts}

To compare points of $\PMF$, we will represent them on a common
horosphere and use the coordinates and estimates of the preceding
sections. We need finitely many charts such that every triple lies
in an inner set and the outer sets have small intrinsic AGY diameter
after this normalization. The nested sets leave room for the perturbations used
later, and a partition of unity will combine the chart distances into
a global metric. The linear coordinates are on the corresponding
cones in $\MF$; their projectivizations give coordinates on $\PMF$.

Choose a $\Mod(S)$--compact set
\index{$K_0$}
$K_0\subset\POT$ whose quotient $\pi(K_0)\subset\POM$ is connected
and has more than $\tfrac23$ of the total Masur--Veech volume.
We will enlarge $K_0$ during this construction and then fix it.

A \emph{chart} $\alpha$ consists of a horosphere $\Hor(F_\alpha)$,
where $F_\alpha$ is a stable foliation of a pseudo-Anosov map
$\phi_\alpha$ with axis in $K_0$, and open sets
\[
U_\alpha\subset V_\alpha^{(1)}\subset V_\alpha^{(2)}
\subset V_\alpha^{(3)}\subset V_\alpha^{(4)}\subset W_\alpha
\subset\PMF(F_\alpha).
\]
The closure of each set is contained in the next, including
$\overline W_\alpha\subset\PMF(F_\alpha)$. We choose the axes in
the principal stratum. Their stable foliations have the principal
singularity type and no saddle connections, so $\Hor(F_\alpha)$
lies in $\POT$.

For $x\in W_\alpha$ and $t\in\R$, write
\[
x^\alpha(t):=x_{F_\alpha}(t)=Q_1(x,e^{-t}F_\alpha).
\]
For $E\subset W_\alpha$, put
$E_\alpha(t)=\{x^\alpha(t):x\in E\}$, abbreviating
$(U_\alpha)_\alpha(t)$ by $U_\alpha(t)$, and similarly for
$V_\alpha^{(i)}(t)$ and $W_\alpha(t)$. We require
\begin{equation}
\label{eq:dagy:x:F:alpha:y:F:alpha}
d_\AGY\big(x^\alpha(0),y^\alpha(0)\big)\le\frac1{200}
\qquad(x,y\in W_\alpha).
\end{equation}

Here $U_\alpha$ is the covering set, and $V_\alpha^{(1)}$ will
contain the closed-orbit centers constructed in
Proposition~\ref{prop:properties:gamma:n}. The sets
$V_\alpha^{(2)}\subset V_\alpha^{(3)}$ allow small perturbations
while keeping a uniform lower bound for the partition function
$\theta_\alpha$, constructed below. Its support lies in
$V_\alpha^{(4)}$, leaving room inside the full
chart $W_\alpha$ to compare nearby points.

We first choose the covering sets. There is a pseudo-Anosov map
whose axis lies in the principal stratum; this follows, for example,
from \cite[Theorem~1.2]{eskin:mirzakhani:rafi:counting}. Its repelling
foliation has infinitely many distinct images under $\Mod(S)$.
Thus, given $x_1,x_2,x_3\in\PMF$, a conjugate $\phi$ has repelling
foliation outside this triple. Choose open sets $U_\phi,Z_\phi$
with disjoint closures such that
\[
x_1,x_2,x_3,[F_+(\phi)]\in U_\phi,
\qquad [F_-(\phi)]\in Z_\phi.
\]
The sets $U_\phi^3$ cover the compact space $\PMF^3$. Take a finite
subcover, indexed by $\cA$, and denote its maps and sets by
$\phi_\alpha$ and $U_\alpha$. In particular, every triple lies
in some $U_\alpha$.

We now choose the normalizations so that the outer sets are small
in the intrinsic AGY metric. For each $\alpha$, choose $q_\alpha'$
on the axis of $\phi_\alpha$ and put $F_\alpha'=F_-(q_\alpha')$.
Choose an open neighborhood $W_\alpha'$ of $[F_+(\phi_\alpha)]$
such that
\[
\overline W_\alpha'\subset
\left\{x\in\PMF(F_\alpha'):
d_\AGY\big(x_{F_\alpha'}(0),q_\alpha'\big)<\frac1{400}\right\},
\qquad
\R_+\overline W_\alpha'\subset\R_+\cV_+(q_\alpha'),
\]
where $\cV_+$ is the neighborhood from
Proposition~\ref{prop:Local-Product}. For $E\subset\PMF$, the cone
$\R_+E$ consists of the measured foliations whose projective classes
lie in $E$; for $E\subset\MF$, it consists of positive multiples
of elements of $E$.
Since $\overline U_\alpha$ avoids the repelling fixed point,
north--south dynamics gives a power $k$ such that
$\phi_\alpha^k(\overline U_\alpha)\subset W_\alpha'$. Set
\[
q_\alpha=\phi_\alpha^{-k}q_\alpha',\qquad
F_\alpha=\phi_\alpha^{-k}F_\alpha',\qquad
W_\alpha=\phi_\alpha^{-k}W_\alpha',
\]
and take the product box at $q_\alpha$ to be the translate of the
one at $q_\alpha'$. Equivariance of $Q_1$ gives
\[
x_{F_\alpha}(0)
=\phi_\alpha^{-k}\big((\phi_\alpha^kx)_{F_\alpha'}(0)\big).
\]
The mapping class action is an AGY isometry. Hence
$\overline U_\alpha\subset W_\alpha$, the cone containment is
preserved, and every point of $Q_1(\overline W_\alpha,F_\alpha)$ is within
$1/400$ of $q_\alpha$ in the intrinsic metric. The triangle
inequality proves \eqref{eq:dagy:x:F:alpha:y:F:alpha}.

Insert the four open sets $V_\alpha^{(1)},\ldots,V_\alpha^{(4)}$
between $\overline U_\alpha$ and $W_\alpha$, with the compact
containments specified above. Choose continuous functions
$\chi_\alpha\colon\PMF\to[0,1]$ equal to one on
$\overline V_\alpha^{(3)}$ and supported in $V_\alpha^{(4)}$, and set
\[
\theta_\alpha(x)=
\frac{\chi_\alpha(x)}{\sum_{\beta\in\cA}\chi_\beta(x)}.
\]
Since the $U_\alpha$ cover $\PMF$, the denominator lies between
$1$ and $|\cA|$. Thus the $\theta_\alpha$ form a partition of unity,
are supported in $V_\alpha^{(4)}$, and are at least $1/|\cA|$
on $V_\alpha^{(3)}$.

By backward nonexpansion (Theorem~\ref{theorem:properties:AGY}),
$W_\alpha(-s)$ remains within intrinsic distance $1/400$ of
$q_\alpha(-s)$ for every $s\ge0$. The finitely many axes have
compact images in moduli space. Properness of the ambient quotient
metric therefore gives a $\Mod(S)$--compact neighborhood containing
all these sets. Enlarge $K_0$ to contain this neighborhood, adjoining
finitely many connecting paths in the quotient if necessary to keep
$\pi(K_0)$ connected. This finishes the choice of $K_0$.

The finite union of the lifted closures
$Q_1(\overline W_\alpha,F_\alpha)$ is compact in the marked
principal stratum. Its projection to $\cT$ is therefore bounded
relative to $P$, giving a constant $\sD_{\rm chart}$ for both
inequalities below. For each consecutive pair between $U_\alpha$
and $W_\alpha$, the
lifted inner closure is compact inside the next lifted open set.
Its intrinsic AGY distance from the complement is positive. Choose
$\rr_{\rm chart}>0$ strictly smaller than the minimum of these
finitely many distances. This gives the following properties.

\begin{theorem}
\label{theorem:charts}
There is a finite set of charts $\cA$, with data
\[
\alpha=\big\{(X_\alpha,q_\alpha),U_\alpha,V_\alpha^{(1)},
V_\alpha^{(2)},V_\alpha^{(3)},V_\alpha^{(4)},W_\alpha,
\theta_\alpha\big\},
\]
and constants $\sD_{\rm chart}>0$ and
\index{$\rr_{\rm chart}>0$}$\rr_{\rm chart}>0$ satisfying:
\begin{enumerate}[(a)]
\item For $\alpha,\alpha'\in\cA$ and
$x\in W_\alpha\cap W_{\alpha'}$,
\[
d_\cT\big(x^\alpha(0),x^{\alpha'}(0)\big)
\le\sD_{\rm chart},
\qquad
d_\cT\big(x^\alpha(0),P\big)\le\sD_{\rm chart}.
\]
\item The sets
\[
U_\alpha\subset V_\alpha^{(1)}\subset V_\alpha^{(2)}
\subset V_\alpha^{(3)}\subset V_\alpha^{(4)}\subset W_\alpha
\subset\PMF(F_\alpha)
\]
are open, with the closure of each contained in the next.
Every three points of $\PMF$ lie in some $U_\alpha$.
\item For any $x,y\in W_\alpha$,
\eqref{eq:dagy:x:F:alpha:y:F:alpha} holds.
\item We have
\begin{equation}
\label{eq:Walpha:in:U:plus}
\R_+\overline W_\alpha\subset\R_+\cV_+(q_\alpha),
\end{equation}
where $\cV_+$ is as in Proposition~\ref{prop:Local-Product}.
\item The continuous functions $\theta_\alpha\colon\PMF\to[0,1]$
form a partition of unity with
$\supp(\theta_\alpha)\subset V_\alpha^{(4)}$, and
\[
\theta_\alpha(x)\ge\frac1{|\cA|}
\qquad(x\in V_\alpha^{(3)}).
\]
\item For any consecutive pair $A,B$ in the chain
$U_\alpha,V_\alpha^{(1)},\ldots,V_\alpha^{(4)},W_\alpha$,
if $x\in A$, $y\in\PMF(F_\alpha)$, and
\[
d_\AGY\big(x_{F_\alpha}(0),y_{F_\alpha}(0)\big)
\le\rr_{\rm chart},
\]
then $y\in B$.
\item For every $\alpha\in\cA$ and $s>0$,
$W_\alpha(-s)\subset K_0$.
\end{enumerate}
\end{theorem}

We fix these charts and $K_0$ for the rest of the paper. A
\emph{uniform constant} is one depending only on the topology of
$S$, the chosen charts, and the fixed basepoint $P$.

Changing charts amounts to changing the horizontal foliation and
renormalizing the area. When both points already belong to both
charts, the projection theorem gives the comparison directly.
Indeed, let $x,y\in W_\alpha\cap W_{\alpha'}$ and $s,s'\in\R$,
and assume
\[
D_\AGY\big(x^\alpha(s),x^{\alpha'}(s')\big)\le\frac1{100},
\qquad
d_\AGY\big(x^\alpha(s),y^\alpha(s)\big)\le\frac1{100}.
\]
Then
\[
d_\AGY\big(x^{\alpha'}(s'),y^{\alpha'}(s')\big)
\le4\,d_\AGY\big(x^\alpha(s),y^\alpha(s)\big).
\]
To see this, apply Theorem~\ref{theorem:cs-projection} with target
$q=x^{\alpha'}(s')$ and input-ball center $x^\alpha(s)$. The target
is within $1/100$ of the center, and both inputs lie strictly inside
the center's ambient $1/50$ ball. Their projections are
$x^{\alpha'}(s')$ and $y^{\alpha'}(s')$, so the estimate follows
from the Lipschitz constant $4$.

The next lemma also allows us to put a nearby point into the second
chart. We first compare the two representatives of a common
projective foliation and align their measured vertical foliations.

\begin{lemma}
\label{lem:uniform-change-of-charts}
There are uniform constants $\kk$ and $\sD_{\rm temp}$ such that,
for $\alpha,\alpha'\in\cA$, the following hold.
\begin{enumerate}[(a)]
\item If $x\in W_\alpha\cap W_{\alpha'}$, then for every $t\ge0$,
\begin{equation}
\label{eq:bounded:for:ever}
d_\cT\big(x^\alpha(t),x^{\alpha'}(t)\big)
\le\sD_{\rm temp}.
\end{equation}
\item For $x\in W_\alpha\cap W_{\alpha'}$, there is a time
$t_0=t_0(x,\alpha,\alpha')$ with $|t_0|\le\sD_{\rm chart}$ such that
\[
F_+\big(x^{\alpha'}(t+t_0)\big)=F_+\big(x^\alpha(t)\big)
\qquad(t\in\R).
\]
\item Let $x\in W_\alpha\cap V_{\alpha'}^{(4)}$, $y\in W_\alpha$,
and let $t_0$ be as in {\rm (b)}. Suppose that, for some $t\ge0$,
\[
d_\AGY\big(x^\alpha(t),x^{\alpha'}(t+t_0)\big)\le\frac1{100},
\qquad
d_\AGY\big(x^\alpha(t),y^\alpha(t)\big)\le\frac1{100\kk}.
\]
Then $y\in W_{\alpha'}$ and
\[
d_\AGY\big(x^{\alpha'}(t),y^{\alpha'}(t)\big)<\frac1{100}.
\]
Moreover, if
$d_\AGY\big(x^\alpha(t),x^{\alpha'}(t+t_0)\big)\le1/500$,
then, for every $z,w\in\{x,y\}$,
\[
D_\AGY\big(z^\alpha(t),w^{\alpha'}(t+t_0)\big)<\frac1{100}.
\]
\end{enumerate}
\end{lemma}

In part {\rm (c)}, the first intrinsic distance is in the common
stable horosphere supplied by {\rm (b)}. The other intrinsic
distances are in unstable horospheres; the final mixed distance
is ambient.

\begin{proof}
For part {\rm (a)}, the initial points lie in a fixed compact
subset of $\cT$ and have distance at most $\sD_{\rm chart}$.
The rays have the same projective vertical foliation. Theorem~7.1
and Remark~7.2 of \cite{rafi:hyperbolicity} give uniform fellow
travelling. To pass from initial distance at most $1$ to
$\sD_{\rm chart}$, subdivide a connecting geodesic into boundedly
many segments of length at most $1$; the subdivision points remain
in a fixed thick part. This gives $\sD_{\rm temp}$.

For part {\rm (b)}, put $F=F_+(x^\alpha(0))$ and
$F'=F_+(x^{\alpha'}(0))$, and write $F=e^{t_0}F'$.
Let $X,X'$ be the underlying surfaces. Unit area gives
$\Ext_X(F)=\Ext_{X'}(F')=1$, so Kerckhoff's formula gives
\[
e^{-2\sD_{\rm chart}}
\le\Ext_X(F')=e^{-2t_0}\le e^{2\sD_{\rm chart}}.
\]
Thus $|t_0|\le\sD_{\rm chart}$. Since the flow multiplies the
vertical foliation by $e^t$, the asserted identity holds for all
real $t$.

For part {\rm (c)}, choose $\kk\ge8$ so that
\[
\frac{4e^{2\sD_{\rm chart}}}{100\kk}
<\min\left\{\rr_{\rm chart},\frac1{100}\right\}.
\]
We project the nearby point to the second horosphere, check its
chart membership at time zero, and then estimate at time $t$.
Set $u=t+t_0$ and $q=x^{\alpha'}(u)$. Apply
Theorem~\ref{theorem:cs-projection} with target $q$ and
input-ball center $x^\alpha(t)$. The hypotheses put the target
within $1/100$ of the center and both inputs strictly inside the
center's ambient $1/50$ ball. Thus $p=\pi_q^{\rm cs}(y^\alpha(t))$ exists,
and $\pi_q^{\rm cs}(x^\alpha(t))=q$ gives
\[
d_\AGY(q,p)\le4\,d_\AGY\big(x^\alpha(t),y^\alpha(t)\big)
\le\frac4{100\kk}.
\]
The projected point has $[F_+(p)]=y$ and
$F_-(p)=e^{-u}F_{\alpha'}$. In particular,
$y\in\PMF(F_{\alpha'})$. Flowing $q,p$ by time $-u$ gives
$q_0=x^{\alpha'}(0)$ and $p_0=Q_1(y,F_{\alpha'})$.
Since $u\ge-\sD_{\rm chart}$, backward nonexpansion and the
forward bound in \eqref{eq:Monotone} give
\[
d_\AGY(q_0,p_0)
\le e^{2\sD_{\rm chart}}d_\AGY(q,p)
\le\frac{4e^{2\sD_{\rm chart}}}{100\kk}<\rr_{\rm chart}.
\]
Now $x\in V_{\alpha'}^{(4)}$, so
Theorem~\ref{theorem:charts}{\rm (f)} puts $y$ in $W_{\alpha'}$.
We may therefore write $p=y^{\alpha'}(u)$. To estimate at time $t$,
apply \eqref{eq:Monotone} between $u$ and $t$:
\[
d_\AGY\big(x^{\alpha'}(t),y^{\alpha'}(t)\big)
\le e^{2|t_0|}d_\AGY(q,p)
\le\frac{4e^{2\sD_{\rm chart}}}{100\kk}<\frac1{100}.
\]
For the additional assertion, use the triangle inequality through
$x^\alpha(t)$ and $q$. For every $z,w\in\{x,y\}$,
\[
D_\AGY\big(z^\alpha(t),w^{\alpha'}(u)\big)
\le\frac1{100\kk}+\frac1{500}+\frac4{100\kk}<\frac1{100},
\]
since $\kk\ge8$.
\end{proof}

\section{Norm estimates for the Gauss--Manin map}
\label{sec:GM-new}

The idea of this section is to measure a change of foliation by the
change it makes in intersection numbers. Extremal length makes this
quantity vary at most exponentially in \Teich distance. On compact
subsets of the principal stratum, it is comparable to the Hodge norm
of the corresponding change in periods. Comparing at the two endpoints
will give the required bound on Gauss--Manin transport.

This is the same principle as Forni's estimate for the variation of the
Hodge norm along the \Teich flow
\cite[Lemma~2.1$'$]{forni:deviation}: control the logarithmic change
of a norm to bound the growth of transported vectors. Here we use a
norm of intersection functions, whose variation follows directly from
extremal length. A fixed projective foliation along the path will
identify the vectors being compared at its endpoints.

\subsection*{The estimate}

For a compact Riemann surface $Z$, the Hodge norm of
$\alpha\in H^1(Z;\C)$ is
\begin{equation}
\label{eq:GM-H}
\|\alpha\|_{H,Z}^2
=\int_Z\alpha_{\rm harm}\wedge\star_Z\overline{\alpha_{\rm harm}},
\end{equation}
where $\alpha_{\rm harm}$ is its harmonic representative. For
$q=(X,q)\in\PQT(S)$, we write $\|\cdot\|_{H,q}$ for this norm
on the odd cohomology of the compactified orienting cover
$\widetilde X$, and on its real subspace $L(q)$.

On marked compact surfaces, the marking-induced map
$\GM_{Z,W}:H^1(Z;\C)\to H^1(W;\C)$ satisfies
\[
\max\{\|\GM_{Z,W}\|_H,\|\GM_{Z,W}^{-1}\|_H\}
\le e^{d_\cT(Z,W)}.
\]
Indeed, pulling a harmonic form across an extremal quasiconformal map
increases its $L^2$ norm by at most the square root of the dilatation,
$e^{d_\cT(Z,W)}$, and taking the harmonic representative can only
decrease that norm. Apply this also to the inverse map.
For orienting covers, this bounds transport by the distance between
the covers. We need the distance between the base surfaces $X,Y$.
The restriction to a fixed projective foliation is essential: a
general path can have nontrivial monodromy even with equal endpoints.

A measured foliation is \emph{generic} if it has the singularity type
of the principal stratum and every arc between singularities has
positive intersection with it, unless the arc is homotopic rel endpoints
into a singularity. A differential with a generic vertical or horizontal
foliation lies in the principal stratum.

\begin{proposition}
\label{prop:GM-new}
For every $\Mod(S)$--compact set $K\subset\POT(S)$ there is
$B_K\ge0$ with the following property. Let $F$ be generic and let
\[
q=(X,q),\quad p=(Y,p)\in K\cap D\Hor_\pm(F).
\]
For every path $\xi$ in $D\Hor_\pm(F)$ from $q$ to $p$,
\begin{equation}
\label{eq:GM-bound}
\log\max\{\|\GM_\xi\|_H,\|\GM_\xi^{-1}\|_H\}
\le d_\cT(X,Y)+B_K.
\end{equation}
The constant is independent of $F$ and $\xi$. Moreover, $\GM_\xi$
is independent of the path in $D\Hor_\pm(F)$ with these endpoints.
\end{proposition}

Only the endpoints must lie in $K$. The coefficient of distance is
$1$; all dependence on the compact set is in the additive constant.

\subsection*{A norm from intersection numbers}

View a measured foliation $G$ through its intersection function
$J\mapsto\I(G,J)$. For $X\in\cT(S)$, measure the difference of
two such functions by
\begin{equation}
\label{eq:GM-delta}
\delta_X(G,H)
=\sup_{J\in\MF\setminus\{0\}}
\frac{|\I(G,J)-\I(H,J)|}{\sqrt{\Ext_X(J)}}.
\end{equation}
This is a weighted supremum norm of the difference of the intersection
functions. Equivalently, we take the supremum over the compact set
$\{J:\Ext_X(J)=1\}$. Intersection numbers separate measured
foliations, so $\delta_X$ is a metric on $\MF$, with
$\delta_X(cG,cH)=c\delta_X(G,H)$ for $c>0$.

The advantage of this norm is its simple dependence on $X$.
Kerckhoff's formula bounds the change in every denominator and gives
\begin{equation}
\label{eq:GM-change}
e^{-d_\cT(X,Y)}\delta_X(G,H)
\le\delta_Y(G,H)
\le e^{d_\cT(X,Y)}\delta_X(G,H).
\end{equation}
Thus, for fixed $G\ne H$, the function $X\mapsto\log\delta_X(G,H)$
is $1$--Lipschitz. Along a unit-speed \Teich geodesic $X(t)$,
this says $|\frac{d}{dt}\log\delta_{X(t)}(G,H)|\le1$ almost
everywhere. We use the finite comparison \eqref{eq:GM-change},
which requires no differentiation in the foliation variable.
We also have mapping class equivariance:
$\delta_{\gamma X}(\gamma G,\gamma H)=\delta_X(G,H)$.

It remains to compare this distance with the size of a change in
periods. Recall from Proposition~\ref{prop:Local-Product} that a
triangulation transverse to both foliations gives a product box
$\mathcal B=\Phi^{-1}(U+iV)$ with fixed train-track maps $f,h$:
\begin{equation}
\label{eq:GM-box}
F_+\bigl(\Phi^{-1}(u+iv)\bigr)=f(u),\qquad
F_-\bigl(\Phi^{-1}(u+iv)\bigr)=h(v).
\end{equation}
Here $U,V$ are open convex subsets of a fixed real odd cohomology
space $L$, and the track weights are linear in $u,v$. In this
trivialization, $\Psi_q=f^{-1}$ near $F_+(q)$, with consistent signs.
These boxes exist even when the foliations have saddle connections.

\begin{lemma}
\label{lem:GM-norm}
For every $\Mod(S)$--compact $K\subset\POT(S)$ there is $C_K\ge1$
such that, at each $q=(X,q)\in K$, one can choose a neighborhood
$\mathcal U_q$ of $F_+(q)$ on which
\begin{equation}
\label{eq:GM-local}
C_K^{-1}\|\Psi_q(G)-\Psi_q(H)\|_{H,q}
\le\delta_X(G,H)
\le C_K\|\Psi_q(G)-\Psi_q(H)\|_{H,q}
\end{equation}
for all $G,H\in\mathcal U_q$. The same constant works near any
positive multiple of $F_+(q)$, using the homogeneous extension of
$\Psi_q$.
\end{lemma}

\begin{proof}
Work on a smaller product box with compact closure in $\mathcal B$,
and choose a Euclidean norm $|\cdot|$ on $L$. In this proof,
$\prec$ and $\asymp$ allow constants depending on the chosen box,
but uniform on its smaller closure. Continuity of the Hodge norm gives
\begin{equation}
\label{eq:GM-norms}
\|z\|_{H,q}\asymp |z|.
\end{equation}
We prove the lower bound using finitely many fixed test foliations,
and the upper bound by controlling all normalized tests.

For the lower bound, consider the constant pairing
$\Omega(z,w)=\tfrac12\int_{\widetilde X}z\wedge w$ on $L$.
It is nondegenerate: the even and odd eigenspaces of the deck
involution are orthogonal for cup product, which is nondegenerate
on the full cohomology. The area identity is
\begin{equation}
\label{eq:GM-area}
\I(f(u),h(v))=\area\bigl(\Phi^{-1}(u+iv)\bigr)=\Omega(u,v).
\end{equation}
The factor $1/2$ accounts for the degree of the orienting cover.
Choose a basis $v_1,\ldots,v_N$ of $L$ in the open set $V$, and
put $H_j=h(v_j)$. These foliations detect every real-period
direction. For $G=f(u)$ and $H=f(u')$ in the smaller box,
\[
|u-u'|\prec\max_j|\Omega(u-u',v_j)|
=\max_j|\I(G,H_j)-\I(H,H_j)|\prec\delta_X(G,H).
\]
The last comparison is uniform because the finitely many
$\Ext_X(H_j)$ are bounded as $X$ ranges over the compact projection
of the smaller box. Together with \eqref{eq:GM-norms}, this gives
the lower bound. The mixed differentials with foliations $(G,H_j)$
are used only for the area identity; extremal lengths are evaluated
on $X$.

For the upper bound, use the Dehn--Thurston metric $d_{\rm DT}$
of Luo and Stong, and put $N(J)=d_{\rm DT}(J,0)$. This metric
is the sum, over coordinate pairs, of
\[
d([x,t],[x',t'])
=\min\{|x-x'|+|t-t'|,\ |x+x'|+|t+t'|\},
\]
where the brackets mean simultaneous identification of both signs.
Their intersection estimate, extended to measured foliations, is
\begin{equation}
\label{eq:GM-LS}
|\I(G,J)-\I(H,J)|\le3N(J)\,d_{\rm DT}(G,H).
\end{equation}
See \cite[p.~2, inequality~(1.1), and Section~6]{luo:stong:DT}.
For punctured surfaces, restrict their surface-with-boundary
coordinates to the closed-curve subspace, with zero boundary
coordinates. Thus no tests by proper arcs are needed.

Changes from train-track to Dehn--Thurston coordinates are piecewise
linear \cite[Section~2.6]{ph:TT}. Finitely many linear pieces suffice
on our smaller compact coordinate closure, so
$d_{\rm DT}(f(u),f(u'))\prec|u-u'|$.
The tests with $\Ext_X(J)=1$ also lie in a fixed compact subset of
$\MF$: choose $X_*$ and $R$ with $d_\cT(X,X_*)\le R$ on the
smaller box, and use $\Ext_{X_*}(J)\le e^{2R}$. Their $N(J)$ are therefore
uniformly bounded. Taking the supremum in \eqref{eq:GM-LS} gives
\[
\delta_X(G,H)\prec |u-u'|
\asymp\|\Psi_q(G)-\Psi_q(H)\|_{H,q},
\]
as required.

Choose finitely many such smaller boxes around marked lifts covering
$K/\Mod(S)$, and take the maximum of their comparison constants
and $1$. Equivariance of all the constructions gives $C_K$ on $K$.
This also explains uniformity near saddle connections: at each
limiting differential the triangulation is chosen transverse to both
foliations, so its edge signs and linear coordinate maps stay fixed
nearby. None of the bounds divides by a transverse measure tending
to zero. The compact set lies inside the principal stratum, so zeros
do not collide. Finally, homogeneity of $\delta_X$ and
Lemma~\ref{lem:Homogeneous} give the same constant near each positive
multiple of $F_+(q)$. The sizes of these neighborhoods need not be uniform.
\end{proof}

\subsection*{From foliation coordinates to transport}

The comparison just proved contains the analytic part of the argument.
We only need to identify the same small change of foliation at the
two endpoints. This is where the path's fixed projective foliation
is used.

\begin{lemma}
\label{lem:GM-leaf}
Let $\xi:[0,1]\to\PQT(S)$ be a path with
$F_+(\xi(t))=a_tF$, where $a_t>0$ and $F\ne0$ is fixed. Put
$q=\xi(0)$ and $p=\xi(1)$. With endpoint period signs chosen
consistently along $\xi$, there is an open neighborhood $\mathcal W$
of $F$ in $\MF$ such that
\begin{equation}
\label{eq:GM-leaf}
\GM_\xi\bigl(\Psi_q(G)\bigr)=\Psi_p(G)
\qquad(G\in\mathcal W).
\end{equation}
The endpoint maps are the homogeneous local coordinates near the ray of $F$.
\end{lemma}

\begin{proof}
Replace $\xi(t)$ by $a_t^{-2}\xi(t)$ to make its vertical
foliation equal to $F$. Scalar multiplication leaves the orienting
cover canonically unchanged, and Lemma~\ref{lem:Homogeneous}
identifies the homogeneous coordinate maps. It suffices to prove
the assertion on this normalized path.

Subdivide the path into finitely many subpaths lying in product boxes
from Proposition~\ref{prop:Local-Product}. For $G$ near $F$,
\eqref{eq:GM-box} makes its coordinate $f^{-1}(G)$ constant in
the flat trivialization along each subpath.
At a junction $r$, the section $G\mapsto Q(G,F_-(r))$ lies in both
boxes for $G$ near $F$; hence their coordinate maps agree under the
Gauss--Manin change of trivialization. Choose signs successively
and intersect these finitely many neighborhoods of $F$.
Composing the local identities proves \eqref{eq:GM-leaf}.
\end{proof}

The neighborhood may depend on the whole path. Its size will not
enter the estimate, since a neighborhood of zero suffices to test a
linear map. Without a fixed projective foliation, these local
identities need not compare the same foliation at both endpoints.

\begin{proof}[Proof of Proposition~\ref{prop:GM-new}]
First suppose $\xi\subset D\Hor_+(F)$. Let $C$ be the comparison
constant from Lemma~\ref{lem:GM-norm} for $K\cup(-K)$; using this
compact set also covers the horizontal case below. Shrink
$\mathcal W$ from Lemma~\ref{lem:GM-leaf} so that both endpoint
comparisons hold there. The homogeneous assertion allows this even
when $F_+(q)$ and $F_+(p)$ are different multiples of $F$.

For $G\in\mathcal W$, put $z=\Psi_q(G)-\Psi_q(F)$. Then
\begin{align*}
\|\GM_\xi z\|_{H,p}
&=\|\Psi_p(G)-\Psi_p(F)\|_{H,p}\\
&\le C\,\delta_Y(G,F)
\le C e^{d_\cT(X,Y)}\delta_X(G,F)\\
&\le C^2e^{d_\cT(X,Y)}\|z\|_{H,q}.
\end{align*}
These $z$ fill a neighborhood of zero in $L(q)$, so linearity gives
the bound for every vector. Reversing $\xi$ gives the inverse bound.
Both extend to complex cohomology because transport is real and
$\|z+iw\|_H^2=\|z\|_H^2+\|w\|_H^2$ for real $z,w$.

For $D\Hor_-(F)$, apply this argument to $-\xi$. The rotation
$q\mapsto-q$ interchanges the two foliations and identifies the
conformal orienting covers by Hodge isometries that intertwine
transport. Thus \eqref{eq:GM-bound} holds with $B_K=2\log C$.
Finally, Hubbard--Masur identifies $D\Hor_+(F)$ with
$\cT(S)\times\R_{>0}$ by $(X,c)\mapsto c^2Q(F,X)$.
Genericity puts this contractible leaf entirely in the principal
stratum, so homotopy invariance (Lemma~\ref{lem:GM}) gives path
independence. The horizontal case is the same.
\end{proof}

\subsection*{The estimate in the AGY norm}

The Hodge and AGY norms on the real spaces $L(q)$ are continuous
and $\Mod(S)$--invariant, so they are uniformly comparable on a $\Mod(S)$--compact subset
$K\subset\POT$. Choose $A_K\ge1$ with
\[
A_K^{-1}\|z\|_{H,q}\le\|z\|_{\AGY,q}
\le A_K\|z\|_{H,q}\qquad(q\in K).
\]
Under the hypotheses of Proposition~\ref{prop:GM-new},
\begin{equation}
\label{eq:GM-AGY}
\max\{\|\GM_\xi\|_{\AGY},\|\GM_\xi^{-1}\|_{\AGY}\}
\le A_K^2e^{B_K}e^{d_\cT(X,Y)}.
\end{equation}
Changing norms affects only the prefactor. In particular,
$d_\cT(X,Y)\le ct+L$ gives the bound
$A_K^2e^{B_K+L}e^{ct}$, with the coefficient $c$ unchanged.

\section{Good geodesics}
\label{sec:good-geodesics}

We construct pseudo-Anosov maps whose attracting foliations will be
the centers of the spikes. Their periods lie near prescribed times
$R_n=n\tau'$. The appendix supplies many times with good backward
statistics. We add a bounded return segment and apply the closing
lemma, then project the resulting periodic orbit to the original
chart. We keep track of the statistics and the approximation error
because both enter the spike estimates.

Recall that $K\subset\POT$ is \emph{$\Mod$--compact} if it is
$\Mod(S)$--invariant and its moduli quotient is compact; in particular,
$K$ is closed. We use $\axis(\gamma)$ for both the \Teich axis of a
pseudo-Anosov map and its oriented lift to the space of unit-area
differentials. On this lift, $g_{\ell(\gamma)}z=\gamma z$.

\begin{proposition}[Closing Lemma]
\label{prop:closing:lemma}
There are constants $\ep_{\rm clos}>0$ and $C_{\rm clos}>1$ depending
on $K_0$ with the following property. For every $\Mod$--compact set
$K_1\subset\POT$, there is $R_1>0$ depending on $K_0$,
$\ep_{\rm clos}$, and $K_1$ such that, if $0<\ep<\ep_{\rm clos}$,
$q\in K_0$, and $t>0$ satisfy
\[
D_\AGY\big(\pi(q),\pi(q(t))\big)<\ep,
\qquad
\big|\{s\in[0,t]:q(s)\in K_1\}\big|\ge R_1,
\]
then there is a pseudo-Anosov map $\gamma$ with axis in the principal
stratum such that
\[
|\ell(\gamma)-t|\le C_{\rm clos}\ep,
\qquad
D_\AGY\big(q(s),\axis(\gamma)\big)\le C_{\rm clos}\ep
\quad(0\le s\le t).
\]
\end{proposition}

\begin{proof}
We adapt the closing argument of
\cite[Lemma~8.1 and Remark~8.2]{eskin:mirzakhani:rafi:counting}
to the AGY metric. Fix local product neighborhoods and their comparison
constants over $K_0$; these determine $\ep_{\rm clos}$ and
$C_{\rm clos}$. Stable forward and unstable backward nonexpansion
(Theorem~\ref{theorem:properties:AGY}) keep nearby comparison paths
in a fixed-radius compact thickening of $K_1$ whenever the original
orbit visits $K_1$. By
Proposition~\ref{prop:AG43}, sufficiently large total time there
gives the contraction required by the fixed-point argument. This
determines $R_1$. The construction bounds the displacement and period
correction by $C_{\rm clos}\ep$, and nonexpansion gives the tracking
bound. Thus $K_1$ affects the occupation threshold, not the closing
accuracy or the tracking constant.
\end{proof}

Good statistics requires control on every interval ending at the
point under consideration, however far into the past it begins.

\begin{definition}
\label{def:good-statistics}
For a $\Mod$--compact set $K$ and $0<\eta<1/2$, a differential
$q\in\POT$ has $(K,\eta)$--good statistics if $q\in K$ and,
for every $s<0$,
\begin{equation}
\label{eq:def:good-statistics}
\big|\{u\in[s,0]:q(u)\in K\}\big|\ge(1-\eta)|s|.
\end{equation}
\end{definition}
We use the same definition for open neighborhoods of $\Mod$--compact
sets. For rays starting in our charts, the entire negative ray lies
in $K_0$, so only the positive-time part needs checking.

\begin{lemma}
\label{lemma:good:statistics:beyond:zero}
Suppose $0<\eta<1/2$, $K\supset K_0$ is $\Mod$--compact,
$q\in U_\alpha(0)$, and $R>0$. If, for every $0<s<R$,
\begin{equation}
\label{eq:alt:good-statistics}
\big|\{u\in[s,R]:q(u)\in K\}\big|\ge(1-\eta)(R-s),
\end{equation}
then $q(R)$ has $(K,\eta)$--good statistics.
\end{lemma}

\begin{proof}
Letting $s\downarrow0$ gives \eqref{eq:alt:good-statistics} at
$s=0$. Letting $s\uparrow R$ gives points of $K$ approaching
$q(R)$, so closedness gives $q(R)\in K$. For $u<0$,
Theorem~\ref{theorem:charts}{\rm (g)} gives $q(u)\in K_0$.
Adding this negative-time interval proves
\eqref{eq:alt:good-statistics} for every $s<0$, hence the assertion.
\end{proof}

The next proposition records both the geometry of the periodic orbits
and how often they approximate rays from each chart. All uniform
constants depend only on the fixed geometric data of
Theorem~\ref{theorem:charts}.

\begin{proposition}
\label{prop:properties:gamma:n}
There are uniform constants $\sD_0,\sD_1$ such that, for every
$0<\eta_1<1/2$, there are a $\Mod$--compact set $K_1\supset K_0$,
constants $\rho_1,\tau>0$, $\tau'_0>\rho_1$, and $0<\kappa_0<1$
depending on $\eta_1$ with the following properties.
Fix $\tau'\in\tau\natls$ with $\tau'>\tau'_0$, and set
$R_n=n\tau'$. For each chart $\alpha$ and all sufficiently large
$n$ depending on $\eta_1$, there is a set
$\Gamma_\alpha(R_n)\subset\Mod(S)$ satisfying {\rm (a)--(f)} below.
Parts {\rm (a)--(e)} hold for each $\gamma\in\Gamma_\alpha(R_n)$;
part {\rm (f)} concerns the whole family.
\begin{itemize}
\item[{\rm (a)}] The map $\gamma$ is pseudo-Anosov, and its axis
passes within $d_\cT$-distance $\sD_0$ of $P$.
\item[{\rm (b)}] Its translation length satisfies
\[
R_n-1\le\ell(\gamma)\le R_n+\rho_1+1.
\]
\item[{\rm (c)}] Write $a_\gamma=[F_+(\gamma)]$. Then
\[
D_\AGY\big(a_\gamma^\alpha(t),\axis(\gamma)\big)\le\frac1{100}
\qquad(t\ge0).
\]
\item[{\rm (c')}] We have
\[
d_\cT\big(a_\gamma^\alpha(\ell(\gamma)),\gamma(P)\big)\le\sD_1.
\]
\item[{\rm (d)}] We have $a_\gamma\in V^{(1)}_\alpha$, and
$a_\gamma^\alpha(\ell(\gamma))$ has $(K_1,\eta_1)$--good statistics.
\item[{\rm (e)}] Set
\[
G_\gamma=\{u\in[0,\ell(\gamma)]:a_\gamma^\alpha(u)
\text{ has }(K_1,\eta_1)\text{--good statistics}\}.
\]
For every $0\le s<\ell(\gamma)$,
\[
\big|[s,\ell(\gamma)]\setminus G_\gamma\big|
\le\eta_1\big(\ell(\gamma)-s\big).
\]
\item[{\rm (f)}] Let $\kk$ be as in
Lemma~\ref{lem:uniform-change-of-charts}, and define
\[
\begin{aligned}
E_\alpha(R_n)=\bigl\{q\in U_\alpha(0):\ &
\text{there exists }\gamma\in\Gamma_\alpha(R_n)\text{ with}\\
&d_\AGY\big(q(\ell(\gamma)),a_\gamma^\alpha(\ell(\gamma))\big)
<\tfrac1{100\kk}\bigr\}.
\end{aligned}
\]
At the excluded initial indices, take both $\Gamma_\alpha(R_j)$
and $E_\alpha(R_j)$ to be empty. Put
\[
k_n(q)=\#\{1\le j\le n:q\in E_\alpha(R_j)\}.
\]
For all sufficiently large $n$ depending on $\eta_1$ and $\tau'$,
\[
\mu\bigl(\{q\in U_\alpha(0):k_n(q)\le(1-\eta_1)n\}\bigr)
\le\kappa_0^n\mu(U_\alpha(0)).
\]
Here $\mu$ is the conditional Masur--Veech volume on
$\Hor(F_\alpha)$, normalized as in Section~\ref{sec:volume-balls}.
\end{itemize}
\end{proposition}

\begin{proof}
\emph{Closing accuracy and return data.}
First choose the accuracy needed for the projections and the radius
in part {\rm (f)}. With the constants of
Proposition~\ref{prop:closing:lemma}, Theorem~\ref{theorem:charts},
and Lemma~\ref{lem:uniform-change-of-charts}, fix $\epsilon>0$ such that
\[
\epsilon<\frac{\epsilon_{\rm clos}}2,
\qquad
C_{\rm clos}\epsilon<
\min\left\{\frac1{2000},\frac1{800\kk},\frac{\rr_{\rm chart}}4\right\}.
\]
This choice uses only the fixed geometric data. Choose finite sets
$\widehat{\cB}(\alpha)\subset U_\alpha(0)$ whose ambient
$\epsilon/4$--balls cover $U_\alpha(0)$, and put
$\widehat{\cB}=\bigcup_\alpha\widehat{\cB}(\alpha)$.

Now fix $\eta=\eta_1/2$. Proposition~\ref{prop:appendix-final}
gives $K\supset K_0$, $\tau$, $\kappa$, and $\theta_0$; set
$\theta=\theta_0/2$. Write $\lambda=\muMV$ for the normalized
Masur--Veech probability measure on $\POM$.

\begin{claim}
\label{claim:ergodic:small:set}
For every $\theta>0$, there are $\rho_\theta>0$ and a compact set
$\Omega_\theta\subset\POM$ with $\lambda(\Omega_\theta)\ge1-\theta$
such that, for every $p\in\Omega_\theta$ and $z\in\widehat{\cB}$,
the segment $\{p(t):0\le t\le\rho_\theta\}$ meets
$\pi(B_\AGY(z,\epsilon/4))$.
\end{claim}

\begin{proof}
Each target is open and has positive measure. Ergodicity implies that
almost every forward orbit meets all these finitely many targets.
Thus the measure of the set that does so by time $R$ tends to one
as $R\to\infty$. Choose $\rho_\theta$ so that this measure exceeds
$1-\theta$, then take a compact subset of measure at least $1-\theta$.
\end{proof}

Put $\rho_1=\rho_\theta$ and use $\Omega_\theta$ also for its full
lift to $\POT$. Write $N_r$ for ambient AGY neighborhoods. Define
\[
\begin{aligned}
\mathcal K_{\rm aux}
&=N_2\bigl(\{p(t):p\in K,\ 0\le t\le\rho_1\}\bigr),\\
K_1&=\left\{q\in\POT:\sys(q)\ge
\inf_{p\in\mathcal K_{\rm aux}}\sys(p)\right\},
\end{aligned}
\]
where $\sys$ is the shortest saddle-connection length. The swept set
is compact modulo $\Mod$, and its closed AGY neighborhood is also
compact modulo $\Mod$. Hence the infimum is positive and $K_1$
is $\Mod$--compact. By Theorem~\ref{theorem:charts}{\rm (g)} and
closedness, $U_\alpha(0)\subset K_0\subset K$. Thus $K_1$ contains
the ambient $2$--neighborhoods of $K_0$ and of every $U_\alpha(0)$.

The AGY norm gives
$D_\AGY(g_tp,g_tp')\le e^{2|t|}D_\AGY(p,p')$: the numerator
in each norm ratio increases by at most $e^{|t|}$, and its
denominator decreases by at most that factor; integrate along
ambient paths. Choose
$0<\delta_\theta<\min\{{1/400},\epsilon e^{-2\rho_1}/8\}$.
Then
\begin{equation}
\label{eq:delta-epsilon}
D_\AGY(p,p')\le2\delta_\theta,\quad |t|\le\rho_1
\quad\Longrightarrow\quad
D_\AGY(p(t),p'(t))\le\epsilon/4.
\end{equation}
Only now choose $\tau'_0$ from Proposition~\ref{prop:appendix-final},
enlarged so that $\tau'_0>\rho_1$. Its measure hypothesis holds
because $\lambda(\Omega_\theta/\Mod)\ge1-\theta$.
Thus $K_1$ precedes $\delta_\theta$, and the later data affect
only time thresholds. In particular, $\epsilon$ is unchanged.

\emph{Closing the orbit segments.}
A point $q\in U_\alpha(0)$ is $(K,\eta,\delta,\Omega)$--generic
at $R$ (Definition~\ref{def:generic}) if
\begin{equation}
\label{eq:i}
q(R)\in N_{2\delta}(\Omega)
\quad\text{and has }(K,\eta)\text{--good statistics},
\end{equation}
and, for every $0\le s<R$,
\begin{equation}
\label{eq:ii}
\big|\{u\in(s,R):q(u)\text{ has }(K,\eta)\text{--good statistics}\}\big|
\ge(1-\eta)(R-s).
\end{equation}
Let $\cG_\alpha(R)$ be the set of points in $U_\alpha(0)$ that
are $(N_{\delta_\theta}(K),\eta,\delta_\theta,\Omega_\theta)$--generic
at $R$. Take the closing threshold $R_1$ for $K_1$, and choose an
integer $j_0\ge1$ with $(1-\eta)j_0\tau'_0>R_1$.
This choice works for every $\tau'\in\tau\natls$ above $\tau'_0$.

Fix $j\ge j_0$ and $q\in\cG_\alpha(R_j)$. Choose
$z_q\in\widehat{\cB}$ with $D_\AGY(q,z_q)<\epsilon/4$ and,
by \eqref{eq:i}, $p'\in\Omega_\theta$ with
$D_\AGY(q(R_j),p')<2\delta_\theta$.
Claim~\ref{claim:ergodic:small:set} gives a return time
$0\le t_{\rm ret}\le\rho_1$ with
$\pi(p'(t_{\rm ret}))\in\pi(B_\AGY(z_q,\epsilon/4))$.
Together with \eqref{eq:delta-epsilon}, the three errors give
\[
D_\AGY\big(\pi(q),\pi(q(R_j+t_{\rm ret}))\big)
\le\epsilon/4+\epsilon/4+\epsilon/4<\epsilon.
\]
Good statistics at $R_j$ supply at least
$(1-\eta)R_j>R_1$ time in $K_1$. The closing lemma therefore
gives a pseudo-Anosov map $\gamma_q$. Set
\[
\Gamma_\alpha(R_j)=\{\gamma_q:q\in\cG_\alpha(R_j)\},
\]
and take $\Gamma_\alpha(R_j)=E_\alpha(R_j)=\varnothing$ for $j<j_0$.

\emph{Geometry of the closed orbit.}
Write $\gamma=\gamma_q$, $\ell=\ell(\gamma)$, and
$T=R_j+t_{\rm ret}$. The initial points lie in fixed marked
compact chart closures. An axis point is within
$C_{\rm clos}\epsilon$ of $q$, so its base distance from $P$
is bounded by a uniform $\sD_0$, proving {\rm (a)}.
The period estimate gives
\[
|\ell-T|\le C_{\rm clos}\epsilon<1,
\qquad R_j-1\le\ell\le R_j+\rho_1+1,
\]
which proves {\rm (b)}. The tracking estimate is
\begin{equation}
\label{eq:axis-close}
D_\AGY\big(q(t),\axis(\gamma)\big)\le C_{\rm clos}\epsilon
\qquad(0\le t\le T).
\end{equation}

Put $a_\gamma=[F_+(\gamma)]$ and choose $z\in\axis(\gamma)$
within $C_{\rm clos}\epsilon$ of $q$. Apply
Theorem~\ref{theorem:cs-projection} with target and ball center $q$.
The projection of $z$ is $Q_1(a_\gamma,F_\alpha)$; the theorem
establishes that it is defined, principal, and at intrinsic distance
at most $4C_{\rm clos}\epsilon<\rr_{\rm chart}$ from $q$.
Theorem~\ref{theorem:charts}{\rm (f)} now gives
$a_\gamma\in V^{(1)}_\alpha$, so this projection is
$a_\gamma^\alpha(0)$. At each time $0\le t\le T$, the same
argument with target and center $q(t)$ gives
\begin{equation}
\label{eq:axis-projection-close}
d_\AGY\big(q(t),a_\gamma^\alpha(t)\big)
\le4C_{\rm clos}\epsilon.
\end{equation}
All inputs lie inside the ambient $1/50$--balls required by
Theorem~\ref{theorem:cs-projection}.

To follow the axis for all future time, observe that
$D_\AGY(a_\gamma^\alpha(0),z)\le5C_{\rm clos}\epsilon<1/200$.
Apply the time-reversed projection theorem to project $z$ onto
$\Hor_+(a_\gamma^\alpha(0))$. Its image $z'$ has the same two
projective foliations as $z$, so it lies on the same axis, and
\[
d_\AGY\big(a_\gamma^\alpha(0),z'\big)
\le20C_{\rm clos}\epsilon.
\]
Forward stable nonexpansion, the time-reversed form of
\eqref{eq:Monotone}, gives
\begin{equation}
\label{eq:part:c}
D_\AGY\big(a_\gamma^\alpha(t),\axis(\gamma)\big)
\le20C_{\rm clos}\epsilon<\frac1{100}
\qquad(t\ge0).
\end{equation}
This proves {\rm (c)}. Since $g_\ell z'=\gamma z'$, it also gives
\[
\begin{aligned}
D_\AGY\big(\gamma^{-1}a_\gamma^\alpha(\ell),a_\gamma^\alpha(0)\big)
&\le D_\AGY\big(a_\gamma^\alpha(\ell),g_\ell z'\big)
    +D_\AGY\big(z',a_\gamma^\alpha(0)\big)\\
&\le40C_{\rm clos}\epsilon<1/50.
\end{aligned}
\]
The closed ambient $1/50$--neighborhood of
$\bigcup_\alpha\overline{V^{(1)}_\alpha(0)}$ is compact in marked
space. Its base distances from $P$ are bounded by a uniform $\sD_1$,
which proves {\rm (c')} by $\Mod$--invariance.

\emph{Transferring the statistics.}
Write $a(t)=a_\gamma^\alpha(t)$, and set
\[
L=\overline{N_1\bigl(\{p(u):p\in K,\ 0\le u\le\rho_1\}\bigr)}.
\]
The entire negative ray $a(t)$ lies in $K_0\subset L$.
On $[0,R_j]$, \eqref{eq:axis-projection-close} puts $a(t)$ in
$L$ whenever $q(t)\in N_{\delta_\theta}(K)$.
Finally, since $q(R_j)\in N_{\delta_\theta}(K)$,
\eqref{eq:delta-epsilon} and \eqref{eq:axis-projection-close}
put the whole appended segment $a([R_j,T])$ in $L$.
These three observations transfer good statistics from $q(u)$ to
$a(u)$ at every good time $u\in[0,R_j]$, and give good statistics
also for $u\le0$ and $R_j\le u\le T$. Indeed, every backward
interval gains all its time before zero and after $R_j$, while the
remaining part satisfies the original occupation bound.
In particular, $a(T)$ has $(L,\eta)$--good statistics, and
\eqref{eq:ii} gives, for every $v<T$,
\[
\big|\{u\in[v,T]:a(u)\text{ fails }(L,\eta)\text{--good statistics}\}\big|
\le\eta(T-v).
\]

The period differs slightly from $T$. Put $h=\ell-T$, so
$|h|\le C_{\rm clos}\epsilon<1/300$.
By \eqref{eq:gt-Lip}, $D_\AGY(g_hr,r)\le|h|$; since
$1+|h|<2$, our definition gives $g_hL\subset K_1$.
Thus translating the entire past by $h$ sends every
$(L,\eta)$--good time $u$ to a $(K_1,\eta)$--good time $u+h$.
It follows that $a(\ell)$ has these statistics. Applying the
preceding frequency estimate with $v=s-h$, also when this is
negative, gives
\[
\big|[s,\ell]\setminus G_\gamma\big|
\le\eta(\ell-s)\qquad(0\le s<\ell).
\]
Since $\eta=\eta_1/2<\eta_1$ and $a_\gamma\in V^{(1)}_\alpha$,
this proves {\rm (d)} and {\rm (e)}.

\emph{Approximation and frequency.}
For part {\rm (f)}, use \eqref{eq:axis-projection-close} at time
$T$, then \eqref{eq:Monotone} over the period correction:
\[
\begin{aligned}
d_\AGY\big(q(\ell),a_\gamma^\alpha(\ell)\big)
&\le e^{2|\ell-T|}d_\AGY\big(q(T),a_\gamma^\alpha(T)\big)\\
&\le4e^{2C_{\rm clos}\epsilon}C_{\rm clos}\epsilon
 <8C_{\rm clos}\epsilon<\frac1{100\kk}.
\end{aligned}
\]
Here the initial bounds $C_{\rm clos}\epsilon<1/2000$ and
$C_{\rm clos}\epsilon<1/(800\kk)$ give the last two inequalities.
Thus $\cG_\alpha(R_j)\subset E_\alpha(R_j)$ for $j\ge j_0$.
Put $g_n(q)=\#\{1\le j\le n:q\in\cG_\alpha(R_j)\}$.
The omitted initial indices give $g_n(q)\le k_n(q)+j_0-1$, so,
for $n\ge2(j_0-1)/\eta_1$,
\[
\{q:k_n(q)\le(1-\eta_1)n\}
\subset\{q:g_n(q)\le(1-\eta_1/2)n\}.
\]
Proposition~\ref{prop:appendix-final}, applied with $\eta=\eta_1/2$,
bounds the measure of the latter set by
$\kappa(\eta_1/2)^n\mu(U_\alpha(0))$ for all sufficiently large
$n$. Taking $\kappa_0=\kappa(\eta_1/2)$ proves {\rm (f)}.
Only the threshold for $n$, not this rate, depends on $\tau'$.
\end{proof}

Define the set of good geodesics by
\begin{equation}
\label{eq:Gamma}
\Gamma(R)=\bigcup_\alpha\Gamma_\alpha(R).
\end{equation}

\section{Convergence of rays with a common foliation}
\label{sec:rays-new}

We prove a convergence estimate for two rays $q,q'$ with the same
projective vertical foliation. If $q(T)$ has $(K,\eta)$--good statistics,
then $q(T)$ and $q'(T+a)$ are close in the intrinsic AGY metric once $T$ exceeds the
\Teich distance between their initial surfaces by an additive constant.
The shift $a$ makes their measured vertical foliations agree.
The additive constant depends on $K$, $\eta$, and the desired accuracy;
the coefficient of the initial distance is one. We allow $q'$ to lie in any stratum, and prove
that it must be principal.

We first use Hodge-norm contraction to prove convergence from a bounded
initial distance once $q$ has spent enough time in $K$. The thin-triangle
theorem of Dowdall--Duchin--Masur (DDM) then reduces the general case
to this one.

All differentials in this section have area one, and $X_p$ denotes the
base surface of $p$. Recall that a $\Mod$--compact subset of $\POT$
is $\Mod(S)$--invariant and has compact image in moduli space
(Definition~\ref{def:Mod-compact}).

\begin{theorem}
\label{thm:rays-linear}
Fix a $\Mod$--compact set $K\subset\POT$, $0<\eta<1/2$, and $\ep>0$.
There exists $B=B(K,\eta,\ep)>0$ with the following property.
Let $q\in\POT$ and $q'\in\OT$ have the same projective vertical
foliation $[F]$, and put
\[
D=d_\cT(X_q,X_{q'}),\qquad
a=\frac12\log
\frac{\Ext_{X_{q'}}(F)}{\Ext_{X_q}(F)}.
\]
Suppose $T\ge D+B$, $q(T)\in K$, and
\begin{equation}
\label{eq:rays-terminal-statistics}
\big|\{t\in[u,T]:q(t)\in K\}\big|
\ge(1-\eta)(T-u)\qquad(0\le u<T).
\end{equation}
Then $q'\in\POT$, $|a|\le D$, $T+a>0$, and
\begin{equation}
\label{eq:rays-linear-conclusion}
d_\AGY\big(q(T),q'(T+a)\big)<\ep.
\end{equation}
If $F_+(q)=F_+(q')$, then $a=0$.
\end{theorem}

If $q(T)$ has $(K,\eta)$--good statistics
(Definition~\ref{def:good-statistics}), then $q(T)\in K$ and
\eqref{eq:rays-terminal-statistics} holds.

The time adjustment is necessary because equality of projective
foliations does not imply equality of their transverse measures.
For a representative $F$ of the common projective foliation,
\[
F_+(q)=\frac{F}{\sqrt{\Ext_{X_q}(F)}},\qquad
F_+(q')=\frac{F}{\sqrt{\Ext_{X_{q'}}(F)}}.
\]
Since $F_+(g_t p)=e^t F_+(p)$, the value of $a$ in the theorem is
exactly the adjustment for which $F_+(q'(a))=F_+(q)$.
Kerckhoff's extremal-length inequality gives $|a|\le D$.

\begin{samepage}
The following proposition treats a bounded initial distance. It uses
the total time spent by $q$ in $K$ during $[0,t]$.

\begin{proposition}
\label{cor:rays-converge-new}
For every $\Mod$--compact $K\subset\POT$, $D_0>0$, and $\ep>0$, there
is $R=R(K,D_0,\ep)>0$ such that the following holds.
Suppose $q\in K$, $q'\in\OT$,
\[
F_+(q)=F_+(q'),\qquad d_\cT(X_q,X_{q'})\le D_0.
\]
If
\[
\big|\{u\in[0,t]:q(u)\in K\}\big|\ge R,
\]
then $q'\in\POT$ and
\[
d_\AGY\big(q(t),q'(t)\big)<\ep.
\]
The threshold is uniform over the common foliation, including
foliations with saddle connections.
\end{proposition}
\end{samepage}

\subsection*{Convergence from a bounded initial distance}

We first work on a closed surface of genus at least two. On the
principal part of a stable leaf $\Hor_+(q)$, let $d_{\rm H}$ be the
path distance induced by the ordinary Hodge norm from
\eqref{eq:GM-H}. The stable tangent space is
\[
E^s(q)=\left\{\alpha\in H^1_\odd(\widetilde X;\R):
\int_{\widetilde X}\Re(\omega_q)\wedge\alpha=0\right\}.
\]
Recall that a generic foliation has principal singularity type and
no saddle connections. Every differential realizing it is principal,
so we can use Hodge length on its entire stable leaf. We need an
initial length bound that remains uniform as the foliation approaches
one with saddle connections.

\begin{lemma}
\label{lem:rays-Hodge-bound}
Suppose $S$ is closed of genus at least two. For every
$\Mod$--compact $K\subset\POT$ and $D_0>0$ there is
$M=M(K,D_0)$ such that
\[
d_{\rm H}(p,p')\le M
\]
whenever $p\in K$, $F_+(p)=F_+(p')$ is generic, and
$d_\cT(X_p,X_{p'})\le D_0$.
\end{lemma}

\begin{proof}
Put $F=F_+(p)=F_+(p')$, so that
$\Ext_{X_p}(F)=\Ext_{X_{p'}}(F)=1$.
For the \Teich segment $Z_r$, $0\le r\le1$, joining the two bases, set
\begin{equation}
\label{eq:rays-HM-path}
h_r=g_{\frac12\log\Ext_{Z_r}(F)}Q_1([F],Z_r).
\end{equation}
This is a continuous unit-area path from $p$ to $p'$ with measured
vertical foliation $F$. Extremal-length comparison bounds its time
adjustment by $D_0$, so its bases lie within $2D_0$ of $X_p$.
After applying a common mapping class to $p,p'$, we may assume that
$p$ lies in a fixed compact set of marked representatives for $K$.
By properness of \Teich space, the pairs $(p,X_{p'})$ in the lemma
then range over a compact parameter set after including nongeneric
limits. Continuity of the Hubbard--Masur map and extremal length
shows that $h_r$ depends continuously on this pair and on $r$.
Thus all these paths, and their base surfaces, lie in fixed compact sets.

For points $b,c$ on the same leaf with bases in this compact set,
the Euclidean and modified Hodge distances of
\cite[Section~3.5]{athreya:bufetov:eskin:mirzakhani} satisfy
\[
d_{\rm H}^{\,\rm mod}(b,c)
\le C\,d_{\rm E}(b,c)|\log d_{\rm E}(b,c)|^{1/2}
\qquad(0<d_{\rm E}(b,c)<1/2)
\]
by \cite[Theorem~3.10]{athreya:bufetov:eskin:mirzakhani}; the stable
version follows by $q\mapsto-q$.
Uniform continuity of the path family gives one subdivision of $[0,1]$ on which
successive endpoints of every path have $d_{\rm E}<1/2$.
Summing over this subdivision bounds $d_{\rm H}^{\,\rm mod}(p,p')$ uniformly.
The modified norm dominates the ordinary Hodge norm
\cite[Section~3.3.2]{athreya:bufetov:eskin:mirzakhani}, and a generic
leaf is entirely principal. Thus the same bound applies to
$d_{\rm H}(p,p')$, uniformly even near a nongeneric foliation.
\end{proof}

\begin{proof}[Proof of Proposition~\ref{cor:rays-converge-new} for closed surfaces]
First suppose the common foliation is generic. Choose a piecewise
$C^1$ stable path $\gamma$ from $q$ to $q'$ of Hodge length at most
$M(K,D_0)+1$, and put
$\mathcal L(u)=\length_{\rm H}(g_u\gamma)$.
Hodge norms of stable vectors are non-increasing, and contract
at a uniform positive rate on each $\Mod$--compact subset of $\POT$
\cite[Theorem~1.3]{frankel:meromorphic}.
Here the flow multiplies stable period vectors by $e^{-u}$.
At any principal differential $p$, Frankel's estimate applies on
the Hodge-orthogonal complement of
$\operatorname{span}\{\Re\omega_p,\Im\omega_p\}$.
The remaining stable direction is $\R\Re\omega_p$, which contracts
by $e^{-2u}$. These two components remain Hodge orthogonal under the flow.

Fix a $\Mod$--compact neighborhood $N\subset\POT$ of $K$.
Nested period charts around compact marked representatives of $K$,
and uniform comparison of Hodge and Euclidean norms, give
$\rho>0$ such that every stable path starting in $K$ and leaving
$N$ has Hodge length greater than $\rho$.
Let $\lambda>0$ be a strict contraction constant on $N$, and choose
$C_N$ so that the AGY norm is at most $C_N$ times the Hodge norm there.
At each visit $q(u)\in K$, the initial portion of $g_u\gamma$ of
length $\min\{\mathcal L(u),\rho\}$ stays in $N$. Contraction on
this portion and nonexpansion on the rest give
\begin{equation}
\label{eq:rays-path-contraction}
\mathcal L'(u)\le-\lambda\,{\boldsymbol 1}_{K}(q(u))
\min\{\mathcal L(u),\rho\}
\end{equation}
for almost every $u$. On each interval $[0,U]$, the transported paths
stay in a compact subset of the principal stratum, so $\mathcal L$
is locally absolutely continuous.

Choose $0<h<\min\{\rho/2,\ep/C_N\}$. Inequality
\eqref{eq:rays-path-contraction} shows that time spent in $K$
decreases $\mathcal L$ at rate at least $\lambda\rho$ while
$\mathcal L\ge\rho$, and gives exponential decay once
$\mathcal L\le\rho$. Thus after
$R_0=R_0(K,D_0,\ep)$ units of time in $K$, we have $\mathcal L\le h$.
Use $R_{\rm gen}=R_0+1$ as the generic threshold. The additional
unit of time in $K$ lets us choose $s\le t$ with $q(s)\in K$ and
$\mathcal L(s)\le h$. Since $h<\rho$, the entire path $g_s\gamma$ lies in $N$,
and hence
\[
d_\AGY(q(s),q'(s))\le C_N\mathcal L(s)<\ep.
\]
Stable AGY nonexpansion gives the same bound at $t$.

We now pass to an arbitrary common foliation. Fix a $\Mod$--compact
neighborhood $K^+\subset\POT$ of $K$.
Choose $0<\delta<1/100$ with $4\delta<\ep$ so that every ambient
$2\delta$--ball centered at a point of $K^+$ lies in a period chart
whose closure is compact and contained in $\POT$.
This choice is uniform by local norm comparison in finitely many
nested charts around compact marked representatives of $K^+$. Set
$R_*=R_{\rm gen}(K^+,D_0+1,\delta)$ and use $R=R_*+1$.
The additional unit of time in $K$ again lets us choose $s\le t$
with $q(s)\in K$ and at least $R_*$ units of time in $K$ during $[0,s]$.

Write $F=F_+(q)=F_+(q')$. Choose generic $F_n\to F$ with
$\Ext_{X_q}(F_n)=1$, and put
\[
q_n=Q_1([F_n],X_q),\qquad
a_n=\frac12\log\Ext_{X_{q'}}(F_n),\qquad
p_n=g_{a_n}Q_1([F_n],X_{q'}).
\]
We have $F_+(q_n)=F_+(p_n)=F_n$, and
$a_n\to0$, $q_n\to q$, $p_n\to q'$ in the full bundle.
For large $n$, the base surfaces of $q_n,p_n$ are at distance at most $D_0+1$.
The rays $g_uq_n$ converge uniformly to $g_uq$ on $[0,s]$, so
\begin{equation}
\label{eq:rays-preserve-visits}
q(u)\in K\quad\Longrightarrow\quad q_n(u)\in K^+
\qquad(0\le u\le s)
\end{equation}
for all sufficiently large $n$. The generic case therefore gives
$d_\AGY(q_n(s),p_n(s))<\delta$ with the same threshold $R_*$.

For large $n$, the ambient triangle inequality gives
$D_\AGY(q(s),p_n(s))<2\delta$. Thus the points $p_n(s)$ lie in a
fixed neighborhood whose closure is compact and contained in $\POT$.
Since $p_n(s)\to q'(s)$ in $\OT$, we have $q'(s)\in\POT$;
flow invariance then gives $q'\in\POT$.
We may now use continuity of the ambient AGY distance to obtain
$D_\AGY(q(s),q'(s))\le\delta$. Their measured vertical foliations
agree, so the stable version of Lemma~\ref{lemma:DAGY:vs:dAGY}
gives $d_\AGY(q(s),q'(s))\le4\delta<\ep$.
Stable nonexpansion finishes the proof at time $t$.
\end{proof}

\subsection*{A delay linear in the initial distance}

\begin{proof}[Proof of Theorem~\ref{thm:rays-linear} for closed surfaces]
We apply DDM far enough along the reference ray that the nearby point
it gives cannot lie on $[X_q,X_{q'}]$. Passing to a limit gives a
nearby point on the comparison ray. After matching measured foliations,
we apply Proposition~\ref{cor:rays-converge-new} over the remaining time.

Fix a $\Mod$--compact neighborhood $K^+\subset\POT$ of $K$ and choose
$\delta<1/100$ with $4\delta<\ep$ as in the preceding proof.
By \cite[Theorem~A]{masur:SH}, there are $C>0$ and $L\ge1$ such
that whenever an interval of length $L$ on one side of a finite
\Teich triangle spends at least half its time over the projection
of $K^+$, it comes within $C$ of one of the other two sides. Set
\begin{equation}
\label{eq:rays-constant-order}
\begin{split}
D_*&=2(C+L),\qquad R_*=R(K^+,D_*,\delta),\\
H&=\frac{R_*+1}{1-\eta},\qquad
W=\frac{H+2L+1}{1-2\eta},\qquad B=W+C+1.
\end{split}
\end{equation}
Here $R_*$ comes from Proposition~\ref{cor:rays-converge-new} with
the fixed distance bound $D_*$. We reserve the last $H$ units of
time for convergence and find an interval to which DDM applies in
$[T-W,T-H]$. All these constants are independent of $D$.

With these constants fixed, suppose the common foliation is generic
and uniquely ergodic. We claim that $T\ge D+B$ and
\eqref{eq:rays-terminal-statistics} with $K^+$ in place of $K$ already
imply $d_\AGY(q(T),q'(T+a))<\delta$.
Partition $[T-W,T-H]$ into intervals of length $L$ and a remainder
of length $r<L$. If each full interval spent less than half its
time in $K^+$, the $K^+$--time on $[T-W,T]$ would be less than
\[
\frac{W-H-r}{2}+r+H
<\frac{W+H+L}{2}<(1-\eta)W.
\]
This contradicts \eqref{eq:rays-terminal-statistics} with $u=T-W$
and $K^+$ in place of $K$.
Thus there is a full interval $I$ during which $q$ spends at least
half its time in $K^+$.
Choose $v\in I$ with $q(v)\in K^+$.

Apply DDM to the segment $X_{q(u)}$, $u\in I$, in the triangle with
vertices $X_q,X_{q'},X_{q(N)}$, for $N>T$.
Every $u\in I$ satisfies $u\ge T-W\ge D+C+1$, so this segment is
more than $C$ from $[X_q,X_{q'}]$. DDM therefore gives a point on
$[X_{q'},X_{q(N)}]$ within $C$ of the reference segment.
We show that these finite sides converge to the ray $q'$ as $N\to\infty$.
Normalize $F$ by $\Ext_{X_q}(F)=1$, and let $G_N$ be the vertical
foliation of the unit differential directing the segment from
$X_{q'}$ to $X_{q(N)}$, of length $t_N\ge N-D$. Then
\[
\I(G_N,F)^2\le
\Ext_{X_{q(N)}}(G_N)\Ext_{X_{q(N)}}(F)
=e^{-2(t_N+N)}.
\]
The $G_N$ lie on the compact $\Ext_{X_{q'}}$--unit sphere.
Since $F$ is filling and uniquely ergodic, every limit is the
representative of $[F]$ with extremal length one at $X_{q'}$.
Hubbard--Masur continuity gives
convergence of the segments to $q'$ on bounded time intervals.
The points supplied by DDM stay within $C$ of the fixed reference
segment, so their distances from $X_{q'}$ are bounded. Passing to a
subsequence gives $b\ge0$ and $u\in I$ with
$d_\cT(X_{q'(b)},X_{q(u)})\le C$.

Since $|u-v|\le L$, replacing $q(u)$ by $q(v)$ increases this
base-distance bound by at most $L$.
Extremal-length comparison then bounds the time needed to match
the measured foliations:
\[
d_\cT(X_{q(v)},X_{q'(b)})\le C+L,
\qquad |b-v-a|\le C+L.
\]
Consequently,
\[
d_\cT(X_{q(v)},X_{q'(v+a)})\le2(C+L)=D_*,
\qquad F_+(q(v))=F_+(q'(v+a)).
\]
Also $v+a\ge v-D\ge C+1>0$.
During $[v,T]$, the ray $q$ spends at least
$(1-\eta)(T-v)\ge R_*+1$ time in $K^+$, since $T-v\ge H$.
Proposition~\ref{cor:rays-converge-new}, applied to $q(v)$ and $q'(v+a)$,
gives the required bound $d_\AGY(q(T),q'(T+a))<\delta$.

For the general case, choose generic uniquely ergodic
$[F_n]\to[F]$. Such foliations are dense: unique ergodicity has
full measure \cite{masur:IE}, and nonprincipal types and
saddle-connection walls have measure zero in foliation coordinates.
Keep both starting surfaces fixed and put
\[
q_n=Q_1([F_n],X_q),\qquad q_n'=Q_1([F_n],X_{q'}),\qquad
a_n=\frac12\log\frac{\Ext_{X_{q'}}(F_n)}{\Ext_{X_q}(F_n)}.
\]
Then $q_n\to q$, $q_n'\to q'$, and $a_n\to a$.
Their initial base distance is exactly $D$, and $|a_n|\le D$.
Uniform convergence on $[0,T]$ gives $q_n(u)\in K^+$ whenever
$q(u)\in K$, for all $u\in[0,T]$ and all sufficiently large $n$.
Thus $q_n$ satisfies \eqref{eq:rays-terminal-statistics} with $K^+$
in place of $K$, and the
preceding argument applies with the same constants, giving
\[
d_\AGY(q_n(T),q_n'(T+a_n))<\delta.
\]
Since $q(T)\in K$ and $q_n(T)\to q(T)$, the points $q_n'(T+a_n)$
eventually lie in the ambient $2\delta$--ball about $q(T)$.
This ball lies in a chart whose closure is compact and contained
in $\POT$. The limit $q'(T+a)$ is therefore principal.
As in the proof of Proposition~\ref{cor:rays-converge-new},
continuity of the ambient AGY distance followed by Lemma~\ref{lemma:DAGY:vs:dAGY}
gives $d_\AGY(q(T),q'(T+a))\le4\delta<\ep$.
Flow invariance gives $q'\in\POT$, and
$T+a\ge T-D\ge B>0$.
\end{proof}

\subsection*{Surfaces with punctures}
\hypertarget{rays-punctures}{}

The Hodge-distance comparison used in Lemma~\ref{lem:rays-Hodge-bound}
and the triangle theorem of \cite[Theorem~A]{masur:SH} are stated for
closed surfaces. We transfer the preceding convergence estimates through
a fixed cover branched only over the punctures. This preserves the
principal quadratic stratum and the initial \Teich distance. The
orienting cover also branches over simple zeros and would not preserve this
stratum.

Suppose $S=S_{g,n}$, where $n>0$, $3g-3+n>0$, and $S\ne S_{0,4}$.
The hyperbolic orbifold obtained by giving the punctures order two has
a finite regular surface cover \cite[Section~2, Theorem~2.5]{scott:geometries}.
Write this cover as $h\colon\widehat S\to\overline S$, where $\overline S$ is obtained by
filling the punctures of $S$ and $\widehat S$ is closed of genus at least two.
It has ramification index two over each puncture and is unbranched
elsewhere. Fix this cover independently of the differential.

Write $d=\deg h$. Pullback gives an isometric embedding
$\iota\colon\cT(S)\to\cT(\widehat S)$
\cite[Section~7]{rafi:schleimer:covers}. Its lift to unit-area differentials is
\[
\mathcal L(q)=\widehat q=d^{-1}h^*q.
\]
The map $\mathcal L$ commutes with the flow and is a proper topological
embedding: pullback is continuous and injective, $\iota$ is isometric,
and the unit-area bundle over a compact subset of \Teich space is compact.
If $q$ has order $k$ at a puncture, its pullback has order $2k+2$
at every point above it. A simple pole becomes a regular point;
if $q$ is holomorphic at the puncture, its pullback has a multiple
zero above that puncture. Thus $q$ is principal if and only if
$\widehat q$ is principal. The normalization also gives
\[
F_+(\widehat q)=d^{-1/2}h^*F_+(q),
\qquad
\Ext_{\iota(X)}(d^{-1/2}h^*F)=\Ext_X(F),
\]
where the second identity follows by pulling back the Hubbard--Masur
differential and comparing areas.

To transfer the estimates, we need closeness of the lifts to imply
closeness downstairs, uniformly for $q$ in a $\Mod$--compact set.
The subgroup
of $\Mod(S)$ that lifts through $h$ has finite index, since there are
only finitely many covers of a fixed degree. Consequently, for every
$\Mod(S)$--compact $K\subset\POT(S)$, the set
\[
\widehat K=\Mod(\widehat S)\cdot\mathcal L(K)
\]
is $\Mod(\widehat S)$--compact. Choose a compact set of representatives
for $K$ under the liftable subgroup. Properness of $\mathcal L$ and
the fact that the ambient AGY metric induces the principal-stratum
topology show that, for each $\ep>0$, there is $\delta>0$ such that
\[
D_\AGY(\widehat q,\widehat p)<\delta
\quad\Longrightarrow\quad
D_\AGY(q,p)<\min\{1/100,\ep/4\}
\]
whenever $q\in K$, $p\in\OT(S)$, and $\widehat p$ is principal.
The same $\delta$ works for every $q\in K$, since the mapping classes
in this subgroup lift through $h$.
If also $F_+(p)=F_+(q)$, Lemma~\ref{lemma:DAGY:vs:dAGY}, applied to $-q,-p$,
gives $d_\AGY(q,p)<\ep$.

For Proposition~\ref{cor:rays-converge-new}, choose $\delta$ as above,
let $\widehat R$ be the closed-surface threshold for
$\widehat K,D_0,\delta$, and set $R=\widehat R+1$.
If $q$ spends at least $R$ time in $K$ during $[0,t]$, there is a return $s\le t$ to $K$
after at least $\widehat R$ time has been spent there: the cumulative
occupation time is continuous, and one further unit of $K$--time
remains after it reaches $\widehat R$. Every $K$--visit lifts to a
$\widehat K$--visit. Apply the closed-surface proposition at $s$ to
obtain $\widehat q'$ principal and
$d_\AGY(\widehat q(s),\widehat q'(s))<\delta$.
The cover preserves principality, and the preceding comparison gives
$d_\AGY(q(s),q'(s))<\ep$; stable nonexpansion gives the same bound at
$t$. The endpoint $q(t)$ need not lie in $K$.

For Theorem~\ref{thm:rays-linear}, the lifted ray satisfies
\eqref{eq:rays-terminal-statistics} with $\widehat K$ in place of $K$
and the same $\eta$.
Apply the closed-surface theorem with target distance $\delta$ and
use the preceding comparison at $q(T)$ and $q'(T+a)$, where $q(T)\in K$.
The extremal-length identity preserves
the time shift $a$, and the isometry preserves the initial distance
$D$. Thus the conclusion still holds for $T\ge D+B$ with coefficient
one on $D$. Only the additive constant and the upstairs accuracy change.

For $S_{0,4}$, every nonzero quadratic differential is principal and its
double cover is a flat torus. Two unit-area differentials $q,p$ with
the same measured vertical foliation differ by a shear with parameter
$c$, where $|c|=2\sinh d_\cT(X_q,X_p)$. Under the flow this parameter
becomes $ce^{-2t}$. Along the straight shear path a relative period
$z=x+iy$ has $dz=ix\,dc$, so $|dz|/|z|\le|dc|$.
This applies to lifts of saddle connections downstairs, and hence
\[
d_\AGY(q(t),p(t))
\le 2\sinh\bigl(d_\cT(X_q,X_p)\bigr)e^{-2t}.
\]
This proves Proposition~\ref{cor:rays-converge-new}, since elapsed time
is at least compact occupation time. For the theorem, $p=q'(a)$ has
the same measured vertical foliation as $q$, and $|a|\le D$ gives
$d_\cT(X_q,X_p)\le2D$. Therefore
\[
d_\AGY(q(T),q'(T+a))
\le2\sinh(2D)e^{-2T}
\le e^{-2(T-D)}.
\]
Choosing $B>0$ sufficiently large proves the theorem with $T\ge D+B$
and $T+a>0$, in this case without an occupation assumption.

\subsection*{Choosing \texorpdfstring{$K_2$}{K2}}

We enlarge $K_1$ to accommodate the bounded time shifts and nearby
differentials used in chart comparisons.

\begin{definition} \label{def:K_2}
Let $K_2$ be a $\Mod$--compact subset of $\POT$ of the form
\[
K_2 = \{q\in\POT:\sys(q)\ge \ss_2\},
\]
where $\ss_2>0$ is chosen small enough so that the following holds. 
If $q_1, q_2$ are points in $\POT$ and $u,v \in \reals$ are such that 
\[
q_1 \in K_1, \quad
\max(|u|,|v|) \leq \max(\sD_0, \sD_{\rm chart}, \sD_{\rm temp})
\qquad\text{and}\quad
D_\AGY(q_1(u), q_2(v)) \leq \tfrac 1{100}, 
\]
then $q_2\in K_2$.
\end{definition}

In later applications of Proposition~\ref{cor:rays-converge-new}, we use
$K=K_2$ and $\ep=1/500$ unless another accuracy is specified.
The auxiliary $\Mod$--compact neighborhoods used in the proofs do not
change this choice of $K_2$.

\section{Gromov distance on $\PMF$}
\label{sec:gromov-distance}
We construct a metric on $\PMF$ using separation times in the charts
of Section~\ref{sec:Charts}. We record the last time at which two rays
are AGY-close and both have good statistics; later separation gives
a smaller distance. We first construct chart metrics and then a global
metric, whose balls we compare with horospherical AGY balls.

Fix $\eta_2=1/4$ for the remainder of the paper, and assume
$0<\eta_1<\eta_2$. The smaller value of $\eta_1$ needed for the spike
estimates will be specified in \eqref{eq:spike-error-choice}.
Thus $\eta_2$ is fixed before the separation times and metric constants
defined below.

For $x,y\in\PMF$ and $F\in\MF$, define the separation time
\begin{multline}
\label{eq:def-tau}
\tau_F(x,y) =
  \sup \Big\{ t \in \R  :  \text{$x_F(t)$ and $y_F(t)$ have $(K_2, \eta_2)$--good statistics}\\
  \text{and} \quad d_\AGY\big(x_F(t),y_F(t)\big) \le \tfrac 1{100}  \Big\}.
\end{multline}
The time \index{$\tau_F(x,y)$}$\tau_F(x,y)$ may be negative or infinite.
We set it equal to $-\infty$ if the admissible set is empty or if either
point lies outside $\PMF(F)$. Let
\begin{equation}
\label{eq:beta} 
\beta = \beta_{K_2},
\end{equation} 
\begin{equation}
\label{eq:CAGY}
C_\AGY = C_{K_2} \ge 1
\end{equation}
as in Corollary~\ref{cor:AGY:half}. For a chart $\alpha\in\cA$, define
\[
\tau_\alpha(x,y)
=
\begin{cases}
\tau_{F_\alpha}(x,y), & x,y\in W_\alpha,\\[2mm]
-1, & \text{otherwise}.
\end{cases}
\]

Define the chart pre-distance by
\[
D_G^\alpha(x,y)
=
e^{-\beta\tau_\alpha(x,y)}.
\]
Here $e^{-\beta(+\infty)}=0$. The value $e^\beta$ outside the chart
ensures that a chain starting in $W_\alpha$ with cost at most one
stays inside it.

%
%

\begin{proposition}
\label{prop:D}
For any chart $\alpha\in\cA$ and distinct $x,y\in W_\alpha$,
the supremum in \eqref{eq:def-tau} is attained, and
\[
0<D_G^\alpha(x,y)\leq 1. 
\]
If either $x$ or $y$ lies outside $W_\alpha$, then
\[
D_G^\alpha(x,y)=e^\beta>1.
\]
\end{proposition}
\begin{proof}
By Theorem~\ref{theorem:charts}{\rm (g)}, both negative rays lie in
$K_0\subset K_2$. Closedness of $K_2$ includes their endpoints at
time zero, which therefore have $(K_2,\eta_2)$--good statistics.
Equation~\eqref{eq:dagy:x:F:alpha:y:F:alpha} gives
$d_\AGY(x^\alpha(0),y^\alpha(0))\le1/100$, so zero is admissible
and $\tau_{F_\alpha}(x,y)\ge0$.

Next we show that $\tau_{F_\alpha}(x,y)$ is finite. Otherwise there is
a sequence $t_i\to+\infty$ such that both $x^\alpha(t_i)$ and
$y^\alpha(t_i)$ have $(K_2,\eta_2)$--good statistics and
\[
d_\AGY\big(x^\alpha(t_i),y^\alpha(t_i)\big)\le\frac1{100}.
\]
For all sufficiently large $i$, Corollary~\ref{cor:AGY:half} applies on
the interval $[0,t_i]$ and gives
\[
d_\AGY\big(x^\alpha(0),y^\alpha(0)\big)
\le \frac1{100}C_\AGY e^{-\beta t_i}.
\]
Letting $i\to\infty$ contradicts $x\ne y$.

Finally, the admissible set is closed. For each fixed backward interval
length, occupation time varies continuously with the endpoint. Together
with closedness of $K_2$, this makes the set of good-statistics times
closed; the AGY distance condition is also closed. The finite supremum
is therefore attained. The value outside $W_\alpha$ follows from the
definition.
\end{proof}

To obtain the triangle inequality, take the infimum of the costs of
finite chains:
\begin{displaymath}
d_G^\alpha(x, y) = \inf \sum_{i=0}^{k-1} D_G^\alpha(x_i, x_{i+1}),
\end{displaymath}
where $x=x_0,\ldots,x_k=y$ is a finite sequence in $\PMF$.
We allow $k=0$ when $x=y$, with cost zero. This defines the chart
metric on $W_\alpha$ (Corollary~\ref{cor:alpha-metric}), even though
the pre-distance need not vanish on the diagonal.

\begin{lemma}
\label{lem:dG:lower:bound}
For any chart $\alpha$ and $x,y \in W_\alpha$ (with $x \ne y$), let $s$ be a time
such that either 
\begin{equation} \label{eq:assumption}
d_\AGY(x^\alpha(s), y^\alpha(s)) \leq \tfrac 1{100} 
\qquad\text{or}\qquad
d_G^\alpha(x, y) \leq e^{-\beta s}/2. 
\end{equation}
Then 
\begin{equation} \label{eq:not-too-close}
100 \, C_\AGY^{-1} d_\AGY(x^\alpha(s),y^\alpha(s)) e^{-\beta s} \le d_G^\alpha(x,y).
\end{equation}
Also, for $t=\tau_{F_\alpha}(x,y)$, we have 
\begin{equation} \label{eq:d-D}
100 \, C_\AGY^{-1} d_\AGY(x^\alpha(t),y^\alpha(t)) D_G^\alpha(x,y)\le d_G^\alpha(x,y) 
  \le D_G^\alpha(x,y).
\end{equation}
\end{lemma}

\begin{proof}
Take chains $x=x_0,\ldots,x_k=y$ whose costs approach
$d_G^\alpha(x,y)$, and delete consecutive repetitions. We may
assume every vertex lies in $W_\alpha$: an edge leaving the chart
costs $e^\beta>1$, whereas the direct chain costs at most one.
If $d_G^\alpha(x,y)\le e^{-\beta s}/2$, we also take the chain costs
less than $e^{-\beta s}$.

Put $t_i=\tau_{F_\alpha}(x_i,x_{i+1})$. These times are finite
and attained by Proposition~\ref{prop:D}. If every $t_i\ge s$,
the good statistics at $t_i$ and Corollary~\ref{cor:AGY:half} give
\[
d_\AGY\big(x_i^\alpha(s),x_{i+1}^\alpha(s)\big)
\le\frac{C_\AGY}{100}e^{-\beta(t_i-s)}.
\]
Since $D_G^\alpha(x_i,x_{i+1})=e^{-\beta t_i}$, rearranging and
using the triangle inequality gives
\begin{equation}
\label{eq:path:lower:bound}
\sum_{i=0}^{k-1}D_G^\alpha(x_i,x_{i+1})
\ge100\,C_\AGY^{-1}
d_\AGY\big(x^\alpha(s),y^\alpha(s)\big)e^{-\beta s}.
\end{equation}

If some $t_i<s$, the chain has cost at least
$e^{-\beta t_i}>e^{-\beta s}$. Our choice of chains excludes
this case under the second alternative in \eqref{eq:assumption}.
The first alternative must therefore hold:
\begin{equation}
\label{eq:s-close}
d_\AGY\big(x^\alpha(s),y^\alpha(s)\big)\le\frac1{100}.
\end{equation}
Since $C_\AGY\ge1$, the right-hand side of
\eqref{eq:path:lower:bound} is at most $e^{-\beta s}$, so the
same bound holds in this case. Taking the infimum proves
\eqref{eq:not-too-close}.

For \eqref{eq:d-D}, the direct chain gives the upper bound.
The lower bound follows by applying \eqref{eq:not-too-close} at
the attained separation time $t=\tau_{F_\alpha}(x,y)$, when
the first alternative in \eqref{eq:assumption} holds and
$e^{-\beta t}=D_G^\alpha(x,y)$.
\end{proof}

\begin{corollary}
\label{cor:alpha-metric}
$d_G^\alpha$ is a metric on $W_\alpha$. 
\end{corollary} 
\begin{proof}
The chain definition gives symmetry, the triangle inequality, and
$d_G^\alpha(x,x)=0$. For $x\ne y$, Proposition~\ref{prop:D} gives
$D_G^\alpha(x,y)>0$, and \eqref{eq:d-D} then gives $d_G^\alpha(x,y)>0$.
\end{proof}

We now combine the chart pre-distances using the partition of unity
$\{\theta_\alpha\}$ from Theorem~\ref{theorem:charts}{\rm (e)}.
For $x,y\in\PMF$, define
\[
D_G(x,y) = \sum_{\alpha \in \cA} 
  \frac{\theta_\alpha(x) + \theta_\alpha(y)}2 \cdot D_G^\alpha(x, y)
\]
and 
\begin{displaymath}
d_G(x, y) = \inf \sum_{i=0}^{k-1} D_G(x_i, x_{i+1}),
\end{displaymath}
where the infimum is over finite chains $x=x_0,\ldots,x_k=y$ in
$\PMF$, again allowing the zero-edge chain when $x=y$.

Choose $\rr_1$ small enough such that 
\[
4 |\cA|\cdot \rr_1\le\frac12
\qquad\text{and}\qquad
\frac{4 C_\AGY |\cA| \cdot\rr_1}{100}\le\rr_{\rm chart}.
\]
This in particular means that 
\begin{equation}
\label{eq:r1:consequence}
\Big( x\in V^{(2)}_\alpha\quad\text{and}\quad 
d_G^\alpha(x,y)\le4|\cA|\rr_1 \Big)
\quad\Longrightarrow\quad 
y\in V^{(3)}_\alpha.
\end{equation}
In fact, assuming the lefthand side of \eqref{eq:r1:consequence},
we have $d_G^\alpha(x,y)\le1/2$ which forces $y\in W_\alpha$, since an
endpoint outside the chart would force an edge of cost $e^\beta$. Then,
Lemma~\ref{lem:dG:lower:bound}, applied with $s=0$, gives
\[
d_\AGY\big(x^\alpha(0),y^\alpha(0)\big)
  \le \frac{4 C_\AGY|\cA|\rr_1}{100}
  \le \rr_{\rm chart}.
\]

Theorem~\ref{theorem:charts}{\rm (f)} therefore gives $y\in V^{(3)}_\alpha$.

\begin{lemma}
\label{lem:dGalpha:vs:dG}
Let $x \in V^{(2)}_\alpha$ and $y \in \PMF$ be such that $d_G(x,y) \leq \rr_1$. Then 
\[
D_G(x,y) \geq d_G(x,y) \geq \frac 1{2|\cA|} d_G^\alpha(x,y). 
\]
\end{lemma}

\begin{proof}
The first inequality is immediate from the definition of $d_G$. We prove the second
inequality. Without loss of generality, assume $x \not = y$. 

Take a chain $x=x_0,\ldots,x_k=y$ with
\begin{displaymath}
 \sum_{i=0}^{k-1} D_G(x_i, x_{i+1}) \leq 2 \rr_1. 
\end{displaymath}
This, in particular, means for every $j \le k$,
\[
 \sum_{i=0}^{j-1} D_G(x_i, x_{i+1}) \leq  2\rr_1. 
\]
Assuming, by induction, that $\theta_\alpha(x_i) \geq \frac 1{|\cA|}$ for $i=0, \dots, (j-1)$, 
we have 
\[
\frac 1{2|\cA|} \sum_{i=0}^{j-1} D_G^\alpha (x_i, x_{i+1}) \leq 2\rr_1,
\]
which means $d_G^\alpha(x, x_j) \leq 4 |\cA| \cdot \rr_1$. Now \eqref{eq:r1:consequence}
implies $x_j \in V^{(3)}_\alpha$ and hence $\theta_\alpha(x_j) \geq \frac 1{|\cA|}$. 
We have shown that, for $i = 0, \dots, k$, $\theta_\alpha(x_i) \geq \tfrac 1{|\cA|}$. 
Therefore, 
\[
\sum_{i=0}^{k-1} D_G(x_i, x_{i+1})
\geq \frac 1{2|\cA|} \sum_{i=0}^{k-1} D_G^\alpha (x_i, x_{i+1})
\geq  \frac 1{2|\cA|} d_G^\alpha(x,y). 
\]
Taking the infimum over all such sequences $x_i$ we get 
\[
d_G(x,y) \geq  \frac 1{2|\cA|} d_G^\alpha(x,y). 
\]
This finishes the proof. 
\end{proof}

Let $R$ be the constant provided by Proposition~\ref{cor:rays-converge-new}
with $K=K_2$, accuracy $1/500$, and initial \Teich distance bound
$D=\sD_{\rm temp}+\sD_{\rm chart}$. Choose
\begin{equation} \label{eq:R1}
\sR_1> \max\left\{ \sD_{\rm chart}, \frac{R}{\eta_2-\eta_1} \right\}.
\end{equation} 

\begin{theorem}
\label{theorem:metric-comparison}
Let $x\in V^{(3)}_\alpha$, $y\in W_\alpha$, and $t>\sR_1$. Assume that
$x^\alpha(t)$ has $(K_1,\eta_1)$--good statistics, and set $r=e^{-\beta t}$. Then:
\begin{enumerate}
\item[{\rm (a)}]
\[
d_\AGY\big(x^\alpha(t),y^\alpha(t)\big) 
\le \frac1{100\kk}
\qquad\Longrightarrow\qquad
d_G(x,y) < r.
\]

\item[{\rm (b)}]
\[
d_G^\alpha(x,y) \le \frac {r}{C_\AGY \kk}
\qquad\Longrightarrow\qquad
d_G(x,y)< r.
\]
\end{enumerate}
\end{theorem}

\begin{proof}
For part {\rm (a)}, we must transfer closeness and good statistics to
every chart contributing to $D_G(x,y)$. The first $R$ units of time in
$K_1$ let us apply Proposition~\ref{cor:rays-converge-new} to bring
representatives of a common foliation within AGY distance $1/500$.
The margin $\eta_2-\eta_1$ absorbs these $R$ units of compact time
when we check good statistics in the new chart.

By \eqref{eq:Monotone},
\[
d_\AGY\big(x^\alpha(0),y^\alpha(0)\big)
\le d_\AGY\big(x^\alpha(t),y^\alpha(t)\big)
\le \frac1{100\kk}.
\]
Let $s_0>0$ be the first time at which
\begin{displaymath}
|\{ u \in [0,s_0] \st x^\alpha(u) \in K_1 \}| = R.
\end{displaymath}
The occupation function is continuous, and the good statistics at $t$
give at least $(1-\eta_1)t>R$ units of time in $K_1$ by
\eqref{eq:R1}. Thus $s_0$ exists and $s_0<t$.

Let $\cA'$ be the set of charts
$\alpha'\in\cA$ such that
\[
\theta_{\alpha'}(x)+\theta_{\alpha'}(y)>0.
\]
Fix $\alpha'\in\cA'$. Since
\[
\supp(\theta_{\alpha'})\subset V^{(4)}_{\alpha'},
\]
at least one of $x,y$ lies in $V^{(4)}_{\alpha'}$. Denote this point by $z$
and the other point by $w$. Thus
\[
z,w\in W_\alpha,
\qquad
z\in V^{(4)}_{\alpha'}\subset W_{\alpha'}.
\]
By part (b) of Lemma~\ref{lem:uniform-change-of-charts}, there is a time $t_0$
with $|t_0|\le\sD_{\rm chart}$ such that 
\[
F_+\big(z^{\alpha'}(t_0)\big)=F_+\big(z^\alpha(0)\big).
\]
Moreover, part (a) of that lemma gives
\[
d_\cT\big(z^\alpha(0),z^{\alpha'}(0)\big)
\le\sD_{\rm temp},
\]
and hence
\[
d_\cT\big(z^\alpha(0),z^{\alpha'}(t_0)\big)
\le \sD_{\rm temp} + \sD_{\rm chart}.
\]

The segment $x^\alpha([0,s_0])$ spends at least $R$ time in $K_1$.
If $z=x$, the same is therefore true for $z^\alpha([0,s_0])$. If $z=y$, then
for every $u\in[0,s_0]$,
\[
d_\AGY\big(x^\alpha(u),y^\alpha(u)\big)
\le
d_\AGY\big(x^\alpha(t),y^\alpha(t)\big)
\le
\frac1{100\kk}.
\]
Thus, whenever $x^\alpha(u)\in K_1$, the choice of $K_2$ implies
$y^\alpha(u)\in K_2$. In either case,
\begin{displaymath}
|\{ u \in [0,s_0]  \st z^\alpha(u) \in K_2 \}| \ge R. 
\end{displaymath}
Proposition~\ref{cor:rays-converge-new}, applied to the two rays
with common vertical foliation, gives
\[
d_\AGY\big(z^\alpha(s_0),z^{\alpha'}(s_0 +  t_0)\big)
\le \frac1{500}.
\]
Stable nonexpansion (Theorem~\ref{theorem:properties:AGY}) gives
\begin{displaymath}
d_\AGY\big(z^\alpha(u),z^{\alpha'}(u +  t_0)\big)
\le \frac1{500}, \qquad u \ge s_0. 
\end{displaymath}
We also have, using \eqref{eq:Monotone}, for every $0 \le u \le  t$, 
\[
d_\AGY\big(z^\alpha(u),w^\alpha(u)\big)
\le \frac1{100\kk}.
\]
Part (c) of Lemma~\ref{lem:uniform-change-of-charts} (applied with
$x=z$, $y=w$, and $s_0\le u\le t$) therefore gives
\[
x,y\in W_{\alpha'}, \qquad
d_\AGY\big(x^{\alpha'}(t),y^{\alpha'}(t)\big)<1/100,  
\]
and,  for $s_0 \le  u \le t$, 
\[
D_\AGY\big(x^\alpha(u),p^{\alpha'}(u+t_0)\big)<1/100
\qquad(p\in\{x,y\}).
\]
This implies that, for $s_0 \le u \le t$, if
\[
x^{\alpha}(u) \in K_1
\]
then in view of the definition of $K_2$ (Definition~\ref{def:K_2}), we have 
\[
y^{\alpha}(u) \in K_2, \qquad x^{\alpha'}(u) \in K_2 \qquad\text{and}\qquad y^{\alpha'}(u) \in K_2. 
\]
We now claim that $x^{\alpha'}(t)$
and $y^{\alpha'}(t)$ have $(K_2, \eta_2)$--good statistics. Let $z  =
x$ or  $y$. For $s_0 \le s \le t$, 
\begin{displaymath}
|\{ u \in [s,t] \st z^{\alpha'}(u)  \in K_2  \}| \ge |\{ u \in [s,t]
\st x^{\alpha}(u)  \in K_1  \}| \ge (1-\eta_1) (t-s),  
\end{displaymath}
and for $0 \le s \le  s_0$, 
\begin{align*}
|\{ u \in [s,t] \st z^{\alpha'}(u)  \in K_2  \}|  & \ge
|\{ u \in [s_0,t] \st z^{\alpha'}(u)  \in K_2  \}|  \\
& \ge |\{ u \in [s_0,t] \st x^{\alpha}(u)  \in K_1  \}|  \\
& = |\{ u \in [0,t] \st x^{\alpha}(u)  \in K_1  \}| - R \\
& \ge (1-\eta_1)t - R \ge (1-\eta_2)t \ge (1-\eta_2)(t-s).   
\end{align*}
For $s<0$, Theorem~\ref{theorem:charts}{\rm (g)} gives
$z^{\alpha'}(s)\in K_0\subset K_2$. Adding the interval $[s,0]$
to the estimate for $[0,t]$ proves the required bound for every
$s<0$ as well.

We can extend these statistics a little past $t$. Indeed, hold
$x^\alpha(t)\in K_1$ fixed. For $p\in\{x,y\}$ and all sufficiently
small $h\ge0$, continuity preserves the strict bound
\[
D_\AGY\big(x^\alpha(t),p^{\alpha'}(t+h+t_0)\big)<1/100.
\]
Definition~\ref{def:K_2}, with shifts $0,t_0$, puts
$p^{\alpha'}(t+h)$ in $K_2$. This additional interval entirely in
$K_2$ preserves good statistics, and continuity also preserves
$d_\AGY(x^{\alpha'}(t+h),y^{\alpha'}(t+h))<1/100$.
Hence
\[
\tau_{F_{\alpha'}}(x,y)> t,
\]
and therefore
\[
D_G^{\alpha'}(x,y)< e^{-\beta t}=r.
\]
It follows that
\[
D_G(x,y)
= \sum_{\alpha'\in\cA'} \frac{\theta_{\alpha'}(x)+\theta_{\alpha'}(y)}2 D_G^{\alpha'}(x,y)
< r \sum_{\alpha'\in\cA'} \frac{\theta_{\alpha'}(x)+\theta_{\alpha'}(y)}2 = r.
\]
Thus
\[
d_G(x,y)\le D_G(x,y)< r,
\]
which proves part (a).

We now prove part (b). If $x=y$, the statement trivially
holds. Suppose $x \ne y$ and
\[
d_G^\alpha(x,y)\le\frac{r}{C_\AGY\kk}.
\]
Since $\kk\ge2$  and $C_\AGY \ge  1$,
the second alternative in \eqref{eq:assumption} holds with
$s=t$. Lemma~\ref{lem:dG:lower:bound} therefore gives 
\[
100\, C_\AGY^{-1} d_\AGY\big(x^\alpha(t),y^\alpha(t)\big)e^{-\beta t}
\le
d_G^\alpha(x,y)
\le
\frac{r}{C_\AGY\kk}.
\]
Since $r=e^{-\beta t}$, we obtain
\[
d_\AGY\big(x^\alpha(t),y^\alpha(t)\big)
\le
\frac1{100\kk}.
\]
Part (a) now gives
\[
d_G(x,y) < r. \qedhere 
\]
\end{proof}

\begin{corollary}
\label{cor:metric}
$d_G$ is a metric on $\PMF$. Furthermore, $(\PMF, d_G)$ has bounded diameter. 
\end{corollary} 
\begin{proof}
The chain definition gives symmetry, the triangle inequality, and
$d_G(x,x)=0$. Lemma~\ref{lem:dGalpha:vs:dG} and
Corollary~\ref{cor:alpha-metric} give positivity when $x\ne y$.
Since each chart pre-distance is at most $e^\beta$ and the
partition-of-unity weights sum to one,
\[
d_G(x,y)\le D_G(x,y)\le e^\beta. \qedhere
\] 
\end{proof}

The balls of $d_G^\alpha$ and $d_G$ are Borel in the usual topology on $\PMF$.
With $F=F_\alpha$, the conditions in \eqref{eq:def-tau} are closed on
$W_\alpha^2\times\R$, since the leafwise distance is continuous and the
time spent in the closed set $K_2$ on each fixed backward interval is
upper semicontinuous. For $r>0$, exhaust the open subset where
$e^{-\beta t}<r$ by compact sets and project the admissible triples
to $(x,y)$. This writes $\{D_G^\alpha<r\}$ inside the chart as a
countable union of compact sets, and the cost is constant off $W_\alpha^2$.
Taking rational upper bounds for the individual costs gives the same
property for finite sums with continuous nonnegative weights, hence for
the cost of a chain with a fixed number of edges. Projecting to the
endpoints and taking the union over chain lengths gives the same property
for $\{d_G^\alpha<r\}$ and $\{d_G<r\}$.

\subsection{Volumes of balls}
\label{sec:volume-balls}

We use the horospherical volume forms of
\cite[Example~4.2]{AG}. We first specify their normalization and its
relation to Thurston measure; the factor relating a hypersurface measure
to the volume of its cone will be needed below. For $q\in\POT$, recall that
\[
\Hor(q)=\Hor_-(F_-(q)) = \big\{Q_1(x,F_-(q)) : x\in\PMF(F_-(q)) \big\}.
\]
Write the period of $q$ as $a+ib\in L(q)+iL(q)$, and let $\lambda_q$
be the real period volume form on $L(q)$ for which the integral period
lattice has covolume one. Here the lattice is
$L(q)\cap H^1(\widetilde X,\cZ(q);\mathbb Z)$, using the relative
realization of $L(q)$. We use the same notation for its associated
Lebesgue measure. Gauss--Manin transport preserves this measure. Let
$A_b(u)=A(u,b)$ be the area functional with imaginary period $b$
fixed, so that $A_b(a)=1$. In local period coordinates,
\[
T_q\Hor(q)=\ker A_b.
\]
Following Avila--Gou\"ezel, define the positive horospherical measure by
the volume density
\begin{equation}
\label{eq:AG-hor-volume}
d\mu_q(v_1,\ldots,v_{w-1})
=\big|\lambda_q(a,v_1,\ldots,v_{w-1})\big|,
\qquad v_j\in\ker A_b.
\end{equation}
This density is constant along the affine horosphere: replacing $a$ by
$a+u$, with $u\in\ker A_b$, does not change its restriction to
$\ker A_b$. Thus $\mu_q$ depends only on the horosphere containing $q$.

We retain the Thurston measure $\muT$ fixed in Section~\ref{sec:prelims}.
In real period coordinates it is $c_S\lambda_q$, where $c_S>0$ is a
fixed constant depending only on the surface and the two lattice
normalizations. To see that this factor is independent of $q$, use the
adapted train tracks of Proposition~\ref{prop:Local-Product}. Their transverse weights
are signed real edge periods. The Thurston symplectic form in these
coordinates is one half of the cup-product form on odd cohomology;
this is the train-track calculation in
\cite[proof of Theorem~5.8]{dumas:skinning}, applied to real periods.
The same calculation includes punctured monogons at simple poles.
Taking top exterior powers gives the stated constant factor. The local
identities agree on overlaps and extend along positive rays by
homogeneity.

For a measurable set $E$ in a small horospherical period chart, write
$a_p+ib$ for the period of $p\in E$ and set
\[
\mathcal C_q(E)=\{t a_p:p\in E,\ 0<t\le1\}\subset L(q).
\]
This is a cone in the vector space; no extension of the period chart
along the whole cone is required. Since $A_b(a_p)=1$, the map
$(t,p)\mapsto t a_p$ is injective and has radial Jacobian $t^{w-1}$
relative to \eqref{eq:AG-hor-volume}. Therefore
\begin{equation}
\label{eq:AG-cone-factor}
\lambda_q(\mathcal C_q(E))=\frac1w\mu_q(E),
\qquad
\muT\bigl(\{tF_+(p):p\in E,\ 0<t\le1\}\bigr)
=\frac{c_S}{w}\mu_q(E).
\end{equation}
By subdividing into period charts, the second equality holds for every
measurable $E$ in the principal-stratum part of $\Hor(q)$. Equivalently,
\begin{equation}
\label{eq:AG-Thurston-cone}
\mu_q(E)=\frac{w}{c_S}\,
\muT\Bigl(\bigl\{F\in\MF:
Q_1([F],F_-(q))\in E,\ \I(F,F_-(q))\le1\bigr\}\Bigr).
\end{equation}

For any $E\subset W_\alpha$ and $t\ge0$, define
\[
E_\alpha(t)=\{x^\alpha(t):x\in E\}.
\]
Then
\[
E_\alpha(0)\subset \Hor(q_\alpha)
\qquad\text{and}\qquad
E_\alpha(t)\subset \Hor(g_tq_\alpha).
\]
Since $F_-(g_tq_\alpha)=e^{-t}F_\alpha$, we have
\begin{align*}
\mu_{g_tq_\alpha}(E_\alpha(t))
&= \frac{w}{c_S}\,
\muT\Big( \big\{ F\in\MF: [F]\in E,\  \I(F,e^{-t}F_\alpha)\le 1 \big\} \Big) \\
&= \frac{w}{c_S}\,
\muT\Big( \big\{ F\in\MF: [F]\in E,\  \I(F,F_\alpha)\le e^t \big\} \Big).
\end{align*}
By homogeneity of Thurston measure of degree $w=\dim(\MF)$, this gives
\begin{equation}
\label{eq:Volume-Growth}
\mu_{g_tq_\alpha}(E_\alpha(t))
=
e^{wt}\mu_{q_\alpha}(E_\alpha(0)).
\end{equation}
In what follows, when the base point is clear, we write $\mu$ for the
corresponding horospherical measure $\mu_q$.

\begin{lemma} \label{lem:MV-Thurston}
For any measurable set $E \subset W_\alpha$, we have 
\[
\nu_\rho(E) \asymp \mu (E_\alpha(0)).
\]
Here $\mu$ is the horospherical measure on $\Hor(q_\alpha)$ defined above.
\end{lemma}

\begin{proof}
For $x\in W_\alpha$, set $F_x=F_+(x^\alpha(0))$, so that
$\I(F_x,F_\alpha)=1$. The cone defining $\nu_\rho$ cuts the ray
through $F_x$ at $\rho(F_x)^{-1}F_x$, whereas the cone in
\eqref{eq:AG-Thurston-cone} cuts it at $F_x$. Homogeneity of degree
$w$ and \eqref{eq:AG-cone-factor} therefore give the exact density
\begin{equation}
\label{eq:Thurston-hor-density}
d\nu_\rho(x)
=\frac{c_S}{w\,m_\rho}\,
\rho(F_+(x^\alpha(0)))^{-w}\,d\mu_{q_\alpha}(x^\alpha(0)).
\end{equation}
Here $m_\rho=\muT(\{\rho\le1\})$ is the normalizing mass from
Section~\ref{sec:prelims}. The displayed density is positive and
continuous on the compact chart closure, so it is bounded above and
below there.
\end{proof}

We next bound the volume of a small intrinsic horospherical ball, uniformly
even when its center leaves every compact set. Shuvaeva proves an ambient
small-ball estimate for the holomorphic principal stratum
\cite[Theorem~2.1.3]{shuvaeva:AGY-volume}. For the horospherical measure
we give a direct argument that also includes simple poles. The key is
that the real AGY norm controls the Hodge norm, whose balls have fixed
volume in odd cohomology.

\begin{lemma} \label{lem:AGY-volume}
There is a constant $\sV_0>0$, depending only on $S$, such that for every
$q\in\POT$,
\[
\mu_q\bigl(B(q,\tfrac1{100})\bigr)\leq\sV_0.
\]
Here $B(q,r)$ is the ball in $\Hor(q)$ for the intrinsic metric
$d_\AGY$. Furthermore, for every $\Mod$--compact set $K\subset\POT$
there is a constant $\sV_1=\sV_1(S,K)>0$ such that, for every $q\in K$,
\[
\mu_q\bigl(B(q,\tfrac1{100})\bigr)\geq\sV_1.
\]
\end{lemma}

\begin{proof}
We first bound the volume of the real norm ball
\[
\mathbb B_q=\{u\in L(q):\|u\|_q\le1\}.
\]
We then place the cone over a small horospherical ball inside a fixed
multiple of $\mathbb B_q$ and use \eqref{eq:AG-cone-factor}.

Lemma~\ref{lem:AGY-area} gives
$\|u\|_{H,q}\le\sqrt2\,\|u\|_q$ for $u\in L(q)$,
including at nongeneric centers and in the presence of simple poles.

The $\lambda_q$--volume of a Hodge unit ball in $L(q)$ depends only on
$S$. To check this,
let $\Omega$ be the cup-product form restricted to $L(q)$, and write
$w=2m$. The form $\Omega$ is nondegenerate, and the Hodge star $J_q$
preserves $L(q)$ because it commutes with the covering involution.
The Hodge inner product is
\[
G_q(u,v)=\Omega(u,J_qv).
\]
Since $J_q^2=-1$ and $\det J_q=1$, the volume form of $G_q$ is
$|\Omega^m|/m!$. Its ratio to $\lambda_q$ depends only on the
topological orienting cover and its period lattice. In particular,
\eqref{eq:AGY-Hodge-upper} implies that there is a constant
$C_{\rm norm}=C_{\rm norm}(S)$ such that
\begin{equation}
\label{eq:real-AGY-ball-volume}
\lambda_q(\mathbb B_q)
\le\lambda_q\bigl(\{u:\|u\|_{H,q}\le\sqrt2\}\bigr)
\le C_{\rm norm}
\qquad(q\in\POT).
\end{equation}

The norm comparison used all saddle connections, without excluding
any direction. The identification of odd relative and absolute
cohomology and the volume of the Hodge ball are topological. Thus
the bound is uniform also when the center has horizontal or vertical
saddle connections.

Now put $r_0=1/200$ and $E=B(q,r_0)$. The uniform leaf chart of
Lemma~\ref{lem:AGY:vs:cV}, valid also when the fixed foliation has
saddle connections, identifies $E$ with a subset of $\ker A_b$.
Its points have periods $a+u_p+ib$, and
Proposition~\ref{prop:AGY-L(q)} gives
\[
\|u_p\|_q\le4r_0=1/50.
\]
Also $\|a\|_q\le1$,
since the absolute real part of a saddle-connection period is at most
its length. It follows that
\[
\mathcal C_q(E)
=\{t(a+u_p):p\in E,\ 0<t\le1\}
\subset\frac{51}{50}\mathbb B_q.
\]
Consequently, \eqref{eq:AG-cone-factor} and
\eqref{eq:real-AGY-ball-volume} give
\[
\mu_q(B(q,1/200))
=w\lambda_q(\mathcal C_q(E))
\le w\left(\frac{51}{50}\right)^w C_{\rm norm}.
\]
The doubling estimate \cite[Corollary~5.4]{AG}, recalled in
Lemma~\ref{prop:Doubling}, gives
\[
\mu_q(B(q,1/100))
\le c_D^{-1}\mu_q(B(q,1/200)).
\]
Thus we may take
$\sV_0=c_D^{-1}w(51/50)^wC_{\rm norm}$, which depends only on $S$.

For the lower bound, fix $q_0\in K$ and choose a small period box
around it. The spaces $\ker A_b$ admit a continuous local choice of
basis. A fixed sufficiently small Euclidean box in these bases maps
into $B(q,1/100)$ for every center $q$ in a smaller neighborhood of
$q_0$, by continuity of the AGY norm and the affine chart estimates.
The density in \eqref{eq:AG-hor-volume} is positive and continuous,
so these boxes have a common positive lower bound for their measures
after shrinking the neighborhood. A finite cover of $K/\Mod$ by such
neighborhoods gives $\sV_1(S,K)>0$.
\end{proof}

For $x\in\PMF$ and $r>0$, write
\[
B_G(x,r) = \{y\in\PMF:d_G(x,y)<r\}.
\]
For $x\in W_\alpha$, write
\[
B_G^\alpha(x,r) = \{y\in W_\alpha:d_G^\alpha(x,y)<r\}.
\]

\begin{proposition}
\label{prop:volume}
Let $W=w/\beta$. For every
\[
0<r\le \frac{\rr_1}{4|\cA|C_\AGY}
\qquad\text{and}\qquad
x\in\PMF,
\]
we have
\[
\nu_\rho\big(B_G(x,r)\big)\prec r^W.
\]
Furthermore, suppose
\[
x\in V^{(1)}_\alpha,
\qquad
t>\sR_1,
\qquad
r=e^{-\beta t},
\]
and $x^\alpha(t)$ has $(K_1,\eta_1)$--good statistics. Then
\[
\nu_\rho\big(B_G(x,r)\big)\succ r^W.
\]
\end{proposition}

\begin{proof}
We first prove the upper bound. Choose a chart $\alpha$ such that
$x\in U_\alpha\subset V^{(1)}_\alpha$. Set
\[
\widehat r=2|\cA| \cdot r 
\qquad\text{which implies}\qquad 
\widehat r\le  C_\AGY^{-1}\rr_1/2.
\]

If $y\in B_G(x,r)$, then Lemma~\ref{lem:dGalpha:vs:dG} gives
\[
d_G^\alpha(x,y)
\le 2|\cA|\,d_G(x,y)
< \widehat r.
\]
Thus
\[
B_G(x,r)\subset B_G^\alpha(x,\widehat r).
\]
Note  that $2  C_\AGY \widehat{r} \le  \rr_1 < 1$. 
Choose $s\ge0$ so that
\[
e^{-\beta s}=2\, C_\AGY\widehat r.
\]
If $y\in B_G^\alpha(x,\widehat r)$, then
\[
d_G^\alpha(x,y)
<
\widehat r
=
\frac12 C_\AGY^{-1} e^{-\beta s}.
\]
The second alternative in \eqref{eq:assumption} therefore holds, and
Lemma~\ref{lem:dG:lower:bound} gives
\[
100\,C_\AGY^{-1} d_\AGY\big(x^\alpha(s),y^\alpha(s)\big)e^{-\beta s}
\le
d_G^\alpha(x,y)
<
\widehat r.
\]
Since $e^{-\beta s}=2\, C_\AGY \widehat r$, we obtain
\[
d_\AGY\big(x^\alpha(s),y^\alpha(s)\big)
< \tfrac 1{200}.
\]

Let $E=B_G^\alpha(x,\widehat r)$. Then
\[
E_\alpha(s)
\subset
B\big(x^\alpha(s),\ooth\big).
\]
Lemmas~\ref{lem:MV-Thurston} and~\ref{lem:AGY-volume}, together with
the volume growth formula \eqref{eq:Volume-Growth}, now give
\begin{align*}
\nu_\rho(B_G(x,r))
&\le\nu_\rho(E)\asymp\mu(E_\alpha(0))
=e^{-ws}\mu(E_\alpha(s))\\
&\prec e^{-ws}=(2C_\AGY\widehat r)^W\asymp r^W.
\end{align*}

We now prove the lower bound. Define
\[
\mathcal E_t=B\left(x^\alpha(t),\frac1{100\kk}\right),
\qquad
E=\{y\in W_\alpha:y^\alpha(t)\in\mathcal E_t\}.
\]
Every point of $\mathcal E_t$ has the form $y^\alpha(t)$ for some
$y\in\PMF(F_\alpha)$. Moreover, by \eqref{eq:Monotone},
\[
d_\AGY\big(x^\alpha(0),y^\alpha(0)\big)
\le
d_\AGY\big(x^\alpha(t),y^\alpha(t)\big)
<
\frac1{100\kk}
\le
\rr_{\rm chart}.
\]
Since $x\in V^{(1)}_\alpha$, Theorem~\ref{theorem:charts}{\rm (f)}
puts $y$ in $V^{(3)}_\alpha\subset W_\alpha$. Thus
$E_\alpha(t)=\mathcal E_t$. Theorem~\ref{theorem:metric-comparison}{\rm (a)} gives
\[
d_G(x,y)< e^{-\beta t}=r.
\]
That is,
\[
E\subset B_G(x,r).
\]

Since $x^\alpha(t)\in K_1$, applying Lemma~\ref{lem:AGY-volume} and
Lemma~\ref{prop:Doubling} at the radius
$1/(100\kk)$ gives
\[
\mu\big(E_\alpha(t)\big)\succ1.
\]
Using \eqref{eq:Volume-Growth} and Lemma~\ref{lem:MV-Thurston}, we obtain
\begin{align*}
\nu_\rho(B_G(x,r))
&\ge\nu_\rho(E)\asymp\mu(E_\alpha(0))
=e^{-wt}\mu(E_\alpha(t))\\
&\succ e^{-wt}=r^W. \qedhere
\end{align*}
\end{proof}

\begin{corollary} \label{cor:measure_decay}
For every $W'>W$, there is a constant $D_\rho>0$, depending only on
$\nu_\rho$, $W$, and $W'$, such that for every $x\in\PMF$ and $r>0$,
\begin{gather}\label{eq:measure_decay}
\int_{\PMF-B_G(x,r)}\frac{1}{d_G(y,x)^{W'}}\,d\nu_\rho(y)
\le \frac{D_\rho}{r^{W'-W}}.
\end{gather}
\end{corollary}

\begin{proof}
Since $\nu_\rho$ is a probability measure, enlarging the constant in
Proposition~\ref{prop:volume} extends its upper bound to all radii.
Thus there is $C_{\rm vol}>0$ such that
$\nu_\rho(B_G(x,s))\le C_{\rm vol}s^W$ for all $x\in\PMF$ and $s>0$.
Decompose the domain of integration into the shells
$2^jr\le d_G(x,y)<2^{j+1}r$, $j\ge0$. Then
\begin{align*}
\int_{\PMF-B_G(x,r)}\frac{d\nu_\rho(y)}{d_G(y,x)^{W'}}
&\le C_{\rm vol}\sum_{j\ge0}(2^jr)^{-W'}(2^{j+1}r)^W\\
&=\frac{2^W C_{\rm vol}}{1-2^{W-W'}}\,r^{W-W'}.
\end{align*}
This proves the claim with
$D_\rho=2^W C_{\rm vol}/(1-2^{W-W'})$.
\end{proof}

\section{Good geodesics, revisited}
\label{sec:good-geodesics-revisited}

We associate a Gromov ball to each mapping class in $\Gamma_\alpha(R)$.
At each sufficiently large scale, these balls have uniformly bounded
overlap. We also show that the set of points covered at no more than
$(1-\eta_1)n$ of the first $n$ scales has exponentially small measure.
These are the covering estimates used in Section~\ref{sec:first_moment}.
For $\gamma\in\Gamma_\alpha(R)$, write
\[
a_\gamma=[F_+(\gamma)],
\]
and set
\[
B_G(\gamma)=B_G\big(a_\gamma,e^{-\beta\ell(\gamma)}\big).
\]

Let $\sR_1$ be as in \autoref{theorem:metric-comparison}. 
\begin{proposition}
\label{prop:extra:properties:gamma:n}
Let  $\tau'$ be as in Proposition~\ref{prop:properties:gamma:n}  with
$\tau' >  \sR_1+1$  and
let $R_j = j \tau'$. 
We have:
\begin{enumerate}[(a)]
\item The set $\Gamma(R_j)$ is finite. 
  There is a constant $\sB=\sB(\eta_1)$ such that, for any
  sufficiently large  $j$, a point in
$\PMF$ belongs to at most $\sB$ of the balls
\[
B_G(\gamma),
\qquad
\gamma\in\Gamma(R_j).
\]

\item For each chart $\alpha$ and $x\in U_\alpha$, let
\[
k_n'(x)
=
\left|
\left\{
j\in[1,n]\,\Big|\,
x\in\bigcup_{\gamma\in\Gamma_\alpha(R_j)} B_G(\gamma)
\right\}
\right|.
\]
Then, for $n$ large enough, we have
\[
\nu_\rho
\Big(
\{x\in U_\alpha:k_n'(x)\le (1-\eta_1) n\}
\Big)
\le
\kappa^n, 
\]
\end{enumerate}
where $0 <  \kappa = \kappa(\eta_1) < 1$. 
\end{proposition}

\begin{proof}
By Proposition~\ref{prop:properties:gamma:n}{\rm (a),(b)},
every $\gamma\in\Gamma(R_j)$ satisfies
\[
d_\cT(P,\gamma(P))
\le 2\sD_0+\ell(\gamma)
\le 2\sD_0+R_j+\rho_1+1.
\]
Since closed Teichm\"uller balls are compact and the mapping class
group acts properly discontinuously, $\Gamma(R_j)$ is finite.

To bound the overlap, we work one chart at a time. Balls containing
a common point will force the orbit points $\gamma_i(P)$ into a
uniformly bounded \Teich neighborhood of one another. Proper discontinuity
will then bound their number.

Fix $x_0\in\PMF$ and $\alpha\in\cA$, and list all elements with
\[
x_0\in B_G(\gamma_i),
\qquad
\gamma_i\in\Gamma_\alpha(R_j),
\qquad
i=1,\ldots,m_\alpha.
\]
If the list is empty, its contribution to the overlap is zero.
Otherwise, fix one element and call it $\gamma$. By
Proposition~\ref{prop:properties:gamma:n}{\rm (b)}, for all $i$ we have
\[
|\ell(\gamma_i)-\ell(\gamma)|\le \rho_1+2.
\]
Hence
\[
d_G(a_{\gamma_i},a_\gamma)
\le
d_G(a_{\gamma_i},x_0)+d_G(x_0,a_\gamma)
\le
C_0e^{-\beta\ell(\gamma)}
\]
for a constant $C_0$ depending only on $\rho_1$.

By Proposition~\ref{prop:properties:gamma:n}{\rm (d)}, both
$a_\gamma$ and $a_{\gamma_i}$ lie in $V^{(1)}_\alpha  \subset V^{(2)}_\alpha$. For $R_j$ large enough, the
preceding inequality is smaller than $\rr_1$. Therefore
Lemma~\ref{lem:dGalpha:vs:dG} gives
\[
d_G^\alpha(a_{\gamma_i},a_\gamma)
\le
C_1e^{-\beta\ell(\gamma)}.
\]
Choose a constant $L_0$ large enough, depending only on $C_1$ and
$\beta$ so that
\[
C_1e^{-\beta\ell(\gamma)}
\le
\frac12 e^{-\beta(\ell(\gamma)-L_0)}.
\]
If $a_{\gamma_i}\ne a_\gamma$, apply Lemma~\ref{lem:dG:lower:bound} with
\[
s=\ell(\gamma)-L_0
\]
to obtain a uniform bound
\[
d_\AGY
\big(
a_{\gamma_i}^\alpha(s),
a_\gamma^\alpha(s)
\big)
\le C_2.
\]
If the centers coincide, the two chart representatives coincide,
so the same bound holds with zero on the left.
Using \eqref{eq:Monotone} over the bounded interval of length $L_0$, we obtain
\[
d_\AGY
\big(
a_{\gamma_i}^\alpha(\ell(\gamma)),
a_\gamma^\alpha(\ell(\gamma))
\big)
\le C_3.
\]
By Proposition~\ref{prop:properties:gamma:n}{\rm (d)} and
Definition~\ref{def:good-statistics},
$a_\gamma^\alpha(\ell(\gamma))\in K_1$. Apply a common mapping class
so that this point lies in a fixed compact set of marked representatives
of $K_1$. Since $D_\AGY\le d_\AGY$, the other differential lies in
the closed ambient $C_3$--neighborhood of this compact set. This
neighborhood is compact: the marked ambient AGY metric is a complete,
locally compact length metric, as discussed in Section~\ref{sec:AGY}.
Continuity of the \Teich distance therefore gives
\[
d_\cT
\big(
a_{\gamma_i}^\alpha(\ell(\gamma)),
a_\gamma^\alpha(\ell(\gamma))
\big)
\le C_4.
\]
Also, by Proposition~\ref{prop:properties:gamma:n}{\rm (b)},
\[
d_\cT
\big(
a_{\gamma_i}^\alpha(\ell(\gamma)),
a_{\gamma_i}^\alpha(\ell(\gamma_i))
\big)
\le \rho_1+2.
\]
By Proposition~\ref{prop:properties:gamma:n}{\rm (c')}, 
the points
$a_{\gamma_i}^\alpha(\ell(\gamma_i))$ and $\gamma_i(P)$ are at uniformly
bounded \Teich distance, and similarly
$a_\gamma^\alpha(\ell(\gamma))$ is uniformly close to $\gamma(P)$. Therefore
there is a constant $L_1=L_1(\eta_1)$ such that
\[
d_\cT\big(\gamma_i(P),\gamma(P)\big)\le L_1
\]
for every $i$. Equivalently,
\[
d_\cT\big(\gamma^{-1}\gamma_i(P),P\big)\le L_1.
\]
Since the action of $\Mod(S)$ on $\cT(S)$ is properly discontinuous, there are
only finitely many mapping classes $h$ with
\[
d_\cT(h(P),P)\le L_1.
\]
Thus $m_\alpha$
is bounded by a constant depending only on $L_1$, and hence only on
$\eta_1$. The  total multiplicity $\sB$ is bounded by $\sum_{\alpha \in
  \cA} m_\alpha$. This proves part {\rm (a)}.

We now prove part {\rm (b)} by showing that the approximations in
Proposition~\ref{prop:properties:gamma:n}{\rm (f)} imply membership
in the balls $B_G(\gamma)$. Let $E_\alpha(R_j)$ be the sets
defined there, and suppose
$q\in E_\alpha(R_j)$. Then there is $\gamma\in\Gamma_\alpha(R_j)$ such that
\[
d_\AGY
\big(
q(\ell(\gamma)),
a_\gamma^\alpha(\ell(\gamma))
\big)
<
\frac{1}{100\kk}.
\]
Let
\[
x=[F_+(q)].
\]
Since $q\in U_\alpha(0)$, we have $x\in U_\alpha\subset W_\alpha$ and
\[
q=x^\alpha(0),
\qquad
q(t)=x^\alpha(t).
\]
By Proposition~\ref{prop:properties:gamma:n}{\rm (d)}, $a_\gamma\in
V^{(1)}_\alpha\subset V^{(3)}_\alpha$ and
$a_\gamma^\alpha(\ell(\gamma))$ has $(K_1,\eta_1)$--good statistics.
Also, part {\rm (b)} of that proposition and
$\tau'>\sR_1+1$ give
\[
\ell(\gamma)\ge R_j-1=j\tau'-1>\sR_1.
\]
We may therefore apply
Theorem~\ref{theorem:metric-comparison}{\rm (a)} with base point $a_\gamma$, time
$\ell(\gamma)$, and point $x$, to obtain
\[
d_G(a_\gamma,x) < e^{-\beta\ell(\gamma)}.
\]
Thus $x\in B_G(\gamma)$: every index counted by $k_n(x^\alpha(0))$
is counted by $k_n'(x)$. At the excluded initial indices,
$\Gamma_\alpha(R_j)=E_\alpha(R_j)=\varnothing$, so
\[
k_n'(x)\ge k_n(x^\alpha(0))\qquad(x\in U_\alpha).
\]

Under $x\mapsto x^\alpha(0)$, Lemma~\ref{lem:MV-Thurston} compares
$\nu_\rho$ with the horospherical measure $\mu$ on $\Hor(F_\alpha)$.
Together with Proposition~\ref{prop:properties:gamma:n}{\rm (f)},
the counting inequality therefore gives, for all sufficiently large $n$,
\begin{align*}
&\nu_\rho\bigl(\{x\in U_\alpha:k_n'(x)\le(1-\eta_1)n\}\bigr)\\
&\qquad\le C_\alpha\,
\mu\bigl(\{q\in U_\alpha(0):k_n(q)\le(1-\eta_1)n\}\bigr)\\
&\qquad\le C_\alpha\kappa_0^n\mu(U_\alpha(0)),
\end{align*}
where $C_\alpha>0$ is independent of $n$ and
$\kappa_0=\kappa_0(\eta_1)<1$. Set $\kappa=(1+\kappa_0)/2<1$.
Since $\kappa>\kappa_0$, the last expression is at most $\kappa^n$
for all sufficiently large $n$. There are only finitely many charts,
so the threshold can be chosen uniformly in $\alpha$. It may depend
on $\eta_1$, $\tau'$, and $\rho$, whereas $\kappa$ depends only on
$\eta_1$. This proves part {\rm (b)}.
\end{proof}

\section{Lipschitz Estimates for spike functions}
\label{sec:Lipschitz}
Recall that $\Gamma(R)=\bigcup_\alpha\Gamma_\alpha(R)$, with
$\Gamma_\alpha(R)$ as in Proposition~\ref{prop:properties:gamma:n},
and that $H(\gamma,x)$ is defined by \eqref{eq:H}. We estimate its
size in terms of the distance from the attracting foliation of $\gamma$,
and its variation between points at distance at most
$e^{-\beta\ell(\gamma)}$. These estimates will verify the spike
conditions in Section~\ref{sec:spikes}. As before,
$w=\dim(\MF)=6g-6+2n$.

\begin{theorem}
\label{theorem:h:estimates}
Choose
\begin{equation}
\label{eq:spike-error-choice}
0<\eta_1<\min\left\{\frac1{100},\frac1{512w}\right\},
\qquad \eta=256w\eta_1.
\end{equation}
There are constants $\sC>0$ and $\sR_0>0$, depending on $\eta_1$,
such that the following holds. Let $\gamma\in\Gamma(R)$ for $R\ge\sR_0$
and let $a_\pm = [F_\pm(\gamma)]$ and $r = e^{-\beta \ell(\gamma)}$. Then 
\begin{itemize}
\item[{\rm (a)}] For any point $x \in \PMF$,
\begin{equation} \label{eq:H-bounds}
\frac {e^{w\ell(\gamma)}}{\sC}  \left(1+ \frac{d_G(a_+,x)}{r}\right)^{-\tfrac {2w+\eta}\beta} 
 \!\!  \leq H(\gamma,x)
  \leq \sC\, e^{w\ell(\gamma)}  \left(1+ \frac{d_G(a_+,x)}{r} \right)^{-\tfrac{2w-\eta}\beta}, 
\end{equation}
 In particular, 
if $x,y \in \PMF$ satisfy $d_G(x,y) \le r$, then
\begin{equation}
\label{eq:H:coarse:Lipschitz}
\frac 1{\sC} \left(1+ \frac{d_G(a_+,x)}{r} \right)^{- \frac{2\eta}{\beta}} 
  \le \frac{H(\gamma,y)}{H(\gamma,x)} 
  \le \sC \left(1+ \frac{d_G(a_+,x)}{r} \right)^{\frac{2\eta}{\beta}}, 
\end{equation}

\item[{\rm (b)}] For any $x,y \in \PMF$ with $d_G(x,y) \le r$,
\begin{equation}
\label{eq:H:Lipschitz}
\Big| H(\gamma,x)-H(\gamma,y)\Big| 
  \le \frac{\sC}{r} H(\gamma,x) \left(1+\frac{d_G(a_+,x)}{r} \right)^{\frac{\eta}{\beta}} \, d_G(x,y).
\end{equation}
\end{itemize}
\end{theorem}

\subsection*{Thin triangles}
\label{subsec:thin-triangles}

Fix $\gamma\in\Gamma(R)$, and write
\[
a_+=[F_+(\gamma)],\qquad a_-=[F_-(\gamma)].
\]
For $x\in\PMF$, we compare two rays toward $x$, starting near $P$
and $\gamma(P)$. We will choose times $\tau,\tau'$ whose sum is
close to $\ell(\gamma)$ and whose corresponding points can be joined
by a short geodesic segment. The time $\tau$ also controls
$d_G(a_+,x)$. These comparisons will relate the change in the size
of a foliation to its distance from $a_+$.

We use the following form of the thin-triangle theorem
\cite[Theorem~A]{masur:SH}: a long interval that spends at least half
its time in a fixed thick part comes within bounded distance of one
of the other two sides.

\begin{theorem}[Dowdall--Duchin--Masur]
\label{theorem:DDM}
For any $\Mod$--compact set $\mathcal K\subset\cT$ there are constants $D_{\DDM}$ and $L_{\DDM}$ such that 
the following holds. Let $[x,y]$ be a geodesic segment in $\cT$ and let $I$ be a geodesic subinterval 
of $[x,y]$, with
\[
|I|\ge L_{\DDM},
\qquad
\big|\{u\in I:u\in\mathcal K\}\big|\ge\frac12|I|.
\]
Then for all $z \in \cT$, we have 
\[
I \cap N_{D_{\DDM}} ([x, z] \cup [y, z]) \not = \emptyset.
\]
\end{theorem} 

Below we choose $D_{\DDM}$ and $L_{\DDM}$ for the projection of $K_2$
to $\cT$.

\begin{figure}[ht]
\setlength{\unitlength}{0.01\linewidth}
\begin{picture}(100, 26) 

\put(9,0){
\begin{tikzpicture} [>=stealth]

\tikzstyle{vertex} =[circle,draw,fill=black,thick, inner sep=0pt,minimum size= .7 mm]


\draw[thick] (-5.5,0) -- (5.5,0);
\draw[thick] (1,2) -- (1,2.7);

\draw[thick] (-1,0) arc (-90: 0: 2);
\draw[thick] (3, 0)  arc(-90: -180: 2);

\draw[thick] (-5,-.5) arc (150: 90: 1);
\draw[thick] (5,-.5) arc (30: 90: 1);

\node (a-) at (-5.5,0)  [label=left:$a_-$]{};
\node (a+) at (5.5,0)  [label=right:$a_+$]{};
\node (x) at (1,2.7)  [label=above:$x$]{};

\node [vertex] (xa) at (-5,-.5)  [label=below:$x^{\alpha_+}(0)$]{};
\draw [thick, ->] (xa) -- (-4.8, -.1) ;

\node [vertex] (xa') at (5,-.5)  [label=below:$x_{F'}(0)$]{};
\draw [thick, ->] (xa') -- (4.8, -.1) ;

\node [vertex] (xt) at (-.4,.1)  [label=above left:$x^{\alpha_+}(\tau)$]{};
\draw [thick, ->] (xt) -- (0.1, .2) ;

\node [vertex] (xt') at (2.4,.1)  [label=above right:$x_{F'}(\tau')$]{};
\draw [thick, ->] (xt') -- (1.9, .2) ;

\node [vertex] (xT) at (1, 2.1)  [label=left:$x^{\alpha_+}(T)$]{};
\draw [thick, ->] (xT) -- (1, 2.5) ;

\draw[thick] [decorate,decoration={brace,amplitude=10pt,mirror}] 
        (-.4,-.1) -- (2.4, -.1) node[midway,yshift=-20pt]{$\leq 16 \eta_1 \tau' + \sL_2$};

\end{tikzpicture}}
\end{picture}
\caption{The two rays toward $x$ start near $P$ and $\gamma(P)$.
The selected points at times $\tau$ and $\tau'$ are joined by a
geodesic segment of length at most $16\eta_1\tau'+\sL_2$.
Part~(d) of Proposition~\ref{prop:thin} compares rays from this
segment with $x^{\alpha_+}(T)$.}
\label{fig:Thin} 
\end{figure}
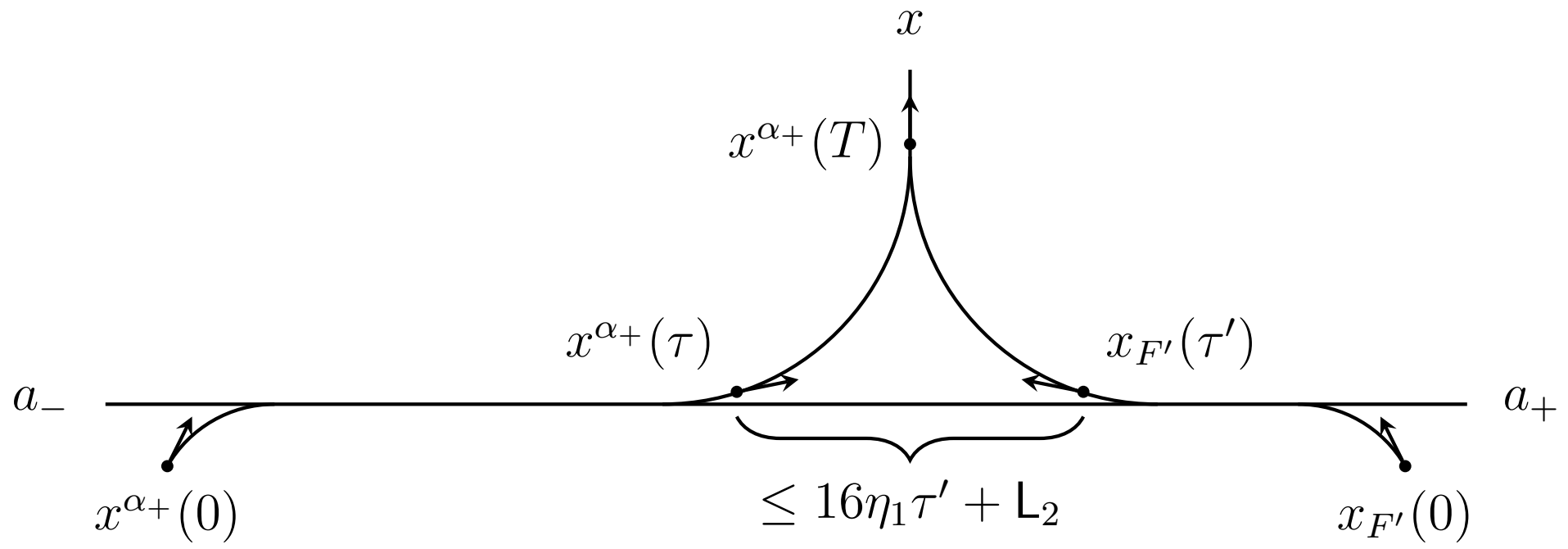

For $x\in\PMF$, set
\[
\bar x=\gamma^{-1}x.
\]
If $\alpha_-\in\cA$ is a chart with $\bar x,a_-\in W_{\alpha_-}$, we set
\[
F'=\gamma(F_{\alpha_-}),
\]
so that
\[
x_{F'}(t)=\gamma\big(\bar x^{\alpha_-}(t)\big).
\]
Proposition~\ref{prop:thin} gives the times shown schematically in
Figure~\ref{fig:Thin}. We call the geodesic segment joining the selected
points the \emph{bridge}. Its length and the deficit
$\ell(\gamma)-\tau-\tau'$ are bounded by fixed numerical multiples
of $\eta_1\tau'$, up to additive errors depending on the compact sets.

These times explain the decay in \eqref{eq:H-bounds}. Ignoring the
controlled errors and fixed multiplicative constants, the change of
normalization gives $H(\gamma,x)^{1/w}\approx e^{\tau-\tau'}$, while
$1+d_G(a_+,x)/r\approx e^{\beta\tau'}$, where
$r=e^{-\beta\ell(\gamma)}$. Together with
$\tau+\tau'\approx\ell(\gamma)$, these relations give the decay
exponent $-2w/\beta$. The estimates below quantify the errors in
this picture.

\begin{proposition}
\label{prop:thin}
Assume $0<\eta_1\le1/100$, and recall that $\eta_2=1/4$.
There are constants
\[
\sC_2\ge16,\qquad \sL_2,\sT_2,\sR_2>0
\]
depending only on $K_1,K_2,\eta_1$ and the fixed chart system, such that if
$\gamma\in\Gamma(R)$ and $\ell(\gamma)\ge \sR_2$, then the following holds.

For every $x\in\PMF$, there exist charts $\alpha_+,\alpha_-\in\cA$, with
\[
x,a_+\in V^{(2)}_{\alpha_+},
\qquad
\bar x,a_-\in V^{(2)}_{\alpha_-},
\]
and times
\[
0\le \tau,\tau'\le \ell(\gamma),
\]
such that, setting $F'=\gamma(F_{\alpha_-})$, we have:

\begin{enumerate}[(a)]
\item $x^{\alpha_+}(\tau)\in K_2$ and
\begin{equation}
\label{eq:tau}
d_G(a_+,x)\le \sC_2 e^{-\beta\tau}.
\end{equation}

\item Either $\tau'\le \sL_2$, or
\begin{equation}
\label{eq:tau2}
\frac1{\sC_2}
e^{-\beta(\tau+8\eta_1\tau')}
\le d_G(a_+,x).
\end{equation}

\item $x_{F'}(\tau')\in K_2$, and
\begin{equation}
\label{eq:tau'}
d_\cT\big(x^{\alpha_+}(\tau),x_{F'}(\tau')\big)
\le 16\eta_1\tau'+\sL_2.
\end{equation}
Moreover,
\begin{equation}
\label{eq:ell}
\ell(\gamma)-8\eta_1\tau'-\sL_2
\le \tau+\tau'
\le \ell(\gamma).
\end{equation}

\item Assume that there is a time
\[
T\ge \ell(\gamma)+\sT_2
\]
such that $x^{\alpha_+}(T)$ has $(K_2,\eta_2)$--good statistics. Let $\zeta$
be a \Teich geodesic segment joining $x^{\alpha_+}(\tau)$ to
$x_{F'}(\tau')$. For every point $Z$ in the image of $\zeta$, let
\[
q_Z=Q_1(x,Z).
\]
Then $q_Z\in\POT$, and there exists a time $s_Z>0$ such that
\begin{equation}
\label{eq:s_Z}
d_\AGY\big(x^{\alpha_+}(T),q_Z(s_Z)\big)\le \tfrac1{100}.
\end{equation}
\end{enumerate}
\end{proposition}

The coefficients $8$ and $16$ are independent of the compact sets
and of $\eta_1$. This will allow us to choose $\eta_1$ before the
compact sets. We first fix the auxiliary constants and prove a
comparison for ideal triangles; we then construct the two times in
three cases.

Fix $\alpha$ such that $\gamma\in\Gamma_\alpha(R)$.
Proposition~\ref{prop:properties:gamma:n}{\rm(d)} gives
$a_+\in V^{(1)}_\alpha$. The ideal endpoints $F_\alpha$, $F'$, and
$a_\pm$ used below are generic and uniquely ergodic: they are the
invariant foliations of pseudo-Anosov axes in the principal stratum. Write
$L=\ell(\gamma)$ and $q(u)=a_+^\alpha(u)$, and put
\[
G=\{u\in[0,L]:q(u)\text{ has }(K_1,\eta_1)
\text{--good statistics}\}.
\]
Proposition~\ref{prop:properties:gamma:n}{\rm(e)} says precisely that
\begin{equation}
\label{eq:bridge-good-time-budget}
|[s,L]\setminus G|\le\eta_1(L-s)\qquad(0\le s<L).
\end{equation}
All constants in the proof may depend on $K_1,K_2,\eta_1$
and the fixed chart system. They enter only additive errors or
multiplicative bounds; the numerical coefficients displayed below
do not depend on these data.

Choose $\delta>0$ such that
\[
8\delta<\rr_{\rm chart},\qquad
3\delta<\frac1{100},\qquad
\delta<\frac1{100\kk},\qquad 2\delta<\frac1{50}.
\]
We may decrease $\delta$ further to work inside fixed local neighborhoods
of $K_1$ and $K_2$. In particular, an AGY comparison of size at most
$3\delta$ with a point of $K_1$ implies membership in $K_2$, by
Definition~\ref{def:K_2}. We use this also after reversing both
differentials: $K_1$ and $K_2$ are invariant under $q\mapsto-q$.

Choose $D_{\rm aux}$ to bound the following base distances:
between any two chart origins $z^\alpha(0)$, for $z\in W_\alpha$;
between a point within ambient AGY distance $\delta$ of $\pm q(L)$ and any
$\gamma(z^{\alpha_-}(0))$, for $z\in W_{\alpha_-}$;
and between $q(0)$ and an axis point within ambient AGY distance
$1/80$ of it. These bounds are uniform.
The chart slices have compact closures, and
Proposition~\ref{prop:properties:gamma:n}{\rm(c')} puts $q(L)$
at bounded base distance from $\gamma(P)$. The fixed AGY
neighborhoods also have compact closures modulo $\Mod(S)$.
If two points have the same projective vertical foliation,
extremal-length distortion bounds the time shift needed to match
their measured foliations by their base distance. Thus
$D_{\rm aux}$ also bounds each such shift, and $2D_{\rm aux}$
bounds the base distance after matching. We choose this constant
before the occupation thresholds and buffer lengths.

Let $h_c>D_{\rm aux}$ be large enough that the backward contraction from a
separation time, over an interval of length $h_c$, reduces distance
$1/100$ to $\delta$. This follows from Corollary~\ref{cor:AGY:half},
since good statistics with parameter $\eta_2<1/2$ supply the required
half occupation on every backward interval.

We will sometimes replace a time by zero. Fix a threshold $R_*$ larger
than $\sR_1$, $D_{\rm aux}$, and
\[
D_{\rm aux}+B(K_2,\eta_2,\delta),
\]
where $B$ is supplied by Theorem~\ref{thm:rays-linear}. A change from a
time at most $R_*$ to zero costs only a fixed additive amount.

The next lemma turns DDM's bounded base-distance comparison into
AGY closeness at one end of an interval. The intervals on either
side of the DDM interval provide the compact time needed for
Proposition~\ref{cor:rays-converge-new}.

\begin{lemma}
\label{lem:thin:DDM-side-switching}
Let $A,B$ be distinct generic uniquely ergodic projective foliations, and
let $C$ be a projective foliation. When $C$ differs from $A$ and $B$,
suppose the pairs defining the three sides fill. Orient the reference side $\omega$ from
$A$ to $B$. Set $l=\max(L_{\DDM},1)$. There is an occupation threshold
$R_c$ with the following property.

Suppose
\[
I_-=[u_0,u_1],\qquad I=[u_1,u_2],\qquad I_+=[u_2,u_3]
\]
are consecutive intervals, $|I|=l$, at least half of $I$ lies in $K_2$,
and each of $I_-,I_+$ contains at least $R_c$ time in $K_2$. Then one of
the following holds:
\begin{enumerate}[(i)]
\item There is a differential $p_A$ on the side from $A$ to $C$, oriented
toward $C$, such that
\[
F_-(p_A)=F_-(\omega(u_0)),\qquad
d_\AGY(\omega(u_0),p_A)\le\delta.
\]
\item There is a differential $p_B$ on the side from $C$ to $B$, oriented
toward $B$, such that
\[
F_+(p_B)=F_+(\omega(u_3)),\qquad
d_\AGY(\omega(u_3),p_B)\le\delta.
\]
\end{enumerate}
If $C=A$ or $C=B$, the corresponding comparison with the identical side
holds directly.
\end{lemma}

\begin{proof}
We first show that one of the other two ideal sides has a point whose
base is within $D_{\DDM}$ of $I$. We use the finite-triangle
conclusion of Theorem~\ref{theorem:DDM}. For punctured surfaces, this
follows by lifting the triangle to the closed-surface cover described in
\hyperlink{rays-punctures}{Section~\ref*{sec:rays-new},
\textit{Surfaces with punctures}}; the lift is isometric and
the lifted set has compact image in moduli space. For the four-punctured
sphere it is the thin-triangle property of its Teichm\"uller disk.
The cases $C=A$ and $C=B$ use the identical side, so assume $C$ differs
from both. Suppose first that $C$ is uniquely ergodic.
Approximate $A$ and $B$ by $\omega(-N)$ and
$\omega(N)$, and approximate $C$ by time $N$ on a ray toward $C$. The
fixed interval $I$ lies on the first finite side. DDM gives a point
on one of the other two finite sides within distance $D_{\DDM}$ of
$I$. These points remain in a compact base set. Pass to a subsequence
on which the side is fixed and the oriented unit-area differentials
at these points converge.

For example, suppose the nearby point is on the side from $\omega(-N)$
to the vertex tending to $C$. Let $p_N$ be its differential, oriented
toward $C$.
Its distances to both vertices tend to infinity. The extremal lengths at
the first vertex of a fixed representative of $A$ and of $F_-(p_N)$
both tend to zero. The intersection inequality therefore implies that
the limiting horizontal foliation has zero intersection with $A$.
Similarly, the limiting vertical foliation has zero intersection with
$C$. Unique ergodicity identifies these projective foliations with $A$
and $C$. Thus the limiting differential lies on the actual side from
$A$ to $C$. The other case gives a point on the side from $B$ to $C$.

For general $C$, approximate it by uniquely ergodic foliations $C_j$;
the pairs with $A$ and $B$ still fill for all sufficiently large $j$.
The nearby points again lie over a compact base set. After passing
to a subsequence, their measured foliations are positive multiples
of $A,C_j$, or of $B,C_j$. Unit-area normalization and continuity of extremal length
bound these multiples above and below. The continuity of $Q(F,H)$
therefore identifies their limit with a point of the required side.
The limiting differential is principal because its other foliation is
generic. This proves the ideal-side comparison.

Choose a point $v\in I$ with $\omega(v)\in K_2$. Its base is within
$D_{\DDM}+l$ of the point just found on the other side. Change
the time on that side so that the foliation associated with
the common ideal endpoint agrees as a measured foliation. The required
shift has absolute value at most $D_{\DDM}+l$, by extremal-length
distortion. The new base-distance bound is therefore at most
$2(D_{\DDM}+l)$.

Choose $R_c$ using Proposition~\ref{cor:rays-converge-new} for this bound,
the compact set $K_2$, and tolerance $\delta$. For a point on the
$A$-side, reverse both differentials and apply that proposition while
moving back from $v$ to $u_0$. For a point on the $B$-side, orient its
differential toward $B$ and apply the proposition from $v$ to $u_3$.
The respective buffers supply the required compact occupation. These
are exactly the two stated comparisons, with the indicated measured
foliations and orientations.
\end{proof}

\begin{proof}[Proof of Proposition~\ref{prop:thin}]
We first construct the charts and times for (a)--(c), using the
constants fixed above. Let $R_t$ be the occupation threshold in
Proposition~\ref{cor:rays-converge-new} for $K=K_2$ and initial base distance $2D_{\rm aux}$
and target $\delta$. Fix
\[
h\ge\max(2R_c,2R_t,8l,2D_{\rm aux},1),\qquad M=8h+4.
\]
These constants are chosen before dividing the proof into cases.

Before treating the cases, we record how to obtain the chart inclusion
in Proposition~\ref{prop:thin}. A preliminary positive time $v>R_*$ below
will satisfy
\[
q(v)\text{ has }(K_1,\eta_1)\text{--good statistics},\qquad
d_\AGY(q(v),x^\alpha(v))\le\delta.
\]
Backward nonexpansion transfers every negative-time suffix, so
$x^\alpha(v)$ has $(K_2,\eta_2)$--good statistics. Choose $\alpha_+$
with $x,a_+\in U_{\alpha_+}$. Applying
Theorem~\ref{thm:rays-linear} from the two chart origins gives a shift
$s_x$, with $|s_x|\le D_{\rm aux}$, for which
\begin{equation}
\label{eq:bridge-final-chart}
d_\AGY(x^\alpha(v),x^{\alpha_+}(v+s_x))\le\delta.
\end{equation}
Set $\tau=v+s_x$. Its endpoint belongs to $K_2$ by comparison with
$q(v)\in K_1$. We reserve a buffer $D_{\rm aux}$ below so that
$\tau+\tau'\le L$. If the preliminary choice is at most $R_*$, use
$\tau=0$ in the new chart instead; its endpoint is in $K_2$ by
Proposition~\ref{prop:D}. The bounded chart change affects only additive
errors and the multiplicative constants in (a) and (b).

For $x\in V^{(3)}_\alpha$, set $\sigma=\tau_\alpha(a_+,x)$.
We distinguish whether this time is near $L$, earlier than $L-M$,
or whether $x$ lies outside $V^{(3)}_\alpha$.

\smallskip
\noindent\emph{Case 1: $x\in V^{(3)}_\alpha$ and $\sigma\ge L-M$.}
Take $t_* =\min(\sigma,L)$, with $t_*=L$ if $\sigma=+\infty$.
In the latter case $x=a_+$ by Proposition~\ref{prop:D}, so the two rays
coincide.
At an attained finite separation time the two rays are within $1/100$;
backward nonexpansion gives the same bound at $t_*$ if $\sigma>L$.
When $t_*=L$, this comparison and the good statistics of $q(L)$ also
give good statistics for the companion ray at $L$. Thus backward
contraction gives distance at most $\delta$ at $t_*-h_c$ and all
earlier times.

By \eqref{eq:bridge-good-time-budget}, an interval ending at
$t_*-h_c$ and having length greater than
$\eta_1(L-t_*+h_c)/(1-\eta_1)$ cannot consist entirely of bad times.
Since $L-t_*\le M$, choose a preliminary good time $v$ within a fixed
distance before $t_*-h_c$. For sufficiently large $L$, we have $v>R_*$.
Make the final chart change \eqref{eq:bridge-final-chart}. The choice
$h_c>D_{\rm aux}$ gives $0\le\tau\le L$, and $L-\tau$ is bounded.

Choose $\alpha_-$ with $\bar x,a_-\in U_{\alpha_-}$ and put $\tau'=0$.
The point $x_{F'}(0)$ is in $K_2$, and its base point is a bounded
distance from $\gamma(P)$. The same holds for $q(L)$. Consequently
the distance between the two selected endpoints is bounded.
The upper metric comparison at $v>\sR_1$
gives (a), after absorbing the bounded chart shift into its prefactor.
Part (b) holds by its first alternative. The bounded distance and
deficit give (c), once $\sL_2$ is chosen.

\smallskip
\noindent\emph{Case 2: $x\in V^{(3)}_\alpha$ and $\sigma<L-M$.}
Here $x\ne a_+$, since time $L$ is admissible when the two rays
coincide. By Proposition~\ref{prop:D}, $\sigma\ge0$ is finite and
attained. Put $H=L-\sigma$ and $m=\eta_1H$.
We first find an interval after $\sigma$ to which
Lemma~\ref{lem:thin:DDM-side-switching} applies, leaving enough
compact time afterward to compare with the ray from $\gamma(P)$.

The bad-time budget on $[\sigma,L]$ is at most $m$, so choose
\[
t_0\in G\cap(\sigma,\sigma+m+1),\qquad
\bar t\in G\cap(t_0+4h,t_0+4h+m+1).
\]
These intervals lie before $L$, and
\begin{equation}
\label{eq:bridge-main-forward}
\bar t-\sigma\le2\eta_1H+4h+2.
\end{equation}
In particular,
\[
L-\bar t\ge(1-2\eta_1)H-4h-2\ge3h.
\]
Since $\bar t$ is good, at most $2\eta_1h$ of
$[\bar t-2h,\bar t-h]$ lies outside $K_1$. Partition it into full
$l$-blocks, discarding a remainder shorter than $l$. At least one block
$I$ has half its time in $K_1$. The interval from $t_0$ to the beginning
of $I$ contains $[\bar t-3h,\bar t-2h]$, with at least
$(1-3\eta_1)h\ge R_c$ compact time. The interval from the end of $I$
to $\bar t$ contains $[\bar t-h,\bar t]$, with at least
$(1-\eta_1)h\ge R_c$ compact time.

Apply Lemma~\ref{lem:thin:DDM-side-switching} to $F_\alpha,a_+,x$.
Its first alternative gives the comparison with $x^\alpha(t_0)$ at
exactly time $t_0$, because the measured horizontal foliations agree.
Backward nonexpansion would transfer every negative-time suffix of the
good statistics at $t_0$ to $x^\alpha(t_0)$. This contradicts
$t_0>\sigma$. Hence the second alternative gives a differential $p$
oriented toward $a_+$, with
\[
F_+(p)=F_+(q(\bar t)),\qquad d_\AGY(p,q(\bar t))\le\delta.
\]

Put $b=L-\bar t$ and define the $x$-ray $r(u)=-g_{b-u}p$.
Forward stable nonexpansion gives
\[
D_\AGY(r(u),-q(L-u))\le\delta\qquad(0\le u\le b).
\]
Thus $r(0)\in K_2$, its base point is a bounded distance from
$\gamma(P)$, and it spends at least $(1-\eta_1)b$ time in $K_2$.
Choose $\alpha_-$ with $\bar x,a_-\in U_{\alpha_-}$. Match the measured
vertical foliation of $r(0)$ to that of $x_{F'}(s_0)$, with
$|s_0|\le D_{\rm aux}$. Since $b\ge3h$, we have $b>D_{\rm aux}$ and
$(1-\eta_1)b\ge R_t$. Proposition~\ref{cor:rays-converge-new}
therefore gives
\[
\tau'=L-\bar t+s_0>0,\qquad
D_\AGY(x_{F'}(\tau'),-q(\bar t))\le2\delta.
\]
In particular $x_{F'}(\tau')\in K_2$. This terminal comparison has
only a bounded time error.

We next choose a preliminary time $v\le\sigma$. Backward contraction
from $\sigma$ gives distance at most $\delta$ by time $\sigma-h_c$.
If the interval
\[
J=\big(\sigma-h_c-[2\eta_1(H+h_c)+1],\ \sigma-h_c\big)
\]
lies in $[0,L]$, it contains a point of $G$. Indeed, otherwise its
length $\Delta$ would satisfy
$\Delta\le\eta_1(H+h_c+\Delta)$ by
\eqref{eq:bridge-good-time-budget}, which contradicts the displayed
choice of $\Delta$. Choose such a point if it exceeds $R_*$; otherwise
use $v=0$. Also use $v=0$ if the lower endpoint of $J$ is nonpositive.
This includes the case that $J$ lies entirely before zero. In all cases,
\[
\sigma-v\le2\eta_1H+B_2,\qquad
\bar t-v\le4\eta_1H+B_3.
\]
Make the final chart change described in
\eqref{eq:bridge-final-chart} if $v>0$, and use $\tau=0$ in that chart
otherwise. Since $\bar t-\sigma\ge4h\ge2D_{\rm aux}$, the two bounded time
shifts leave $\tau+\tau'\le L$.

The endpoint comparisons and the segment of $q$ between $v$ and
$\bar t$ give bridge length at most $4\eta_1H+B_4$; the same bound
holds for $L-\tau-\tau'$. Moreover
\[
H=\tau'+\bar t-\sigma-s_0\le\tau'+2\eta_1H+B_5.
\]
Hence both quantities are bounded by
\[
\frac{4\eta_1}{1-2\eta_1}\tau'+B_6
\le5\eta_1\tau'+B_6.
\]
Thus the distance and deficit satisfy (c), with the additive
error to be included in $\sL_2$.

It remains to relate these choices to $d_G(a_+,x)$.
If $v>0$, the upper metric comparison applies at $v>\sR_1$;
if $v=0$, bounded diameter applies. The final chart shift is bounded,
so we obtain the upper bound in (a) with a fixed prefactor.
At $\bar t>\sigma$ the
distance between $q(\bar t)$ and $x^\alpha(\bar t)$ exceeds $1/100$:
otherwise its good statistics would transfer backwards and contradict
the definition of $\sigma$. At $v$ the distance is smaller than
$1/100$, including when $v=0$ by the chart diameter. Thus there is
$s\in[v,\bar t]$ where that distance equals $1/100$.
Lemma~\ref{lem:dG:lower:bound} gives
\[
d_G^\alpha(a_+,x)\ge C_\AGY^{-1}e^{-\beta s}.
\]
Since $s\le\bar t\le\tau+8\eta_1\tau'+B_7$,
Lemma~\ref{lem:dGalpha:vs:dG} gives, when $d_G(a_+,x)\le\rr_1$,
\[
d_G(a_+,x)\ge
\frac{e^{-\beta(\tau+8\eta_1\tau')}}
     {2|\cA|C_\AGY e^{\beta B_7}}.
\]
When $d_G(a_+,x)>\rr_1$, use $\rr_1$ as the lower bound in (b).

\smallskip
\noindent\emph{Case 3: $x\notin V^{(3)}_\alpha$.}
Since $a_+\in V^{(1)}_\alpha$, Lemma~\ref{lem:dGalpha:vs:dG} and
\eqref{eq:r1:consequence} would give $x\in V^{(3)}_\alpha$ if
$d_G(a_+,x)\le\rr_1$. Thus $d_G(a_+,x)>\rr_1$.
Choose $\alpha_+,\alpha_-$ with
$x,a_+\in U_{\alpha_+}$ and $\bar x,a_-\in U_{\alpha_-}$.
Set $\tau=0$, so $x^{\alpha_+}(\tau)\in K_2$, and put
\[
\omega(t)=\gamma(a_-^{\alpha_-}(t)),\qquad F'=\gamma F_{\alpha_-}.
\]

Choose a point $z_0$ of the oriented axis with
$D_\AGY(z_0,q(0))<1/80$. Such a point exists by
Proposition~\ref{prop:properties:gamma:n}{\rm(c)}; no attainment of
the distance to the axis is needed. Choose $|h_0|\le D_{\rm aux}$ so that
$a(u)=g_{u+h_0}z_0$ and $q(u)$ have the same measured vertical
foliation. Put $m=\eta_1L$ and $U=m+R_{\rm ax}$, where
$R_{\rm ax}\ge\max(R_t,2D_{\rm aux}+2)$ is fixed. The interval $[0,U]$ has at least
$R_{\rm ax}$ compact time, by
\eqref{eq:bridge-good-time-budget} with $s=0$. Proposition
\ref{cor:rays-converge-new}, followed by stable nonexpansion, gives
\[
d_\AGY(q(u),a(u))\le\delta\qquad(u\ge U).
\]
Let $A(s)=-a(L-s)$. Match $\omega(c+s)$ to $A(s)$ with equal measured
vertical foliations and $|c|\le D_{\rm aux}$. Their starting base distance is
bounded by $2D_{\rm aux}$, since both starting points lie near $\gamma(P)$.
Set $S_0=2R_t+1$. For sufficiently large $L$, we have $L-S_0\ge U$.
The good statistics of $q(L)$ and the preceding comparison give at
least $(1-\eta_1)S_0>R_t$ time in $K_2$ for $A$ on $[0,S_0]$.
Also $A(0)\in K_2$. A second application of Proposition
\ref{cor:rays-converge-new} therefore gives
\[
d_\AGY(\omega(c+s),A(s))\le\delta\qquad(s\ge S_0).
\]
Consequently,
\begin{equation}
\label{eq:bridge-outside-axis}
D_\AGY(\omega(L-u+c),-q(u))\le2\delta
\quad(U\le u\le L-S_0).
\end{equation}

In this range let $E=\{t:L-t+c\in G\}$; its complement has measure
at most $m$. Every point of $E$ is a $K_2$ point of $\omega$.
Set $A_0=R_c+1$, $\Delta_0=4m+2A_0+2l$, and $b=L-U+c$. Choose
\[
t_0\in E\cap[b-\Delta_0-2m-2,b-\Delta_0-m-1],\qquad
\bar t\in E\cap[t_0+\Delta_0,t_0+\Delta_0+m+1].
\]
Each choice interval has length $m+1$. All these times lie in the range
of \eqref{eq:bridge-outside-axis} for sufficiently large $L$, since
their corresponding $u$-parameters lie between $U$ and
$U+6m+2A_0+2l+2$, and $\eta_1\le1/100$.

The first and last intervals of length $m+A_0$ in $[t_0,\bar t]$
each have at least $R_c$ time in $E$. The middle interval has length
at least $2m+2l$; its full $l$-blocks have total length greater than
$2m$, so one has half its time in $E$. Apply
Lemma~\ref{lem:thin:DDM-side-switching} to this block and the triangle
$F',a_-,x$.

Its second alternative is impossible. Indeed, it gives a point $p$
oriented toward $a_-$, with the same measured vertical foliation as
$\omega(\bar t)$ and within distance $\delta$. Put
$t_*=L+h_0+c$. Then $A(t_*-c)=-z_0$, and $t_*\ge\bar t$ because
$R_{\rm ax}>|h_0|+1$. Flow $p$ and $\omega(\bar t)$ forward to $t_*$,
and put $v=-g_{t_*-\bar t}p$. Then $v$ has vertical class $x$ and
$D_\AGY(v,z_0)\le2\delta$. The projection
$\pi^{\rm cs}_{q(0)}$ is defined on $B_\AGY(z_0,1/50)$, sends $z_0$
to $q(0)$, and is $4$--Lipschitz. Thus $q_x=\pi^{\rm cs}_{q(0)}(v)$ satisfies
\[
[F_+(q_x)]=x,\qquad F_-(q_x)=F_\alpha,\qquad
d_\AGY(q_x,q(0))\le8\delta<\rr_{\rm chart}.
\]
Thus $q_x=x^\alpha(0)$, and Theorem~\ref{theorem:charts}{\rm(f)}
would give $x\in V^{(3)}_\alpha$, a contradiction.

The first alternative consequently holds. Because the measured
horizontal foliations agree, it compares $\omega(t_0)$ to
$x_{F'}(t_0)$ at the exact same time. Set $\tau'=t_0$. This point is
within $3\delta$ of $-q(u_0)$, where $u_0=L-t_0+c\in G$, so it
belongs to $K_2$. Also,
\[
u_0\le7\eta_1L+B_9,
\qquad
d_\cT(x^{\alpha_+}(0),x_{F'}(\tau'))\le u_0+B_{10},
\qquad L-\tau'\le7\eta_1L+B_{11}.
\]
Here the second estimate uses the bounded distance between the two
chart origins and the segment of $q$ of length $u_0$.
The choice $U>|c|$ ensures $\tau'\le L$; positivity holds when $L$
is sufficiently large. Eliminating $L$ by
$\tau'\ge(1-7\eta_1)L-B_{11}$ gives
\[
d_\cT(x^{\alpha_+}(0),x_{F'}(\tau')),\ L-\tau'
\le\frac{7\eta_1}{1-7\eta_1}\tau'+B_{12}
\le8\eta_1\tau'+B_{12}.
\]
Finally, bounded diameter gives the upper bound in (a) with a
fixed prefactor, and $d_G(a_+,x)>\rr_1$ gives the required lower
bound for (b).

\smallskip
We now choose the constants in (a)--(c). Take $\sL_2>0$ to
dominate the bounded distance and deficit in Case 1 and the additive
errors $B_6,B_{12}$ in Cases 2 and 3. Choose $\sC_2\ge16$ to
dominate all upper metric prefactors and
\[
2|\cA|C_\AGY e^{\beta B_7},\qquad \rr_1^{-1}.
\]
These choices give the lower metric bound in (b), with coefficient
$8$, and the distance and deficit bounds in (c), with coefficients
$16$ and $8$. Choose $\sR_2$ larger than all the length thresholds
used in the three cases. Every choice depends only on the data in
the proposition. The charts chosen above have both required points
in $U_{\alpha_\pm}\subset V^{(2)}_{\alpha_\pm}$, completing (a)--(c).

For part (d), set $B_{\rm br}=B(K_2,\eta_2,1/100)$ as in
Theorem~\ref{thm:rays-linear}, and choose
\begin{equation}
\label{eq:bridge-final-delay}
\sT_2\ge\sL_2+B_{\rm br}+1.
\end{equation}
Let $T\ge\ell(\gamma)+\sT_2$ be a time at which
$x^{\alpha_+}(T)$ has $(K_2,\eta_2)$--good statistics, and put
$S=T-\tau$. For every point $Z$ of the indicated geodesic bridge,
\eqref{eq:tau'} gives
\[
D_Z:=d_\cT\big(Z,X_{x^{\alpha_+}(\tau)}\big)
\le16\eta_1\tau'+\sL_2.
\]
Since $\tau+\tau'\le\ell(\gamma)$, we have
$S\ge\tau'+\sT_2$. Consequently,
\[
S-D_Z
\ge(1-16\eta_1)\tau'+\sT_2-\sL_2
\ge B_{\rm br}+1,
\]
where we used the initial restriction $\eta_1\le1/100$.
Apply Theorem~\ref{thm:rays-linear} with reference differential
$p=x^{\alpha_+}(\tau)$, comparison differential $q_Z=Q_1(x,Z)$,
and elapsed time $S$. Their projective vertical foliation is $x$.
The reference endpoint is
$x^{\alpha_+}(T)\in K_2$. Every terminal subinterval of $[0,S]$
corresponds to an interval ending at $T$, so the required occupation
bound follows from its good statistics. The preceding inequality
gives the required time bound.

Since $\Ext_{X_p}(F_+(p))=1$, the theorem gives
\[
s_Z=S+\tfrac12\log\Ext_Z(F_+(p)),\qquad |s_Z-S|\le D_Z.
\]
In particular $s_Z\ge S-D_Z>0$. It also proves that $q_Z$ is
principal and gives the intrinsic AGY estimate \eqref{eq:s_Z}.
The theorem permits nongeneric $x$, and the same constant $B_{\rm br}$
applies to every $Z$ on the segment.

Thus the restriction on $\eta_1$ is imposed before choosing the
statistical compact sets, whereas $\sT_2$ is chosen last from
\eqref{eq:bridge-final-delay}. No smallness condition involving
the compact-dependent prefactor $\sC_2$ is used in (d).
For the later applications to Theorem~\ref{theorem:h:estimates},
take its length threshold $\sR_0\ge\sR_2+1$, using
Proposition~\ref{prop:properties:gamma:n}{\rm(b)}.
\end{proof}

\subsection*{Coarse estimate for the spike function}

We first compare the unit-area normalizations of $x$ in the two charts
chosen by Proposition~\ref{prop:thin}. Extremal-length distortion between
the selected points then expresses $H(\gamma,x)$ in terms of
$\tau,\tau'$. The metric bounds in that proposition will relate these
times to $d_G(a_+,x)$.

\begin{proof}[Proof of Theorem~\ref{theorem:h:estimates}{\rm(a)}]
Fix $x\in\PMF$ and the charts and times from
Proposition~\ref{prop:thin}. Write $F'=\gamma(F_{\alpha_-})$ and put
\[
\begin{gathered}
p_0=x^{\alpha_+}(0),\qquad p_1=x^{\alpha_+}(\tau),\\
p_3=x_{F'}(0),\qquad p_2=x_{F'}(\tau').
\end{gathered}
\]
The vertical foliations $F_0=F_+(p_0)$ and $G_0=F_+(p_3)$ represent
the same projective class $x$, so $F_0=cG_0$ for some $c>0$.
The points $p_0$ and $\gamma^{-1}p_3$ lie in the finitely many fixed
chart slices. These slices have compact closures, so continuity and
positivity of $\rho$ give uniform bounds
\[
\rho(F_0)\asymp1,\qquad \rho(\gamma^{-1}G_0)\asymp1.
\]
By \eqref{eq:H} and homogeneity,
\begin{equation}
\label{eq:H-estimate}
H(\gamma,x)^{1/w}
=\frac{\rho(F_0)}{\rho(\gamma^{-1}F_0)}
=\frac{\rho(F_0)}{c\,\rho(\gamma^{-1}G_0)}
\asymp c^{-1}.
\end{equation}

We estimate $c$ at the two ends of the connecting geodesic. The
vertical foliations of $p_1,p_2$ are $e^\tau F_0,e^{\tau'}G_0$.
Since both differentials have area one, the Hubbard--Masur area
identity gives
\[
\sqrt{\Ext_{X_{p_1}}(F_0)}=e^{-\tau},\qquad
\sqrt{\Ext_{X_{p_2}}(F_0)}=c e^{-\tau'}.
\]
Set $D=d_\cT(X_{p_1},X_{p_2})$. Extremal-length distortion yields
\[
e^{-D}\le c e^{\tau-\tau'}\le e^D,
\qquad
e^{\tau-\tau'-D}\le c^{-1}\le e^{\tau-\tau'+D}.
\]
Using $D\le16\eta_1\tau'+\sL_2$ from \eqref{eq:tau'}, we obtain
\[
e^{\tau-\tau'-16\eta_1\tau'}
\prec H(\gamma,x)^{1/w}
\prec e^{\tau-\tau'+16\eta_1\tau'}.
\]
The prefactors include $e^{\sL_2}$ and the chart comparison constants.
This argument applies directly to every $x$, including nongeneric
foliations; it requires no lift of the connecting geodesic.

Write $L=\ell(\gamma)$ and $d=d_G(a_+,x)$, and abbreviate
$C_+=\sC_2$, $C_-=\sC_2^{-1}$, and $B=\sL_2$.
Combining the deficit estimate \eqref{eq:ell} with
the scalar comparison gives
\begin{equation}
\label{eq:H-final}
e^{L-(2+24\eta_1)\tau'-B}
\prec H(\gamma,x)^{1/w}
\prec e^{L-(2-16\eta_1)\tau'}.
\end{equation}
The coefficient $24$ is the sum of the deficit coefficient $8$ and
the bridge coefficient $16$.

Put $r=e^{-\beta L}$ and $Z=1+d/r\ge1$. If $\tau'\le\sL_2$, then
\begin{equation}
\label{eq:tau'-small}
e^{\tau'}\le e^{\sL_2}\le e^{\sL_2}Z^{1/\beta}.
\end{equation}
Otherwise the lower metric bound and $\tau+\tau'\le L$ give
\[
d/r\ge C_-e^{\beta(L-\tau-8\eta_1\tau')}
\ge C_-e^{\beta(1-8\eta_1)\tau'}.
\]
Combining the two cases, we obtain
\begin{equation}
\label{eq:e^tau'}
e^{(1-8\eta_1)\tau'}\prec Z^{1/\beta},
\end{equation}
and hence
\begin{equation}
\label{eq:upper-tau'}
e^{\tau'}\prec Z^{1/(\beta(1-8\eta_1))}.
\end{equation}
The denominator is positive because $\eta_1\le1/100$.
For the reverse bound, $r\le e^{-\beta\tau}$ and
\eqref{eq:tau} and \eqref{eq:ell} imply
\[
Z\le(1+C_+)e^{\beta(L-\tau)}
\le(1+C_+)e^{\beta B}e^{\beta(1+8\eta_1)\tau'}.
\]
Therefore
\begin{equation}
\label{eq:lower-tau'}
e^{\tau'}\succ Z^{1/(\beta(1+8\eta_1))}.
\end{equation}

Substituting these bounds in \eqref{eq:H-final} and raising to the
power $w$ gives
\[
e^{wL}Z^{-w(2+24\eta_1)/(\beta(1-8\eta_1))}
\prec H(\gamma,x)
\prec
e^{wL}Z^{-w(2-16\eta_1)/(\beta(1+8\eta_1))}.
\]
The two errors in the numerator of the exponent are
\[
\frac{40w\eta_1}{1-8\eta_1},\qquad
\frac{32w\eta_1}{1+8\eta_1}.
\]
Both are at most $44w\eta_1$ for $\eta_1\le1/100$. Thus
\eqref{eq:H-bounds} holds with
\begin{equation}
\label{eq:coarse-eta}
\eta_{\rm c}=44w\eta_1.
\end{equation}
All dependence on $\beta$, the compact sets, and the additive errors
is allowed in the multiplicative constant $\sC$. The coefficient
$44w$ depends only on the surface.

For the ratio estimate, set $Z_u=1+d_G(a_+,u)/r$. If
$d_G(u,v)\le kr$ for a fixed $k\ge1$, the triangle inequality gives
\[
(k+1)^{-1}\le Z_v/Z_u\le k+1.
\]
Taking the ratio of the two coarse bounds therefore yields
\begin{equation}
\label{eq:coarse-ratio-k}
C_k^{-1}Z_u^{-2\eta_{\rm c}/\beta}
\le\frac{H(\gamma,v)}{H(\gamma,u)}
\le C_k Z_u^{2\eta_{\rm c}/\beta}.
\end{equation}
The factor $(k+1)^{(2w+\eta_{\rm c})/\beta}$ is part of $C_k$.
Taking $k=1$ proves \eqref{eq:H:coarse:Lipschitz}.

The conclusions in part (a) remain valid with any
$\eta\ge\eta_{\rm c}$. The choice $\eta=256w\eta_1$ in
\eqref{eq:spike-error-choice} leaves room for the Lipschitz estimate
and the passage from $D_G$ to $d_G$. Its restriction on $\eta_1$
depends only on the surface and is imposed before choosing
$K_1,K_2$ and hence before determining $\beta$. In particular,
$4\eta_{\rm c}\le\eta<1/2$, so part (a) also holds with $\eta/4$
in place of $\eta$. We will check below that the same choice suffices
for part (b), without any comparison between $\eta$ and $\beta$.
\end{proof}

\subsection*{The Lipschitz estimate for the spike function}

To compare two nearby foliations, we transport their coordinate vectors
from the source chart to the target chart along five paths. When
$D_G(x,y)$ is sufficiently small compared with $r$, the two source
rays are close at a time beyond $\ell(\gamma)$. We use this extra
time to place both foliations in compatible coordinates along the
paths. After bounding the transport and checking these domains, we
estimate the relative coordinate errors in the quotient defining $H$.

Fix $x\in\PMF$, choose the charts and times from
Proposition~\ref{prop:thin}, and write $F'=\gamma(F_{\alpha_-})$.
Put
\[
\begin{gathered}
p_0=x^{\alpha_+}(0),\qquad p_1=x^{\alpha_+}(\tau),\qquad
p_2=x_{F'}(\tau'),\\
p_3=x_{F'}(0),\qquad p_4=\gamma(q_{\alpha_-}),\qquad
F_0=F_+(p_0).
\end{gathered}
\]
Thus $[F_0]=x$ and $p_1,p_2\in K_2$. The four paths outside the bridge
are always defined:
\begin{itemize}
\item $\xi_0$ joins $q_{\alpha_+}$ to $p_0$ in
$\Hor(F_{\alpha_+})$;
\item $\xi_1$ is $x^{\alpha_+}([0,\tau])$;
\item $\xi_3$ is the reverse of $x_{F'}([0,\tau'])$;
\item $\xi_4$ joins $p_3$ to $p_4$ in
$\gamma\Hor(F_{\alpha_-})$.
\end{itemize}
Choose $\xi_0$ and $\xi_4$ to be piecewise $C^1$ paths in their
fixed-horizontal, area-one slices, each of AGY length less than $1/300$.
Such paths exist because the construction in Section~\ref{sec:Charts}
gives
\begin{equation}
\label{eq:chart-origin-radius}
d_\AGY\big(q_\alpha,z^\alpha(0)\big)<\frac1{400}
\qquad(z\in W_\alpha).
\end{equation}
Indeed, this bound was imposed on $\overline W_\alpha'$ and is
preserved by the mapping class used to define $W_\alpha$.
The definition of intrinsic distance then gives the required paths;
no minimizing path is needed. Every point of $\xi_0$ has intrinsic
distance less than $1/300$ from $q_{\alpha_+}$, and the analogous
statement holds for $\gamma^{-1}\xi_4$ centered at $q_{\alpha_-}$.
The uniform leaf coordinates of Lemma~\ref{lem:AGY:vs:cV}
therefore compute each short path's transport in one period chart.

Let $\zeta:[0,1]\to\cT$ be the geodesic from $X_{p_1}$ to $X_{p_2}$.
When its normalized lift is principal, define
\[
\xi_2(s)=Q_1(x,\zeta(s)).
\]
This path joins $p_1$ to $p_2$ and fixes the projective vertical foliation.
It is principal whenever $F_0$ is generic, and also whenever the
hypothesis of Proposition~\ref{prop:thin}{\rm(d)} holds. We only use its
transport when this condition has been verified. In that case the
concatenated path gives the diagram
\begin{align*}
L(q_{\alpha_+})&\xrightarrow{\GM_{\xi_0}}L(p_0)
 \xrightarrow{\GM_{\xi_1}}L(p_1)
 \xrightarrow{\GM_{\xi_2}}L(p_2),\\
L(p_2)&\xrightarrow{\GM_{\xi_3}}L(p_3)
 \xrightarrow{\GM_{\xi_4}}L(p_4)
 \xrightarrow{(\gamma^{-1})^*}L(q_{\alpha_-}).
\end{align*}
Whenever all the intervening paths exist, write
\[
\GM_{ij}=\GM_{\xi_j}\circ\cdots\circ\GM_{\xi_i}
\qquad(0\le i\le j\le4).
\]

\begin{claim}
\label{claim:GM:estimates}
There is a prefactor $M_{\rm br}\ge1$, depending only on the fixed compact
data and $\sL_2$, for which the following estimates hold. All norms are
AGY norms.
\begin{enumerate}
\setcounter{enumi}{-1}
\item \label{item:GM0}
$\max\{\|\GM_{\xi_0}\|,\|\GM_{\xi_0}^{-1}\|\}\prec1$, and
\[
\|\Psi_{q_{\alpha_+}}(F_0)\|_{q_{\alpha_+}}
\asymp\|\Psi_{p_0}(F_0)\|_{p_0}.
\]
\item \label{item:GM1}
$\max\{\|\GM_{\xi_1}\|,\|\GM_{\xi_1}^{-1}\|\}\asymp e^\tau$, and
\[
\|\Psi_{p_1}(F_0)\|_{p_1}
=\|\GM_{\xi_1}(\Psi_{p_0}(F_0))\|_{p_1}
\asymp e^{-\tau}\|\Psi_{p_0}(F_0)\|_{p_0}.
\]
\item \label{item:GM2}
If the lift $\xi_2$ is principal, then
\begin{equation}
\label{eq:bridge-operator}
\max\{\|\GM_{\xi_2}\|,\|\GM_{\xi_2}^{-1}\|\}
\le M_{\rm br}e^{16\eta_1\tau'}.
\end{equation}
The coefficient $16$ is independent of the compact sets and statistical
parameters; their dependence is confined to $M_{\rm br}$.
\item \label{item:GM3}
$\max\{\|\GM_{\xi_3}\|,\|\GM_{\xi_3}^{-1}\|\}\asymp e^{\tau'}$, and
\[
\|\Psi_{p_3}(F_0)\|_{p_3}
=\|\GM_{\xi_3}(\Psi_{p_2}(F_0))\|_{p_3}
\asymp e^{\tau'}\|\Psi_{p_2}(F_0)\|_{p_2}.
\]
\item \label{item:GM4}
$\max\{\|\GM_{\xi_4}\|,\|\GM_{\xi_4}^{-1}\|\}\prec1$, and
\[
\|\Psi_{p_3}(F_0)\|_{p_3}\asymp\|\Psi_{p_4}(F_0)\|_{p_4}.
\]
\end{enumerate}
\end{claim}

\begin{proof}
For $\xi_0$ and $\xi_4$, the length bound $1/300$ and
\eqref{eq:AGY-GM-path} bound both operator norms and their inverses
by $e^{1/300}$. Along $\xi_0$, the point with projective vertical
foliation $x$ and the same horizontal foliation is the fixed point
$p_0$; it is at intrinsic distance less than $1/300$ from every
point of $\xi_0$. Along $\xi_4$, the corresponding point is $p_3$.
Lemma~\ref{lem:AGY:vs:cV} and its coordinate extension supply the
domains and continuous leaf paths required by
Lemma~\ref{lem:Psi-is-identity}. Thus transport identifies the
endpoint coordinates of $F_0$, giving the scalar comparisons as well.

For the flow segments $\xi_1,\xi_3$, Lemma~\ref{lem:GM-leaf}
gives the foliation-vector identities and
Lemma~\ref{lem:GM:near:geodesic} gives the operator upper bounds
and the scalar decay or growth. Applying the scalar decay to the
inverse of $\xi_1$, and the scalar growth to $\xi_3$, gives the
matching lower bounds for the maxima of the operator norms.

It remains to bound transport along $\xi_2$. Choose once a
$\Mod$--compact set $\widehat K_2$ in the principal stratum
containing a neighborhood of $K_2$. This is auxiliary endpoint
data; the standing set $K_2$ is unchanged. Put
$D=d_\cT(X_{p_1},X_{p_2})$. If $x$ is generic,
Proposition~\ref{prop:GM-new} and \eqref{eq:GM-AGY} give
\[
\max\{\|\GM_{\xi_2}\|,\|\GM_{\xi_2}^{-1}\|\}
\le A_{\widehat K_2}^{\,2}e^{B_{\widehat K_2}}e^D.
\]
Only the endpoints must lie in $\widehat K_2$.

Now suppose $x$ is nongeneric and $\xi_2$ is principal. Keep the
base geodesic $\zeta$ fixed. Continuity of
$(u,s)\mapsto Q_1(u,\zeta(s))$ and compactness of $[0,1]$ give
a contractible neighborhood $U$ of $x$ such that every lift over
$U\times[0,1]$ is principal and its endpoints lie in
$\widehat K_2$. The pullback cohomology local system on this
contractible product is trivial. In this trivialization all the
lifts have the same transport matrix. Choose generic $x_n\to x$
in $U$ and apply the preceding estimate to
$\xi_{2,n}(s)=Q_1(x_n,\zeta(s))$. Their base distance is exactly
$D$. Continuity of the endpoint AGY norms passes both the operator
and inverse-operator bounds to the limit.

Since $D\le16\eta_1\tau'+\sL_2$, we may therefore take
\[
M_{\rm br}=A_{\widehat K_2}^{\,2}
e^{B_{\widehat K_2}+\sL_2}
\]
in \eqref{eq:bridge-operator}. The neighborhood $U$ may depend on
$\xi_2$, but its size does not enter the bound. The approximants
require no good statistics or new choices of $\tau,\tau'$.
The sets $K_2,\widehat K_2$ and the additive error $\sL_2$ are
fixed before $M_{\rm br}$; this prefactor enters no later
smallness condition on $\eta_1$. The delay $\sT_2$ in
Proposition~\ref{prop:thin}{\rm(d)} remains the one chosen in
\eqref{eq:bridge-final-delay}.
\end{proof}

Whenever $\xi_2$ is principal, Lemma~\ref{lem:GM-leaf} also
identifies the foliation vector across it. Composing with the
four other identities just proved gives
\begin{equation}
\label{eq:F_0-okay}
\GM_{0i}\big(\Psi_{q_{\alpha_+}}(F_0)\big)=\Psi_{p_i}(F_0)
\qquad(0\le i\le4).
\end{equation}

The spike is a quotient of two powers of $\rho$. Nearby points give
nearby vectors in the source chart, but transport may greatly increase
their relative separation in the target chart. We therefore use a
local estimate for the numerator and an estimate for reciprocal powers
that applies to any two nonzero target vectors. The reciprocal-power
estimate in Lemma~\ref{lem:diff:normalized}{\rm(b)} introduces the larger
of the two spike values, which \eqref{eq:coarse-ratio-k} controls.

We now express $\rho$ in the fixed charts. On
$\Psi_{q_\alpha}(\R_+W_\alpha)$, define
\begin{equation}
\label{eq:rho-chart}
\rho_\alpha(v)=\rho\big(\Psi_{q_\alpha}^{-1}(v)\big).
\end{equation}
This function is positive and homogeneous of degree one. We may extend
it to $L(q_\alpha)\setminus\{0\}$ so that its restriction to the AGY
unit sphere is positive and Lipschitz. Here is why the extension can
be made with uniform bounds. In the construction of the charts in
Section~\ref{sec:Charts}, the cone over $\overline W_\alpha$ is contained
in $\R_+\cV_+(q_\alpha)$. Its unit directions form a compact subset of
the coordinate domain. On a neighborhood of this set,
$\rho\circ\Psi_{q_\alpha}^{-1}$ is locally Lipschitz: the foliation
coordinates are locally piecewise linear in period coordinates, as in
Proposition~\ref{prop:Local-Product}, and $\rho$ is locally Lipschitz by our
standing regularity convention. Its restriction
to the compact set of unit directions is thus Lipschitz and has a
positive minimum. Extend this restriction to the whole unit sphere by
a Lipschitz extension, truncate its values to its original minimum and
maximum, and then extend homogeneously. We continue to write
$\rho_\alpha$ for this extension; all evaluations representing
foliations will be in the original cone.

Since there are finitely many charts, there are constants
$0<c_\rho\le M_\rho$ and $L_\rho<\infty$, depending only on the
fixed charts and $\rho$, such that
\[
c_\rho\|v\|_{q_\alpha}\le\rho_\alpha(v)
\le M_\rho\|v\|_{q_\alpha},
\]
and every restriction to the unit sphere has Lipschitz constant at
most $L_\rho$. In particular, these constants are chosen independently
of the statistical compact sets.

\begin{lemma}
\label{lem:diff:normalized}
Let $L$ be a finite-dimensional normed vector space, and let $\rho$ be
the degree-one homogeneous extension of a positive Lipschitz function
$\rho_0$ on its unit sphere. Fix $w\ge1$. For nonzero $\bfv,\bfw\in L$,
put $\delta=\|\bfv-\bfw\|_L/\|\bfv\|_L$. Then:
\begin{itemize}
\item[{\rm (a)}] If $\delta\le1/2$, then
\[
|\rho(\bfv)^w-\rho(\bfw)^w|
\prec_{w,\rho}\|\bfv-\bfw\|_L\,\rho(\bfv)^{w-1}.
\]
\item[{\rm (b)}] Without any restriction on $\delta$,
\begin{equation}
\label{eq:global-reciprocal-power}
|\rho(\bfv)^{-w}-\rho(\bfw)^{-w}|
\prec_{w,\rho}
\delta\max\{\rho(\bfv)^{-w},\rho(\bfw)^{-w}\}.
\end{equation}
\end{itemize}
\end{lemma}

\begin{proof}
Write $u=\bfv/\|\bfv\|_L$ and $u'=\bfw/\|\bfw\|_L$.
The reverse triangle inequality gives $\|u-u'\|_L\le2\delta$.
Choose bounds $0<c_\rho\le\rho_0\le M_\rho$ on the unit sphere,
and let $L_\rho$ be a Lipschitz constant for $\rho_0$.
Expanding $\rho(\bfv)=\|\bfv\|_L\rho_0(u)$ and similarly for
$\bfw$ gives
\[
|\rho(\bfv)-\rho(\bfw)|
\le(M_\rho+2L_\rho)\|\bfv-\bfw\|_L.
\]
If $\delta\le1/2$, the two norms are within a factor of two,
so $\rho(\bfw)\asymp_\rho\rho(\bfv)$. The mean value theorem
therefore gives, for $b\in\{w,-w\}$,
\[
|\rho(\bfv)^b-\rho(\bfw)^b|
\prec_{b,\rho}\|\bfv-\bfw\|_L\,\rho(\bfv)^{b-1}.
\]
This proves (a). Taking $b=-w$ and using
$\rho(\bfv)\ge c_\rho\|\bfv\|_L$ proves (b) when
$\delta\le1/2$. If $\delta>1/2$, the difference of the two
positive reciprocal powers is at most their maximum, which is
at most $2\delta$ times that maximum. This proves (b) in all cases.
\end{proof}

We will estimate $H$ across a single small $D_G$-step. The next
lemma passes this estimate to $d_G$ by summing over chains.

\begin{lemma}[One-step reduction]
\label{lem:one-step-to-path}
Assume that the coarse bounds in part {\rm (a)} of
Theorem~\ref{theorem:h:estimates} hold with $\eta/4$ in place of $\eta$.
Suppose there are constants
$\ep_{\rm step}>0$ and $C_{\rm step}>0$ such that, whenever
\[
D_G(u,v)\le \ep_{\rm step}r,
\]
we have
\begin{equation}
\label{eq:one-step-H-Lipschitz}
\big|H(\gamma,u)-H(\gamma,v)\big|
\le
\frac{C_{\rm step}}{r}\,
H(\gamma,u)
\left(1+\frac{d_G(a_+,u)}{r}\right)^{\eta/(2\beta)}
D_G(u,v).
\end{equation}
Then part {\rm (b)} of Theorem~\ref{theorem:h:estimates} holds, after increasing
the constant $\sC$.
\end{lemma}

\begin{proof}
The coarse bounds with $\eta/4$ in place of $\eta$ give, whenever
$d_G(u,v)\le2r$,
\begin{equation}
\label{eq:H-coarse-two-r}
H(\gamma,v)
\prec
H(\gamma,u)
\left(1+\frac{d_G(a_+,u)}{r}\right)^{\eta/(2\beta)}.
\end{equation}
Indeed, take the ratio of the bounds in \eqref{eq:H-bounds}
and use the triangle inequality for $d_G$.

The conclusion is immediate if $x=y$. Otherwise, let
$0<d_G(x,y)\le r$, and for $0<\varepsilon\le d_G(x,y)$ choose
a chain $x=x_0,\ldots,x_k=y$ with
\[
\sum_{i=0}^{k-1}D_G(x_i,x_{i+1})
\le d_G(x,y)+\varepsilon\le2r.
\]
Every vertex satisfies $d_G(x,x_i)\le2r$. Thus, with
$Z=1+d_G(a_+,x)/r$, the triangle inequality and
\eqref{eq:H-coarse-two-r} give
\[
1+\frac{d_G(a_+,x_i)}r\le3Z,\qquad
H(\gamma,x_i)\prec H(\gamma,x)Z^{\eta/(2\beta)}.
\]
If $D_G(x_i,x_{i+1})\le\ep_{\rm step}r$, the assumed one-step
estimate now implies
\[
|H(\gamma,x_i)-H(\gamma,x_{i+1})|
\prec\frac1r H(\gamma,x)Z^{\eta/\beta}D_G(x_i,x_{i+1}).
\]
For a larger edge, bound the difference by
$H(\gamma,x_i)+H(\gamma,x_{i+1})$. The preceding vertex bound,
$Z\ge1$, and $D_G(x_i,x_{i+1})>\ep_{\rm step}r$ give the
same estimate, with an additional factor $\ep_{\rm step}^{-1}$.
Thus all edges are covered. Summing gives
\[
|H(\gamma,x)-H(\gamma,y)|
\le\frac{\sC}{r}H(\gamma,x)Z^{\eta/\beta}
\big(d_G(x,y)+\varepsilon\big).
\]
Letting $\varepsilon\to0$ proves part (b).
\end{proof}

We now prove the one-step estimate \eqref{eq:one-step-H-Lipschitz}. Thus,
until the end of the one-step argument, we assume
\[
D_G(x,y)\le \ep_{\rm step}r,
\qquad
r=e^{-\beta\ell(\gamma)}.
\]
The inequality is immediate when $y=x$, so assume $y\ne x$.
We will check that $y\in W_{\alpha_+}$ and that the separation time
$T=\tau_{\alpha_+}(x,y)$ is finite and attained. At time $T$, the
source rays $x^{\alpha_+}$ and $y^{\alpha_+}$ are within $1/100$.
To enter the target chart, we need the finer accuracy
$\rr_{\rm chart}$. We will move backward by a fixed time $h$
to bring the source points closer. We then project to horizontal
horospheres at forward points on the rays from $\xi_2$, and flow
the projected points back. Applying Theorem~\ref{thm:rays-linear}
with a smaller error and an additional delay ensures that these
projections are defined after the backward move.

Fix the charts, $\rho$, the statistical compact sets, $\beta$,
$C_\AGY$, and all constants of Proposition~\ref{prop:thin}.
Choose, in this order,
\begin{equation}
\label{eq:domain-buffer-choice}
\delta_{\rm dom}=\min\left\{\frac{\rr_{\rm chart}}4,\frac1{1000}\right\},
\qquad h\ge1\quad\text{with}\quad
\frac{C_\AGY}{100}e^{-\beta h}<\delta_{\rm dom}.
\end{equation}
Set
\[
\ep_{\rm br}=\frac{e^{-2h}}{200},
\qquad B_{\rm dom}=B(K_2,\eta_2,\ep_{\rm br}),
\]
where $B$ is from Theorem~\ref{thm:rays-linear}, and put
\begin{equation}
\label{eq:domain-delay-choice}
\sT_{\rm dom}
=\max\{\sT_2,\sL_2+B_{\rm dom}+h+1\}.
\end{equation}
The error $\ep_{\rm br}$ compensates for the possible factor $e^{2h}$
of backward growth along a stable leaf. The added delay ensures both
this finer approximation and a positive time after moving backward.
These choices leave $K_1,K_2,\eta_1,\eta_2,\beta$, and $\sT_2$
unchanged. We choose $\ep_{\rm step}$ only after these new constants
are fixed.

Let $C_{\rm src}\ge1$ bound the constant in \eqref{eq:diff:Psi} below.
Its proof uses only the chart data and Lemma~\ref{lem:dG:lower:bound},
once $2|\cA|\ep_{\rm step}\le1/2$; it is independent of
$\ep_{\rm step}$. Let $C_{\rm num}\ge1$ be the resulting constant
in \eqref{eq:diff:rho:Psi}, obtained from $C_{\rm src}$,
Lemma~\ref{lem:diff:normalized}{\rm(a)}, and the fixed bounds for
$\rho_\alpha$, under the preliminary restriction
$C_{\rm src}\ep_{\rm step}\le1/2$. Both constants are independent
of $\ep_{\rm step}$. Choose
\begin{equation}
\label{eq:epsilon-step-choice}
0<\ep_{\rm step}\le
\min\left\{1,\frac1{4|\cA|},
\frac{e^{-\beta\sT_{\rm dom}}}{2|\cA|},
\frac1{4C_{\rm src}},\frac1{2C_{\rm num}}\right\}.
\end{equation}
The first three bounds give the required source chart and separation
time; the last two make the source vectors and their numerator values
comparable. There is no smallness requirement on the transported
target vectors. Thus no transport prefactor enters a restriction on
$\eta_1$ or $\eta$.

Since $x\in V^{(2)}_{\alpha_+}$, we have
$\theta_{\alpha_+}(x)\ge1/|\cA|$ and hence
\[
D_G^{\alpha_+}(x,y)
\le2|\cA|D_G(x,y)
\le2|\cA|\ep_{\rm step}r\le\tfrac12.
\]
An endpoint outside $W_{\alpha_+}$ would give
$D_G^{\alpha_+}(x,y)=e^\beta>1$, so $y\in W_{\alpha_+}$.
Set $F=F_+(y^{\alpha_+}(0))$; as before,
$F_0=F_+(x^{\alpha_+}(0))$. Thus $[F]=y$ and $[F_0]=x$.

Before using transport, we verify that the middle lift is principal.
By Proposition~\ref{prop:D}, the separation time
$T=\tau_{\alpha_+}(x,y)$ is finite and attained. Since
$e^{-\beta T}=D_G^{\alpha_+}(x,y)$, the preceding bound and
\eqref{eq:epsilon-step-choice} give
\[
T\ge\ell(\gamma)+\sT_{\rm dom}.
\]
At this time $x^{\alpha_+}(T)$ has $(K_2,\eta_2)$--good statistics
and
\begin{equation}
\label{eq:close-at-T}
d_\AGY\big(x^{\alpha_+}(T),y^{\alpha_+}(T)\big)\le \frac1{100}.
\end{equation}
Since $\sT_{\rm dom}\ge\sT_2$, Proposition~\ref{prop:thin}{\rm(d)}
shows that $\xi_2(s)=Q_1(x,\zeta(s))$ is principal, including
when $x$ is nongeneric. Thus $\GM_{04}$ is defined and
\eqref{eq:bridge-operator} applies. With $x$ and these paths fixed,
define the linear map
\[
\Upsilon_\gamma=(\gamma^{-1})^*\GM_{04}:
L(q_{\alpha_+})\longrightarrow L(q_{\alpha_-}).
\]
We must still show that this map represents $\gamma^{-1}$ on
nearby foliation coordinates. This requires checking the coordinate
domains along all five paths.

\begin{claim}
\label{claim:well-defined}
If $D_G(x,y)\le\ep_{\rm step}r$, then
\[
\bar y=\gamma^{-1}y\in V^{(3)}_{\alpha_-}\subset W_{\alpha_-}.
\]
For every nonzero $F\in\MF$ with $[F]=y$, both coordinates in
the following identity are defined, and
\[
\Upsilon_\gamma\big(\Psi_{q_{\alpha_+}}(F)\big)
=\Psi_{q_{\alpha_-}}(\gamma^{-1}F).
\]
The same identity holds for every positive multiple of $F_0$.
\end{claim}

\begin{proof}
For $F$ proportional to $F_0$, the identity follows from
\eqref{eq:F_0-okay} and equivariance. For $[F]=y$, we construct,
at each point $z$ of the concatenated path
\[
\xi=\xi_4*\xi_3*\xi_2*\xi_1*\xi_0,
\]
a unit-area differential with projective vertical foliation $y$
and horizontal foliation $F_-(z)$. Keeping this differential within
intrinsic distance $1/100$ of $z$ will put $F$ in the leaf-coordinate
domain \eqref{eq:Psi-leaf-domain}. We first move backward from $T$
to make the source points closer. Projection will then give the
required points along $\xi_2$ and the returning ray, close enough
to enter the target chart.

Write
\[
q(t)=x^{\alpha_+}(t),\qquad p(t)=y^{\alpha_+}(t).
\]
The attained separation time $T$ satisfies \eqref{eq:close-at-T},
and $q(T)$ has $(K_2,\eta_2)$--good statistics. In particular,
\[
\big|\{u\in[T-h,T]:q(u)\in K_2\}\big|
\ge(1-\eta_2)h>\tfrac12h.
\]
Corollary~\ref{cor:AGY:half} and \eqref{eq:domain-buffer-choice}
therefore give
\begin{equation}
\label{eq:source-before-separation}
d_\AGY(q(T-h),p(T-h))
\le \frac{C_\AGY}{100}e^{-\beta h}<\delta_{\rm dom}.
\end{equation}
Only good statistics at $T$ are used here; no good-statistics
assertion at $T-h$ is needed. Also,
\[
T-h\ge\ell(\gamma)+\sT_{\rm dom}-h>\ell(\gamma)\ge\tau.
\]

We check the first two paths. For $z\in\xi_0$, the length chosen
for $\xi_0$ and \eqref{eq:chart-origin-radius} give
\[
d_\AGY(z,p(0))
<\frac1{300}+\frac1{400}<\frac1{100}.
\]
Here $z$ and $p(0)$ have the same horizontal foliation
$F_{\alpha_+}$. The required point along $\xi_0$ is therefore
$p(0)$, and \eqref{eq:Psi-leaf-domain} defines $\Psi_z(F)$.
Along $\xi_1$, backward
nonexpansion from \eqref{eq:source-before-separation} gives
\[
d_\AGY(q(u),p(u))<\delta_{\rm dom}<\frac1{100}
\qquad(0\le u\le\tau).
\]
The points $p(u)$ give the same domain conclusion along $\xi_1$.

For $\xi_2$, we need a finer approximation than \eqref{eq:s_Z},
so that it remains useful after moving backward by $h$.
Fix $Z\in\zeta([0,1])$, and put
\[
q_Z=Q_1(x,Z),\qquad
S=T-\tau,\qquad D_Z=d_\cT(Z,X_{q(\tau)}).
\]
By \eqref{eq:tau'} and $\tau+\tau'\le\ell(\gamma)$,
\begin{align*}
S-D_Z
&\ge (1-16\eta_1)\tau'+\sT_{\rm dom}-\sL_2\\
&\ge B_{\rm dom}+h+1.
\end{align*}
The last inequality uses $\eta_1\le1/100$ and
\eqref{eq:domain-delay-choice}. Apply Theorem~\ref{thm:rays-linear}
to $q(\tau)$ and $q_Z$, with elapsed time $S$ and error
$\ep_{\rm br}$. Its terminal-interval hypothesis follows from the
good statistics of $q(T)$. The theorem applies to arbitrary $x$,
including nongeneric foliations, with the same $B_{\rm dom}$
for every $Z$.

For $G=F_+(q(\tau))$, the resulting time is
\[
s_Z=S+\tfrac12\log\Ext_Z(G).
\]
It satisfies
\[
s_Z\ge S-D_Z>h,\qquad
F_+(q_Z(s_Z))=F_+(q(T)),\qquad
d_\AGY(q(T),q_Z(s_Z))<\ep_{\rm br}.
\]
The equality of measured vertical foliations puts these two
differentials in the same stable horosphere. The time-reversed
form of \eqref{eq:Monotone} now gives
\begin{equation}
\label{eq:bridge-before-separation}
d_\AGY(q(T-h),q_Z(s_Z-h))
<e^{2h}\ep_{\rm br}=\frac1{200}.
\end{equation}

We can now project. Set $q_Z^*=q_Z(s_Z-h)$ and define
\[
p_Z^*=\pi^{\rm cs}_{q_Z^*}(p(T-h)).
\]
To justify this definition, apply
Theorem~\ref{theorem:cs-projection} with target $q_Z^*$ and ball
centered at $q(T-h)$. Its assumption on the target follows from
\eqref{eq:bridge-before-separation} and $D_\AGY\le d_\AGY$.
Both inputs $q(T-h)$ and $p(T-h)$ lie strictly inside that
open ball of radius $1/50$, by \eqref{eq:source-before-separation}.
The projection sends $q(T-h)$ to $q_Z^*$, since they have
the same projective vertical foliation. Thus
\[
d_\AGY(q_Z^*,p_Z^*)
\le4D_\AGY(q(T-h),p(T-h))<4\delta_{\rm dom}.
\]
In particular $p_Z^*$ is principal. Define
$p_Z=g_{-(s_Z-h)}p_Z^*$. The projection preserves the projective
vertical foliation and fixes the measured horizontal foliation,
so uniqueness of the differential with prescribed filling
foliations identifies
\[
p_Z=Q_1(y,F_-^Z),\qquad F_-^Z=F_-(q_Z).
\]
This construction establishes that the pair defining $p_Z$ fills;
it does not assume its existence before applying the projection.
Backward nonexpansion in the common horizontal horosphere gives
\begin{equation}
\label{eq:bridge-domain-bound}
d_\AGY(q_Z,p_Z)<4\delta_{\rm dom}
\le\min\left\{\rr_{\rm chart},\frac1{250}\right\}.
\end{equation}
Consequently \eqref{eq:Psi-leaf-domain} defines $\Psi_{q_Z}(F)$.
The formula for $p_Z$ and continuity of the filling-pair
correspondence show that $p_Z$ varies continuously along $\xi_2$.
The constants are independent of $Z$ and of the common foliation.

At the endpoint $q_Z=x_{F'}(\tau')$, the horizontal foliation is
$e^{-\tau'}F'$, so $p_Z=y_{F'}(\tau')$.
Flowing this pair backward proves
\[
d_\AGY(x_{F'}(u),y_{F'}(u))<4\delta_{\rm dom}<\frac1{100}
\qquad(0\le u\le\tau').
\]
Equation~\eqref{eq:Psi-leaf-domain} therefore supplies the domain
inclusion along $\xi_3$. At $u=0$, apply $\gamma^{-1}$ to obtain
\[
d_\AGY\big(\bar x_{F_{\alpha_-}}(0),
                 \bar y_{F_{\alpha_-}}(0)\big)
<4\delta_{\rm dom}\le\rr_{\rm chart},
\qquad \bar y=\gamma^{-1}y.
\]
Both differentials have already been constructed in their common
horizontal horosphere, so $\bar x,\bar y\in\PMF(F_{\alpha_-})$.
Proposition~\ref{prop:thin} gives
$\bar x\in V^{(2)}_{\alpha_-}$. Theorem~\ref{theorem:charts}{\rm(f)}
therefore yields
\[
\bar y\in V^{(3)}_{\alpha_-}\subset W_{\alpha_-}.
\]

Finally consider $\xi_4$. For $z\in\gamma^{-1}\xi_4$, the
chosen length and \eqref{eq:chart-origin-radius} give
\[
d_\AGY\big(z,\bar y^{\alpha_-}(0)\big)
<\frac1{300}+\frac1{400}<\frac1{100}.
\]
The required point along $\xi_4$ is $\gamma\bar y^{\alpha_-}(0)$.
Equation~\eqref{eq:Psi-leaf-domain} and equivariance give the domain.

The points constructed on consecutive pieces agree at the four joints:
they are the unique unit-area differentials with projective vertical
foliation $y$ and the horizontal foliation of the corresponding point
of $\xi$. They therefore vary continuously along $\xi$, everywhere
at intrinsic distance less than $1/100$. The joint leaf-coordinate
construction after Lemma~\ref{lem:AGY:vs:cV} supplies continuous
affine leaf paths joining them to $\xi$. These paths remain principal,
and all open-domain inequalities are strict.
Lemma~\ref{lem:Psi-is-identity} now gives
\[
\GM_{04}\big(\Psi_{q_{\alpha_+}}(F)\big)
=\Psi_{\gamma(q_{\alpha_-})}(F).
\]
Applying $(\gamma^{-1})^*$ and using equivariance proves
\[
\Upsilon_\gamma\big(\Psi_{q_{\alpha_+}}(F)\big)
=\Psi_{q_{\alpha_-}}(\gamma^{-1}F).
\]
Every domain used is a positive cone, so the conclusion holds for
every positive multiple of the chosen representative of $y$.
\end{proof}

The next claim bounds the relative coordinate errors at the source and
target by fixed multiples of $D_G(x,y)$ and
$r^{-1}e^{32\eta_1\tau'}D_G(x,y)$, respectively.

\begin{claim}
\label{claim:diff:psi}
We have
\begin{equation}
\label{eq:diff:Psi}
\big\|
\Psi_{q_{\alpha_+}}(F_0)-\Psi_{q_{\alpha_+}}(F)
\big\|_{q_{\alpha_+}}
\prec
D_G(x,y)\,
\norm{\Psi_{q_{\alpha_+}}(F_0)}_{q_{\alpha_+}},
\end{equation}
and
\begin{equation}
\label{eq:diff:Upsilon:gamma:Psi}
\big\|
\Upsilon_\gamma(\Psi_{q_{\alpha_+}}(F_0))
-
\Upsilon_\gamma(\Psi_{q_{\alpha_+}}(F))
\big\|_{q_{\alpha_-}}
\prec
e^{32\eta_1\tau'}\frac{D_G(x,y)}{r}\,
\norm{\Upsilon_\gamma(\Psi_{q_{\alpha_+}}(F_0))}_{q_{\alpha_-}}.
\end{equation}
\end{claim}

\begin{proof}
Write $D=D_G(x,y)$. For $s=0$ and $s=\tau$, the alternative
hypothesis of Lemma~\ref{lem:dG:lower:bound} follows from
\[
d_G^{\alpha_+}(x,y)\le D_G^{\alpha_+}(x,y)
\le2|\cA|\ep_{\rm step}e^{-\beta\ell(\gamma)}
\le\tfrac12e^{-\beta s},
\]
since $\tau\le\ell(\gamma)$ and
$2|\cA|\ep_{\rm step}\le1/2$. The lemma therefore gives
\[
d_\AGY\big(x^{\alpha_+}(s),y^{\alpha_+}(s)\big)
\prec D_G^{\alpha_+}(x,y)e^{\beta s}
\le2|\cA|D e^{\beta s}.
\]
At $s=0$, both points lie in the intrinsic $1/400$--ball about
$q_{\alpha_+}$. Proposition~\ref{prop:AGY-L(q)} then bounds
the difference of their coordinates by a fixed multiple of $D$.
Since $\|\Psi_{q_{\alpha_+}}(F_0)\|_{q_{\alpha_+}}\asymp1$
on the fixed chart slices, this proves \eqref{eq:diff:Psi}.

At time $\tau$, the vertical foliations are $e^\tau F_0$ and
$e^\tau F$. The construction in Claim~\ref{claim:well-defined}
puts the two differentials at intrinsic distance less than
$\delta_{\rm dom}\le1/1000$. Applying
Proposition~\ref{prop:AGY-L(q)} with center $p_1$, and then using
homogeneity, gives
\[
e^\tau\|\Psi_{p_1}(F)-\Psi_{p_1}(F_0)\|_{p_1}
\prec D e^{\beta\tau}.
\]
The real period has AGY norm one, as proved after
Lemma~\ref{lem:AGY-area}, so
\[
e^\tau\|\Psi_{p_1}(F_0)\|_{p_1}
=\|\Psi_{p_1}(F_+(p_1))\|_{p_1}=1.
\]
Thus the factor $e^\tau$ cancels in the relative error.

To follow that relative error through the remaining paths, put
\[
v_i=\GM_{0i}\big(\Psi_{q_{\alpha_+}}(F_0)\big),\qquad
w_i=\GM_{0i}\big(\Psi_{q_{\alpha_+}}(F)\big)
\quad(1\le i\le4).
\]
The coordinate identities in \eqref{eq:F_0-okay} and the proof of
Claim~\ref{claim:well-defined} identify these with
$\Psi_{p_i}(F_0)$ and $\Psi_{p_i}(F)$, respectively. The preceding
estimates give $\|w_1-v_1\|_{p_1}/\|v_1\|_{p_1}\prec D e^{\beta\tau}$.
Apply \eqref{eq:bridge-operator} to the difference and its inverse
bound to the reference vector. Then
\begin{equation}
\label{eq:GM02:bound}
\begin{aligned}
\frac{\|w_2-v_2\|_{p_2}}{\|v_2\|_{p_2}}
&\le M_{\rm br}^2e^{32\eta_1\tau'}
\frac{\|w_1-v_1\|_{p_1}}{\|v_1\|_{p_1}}\\
&\prec D e^{\beta\tau+32\eta_1\tau'}.
\end{aligned}
\end{equation}
The factor $M_{\rm br}^2$ enters the prefactor, leaving coefficient
$32$ unchanged.

Along the returning flow $\xi_3$, \eqref{eq:GM-exp} multiplies the
numerator by at most $e^{\tau'}$, while \eqref{eq:F-decay}
multiplies the denominator by exactly $e^{\tau'}$. The relative
error therefore does not increase. The final short path changes
it by at most a fixed factor, by
Claim~\ref{claim:GM:estimates}{\rm(}\ref{item:GM4}{\rm)}.
Finally, $(\gamma^{-1})^*$ is an AGY isometry and
$e^{\beta\tau}\le e^{\beta\ell(\gamma)}=r^{-1}$.
Using the definition of $\Upsilon_\gamma$ proves
\eqref{eq:diff:Upsilon:gamma:Psi}.
\end{proof}

We finish by comparing the numerator and the reciprocal denominator
of the spike. Write
\[
\begin{gathered}
u_0=\Psi_{q_{\alpha_+}}(F_0),\qquad
u=\Psi_{q_{\alpha_+}}(F),\qquad
v_0=\Upsilon_\gamma u_0,\qquad v=\Upsilon_\gamma u,\\
D=D_G(x,y),\qquad Z=1+\frac{d_G(a_+,x)}r.
\end{gathered}
\]
All four vectors are nonzero, and Claim~\ref{claim:well-defined}
places them in the respective foliation-coordinate cones.
At the source, \eqref{eq:diff:Psi} and
\eqref{eq:epsilon-step-choice} give relative error at most
$C_{\rm src}D\le1/4$, since $r<1$. Applying
Lemma~\ref{lem:diff:normalized}{\rm(a)} and the fixed comparison
between $\rho_{\alpha_+}$ and the norm gives
\begin{equation}
\label{eq:diff:rho:Psi}
|\rho_{\alpha_+}(u)^w-\rho_{\alpha_+}(u_0)^w|
\le C_{\rm num}D\,\rho_{\alpha_+}(u_0)^w.
\end{equation}
Since $C_{\rm num}D\le1/2$, the two numerator values are within
a factor of two.

At the target, the relative error need not be small. The unrestricted
estimate in Lemma~\ref{lem:diff:normalized}{\rm(b)}, together with
\eqref{eq:diff:Upsilon:gamma:Psi}, gives
\begin{equation}
\label{eq:diff:rho:Upsilon:Psi}
|\rho_{\alpha_-}(v)^{-w}-\rho_{\alpha_-}(v_0)^{-w}|
\prec\frac{D}{r}e^{32\eta_1\tau'}
\max\{\rho_{\alpha_-}(v)^{-w},\rho_{\alpha_-}(v_0)^{-w}\}.
\end{equation}
By \eqref{eq:H} and Claim~\ref{claim:well-defined},
\[
H(\gamma,x)=\frac{\rho_{\alpha_+}(u_0)^w}{\rho_{\alpha_-}(v_0)^w},
\qquad
H(\gamma,y)=\frac{\rho_{\alpha_+}(u)^w}{\rho_{\alpha_-}(v)^w}.
\]
Subtracting the two quotients gives
\begin{align*}
H(\gamma,x)-H(\gamma,y)
&=\frac{\rho_{\alpha_+}(u_0)^w-\rho_{\alpha_+}(u)^w}
        {\rho_{\alpha_-}(v_0)^w}\\
&\quad+\rho_{\alpha_+}(u)^w
\big(\rho_{\alpha_-}(v_0)^{-w}-\rho_{\alpha_-}(v)^{-w}\big).
\end{align*}
Applying \eqref{eq:diff:rho:Psi} and
\eqref{eq:diff:rho:Upsilon:Psi}, and using the comparability of the two
numerators, gives
\[
|H(\gamma,x)-H(\gamma,y)|
\prec
DH(\gamma,x)
+
\frac{D}{r}e^{32\eta_1\tau'}
\max\{H(\gamma,x),H(\gamma,y)\}.
\]
The global reciprocal-power estimate thus reduces the problem to
controlling the larger spike value. Since
$d_G(x,y)\le D\le\ep_{\rm step}r\le r$,
the coarse comparison \eqref{eq:coarse-ratio-k}, with $k=1$, gives
\[
\max\{H(\gamma,x),H(\gamma,y)\}
\prec H(\gamma,x)Z^{2\eta_{\rm c}/\beta}.
\]
Also, \eqref{eq:upper-tau'} gives
\[
e^{32\eta_1\tau'}
\prec
Z^{32\eta_1/(\beta(1-8\eta_1))}.
\]
The first term is absorbed since $r<1$ and $Z\ge1$. Consequently,
\begin{equation}
\label{eq:one-step-exponent}
|H(\gamma,x)-H(\gamma,y)|
\prec \frac{D}{r}H(\gamma,x)
Z^{\big(32\eta_1/(1-8\eta_1)+2\eta_{\rm c}\big)/\beta}.
\end{equation}
We can now see the entire error in the exponent. The bridge and its
inverse contribute $32\eta_1/(1-8\eta_1)$, and the coarse comparison
contributes $2\eta_{\rm c}$. For the initial choices
\eqref{eq:spike-error-choice},
\[
\frac{32\eta_1}{1-8\eta_1}+2\eta_{\rm c}
\le (35+88w)\eta_1
\le128w\eta_1=\frac\eta2,
\]
where $w\ge1$. Thus
\[
|H(\gamma,x)-H(\gamma,y)|
\le
\frac{C_{\rm step}}{r}
H(\gamma,x)
\left(1+\frac{d_G(a_+,x)}{r}\right)^{\eta/(2\beta)}
D_G(x,y).
\]
This is \eqref{eq:one-step-H-Lipschitz}. All compact-dependent
norm comparisons, the factor $M_{\rm br}^2$, and the constants in
\eqref{eq:coarse-ratio-k} and \eqref{eq:upper-tau'} enter
the prefactor $C_{\rm step}$; the exponent remains $\eta/(2\beta)$.
After fixing these constants and $\ep_{\rm step}$, the chain argument
of Lemma~\ref{lem:one-step-to-path} only enlarges $\sC$, allowing
dependence on $\ep_{\rm step}^{-1}$. Thus the order is
$\eta_1,\eta$ first, then the compact sets and $\beta$, then the
prefactors and the additional delay $\sT_{\rm dom}$, and finally
$\ep_{\rm step}$. The compact sets remain fixed throughout the
domain argument in Claim~\ref{claim:well-defined}. No comparison
between $\eta$ and $\beta$ is required. Lemma~\ref{lem:one-step-to-path}
now proves part {\rm (b)} of Theorem~\ref{theorem:h:estimates}.

\section{$Q$-Spikes}
\label{sec:spikes}

We now collect the properties of $H(\gamma,\cdot)$ that we will use in
Section~\ref{sec:basis}. These functions are large on a ball, decay away
from it, and have controlled variation on the scale of its radius.
The next section uses these properties to subtract positive sums of
spikes while keeping the remainder positive and sufficiently regular
to repeat the construction.

\subsection{Defining $Q$-spikes.}
\label{sec:subsec:abstract:setup}
Let $(X,d_G,\nu_\rho)$ be a metric measure space, and let $B_G(a,r)$
denote the open ball of radius $r$ centered at $a$.
Fix $W>0$ and $C_{\rm vol}>0$ such that, for every $a\in X$ and $r>0$,
\begin{equation}
\label{eq:Q:spike:volume:growth}
\nu_\rho(B_G(a,r)) \le C_{\rm vol} r^W. 
\end{equation}
The annular decomposition in the proof of Corollary~\ref{cor:measure_decay}
then gives, for every $W'>W$, $x\in X$, and $r>0$,
\begin{equation}
\label{eq:Qspike:measure:decay}
\int_{X\setminus B_G(x,r)}
\frac{d\nu_\rho(z)}{d_G(x,z)^{W'}}
\le D_\rho\,r^{W-W'},
\end{equation}
where $D_\rho$ depends only on $C_{\rm vol}$, $W$, and $W'$.

For a function $f:X\to\reals$ and $s>0$, define its
\emph{Lipschitz constant at $x$ on scale $s$} by
\[
D_s f(x)=
\sup_{\{y \st 0 < d_G(x,y) \le s\}}
\frac{|f(y)-f(x)|}{d_G(x,y)},
\]
with value zero if the set is empty.

Fix $0<\Xi<W$. We use a variant of Connell--Muchnik's $Q$-spikes
\cite[Definitions~4.1--4.2]{connell:muchnik:QC}, allowing the local
comparison and Lipschitz bounds to deteriorate away from the center.
A triple $(h,a,r)$, where $a\in X$, $r>0$, and $h:X\to\reals$ is
nonnegative with $h(a)>0$, is called a $Q$-spike if the following
estimates hold with implicit constants independent of $x,y$.

First, $h$ decays away from $a$: for every $x\in X$,
\begin{equation}
\label{eq:h-bounds}
h(x) \prec h(a) \left(1+ \frac{d_G(a,x)}{r}\right)^{-2W+\Xi/2}.
\end{equation}
On the ball $B_G(a,r)$, we also require the lower bound
\begin{equation}
\label{eq:h-local-lower}
h(x) \succ h(a).
\end{equation}
For any $x,y\in X$ with $d_G(x,y)\le r$, we require the comparison
\begin{equation}
\label{eq:h-h}
h(y) \prec h(x) \left(1+r^{-1}d_G(a,x)\right)^{\Xi}.
\end{equation}
The selected ball must satisfy the volume lower bound
\begin{equation}
\label{eq:spike:volume-lower}
\nu_\rho(B_G(a,r))\succ r^W.
\end{equation}
Finally, the variation of $h$ on scale $r$ must satisfy, for every $x\in X$,
\begin{equation}
\label{eq:D_r_bound_sharp}
D_r h(x) \prec
\frac{h(x)}{r}
\left(1+\frac{d_G(a,x)}{r}\right)^{\Xi/2}.
\end{equation}
Combining \eqref{eq:h-bounds} and \eqref{eq:D_r_bound_sharp} gives
the following consequence, which controls variation far from the center:
\begin{equation}
\label{eq:D_r_bound_decay}
D_r h(x) \prec
\frac{h(a)}{r}
\left(1+\frac{d_G(a,x)}{r}\right)^{-2W+\Xi}.
\end{equation}
A family has \emph{uniform spike estimates} if the same implicit
constants work for all its members.

\subsection{Verification of the $Q$-spike conditions.}
\label{sec:subsec:Qspike:verification}
Take $X=\PMF$ with the Gromov metric from Section~\ref{sec:gromov-distance}
and the measure $\nu_\rho$ from Section~\ref{sec:prelims}. Set
\[
W = w/\beta
\qquad\text{and}\qquad
\Xi = 2\eta/\beta,
\]
where $\eta$ is the constant from Theorem~\ref{theorem:h:estimates}
and $\beta$ is as in Section~\ref{sec:gromov-distance}.
The initial choice \eqref{eq:spike-error-choice} gives
$2\eta=512w\eta_1<1\le w$, and hence $\Xi<W$.
No further change to $\eta_1$ or the compact sets is needed.
Proposition~\ref{prop:volume} gives \eqref{eq:Q:spike:volume:growth}
for sufficiently small $r$. Since $\nu_\rho$ is a probability measure,
enlarging $C_{\rm vol}$ extends this bound to every $r>0$.

\begin{theorem}
\label{theorem:is:Q:spike}
Suppose $\gamma \in \Gamma(R)$ with $R >  \max(\sR_0, \sR_1)+1$. 
Let
\[
a = [F_+(\gamma)]
\qquad\text{and}\qquad
r = e^{-\beta \ell(\gamma)}.
\]
Let $h(x)=H(\gamma,x)$. Then $(h,a,r)$ is a $Q$-spike, with
constants independent of $\gamma$ and $R$ once the standing data are fixed.
\end{theorem}

\begin{proof}
By \eqref{eq:H}, we have $h(a)=e^{w\ell(\gamma)}$.
The two-sided estimate in Theorem~\ref{theorem:h:estimates}{\rm(a)}
therefore gives \eqref{eq:h-bounds} and \eqref{eq:h-local-lower}.
Its ratio estimate \eqref{eq:H:coarse:Lipschitz} gives \eqref{eq:h-h},
since $\Xi=2\eta/\beta$.

For the volume lower bound, choose $\alpha\in\cA$ such that
$\gamma\in\Gamma_\alpha(R)$. Proposition~\ref{prop:properties:gamma:n}{\rm(b)}
gives $\ell(\gamma)\ge R-1>\sR_1$. By part~{\rm(d)},
$a\in V^{(1)}_\alpha$ and $a^\alpha(\ell(\gamma))$ has
$(K_1,\eta_1)$--good statistics. Proposition~\ref{prop:volume} therefore gives
\eqref{eq:spike:volume-lower}.

Finally, Theorem~\ref{theorem:h:estimates}{\rm(b)} gives
\eqref{eq:D_r_bound_sharp}, since $\eta/\beta=\Xi/2$.
All constants are uniform in $\gamma$ and $R$.
\end{proof}

\section{Approximation by positive sums of spikes}
\label{sec:basis}
We now use the spikes from Section~\ref{sec:spikes} to subtract a
definite proportion of a positive function on the set covered by their
balls. The remainder must stay positive and satisfy the same estimates
at a smaller scale, so that we can repeat the construction.
Proposition~\ref{prop:subfunction} gives this step. In
Section~\ref{sec:first_moment}, we start with $f\equiv1$ and iterate it.

Let $(X,d_G,\nu_\rho)$ be as in \S\ref{sec:subsec:abstract:setup}. 
Let $(h_i,a_i,r_i)_{i=1}^N$ be a collection of $Q$-spikes.
We assume that all members of this collection satisfy the $Q$-spike
estimates with common bounds on their constants. 
We also assume there is a uniform constant $\sW_0$ such that, for
$i,j\in[1,N]$,
\[
r_i\le \sW_0 r_j.
\]
Let $r_{\max} = \max_i  r_i$, $r_{\min} = \min_i  r_i$. 
Let $Y\subset X$ be a set covered by the balls $\{B_G(a_i,r_i)\}$. We further
assume that there is a constant $\sB$ such that each point in $X$ is covered
by at most $\sB$ balls. In this section, we call a constant \emph{uniform} if
it depends on the constants in the $Q$-spike estimates in
\S\ref{sec:subsec:abstract:setup}, $C_{\rm vol}$,
$W$, $\Xi$, $\sW_0$, and $\sB$. Set
\begin{displaymath}
\sigma = \frac{W-\Xi}{2}. 
\end{displaymath}

\begin{proposition} 
\label{prop:subfunction}
There exist uniform constants $\mm \ge 2$ and $\epsilon_0 > 0$
such that the following holds:
Suppose $f \from X \to \R_+$ is a bounded positive function
which is bounded away from 0, and suppose there is a $\delta>0$ such that
\begin{equation}
\label{eq:local:subfunction-ratio}
D_\delta f(x) \le \frac{f(x)}{\delta},
\end{equation}
\begin{equation}
\label{eq:global:subfunction-ratio}
\frac{f(y)}{f(x)}\le  4\left(1+\frac{d_G(x,y)}{\delta}\right)^\sigma \text{ for all $x,y\in X$.}
\end{equation}
and
\begin{equation}
\label{eq:small:enough:balls}
r_{\max} \leq \delta/\mm.
\end{equation}
Then we can find numbers $\muG^i \geq0$ such that for
\[
g=\sum_{i=1}^N \muG^i h_{i}
\] 
and $F = f-g$ we have
\begin{displaymath}
\frac{1}{2}f(x) \le F(x)\leq f(x) \qquad \text{for $x \in X$,}
\end{displaymath}
\begin{displaymath}
F(x) \leq  (1-\epsilon_0)f(x), \qquad \text{for $x \in Y$.}
\end{displaymath}
Also,  $F$
satisfies both \eqref{eq:local:subfunction-ratio} and
\eqref{eq:global:subfunction-ratio} with  $\delta' = \frac{r_{\min}}{4}$ in
place of $\delta$. 
\end{proposition}

The local bound \eqref{eq:local:subfunction-ratio} controls variation
within distance $\delta$. The global bound
\eqref{eq:global:subfunction-ratio} controls contributions from distant
centers, where $f(a_i)$ may be much larger than $f(x)$.
Our choice of $\sigma$ gives $2W-\Xi-\sigma=W+\sigma>W$, so
the integral estimate \eqref{eq:Qspike:measure:decay} still applies
after allowing for this growth.

The next two lemmas express the decay bounds as integrals over the spike
balls. Bounded overlap will then let us estimate sums of spikes.

\begin{lemma}
\label{lemma:spike}
Suppose $(h,a,r)$ is a $Q$-spike.  For each $x \in B_G(a,r)^c$ we have
\[
0 \le h(x)\prec
h(a)\, r^{W-\Xi}
\int_{B_G(a,r)}
\frac{d\nu_\rho(y)}{d_G(x,y)^{2W-\Xi}}.
\]
\end{lemma}

\begin{proof}
Since $d_G(x,a) \ge r$, for any $y \in B_G(a,r)$ we have
\[
d_G(x,y) \le d_G(x,a)+r \le 2d_G(x,a),
\]
and hence
\[
\frac{1}{d_G(x,a)}
\prec
\frac{1}{d_G(x,y)}.
\]
Raising both sides to the power $2W-\Xi$ and integrating over $B_G(a,r)$, we obtain
\[
\frac{\nu_\rho(B_G(a,r))}{d_G(x,a)^{2W-\Xi}}
\prec
\int_{B_G(a,r)}
\frac{d\nu_\rho(y)}{d_G(x,y)^{2W-\Xi}}.
\]
Multiplying both sides by $h(a)r^{W-\Xi}$ and using
\eqref{eq:spike:volume-lower}, we get
\[
h(a)\left(\frac{r}{d_G(a,x)}\right)^{2W-\Xi}
\prec
h(a)r^{W-\Xi}
\int_{B_G(a,r)}
\frac{d\nu_\rho(y)}{d_G(x,y)^{2W-\Xi}}.
\]
On the other hand, since $d_G(a,x)\ge r$, we have
\[
\left(1+\frac{d_G(a,x)}{r}\right)^{-2W+\Xi}
\asymp
\left(\frac{d_G(a,x)}{r}\right)^{-2W+\Xi}
=
\left(\frac{r}{d_G(a,x)}\right)^{2W-\Xi}.
\]
Since \eqref{eq:h-bounds} implies the weaker bound with exponent
$-2W+\Xi$, the lemma follows.
\end{proof}

\begin{lemma}
\label{lemma:D_rh}
Suppose $(h,a,r)$ is a $Q$-spike. For each $x \in B_G(a,r)^c$ we have
\[
D_r h(x)\prec
\frac{h(a)}{r}\, r^{W-\Xi}
\int_{B_G(a,r)}
\frac{d\nu_\rho(y)}{d_G(x,y)^{2W-\Xi}}.
\]
\end{lemma}

\begin{proof}
Apply the integral comparison in the proof of Lemma~\ref{lemma:spike},
with $h(a)/r$ in place of $h(a)$, to the decay bound
\eqref{eq:D_r_bound_decay}.
\end{proof}

\begin{proof}[Proof of Proposition~\ref{prop:subfunction}]
We construct a sum greedily, filling the deficit at each center where
the current sum is below $f$. We bound the sum above on $X$ and below
on $Y$, then rescale it to control both the remainder and its variation.

Choose $\mm\ge2$ with $4\mm\ge2^{1+1/\sigma}$, large enough that
the two tail bounds in the proof of Lemma~\ref{lem:dleqrn} are at most
$f(y)/4$ each. Lemmas~\ref{lem:estimate:I0} and \ref{lem:estimate:I1}
give these bounds with implicit constants independent of $\mm$.
Thus this choice uses only the uniform data; keep $\mm$ fixed throughout.

Order the radii so that $r_i$ is nonincreasing. Let $g_0(x)=0$ and
$$g_n(x)=\begin{cases}
g_{n-1}(x) & \text{ if }\quad g_{n-1}(a_{n})\geq f(a_{n})\\
g_{n-1}(x) + \frac{h_{n}(x)}{h_n(a_n)}(f(a_{n})- g_{n-1}(a_{n})) & \text{ if }
\quad g_{n-1}(a_{n})< f(a_{n})\end{cases}$$

This construction yields a sequence $\{\lambda_i\}_{i=1}^N$ with
\[
g_n(x) =\sum_{i=1}^{n} \la_i h_{i}(x), \qquad\text{ where } 
\lambda_i \, h_i(a_i) = \max(0, f(a_i) - g_{i-1}(a_i)).
\]
In particular, 
\begin{equation} \label{eq:lambda}
\lambda_i h_i(a_i) \leq f(a_i).
\end{equation} 

Fix $y\in X$. We first separate centers beyond distance $\delta$ and
those between $\mm r_i$ and $\delta$. For the remaining centers, we
distinguish whether the ball $B_G(a_i,r_i)$ is farther than $r_n/2$
from $y$. This gives four index sets, whose dependence on $y$ we
suppress in the notation:
\begin{align*}
I_0(n,N) &= \big\{ i \in \{n, . . . , N\}  : d_G(y,a_i) > \delta \big\},\\
I_1(n,N) &= \big\{ i \in \{n, . . . , N\}  : \mm \, r_i < d_G(y,a_i) \le \delta \big\},\\
I_2(n,N) &= \big\{ i \in \{n, . . . , N\} 
    : d_G(y,a_i) \le \mm \, r_i  \text{ and } d_G\big(y,B_G(a_i,r_i)\big)>r_n/2\big\},\\
\intertext{and}
I_3(n,N) &= \big\{ i \in \{n, . . . , N\} 
    : d_G(y,a_i) \le \mm \, r_i  \text{ and } d_G\big(y,B_G(a_i,r_i)\big)\le r_n/2\big\}.
\end{align*}

\begin{lemma}
\label{lem:estimate:I0}
For any $y \in X$, 
\[
\sum_{i \in I_0(1,N)} \lambda_i h_i(y) \prec \frac{f(y)}{\mm^{W-\Xi}}
\qquad\text{and}\qquad
\sum_{i \in I_0(1,N)} \lambda_i D_{r_i} h_i(y) \prec \frac{f(y)}{r_N \, \mm^{W-\Xi}}.
\]
\end{lemma}

\begin{proof} 
From Lemma~\ref{lemma:spike}  and \eqref{eq:lambda} we have
\begin{equation}
\label{eq:I0:step1}
  \sum_{i\in I_0(1,N)} \la_i h_i(y)
\prec
\sum_{i\in I_0(1,N)}
f(a_i) r_i^{W-\Xi}
\int_{B_G(a_i,r_i)}
\frac{d\nu_\rho(x)}{d_G(y,x)^{2W-\Xi}}
\end{equation}
Note that for all $i$, $r_i \le r_1 \le  \delta/\mm \le \delta/2$.   
For $i\in I_0(1,N)$ and
$x\in B_G(a_i,r_i)$, we have $d_G(y,x)>\delta/2$ and
\begin{displaymath}
d_G(y,a_i) \le d_G(y,x)+r_i \le d_G(y,x)+\delta/2, 
\end{displaymath}
so
\begin{displaymath}
\left(1+ \frac{d_G(y,a_i)}{\delta}\right) \le \left(\frac{3}{2}  +
  \frac{d_G(y,x)}{\delta}\right) \le  \frac{4 d_G(y,x)}{\delta}. 
\end{displaymath}
and then, using
\eqref{eq:global:subfunction-ratio},
\begin{displaymath}
f(a_i) \le  4 f(y) \left(1+\frac{d_G(y,a_i)}{\delta}\right)^\sigma
\le 4^{1+\sigma} f(y) \delta^{-\sigma} d_G(y,x)^\sigma.  
\end{displaymath}
Substituting into \eqref{eq:I0:step1}, we get
\begin{displaymath}
\sum_{i\in I_0(1,N)} \la_i h_i(y) \prec \sum_{i \in I_0(1,N)} f(y)
r_i^{W-\Xi} \delta^{-\sigma} \int_{B_G(a_i,r_i)}
\frac{d\nu_\rho(x)}{d_G(y,x)^{2W-\Xi-\sigma}}
\end{displaymath}
Bounded overlap estimates the  right-hand-side as
\[
f(y) \sB  r_1^{W-\Xi}  \delta^{-\sigma} \int_{X\setminus B_G(y,\delta/2)}
\frac{d\nu_\rho(z)}{d_G(y,z)^{2W-\Xi-\sigma}}.
\]
Using \eqref{eq:Qspike:measure:decay}
applied with $W'=2W-\Xi-\sigma=W+\sigma>W$, this is
\[
\prec f(y)r_1^{W-\Xi}\delta^{-\sigma}
\delta^{-(W-\Xi-\sigma)}
=f(y)(r_1/\delta)^{W-\Xi} \le f(y)/\mm^{W-\Xi}.
\]
This proves the first estimate.

The proof of the second estimate is the same, using
Lemma~\ref{lemma:D_rh} instead of Lemma~\ref{lemma:spike}. The
extra factor $1/r_i$ is bounded by $1/r_N$, since $r_i\ge r_N$. Thus
\[
\sum_{i\in I_0(1,N)} \la_i D_{r_i}h_i(y)
\prec \frac{f(y)}{r_N\mm^{W-\Xi}}. 
\]
\end{proof}

\begin{lemma}
\label{lem:estimate:I1}
For any $y \in X$, 
\[
\sum_{i \in I_1(1,N)} \lambda_i h_i(y) \prec \frac{f(y)}{(\mm-1)^{W-\Xi}}
\qquad\text{and}\qquad
\sum_{i \in I_1(1,N)} \lambda_i D_{r_i} h_i(y) \prec \frac{f(y)}{r_N \, (\mm-1)^{W-\Xi}}.
\]
\end{lemma}

\begin{proof} 
Since $\lambda_i h_i(a_i)\le f(a_i)$, Lemma~\ref{lemma:spike} gives
\[
\sum_{i\in I_1(1,N)} \la_i h_i(y)
\prec
\sum_{i\in I_1(1,N)}
f(a_i)r_i^{W-\Xi}
\int_{B_G(a_i,r_i)}
\frac{d\nu_\rho(x)}{d_G(y,x)^{2W-\Xi}}.
\]
For $i\in I_1(1,N)$, we have $d_G(a_i,y)\le\delta$, and hence, by
\eqref{eq:local:subfunction-ratio},
\[
f(a_i)\le 2 f(y).
\]
Also $r_i\ge r_N$, and the balls $B_G(a_i,r_i)$ have multiplicity at most
$\sB$. Therefore
\[
\sum_{i\in I_1(1,N)} \la_i h_i(y)
\prec
f(y)r_1^{W-\Xi}
\int_{X-B_G(y,(\mm-1)r_N)}
\frac{d\nu_\rho(x)}{d_G(y,x)^{2W-\Xi}}.
\]
By \eqref{eq:Qspike:measure:decay} 
applied with $W'=2W-\Xi$, this is
\[
\prec
\frac{f(y)r_1^{W-\Xi}}{((\mm-1)r_N)^{W-\Xi}}
\le 
 \frac{\sW_0^{W-\Xi} f(y)}{(\mm-1)^{W-\Xi}}.
\]
This proves the first estimate.

The proof of the second estimate is identical, using
Lemma~\ref{lemma:D_rh} and the extra factor $1/r_i\le1/r_N$. Hence
\[
\sum_{i\in I_1(1,N)} \lambda_iD_{r_i}h_i(y)
\prec
\frac{\sW_0^{W-\Xi} f(y)}{r_N(\mm-1)^{W-\Xi}}.
\qedhere
\]
\end{proof}

\begin{lemma}
\label{lem:dleqrn} 
There is a uniform constant $\sC$ such that, for $y, y' \in X$ with $d_G(y,y') \le r_n$, 
we have
\begin{displaymath}
g_n(y)\leq \sC\, g_n(y')+ f(y)/2. 
\end{displaymath}
\end{lemma}

\begin{proof} 
By our choice of $\mm$, Lemmas~\ref{lem:estimate:I0} and
\ref{lem:estimate:I1}, applied to the first $n$ spikes, give
\[
\sum_{i \in I_0(1,n)} \lambda_i h_i(y) \le \frac{f(y)}{4}
\qquad\text{and}\qquad
\sum_{i \in I_1(1,n)} \lambda_i h_i(y) \le \frac{f(y)}{4}. 
\]
But, for $i \in I_2(1,n) \cup I_3(1,n)$, we have $d_G(a_i, y) \leq \mm
r_i$. Hence, $d_G(a_i,y') \le (\mm + 1) r_i$. Applying
\eqref{eq:h-h} with $x = y'$ and $r=r_i$, we have
\[
\sum_{i \in I_2(1,n)\cup I_3(1,n)} \lambda_i h_i(y) \prec 
\sum_{i \in I_2(1,n)\cup I_3(1,n)} \lambda_i h_i(y') \left( 1+ \frac{(\mm+1) \, r_i}{r_i} \right)^\Xi
\prec g_n(y')  \left( 2+ \mm \right)^\Xi. 
\]
Since $\mm$ is now fixed, the last factor is uniform. Combining these
estimates proves the lemma.
\end{proof}

\begin{lemma}
\label{lem:gnleqF} 
There is a uniform constant $\sM_1$ such that 
\[
g_N(y) \leq \sM_1 f(y)
\]
 for all $N\in\mathbb{N}$ and $y \in X$.
\end{lemma}

\begin{proof}  
Fix $y \in X$. If $g_{N}(y) < (3\sC + 0.5) f(y)$ then the lemma holds for this $y$ 
as long $\sM_1 \geq (3\sC+0.5)$. Otherwise, let $n$ be the smallest integer such that 
\[
g_{n}(y) \geq (3\sC + 0.5) f(y).
\] 
Then 
\begin{equation} \label{eq:g-n-1}
g_{n-1}(y) < (3\sC+0.5) f(y).
\end{equation}

By Lemma~\ref{lem:dleqrn}, applied to $g_n$, if
$d_G(y,y')\le r_n$, then
\[
g_n(y)\le \sC g_n(y')+\frac12 f(y).
\]
Since $g_n(y)\ge (3\sC+0.5)f(y)$, this implies
\[
g_n(y')\ge 3f(y).
\]
In particular, if $k>n$ and $d_G(y,a_k)\le r_n$, then
\[
g_{k-1}(a_k)\ge g_n(a_k)\ge 3f(y)\ge f(a_k),
\]
where the last inequality follows from \eqref{eq:local:subfunction-ratio}, since
$d_G(y,a_k)\le r_n\le\delta$. Hence, by construction
we have $\lambda_k=0$.

It remains to bound the contribution of indices $n\le i\le N$.
We use the four sets $I_j(n,N)$ defined above. Lemmas~\ref{lem:estimate:I0}
and \ref{lem:estimate:I1} give
\begin{equation} \label{eq:I0I1}
 \sum_{i \in I_0(n,N) \cup I_1(n,N)} \la_i h_{i}(y) \prec f(y). 
\end{equation}

We now estimate the sums over $I_2(n,N)$ and $I_3(n,N)$. 
Since $\lambda_i h_i(a_i)\le f(a_i)$, Lemma~\ref{lemma:spike} gives
\[
\sum_{i\in I_2(n,N)} \la_i h_i(y)
\prec \sum_{i\in I_2(n,N)} f(a_i)r_i^{W-\Xi}
  \int_{B_G(a_i,r_i)} \frac{d\nu_\rho(x)}{d_G(y,x)^{2W-\Xi}}.
\]
For $i\in I_2(n,N)$, we have
\[
d_G(a_i,y)\le \mm r_i\le \mm r_n\le \delta,
\]
so \eqref{eq:local:subfunction-ratio} gives
\[
f(a_i)\le2f(y) \le 3f(y).
\]
Also $r_i\le r_n$, and the balls have multiplicity at most $\sB$. Hence
\[
\sum_{i\in I_2(n,N)} \la_i h_i(y)
\prec f(y)r_n^{W-\Xi}
  \int_{X-B_G(y,r_n/2)} \frac{d\nu_\rho(x)}{d_G(y,x)^{2W-\Xi}}.
\]
By \eqref{eq:Qspike:measure:decay}, applied with $W'=2W-\Xi$,
\begin{equation}
\label{eq:I2}
\sum_{i\in I_2(n,N)} \la_i h_i(y)\prec f(y).
\end{equation}

Now we consider the sum over $I_3(n,N)$. If $i\in I_3(n,N)$ and $r_i<\frac{r_n}{\mm}$,
then
\[
d_G(y,a_i)\le \mm r_i<r_n.
\]
For $i>n$, this implies $\lambda_i=0$ by the observation above. Thus the only
indices in $I_3(n,N)$ which can contribute are contained in
\[
I_4(n,N)
= \left\{ n\le i\le N: r_i\ge \frac{r_n}{\mm}
\text{ and } 
d_G\big(y,B_G(a_i,r_i)\big)\le \frac{r_n}{2} \right\}.
\]
For $i\in I_4(n,N)$, the ball $B_G(a_i,r_i)$ is contained in a fixed
multiple of $B_G(y,r_n)$; say
\[
B_G(a_i,r_i)\subset B_G(y,3r_n).
\]
Using the bounded-overlap assumption and \eqref{eq:spike:volume-lower} and
\eqref{eq:Q:spike:volume:growth} for each $B_G(a_i,r_i)$, we get
\[
\# I_4(n,N)\left(\frac{r_n}{\mm}\right)^{W}
\prec \sum_{i\in I_4(n,N)}\nu_\rho(B_G(a_i,r_i))
\le \sB\,\nu_\rho\big(B_G(y,3r_n)\big)
\prec r_n^{W}.
\]
Hence
\[
\# I_4(n,N)\prec 1,
\]
with an implicit constant depending on $\mm$.

For $i\in I_3(n,N)$ with $\lambda_i\ne0$, we have $i\in I_4(n,N)$. Moreover,
by \eqref{eq:h-bounds} and $d_G(y,a_i)\le\mm r_i$,
\[
h_i(y)\prec h_i(a_i).
\]
Therefore
\[
\sum_{i\in I_3(n,N)}\la_i h_i(y)
\prec \sum_{i\in I_4(n,N)}\la_i h_i(a_i).
\]
For $i\in I_4(n,N)$, the defining distance bound gives
$d_G(y,a_i)\le r_n/2+r_i\le3r_n/2\le\delta$.
Thus \eqref{eq:lambda} and \eqref{eq:local:subfunction-ratio} give
\[
\la_i h_i(a_i)\le f(a_i)\le3f(y).
\]
Thus
\begin{equation} \label{eq:I3}
\sum_{i\in I_3(n,N)}\la_i h_i(y) \prec f(y).
\end{equation}

This argument works for every $y$.  Thus, there exist uniform
constants  $C_j$ such that
\begin{displaymath}
\sum_{i \in  I_j(n,N)} \lambda_i h_i(y) \le C_j f(y), \qquad 0 \le j
\le 3.
\end{displaymath}
Then,
\[
g_N(y)=g_{n-1}(y)+\sum_{i=n}^N\lambda_i h_i(y)
\le(3\sC+1/2+\sum_{j=0}^3 C_j)f(y)
\]
Thus the lemma holds for any $\sM_1$ larger than
$(3\sC+0.5)+ \sum C_j$.
\end{proof}

\begin{lemma}
\label{lem:ggreaterF}  
For every $y \in Y$,
\begin{equation}
\label{eq:ggreaterF}
g_N(y)\ge \frac{f(y)}{6\sC}.
\end{equation}
\end{lemma}

\begin{proof} 
Assume that $n$ is the smallest integer such that $y
\in B_G(a_{n}, r_{n})$. 
By Lemma~\ref{lem:dleqrn},
\[
g_n(y)\ge \frac1{\sC}\left(g_n(a_n)-\frac12 f(a_n)\right).
\]
and, by construction, $g_{n}(a_{n}) \geq f(a_{n})$. Also, as $d_G(y,
a_{n})\leq r_{n}\leq \delta$, we get, by
\eqref{eq:local:subfunction-ratio}, 
$f(y)\leq 2 f(a_{n})$.
Therefore,
$$g_n(y)\geq \frac{g_{n}(a_{n})}{2\sC} \geq
\frac{f(a_{n})}{2\sC} \geq \frac{f(y)}{6\sC}.$$
But $g_N(y) \geq g_n(y)$ hence \eqref{eq:ggreaterF} holds. 
\end{proof}

We first normalize the sum to ensure $g\le f/2$. A further rescaling
below will control its Lipschitz constant. Define
\[
\muG^i = \frac{\lambda_i}{2\sM_1}
\qquad\text{and}\qquad
g = \sum_{i=1}^N \muG^i \, h_i = \frac 1{2\sM_1} g_N. 
\]
For this normalization, set $\ep_0 = \frac{1}{12\sC \sM_1}$ and
$F=f-g$. Then \eqref{eq:ggreaterF} gives, for $y\in Y$,
\[
g(y) \geq \ep_0 f(y) \qquad \Longrightarrow \qquad F(y) \le
(1-\epsilon_0)  f(y)  
\]
and Lemma~\ref{lem:gnleqF} implies that, for $x\in X$,
\[
0 \le g(x)\le \frac12 f(x) \qquad  \Longrightarrow \qquad  \frac{f(x)}{2}
\le F(x) \le f(x). 
\]

It remains to show that $F$ satisfies
\eqref{eq:local:subfunction-ratio} and
\eqref{eq:global:subfunction-ratio}. To do this, we first estimate
$D_{r_{\min}}g$. Fix $y \in X$, and let
\[
g_{N,j}(x) = \sum_{i\in I_j(1,N)}\lambda_i h_i(x).
\]
In the sum above, the sets $I_j(1,N)$ are defined using the point $y$
(and not $x$). For $j=0$ or $j=1$, the second estimates in
Lemmas~\ref{lem:estimate:I0} and \ref{lem:estimate:I1} give
\[
D_{r_N} g_{N,j}(y)\prec \frac{f(y)}{r_N}.
\]
Also, for $j=2,3$, we have $d_G(a_i,y)\le \mm r_i$ for
$i\in I_j(1,N)$. Using \eqref{eq:D_r_bound_sharp}, we obtain
\begin{align*}
D_{r_N} g_{N,j}(y)
&\le \sum_{i\in I_j(1,N)}\la_i D_{r_i}h_i(y) 
\prec \sum_{i\in I_j(1,N)} \lambda_i\frac{h_i(y)}{r_i}
\left(1+\frac{d_G(a_i,y)}{r_i}\right)^\Xi \\
&\prec \sum_{i\in I_j(1,N)} \lambda_i\frac{h_i(y)}{r_i}
\le \frac1{r_N} \sum_{i\in I_j(1,N)}\lambda_i h_i(y) \\
&\le \frac{g_N(y)}{r_N}
\prec \frac{f(y)}{r_N},
\end{align*}
where the last inequality follows from Lemma~\ref{lem:gnleqF}.  Thus, 
\[
D_{r_{\min}}g(y)
=\frac{1}{2\sM_1}D_{r_N}g_N(y)
\le\frac{1}{2\sM_1}\sum_{j=0}^3D_{r_N}g_{N,j}(y)
\le  \frac{\sL f(y)}{r_{\min}}.
\]
where $\sL > 1$ is a uniform constant.

We now divide every coefficient $\muG^i$ and
the lower-bound constant $\ep_0$ by the uniform constant $\sL$
from the Lipschitz estimate.
This  preserves the  inequalities $0\le g\le f/2$ and
$g\ge\ep_0f$ on $Y$. We keep the notation $g,\ep_0$ for the
rescaled quantities, and define $F=f-g$ using this final $g$.
The final removal fraction is therefore
$\ep_0=(12\sC\sM_1\sL)^{-1}$. The constants were chosen in the order
$\mm,\sC,\sM_1,\sL$, and none depends on $f$, $N$, or $\delta$.

Finally, we verify the two comparison estimates for $F=f-g$.
Note that $g$ is  now normalized so that
\begin{displaymath}
D_{r_{\min}}g(x)\le \frac{f(x)}{r_{\min}}.
\end{displaymath}
Recall that $4\mm\ge2^{1+1/\sigma}$ and
$\delta'\le\delta/(4\mm)$. Since $f/2\le F\le f$, the function $F$
is positive, bounded, and bounded away from zero. The relative Lipschitz
hypothesis and the normalized estimate for $g$ give
\[
\begin{aligned}
D_{\delta'}F(x)
&\le D_{\delta'}f(x)+D_{\delta'}g(x)\\
&\le\left(\frac1\delta+\frac1{r_{\min}}\right)f(x)\\
&\le\frac{4F(x)}{r_{\min}}
=\frac{F(x)}{\delta'}.
\end{aligned}
\]
The same estimate applies to pairs at distance exactly $\delta'$, since
$\delta'<\min\{\delta,r_{\min}\}$.

For the polynomial comparison, write $d=d_G(x,y)$. If $d<\delta$, the
relative Lipschitz hypothesis on $f$ gives
\[
\frac{F(y)}{F(x)}
\le2\frac{f(y)}{f(x)}
\le2\left(1+\frac d\delta\right)
\le4
\le4\left(1+\frac d{\delta'}\right)^\sigma.
\]
If $d\ge\delta$, our choice of $\mm$ gives
\[
1+\frac d{\delta'}
\ge4\mm\frac d\delta
\ge2^{1/\sigma}\left(1+\frac d\delta\right).
\]
The polynomial comparison hypothesis on $f$ therefore implies
\[
\frac{F(y)}{F(x)}
\le2\frac{f(y)}{f(x)}
\le8\left(1+\frac d\delta\right)^\sigma
\le4\left(1+\frac d{\delta'}\right)^\sigma.
\]
Thus $F$ satisfies both hypotheses at scale $\delta'$.

\end{proof}

The next lemma explains why a global comparison such as
\eqref{eq:global:subfunction-ratio} is natural in this iteration.
Starting from $f_0\equiv1$, if the scales decrease geometrically and
each step decreases the function by at most half, the local bounds
imply a global power bound.
We do not need this observation in the proof of the main theorem.

\begin{lemma}
Let $A>1$, and suppose that positive scales
$\delta_0,\ldots,\delta_n$ satisfy
\[
\delta_j\le A^{-1}\delta_{j-1}\qquad (1\le j\le n).
\]
Suppose $f_0\equiv1$, and that the positive functions
$f_0,\ldots,f_n$ satisfy
\begin{equation}
\label{eq:remark:doubling}
\frac12 f_{j-1}\le f_j\le f_{j-1},
\end{equation}
and
\begin{equation}
\label{eq:remark:local:ratio}
D_{\delta_j}f_j(x)\le\frac{f_j(x)}{\delta_j}
\quad\text{for all $x\in X$ and $1\le j\le n$.}
\end{equation}
Then, setting $\alpha=\log 2/\log A$, we have
\begin{equation}
\label{eq:remark:global:ratio}
\frac{f_n(y)}{f_n(x)}
\le 4\left(1+\frac{d_G(x,y)}{\delta_n}\right)^\alpha
\qquad\text{for all $x,y\in X$.}
\end{equation}
In particular, if $A\ge2^{1/\sigma}$, then $f_n$ satisfies
\eqref{eq:global:subfunction-ratio} with $\delta=\delta_n$.
\end{lemma}

\begin{proof}
We record the following immediate consequence of
\eqref{eq:remark:local:ratio}:
\begin{equation}
\label{eq:remark:proof:ratio}
f_j(y) \le  2 f_j(x) \qquad\text{if $d_G(x,y) \le \delta_j$}.
\end{equation}
Put $d=d_G(x,y)$.

  For $d\le\delta_n$, \eqref{eq:remark:global:ratio} follows
directly from \eqref{eq:remark:proof:ratio}. If
\[
\delta_{n-k+1}<d\le\delta_{n-k}
\qquad (1\le k\le n),
\]
the estimate \eqref{eq:remark:proof:ratio} for
$f_{n-k}$ and repeated use of \eqref{eq:remark:doubling} give
\[
\frac{f_n(y)}{f_n(x)}
\le 2^k\frac{f_{n-k}(y)}{f_{n-k}(x)}
\le 2^{k+1}.
\]
Since $d>\delta_{n-k+1}\ge A^{k-1}\delta_n$, this is bounded by
$4(d/\delta_n)^\alpha$. Finally, if $d>\delta_0$, then by
\eqref{eq:remark:doubling}
\[
\frac{f_n(y)}{f_n(x)}\le2^n
\le\left(\frac d{\delta_n}\right)^\alpha.
\]
\end{proof}

Including \eqref{eq:global:subfunction-ratio} as one of the
assumptions of Proposition~\ref{prop:subfunction} allows us to
formulate Proposition~\ref{prop:subfunction} as a single inductive
step, instead of keeping track of the entire history of earlier
scales.   

\section{Construction of the measure \texorpdfstring{$\muG$}{mu\_Gamma}}
\label{sec:first_moment}
We apply Proposition~\ref{prop:subfunction} successively to the spikes
with translation lengths near $n\tau'$, starting with the constant
function $1$. Proposition~\ref{prop:extra:properties:gamma:n} ensures
that most points are covered at many scales, so the remainders tend to
zero exponentially in $L^1$. The resulting coefficients satisfy the
stationarity equation \eqref{eq:main:goal:spikes} and have an exponential
\Teich moment.

Let $\beta$, $W=w/\beta$ and $\Xi = 2 \eta/\beta$
be as in \S\ref{sec:spikes}, with the initial choice
\eqref{eq:spike-error-choice}, which already ensures $\Xi<W$.
Let $\rho_1$, $\tau$, $\tau_0'$ and $\sB$ be as in
Proposition~\ref{prop:properties:gamma:n} and
Proposition~\ref{prop:extra:properties:gamma:n}.  

Set $\sW_0=e^{\beta(\rho_1+2)}$, and choose $\mm$ and $\epsilon_0$ as
in Proposition~\ref{prop:subfunction}, decreasing $\epsilon_0$ to at most
$1/2$ if necessary. The uniform spike estimates, $\sW_0$, and $\sB$ are
independent of $\tau'$, so $\mm$ and $\epsilon_0$ are now fixed.
Choose $\tau'\in\tau\natls$ with
\[
\tau'>\max\left\{\tau'_0,\ \rho_1+2+\frac{\log(4\mm)}\beta,\
\sR_0+1,\ \sR_1+1\right\}.
\]
Here $\sR_0$ and $\sR_1$ are the thresholds in
Theorem~\ref{theorem:is:Q:spike}. The second bound will separate
successive spike radii enough to apply Proposition~\ref{prop:subfunction}.

For $j\ge0$, let
\[
R_j = j \tau'.  
\]
Choose $n_0\ge1$ so that Proposition~\ref{prop:properties:gamma:n} and
the bounded-overlap conclusion of
Proposition~\ref{prop:extra:properties:gamma:n}{\rm(a)} apply for every
$n\ge n_0$. Set $\Gamma_\alpha(R_n)=\varnothing$ for every $\alpha$
when $n<n_0$. For $n\ge1$, let $\Gamma_n\subset\Mod$ be the finite set
\[
\Gamma_n=\Gamma(R_n)=\bigcup_{\alpha\in\cA}\Gamma_\alpha(R_n),
\]
using the families from Proposition~\ref{prop:properties:gamma:n}
when $n\ge n_0$. Their translation-length windows are disjoint by
part~{\rm(b)} and $\tau'>\rho_1+2$, so the $\Gamma_n$ are pairwise disjoint.
For $\alpha\in\cA$, define
\[
Y_n^\alpha
= U_\alpha\cap
\bigcup_{\gamma\in\Gamma_\alpha(R_n)}B_G(\gamma),
\qquad\text{and set}\qquad
Y_n=\bigcup_{\alpha\in\cA}Y_n^\alpha.
\]
By Proposition~\ref{prop:extra:properties:gamma:n}{\rm (a)}, the balls
$\{B_G(\gamma):\gamma\in\Gamma_n\}$ cover each point at most $\sB$ times.

Let $\cZ_n$ denote the collection of spikes
\begin{displaymath}
\cZ_n = \{ (H(\gamma,\cdot), [F_+(\gamma)], e^{-\beta \ell(\gamma)}) \st
\gamma  \in \Gamma(R_n)  \}. 
\end{displaymath}
Theorem~\ref{theorem:is:Q:spike} and the choice of $\tau'$ give uniform
spike estimates for these families. Set
\[
\delta_n=\tfrac14 e^{-\beta(R_n+\rho_1+1)}\qquad(n\ge0),
\qquad \sigma=(W-\Xi)/2.
\]

We construct positive functions $f_n:\PMF\to\reals$ inductively, with
$2^{-n}\le f_n\le1$ and satisfying
\eqref{eq:local:subfunction-ratio} and
\eqref{eq:global:subfunction-ratio} at scale $\delta_n$. These bounds
supply the boundedness and positive lower bound required by
Proposition~\ref{prop:subfunction}; its constants do not depend on this
lower bound. Start with
\[
f_0(x)=1,
\]
which satisfies all these conditions at scale $\delta_0$.

Suppose $f_{n-1}$ has been constructed. If $\Gamma_n=\varnothing$,
set $g_n=0$ and $f_n=f_{n-1}$. Otherwise, the radii in $\cZ_n$ are
$e^{-\beta\ell(\gamma)}$, so Proposition~\ref{prop:properties:gamma:n}{\rm(b)}
gives
\begin{equation}
\label{eq:delta:bookkeeping}
4\delta_n = e^{-\beta (R_n+\rho_1+1)} \le r_{{\rm  min},n}  \le r_{{\rm max},n}
\le e^{-\beta(R_n-1)}. 
\end{equation}
In particular, $r_{{\rm max},n}/r_{{\rm min},n} \le \sW_0$. 
By the choice of $\tau'$, 
\begin{displaymath}
  r_{{\rm max},n} \le  e^{-\beta(R_n-1)} = 4 \delta_{n-1} e^{-\beta(
  \tau'-(\rho_1+2))} \le  \frac{\delta_{n-1}}{\mm}. 
\end{displaymath}
We can therefore apply
Proposition~\ref{prop:subfunction} with $Y=Y_n$, $f=f_{n-1}$, $\delta =
\delta_{n-1}$ and the collection $\cZ_n$. Denote its nonnegative
coefficients by $\muG(\gamma)$, $\gamma\in\Gamma_n$, and put
\[
g_n(x)=\sum_{\gamma\in\Gamma_n}\muG(\gamma)H_\rho(\gamma,x).
\]
Let $f_n=f_{n-1}-g_n$ be the remainder supplied by the proposition.
Since the $\Gamma_n$ are disjoint, the coefficients chosen at different
steps do not conflict.

The two comparison estimates hold at scale
$\delta_n'=r_{{\rm min},n}/4\ge\delta_n$, by
\eqref{eq:delta:bookkeeping}. They persist at the smaller scale $\delta_n$:
\begin{displaymath}
D_{\delta_n} f_n(x) \le D_{\delta_n'} f_n(x) \le
\frac{f_n(x)}{\delta_n'} \le \frac{f_n(x)}{\delta_n},
\end{displaymath}
and, for $x,y\in\PMF$,
\begin{displaymath}
\frac{f_n(y)}{f_n(x)}  \le 4 \left( 1+
  \frac{d_G(x,y)}{\delta_n'}\right)^\sigma \le 4 \left( 1+
  \frac{d_G(x,y)}{\delta_n}\right)^\sigma.
\end{displaymath}
If $\Gamma_n$ is empty, the same comparisons with $\delta_{n-1}$ in
place of $\delta_n'$ show that $f_n=f_{n-1}$ satisfies both estimates.
In either case,
\begin{displaymath}
\frac{f_{n-1}(x)}{2} \le f_n(x) \le f_{n-1}(x),
\end{displaymath}
which completes the induction and gives, for all $n\ge0$,
\begin{equation}
\label{eq:fn:bounds}
2^{-n}\le f_n(x)\le 1.
\end{equation}
Also, by Proposition~\ref{prop:subfunction}, for $x\in Y_n$,
\begin{equation}
\label{eq:fn:one:step}
f_n(x)\le(1-\ep_0)f_{n-1}(x).
\end{equation}

Each time $x$ belongs to $Y_j$, the remainder is multiplied by a factor
at most $1-\ep_0$. Thus we count these visits by setting
\[
k_n(x)
=
\#\{j\in[1,n]:x\in Y_j\}.
\]
It follows from \eqref{eq:fn:one:step} that
\begin{equation}
\label{eq:fn:multistep}
f_n(x)\le(1-\ep_0)^{k_n(x)}.
\end{equation}

At points with $k_n(x)>n/2$, the remainder is exponentially small.
We bound the measure of the complement using
Proposition~\ref{prop:extra:properties:gamma:n}{\rm(b)}.
For a fixed chart $\alpha\in\cA$, let
\[
k_n^\alpha(x)
=
\#\left\{
j\in[1,n]:
x\in\bigcup_{\gamma\in\Gamma_\alpha(R_j)}B_G(\gamma)
\right\}.
\]
The initial cutoff discards at most $n_0$ coverage times. Since
$\eta_1<1/2$, we have $n/2+n_0\le(1-\eta_1)n$ for all sufficiently
large $n$. Proposition~\ref{prop:extra:properties:gamma:n}{\rm(b)}, with
its threshold $(1-\eta_1)n$, therefore gives
\[
\nu_\rho\big(\{x\in U_\alpha:k_n^\alpha(x)\le n/2\}\big)\le \kappa^n,
\]
where $\kappa=\kappa(\eta_1)<1$ is the constant in that proposition.
Since the sets $U_\alpha$ cover $\PMF$, if $k_n(x)\le n/2$, then for every
chart $\alpha$ with $x\in U_\alpha$ we have $k_n^\alpha(x)\le n/2$. Hence
\[
\{x\in\PMF:k_n(x)\le n/2\}
\subset
\bigcup_{\alpha\in\cA}
\{x\in U_\alpha:k_n^\alpha(x)\le n/2\}.
\]
Therefore
\begin{equation}
\label{eq:sets:not:covered:new}
\nu_\rho\big(\{x\in\PMF:k_n(x)\le n/2\}\big)
\le
|\cA|\kappa^n.
\end{equation}
Split the integral according to whether $k_n(x)>n/2$. Using
\eqref{eq:fn:multistep} on this set and $f_n\le1$ on its complement,
\eqref{eq:sets:not:covered:new} gives, for sufficiently large $n$,
\[
\|f_n\|_{L^1(\nu_\rho)}
\le
(1-\ep_0)^{n/2}+|\cA|\kappa^n.
\]
Set $\kappa_1=\max\{\sqrt{1-\ep_0},\kappa\}<1$.
Choosing $C_1\ge1+|\cA|$ large enough to cover the finitely many initial
indices gives, for all $n\ge0$,
\begin{equation}
\label{eq:fn:L1:final}
\|f_n\|_{L^1(\nu_\rho)}\le C_1\kappa_1^n.
\end{equation}

\begin{proof}[Proof of Theorem~\ref{theorem:main}]
By construction, $f_n=f_0-\sum_{j=1}^n g_j$.
Since $f_0\equiv1$
and $f_n\to0$ in $L^1(\nu_\rho)$ by \eqref{eq:fn:L1:final}, we have
\[
\sum_{n=1}^\infty g_n(x)=f_0(x)=1
\]
in $L^1(\nu_\rho)$. Setting $\muG(\gamma)=0$ for
$\gamma\notin\bigcup_{n=1}^\infty\Gamma_n$,
we obtain
\[
\sum_{\gamma\in\Mod(S)}\muG(\gamma)H_\rho(\gamma,x)=1
\]
in $L^1(\nu_\rho)$. The coefficients are nonnegative, so Tonelli's
theorem and the normalization
\[
\int_{\PMF}H_\rho(\gamma,x)\,d\nu_\rho(x)=1
\]
give
\[
\sum_{\gamma\in\Mod(S)}\muG(\gamma)=1.
\]
Thus $\muG$ is a probability measure. Since
\[
d(\gamma_*\nu_\rho)(x)=H_\rho(\gamma,x)\,d\nu_\rho(x),
\]
the preceding identity is exactly the stationarity equation
\[
\nu_\rho=\muG*\nu_\rho.
\]

We now prove the exponential moment estimate. Since
$0\le g_n\le f_{n-1}$, the normalization of $H_\rho$ and
\eqref{eq:fn:L1:final} give
\[
\muG(\Gamma_n)=\|g_n\|_{L^1(\nu_\rho)}
\le\|f_{n-1}\|_{L^1(\nu_\rho)}
\le C_1\kappa_1^{n-1}.
\]
By Proposition~\ref{prop:properties:gamma:n}{\rm (a),(b)},
\[
d_\cT(P,\gamma P)\le n\tau'+C_0
\quad(\gamma\in\Gamma_n),
\qquad C_0=2\sD_0+\rho_1+1.
\]
Thus, for every $0<a< -\log\kappa_1/\tau'$,
\[
\begin{aligned}
\sum_{\gamma\in\Mod(S)}\muG(\gamma)e^{a d_\cT(P,\gamma P)}
&\le C_1e^{aC_0}\sum_{n\ge1}e^{an\tau'}\kappa_1^{n-1}\\
&=\frac{C_1e^{a(C_0+\tau')}}{1-\kappa_1e^{a\tau'}}<\infty.
\end{aligned}
\]
The same exponent works
at every base point $P'$, since
$d_\cT(P',\gamma P')\le d_\cT(P,\gamma P)+2d_\cT(P,P')$.

Finally, the semigroup generated by $\supp(\muG)$ is
non-elementary. It is enough to  show that $\supp(\muG)$ has two
independent pseudo-Anosov elements. We have proved equality with $1$
almost everywhere. The pointwise identities
$\sum_{n=1}^N g_n(x)=1-f_N(x)\le1$ also imply
$\sum_\gamma\muG(\gamma)H_\rho(\gamma,x)\le1$ at every $x$.
If all positive-weight pseudo-Anosovs had the same fixed pair
$a_+,a_-$, then \eqref{eq:H} would give, for each such $\gamma$,
\[
H_\rho(\gamma,a_+)+H_\rho(\gamma,a_-)
=e^{w\ell(\gamma)}+e^{-w\ell(\gamma)}>2.
\]
Summing with weights $\muG(\gamma)$ contradicts the preceding bounds
at $a_+$ and $a_-$ and $\sum_\gamma\muG(\gamma)=1$.
The fixed pairs of two pseudo-Anosovs
are either equal or disjoint, so failure of a common pair supplies
two independent elements already in the support.

Kaimanovich--Masur's theorem
\cite[Theorem~2.2.4]{kaimanovich:masur:rennes} now gives uniqueness of
the non-atomic stationary measure. The extension to punctured surfaces
is recalled in \cite[Remark~3.3]{gadre:singular}.
This completes the proof of
Theorem~\ref{theorem:main}.
\end{proof}

\appendix

\section{Large deviation for Horospheres}
\label{app:large:deviation}
For $q \in \POT$, define 
\[
V(q) = \max \left(1, \frac 1{\sys(q)}\right). 
\]
Along a \Teich geodesic, $V(q)$ changes at most exponentially fast. That is, 
for $t > 0$, we always have 
\[
V(q) \ e^{-t} \leq V(q(t)) \leq V(q)  \ e^t. 
\]
Let $\mu_q$ be the Avila--Gou\"ezel horospherical measure fixed in
Section~\ref{sec:volume-balls}, with density
\eqref{eq:AG-hor-volume}, and let $B(q,r)\subset\Hor(q)$ be the ball
of radius $r$ in the intrinsic $d_\AGY$-metric around $q$.

{
We prove that most initial points in a chart slice have many sampled
times at which the geodesic has good backward statistics and visits
a prescribed neighborhood of large measure. Horospherical estimates
for $V$ control excursions into the thin part. Tracking transfers
these estimates from random walks to geodesics in
Proposition~\ref{prop:good:statistics:final}. Mixing supplies the
returns in Corollary~\ref{cor:geodesic:returns}, and the two estimates
give Proposition~\ref{prop:appendix-final}.
\par}
We first recall a few statements from \cite{AG}. 

\begin{proposition}[Proposition 6.1 in \cite{AG}] \label{prop:AG} 
For  $\alpha \in  (0,1/4)$, there exists $C_\alpha > 0$ such that, 
for any $q \in \POT$ and any  $t\geq 0$,
\[
\frac{1}{\mu_q (B(q,\ooh))}
\int_{B(q,\ooh)}  V(y(t)) \, d\mu_q(y) \leq \tfrac 12 C_\alpha e^{-(1-2\alpha)t} \cdot V(q)+ \tfrac 12 C_\alpha.
\] 
Moreover, the function $\log V$ is $1$-Lipschitz for the $d_\AGY$--metric. 
\end{proposition}

For $\tau>0$  choose $\alpha = \alpha(\tau) \in (0,1/4)$ with
$\alpha(\tau) \to 0$ and  $\log
C_{\alpha(\tau)}/\tau \to 0$ as $\tau \to \infty$.

Let $V_0= e^{(1-2\alpha) \tau}$. Then, Proposition
\ref{prop:AG}  implies
\begin{equation} \label{eq:V_0}
\frac{1}{\mu_q (B(q,\ooh))}
\int_{B(q,\ooh)} \hspace{-.5cm} V(y(\tau)) \, d\mu_q(y) \leq 
\begin{cases} 
C_\alpha e^{-(1-2\alpha)\tau} \cdot V(q) & \text{if} \quad V(q) \geq V_0\\
C_\alpha & \text{if} \quad V(q) \leq V_0
\end{cases} 
\end{equation}

We recall the following Lemma (slightly restated) from \cite{AG}. 

\begin{lemma}[\cite{AG}, Corollary 5.4] \label{prop:Doubling}
There is a constant $0 <c_D<1$ such that, for every $q \in \POT$ and for $0 \leq 2r \leq \ooh$,
if $p \in \Hor(q)$ and $d_\AGY(p, q) \leq r$ then 
\[
\mu_q(B(p,r)) \geq c_D \cdot \mu_q(B(q,2r)). 
\]
\end{lemma}

\subsection{A random walk model}
\label{sec:subsec:first:random:walk}
We now define a random walk as follows. Let $\tau>0$ be a large number and $p_0 \in \POT$.  
For  $n \ge 1$, define
\[
\Omega_{n,p_0}= \Big\{ (q_0,q_1,\dots,q_n)  :
q_0 \in B(p_0,\ooh),\ q_{j+1}\in g_\tau B(q_j,\ooh),
\ 0\leq j \leq n-1 \Big\}. 
\]
We equip $\Omega_{n,p_0}$ with a Markov measure $\PP_{n,p_0}$ as follows. 
For $0\leq j\leq n$, let $\mu_j$ be shorthand for
$\mu_q$, where $q$ is any point of $\Hor(p_0(j\tau))$. This is well-defined
because the conditional measure $\mu_q$ depends only on the horosphere
containing $q$. In particular, $\mu_0=\mu_{p_0}$.

By \eqref{eq:Volume-Growth}, these measures satisfy, for
every measurable $E\subset\Hor(p_0(j\tau))$,
\[
\mu_{j+1}(g_\tau E)=e^{w\tau}\mu_j(E).
\]
Then, define $\PP_{n,p_0}$ to be the probability measure on 
$\Omega_{n,p_0}$ as follows.
The law of $q_0$ is the normalized restriction
$\mu_0\vert_{B(p_0,\ooh)}/\mu_0(B(p_0,\ooh))$.
For $0\leq j\leq n-1$, conditionally on $q_j$, the law of
$q_{j+1}$ is
\[
(g_\tau)_*\left(
\frac{\mu_j\vert_{B(q_j,\ooh)}}{\mu_j(B(q_j,\ooh))}
\right).
\]
Equivalently, one may temporarily choose $q_j'\in B(q_j,\ooh)$ according to
the normalized restriction of $\mu_j$ and then set $q_{j+1}=q_j'(\tau)$.

Define 
\[
f_0 \from \POT \to \reals_+
\qquad\text{by}\qquad 
f_0(q) =
\begin{cases}
1 & V(q) > 3 V_0 e^\tau.\\
0 & V(q) \leq 3 V_0 e^\tau.  
\end{cases}
\]
and 
\[
f_1 \from \POT \to \reals_+
\qquad\text{by}\qquad 
f_1(q) =
\begin{cases}
1 & V(q) > 9 V_0 e^{3\tau}.\\
0 & V(q) \leq 9 V_0 e^{3\tau}.  
\end{cases}
\]

For $\omega \in \Omega_{n,p_0}$ and $0 \leq j \leq n$, define
\[
\tV(\omega,j) = V(q_j) \cdot \exp \left(  \sum_{i=0}^{j-1}  \int_{0}^{\tau} f_0(q_i(t)) \, dt \right),
\]
where the sum is empty when $j=0$.

\begin{proposition}
\label{prop:random-walk}
For all sufficiently large
$\tau > 0$ there exists $\alpha_1 = \alpha_1(\tau)>0$ with
$\alpha_1 \to 0$ as 
$\tau \to \infty$ such that, for $n \geq 1$ and $p_0 \in \POT$, we have 
\[
\EE_{n,p_0}(\tV(\omega, n) ) \leq  e^{\alpha_1  \tau  \, n} \cdot V(p_0). 
\]
\end{proposition} 

\begin{proof}
Let $\cF_j$ be the $\sigma$-algebra generated by $q_0,\dots,q_j$. We claim,
for $0 \leq j<n$,
\begin{equation} \label{eq:inductive-step}
\EE_{n,p_0}(\tV(\omega, {j+1}) | \cF_j )  \leq C_\alpha \cdot e^{2\alpha \tau} \cdot \tV(\omega, j). 
\end{equation}
Note that 
\[
\EE_{n,p_0}(\tV(\omega, {j+1}) \mid \cF_j ) =
 \frac {\exp \left(  \sum_{i=0}^{j}
  \int_{0}^{\tau} f_0(q_i(t)) \, dt \right)}{\mu_j(B(q_j,\ooh))}
  \int_{B(q_j,\ooh)} V(q'(\tau))\, d\mu_j(q').
\]
Dividing both sides by $\tV(\omega,j)$ we get 
\begin{equation}  \label{eq:big}
\frac{\EE_{n,p_0}(\tV(\omega, {j+1}) \mid \cF_j )}{\tV(\omega,j)} =
 \frac {\exp \left(\int_{0}^{\tau} f_0(q_j(t)) \, dt \right)}
 {V(q_j) \cdot \mu_j(B(q_j,\ooh))}
  \int_{B(q_j,\ooh)} V(q'(\tau))\, d\mu_j(q').
\end{equation} 
Since $f_0 \leq 1$, 
\begin{equation} \label{eq:trivial-bound}
\exp \left( \int_0^\tau f_0(q_j(t)) \, dt \right) \leq e^\tau.
\end{equation} 
 
We consider two cases. First assume $V(q_j) \geq V_0$. Then \eqref{eq:V_0}, \eqref{eq:big}
and \eqref{eq:trivial-bound} imply
\begin{align*}
\frac{\EE_{n,p_0}(\tV(\omega, {j+1}) | \cF_j )}{\tV(\omega,j)} 
  &\leq e^\tau\cdot C_\alpha e^{-(1-2\alpha)\tau}
  = C_\alpha e^{2\alpha \tau}.
\end{align*} 
If $V(q_j) \leq V_0$, then for $0\leq t\leq\tau$ we have
$V(q_j(t))\leq e^\tau V_0<3e^\tau V_0$, so $f_0(q_j(t))=0$. Thus
\[
\exp \left( \int_0^\tau f_0(q_j(t)) \, dt \right) =1.
\] 
Then using \eqref{eq:V_0}, \eqref{eq:big} and the fact that $V(q_j) \geq 1$, we have 
\begin{align*}
\frac{\EE_{n,p_0}(\tV(\omega, {j+1}) | \cF_j )}{\tV(\omega,j)}  &\leq C_\alpha.  
\end{align*} 
This finishes the proof of \eqref{eq:inductive-step}. Since $q_0\in
B(p_0,\ooh)$ and $\log V$ is $1$--Lipschitz, we also have
\[
\EE_{n,p_0}(\tV(\omega,0)) = \EE_{n,p_0}(V(q_0))
\leq e^{\ooh} V(p_0).
\]
Iterating \eqref{eq:inductive-step} and using the tower property gives
\[
\EE_{n,p_0}(\tV(\omega,n))
\leq C_\alpha^n e^{2\alpha n\tau}
\EE_{n,p_0}(V(q_0))
\leq e^{\ooh} C_\alpha^n e^{2\alpha n\tau}V(p_0).
\]
We now choose $\alpha_1$ such that 
\[
e^{\ooh} C_\alpha^n e^{2\alpha n \tau}
\leq e^{\alpha_1 n \tau}
\]
and $\alpha_1 \to 0$ as $\tau \to \infty$.
(After enlarging $C_\alpha$ so that $C_\alpha\geq1$, one may take
$\alpha_1(\tau)=2\alpha(\tau)+\bigl(\log C_{\alpha(\tau)}+\ooh\bigr)/\tau$.)
This finishes the proof. 
\end{proof}

For every quadratic differential $p \in B(p_0,\ooth)$, 
define
\[
{\widetilde N}_{\ooth, p_0}(p) = \Big\{ \omega \in \Omega_{n,p_0}  :
   d_{\AGY}( p(j \tau) , q_j) < \ooth, \quad 0\leq j \leq n \Big\}. 
\]
That is, ${\widetilde N}_{\ooth, p_0}(p)$ is the set of random walks that start near $p_0$ and 
fellow travel the \Teich geodesic defined by $p$.

Recall that $c_D$ is defined in Lemma~\ref{prop:Doubling} (and is
independent of $\tau$ or $p_0$). 
\begin{lemma}
\label{lemma:thick-tracking}
For any $p_0 \in \POT$ and 
$p \in B(p_0,\ooth)$, we have 
\[
\PP_{n,p_0}( {\widetilde N}_{\ooth, p_0}(p)) \cdot \mu_0(B(p_0,\ooh)) 
  \geq c_D^{n} e^{-w n \tau} \cdot \mu_{n} (B(p(n\tau),\ooth)). 
\]
\end{lemma}
 
\begin{proof} 
Define 
\[
A_j = \mu_{j} (B(p(j\tau),\ooth)). 
\]
Let $T_j(p)$ be the event that $q_i\in B(p(i\tau),\ooth)$ for
$0\leq i\leq j$. On $T_j(p)$, backward nonexpansion gives
\[
g_{-\tau} B(p((j+1)\tau),\ooth) \subset B(q_j,\ooh).
\]
Thus on $T_j(p)$ we have
\begin{align*}
\PP_{n,p_0} \Big( q_{j+1} \in  B(p((j+1)\tau),\ooth)  \mid \cF_j  \Big) 
  & = \frac{\mu_j (g_{-\tau} B(p((j+1)\tau),\ooth))}{\mu_j(B(q_j,\ooh))} \\
  &= \frac{e^{-w \tau} A_{j+1}}{\mu_j(B(q_j,\ooh))} \geq c_D\cdot e^{-w \tau} \cdot \frac{A_{j+1}}{A_j}.
\end{align*}
where $c_D$ is as in Lemma \ref{prop:Doubling}.
Note that $T_n(p) = {\widetilde N}_{\ooth, p_0}(p)$.
Since
\begin{displaymath}
\PP_{n,p_0}(T_{j+1}(p))=\EE_{n,p_0}[\mathbf 1_{T_j(p)}\PP_{n,p_0}(q_{j+1}\in
B(p((j+1)\tau),\ooth)\mid\cF_j)], 
\end{displaymath}
we have
\begin{displaymath}
\PP(T_{j+1}(p))\geq  c_D\cdot e^{-w \tau} \cdot \frac{A_{j+1}}{A_j} \PP(T_j(p)).
\end{displaymath}
Thus,  
\begin{align*}
\PP_{n,p_0}&({\widetilde N}_{\ooth, p_0}(p)) \geq  
  \PP_{n,p_0}\left(T_0(p) \right) \prod_{j=0}^{n-1}  c_D \cdot e^{-w \tau} \frac{A_{j+1}}{A_j} \\
  &\geq \frac{A_0}{\mu_0(B(p_0,\ooh))} \prod_{j=0}^{n-1}  c_D \cdot e^{-w \tau} \frac{A_{j+1}}{A_j} 
   =  c_D^n \cdot e^{-w n \tau} \frac{A_n}{\mu_0(B(p_0,\ooh))}.
\end{align*}
This finishes the proof. 
\end{proof}

For $p  \in \Hor(p_0)$, let
\begin{equation}
\label{eq:def:tV1}
\tV_1(p,n) = V(p(n\tau)) \exp\left( \int_0^{n\tau} f_1(p(t)) \,
  dt \right). 
\end{equation}

\begin{lemma}
\label{lemma:f0:vs:f1}
For any $p  \in B(p_0,\ooth)$, let $E_{p,n} =
g_{-n\tau}(B(p(n\tau),\ooh))$.  Then, 
\begin{displaymath}
\frac{1}{\mu_0(E_{p,n})} \int_{E_{p,n}} \tV_1(q,n) \, d\mu_0(q)  \le
3\,\frac{1}{\PP_{n,p_0}({\widetilde N}_{\ooth, p_0}(p))} \int_{{\widetilde
  N}_{\ooth, p_0}(p)} \tV(\omega,n) \, d\PP_{n,p_0}(\omega). 
\end{displaymath}
\end{lemma}

\begin{proof}
Fix $q\in E_{p,n}$ and
$\omega\in{\widetilde N}_{\ooth,p_0}(p)$.  Since
$d_\AGY(q(n\tau),p(n\tau))\leq\ooh$, the monotonicity
\eqref{eq:Monotone} of the AGY metric along a horosphere gives, for
$0\leq j\leq n$,
\[
d_\AGY(q(j\tau),p(j\tau))\leq\ooh.
\]
Consequently,
\[
d_\AGY(q(j\tau),q_j)\leq\ooh+\ooth=\frac{3}{200}.
\]
Let $t\in[j\tau,(j+1)\tau]$.  If $f_1(q(t))=1$, then
\[
V(q(t))\geq9V_0e^{3\tau}.
\]
Using the exponential change of $V$ along the flow, the
$1$--Lipschitz property of $\log V$, and $e^{-3/200}>1/3$, we obtain
\begin{align*}
V(q(j\tau))
&\geq e^{-\tau}V(q(t))\geq9V_0e^{2\tau},\\
V(q_j)
&\geq e^{-3/200}V(q(j\tau))\geq3V_0e^{2\tau},\\
V(q_j(t-j\tau))
&\geq e^{-\tau}V(q_j)\geq3V_0e^\tau.
\end{align*}
Thus $f_0(q_j(t-j\tau))=1$, and hence
\[
\int_0^{n\tau}f_1(q(t))\,dt
\leq\sum_{j=0}^{n-1}\int_0^\tau f_0(q_j(s))\,ds.
\]
At time $n\tau$, the triangle inequality and the
$1$--Lipschitz property of $\log V$ also give
\[
V(q(n\tau))\leq e^{3/200}V(q_n)<3V(q_n).
\]
It follows that $\tV_1(q,n)\leq3\tV(\omega,n)$.
Since this holds for every pair
$(q,\omega)\in E_{p,n}\times{\widetilde N}_{\ooth,p_0}(p)$, averaging over
the two factors proves the claim.
\end{proof}

\begin{lemma}
\label{lemma:geodesic:expectation}
There exists $\alpha_2=\alpha_2(\tau)  > 0$ with $\alpha_2 \to 0$ as
$\tau \to  \infty$ such  that
\begin{displaymath}
\frac{1}{\mu_0(B(p_0,\ooth))} \int_{B(p_0,\ooth)} \tV_1(q,n) \, d\mu_0(q)  \le
 e^{\alpha_2 n \tau} \cdot V(p_0). 
\end{displaymath}
\end{lemma}

\begin{proof}
Consider the set $Z= g_{n \tau} B(p_0,\ooth) \subset \Hor(p_0( n
\tau))$. Let  $X'$ be a maximal
$\ooh$-separated subset of $Z$ (in the AGY metric), and let $X =
g_{-n\tau} X'$. Then,  
\begin{itemize} 
\item if $q \in B(p_0,\ooth)$ then $q \in E_{p,n}$ for some $p \in X$ where $E_{p,n}$ is as in Lemma~\ref{lemma:f0:vs:f1}.  
\item for distinct $p, p' \in X$, $B(p(n\tau),\ooth) \cap B(p'(n \tau),\ooth) =
  \emptyset$.
\end{itemize}
Then,
\begin{equation}
\label{eq:sum:over:Epn}
\int_{B(p_0,\ooth)} \tV_1(q,n) \, d\mu_0(q) \le \sum_{p \in X} \int_{E_{p,n}} \tV_1(q,n) \, d\mu_0(q), 
\end{equation}
By Lemma~\ref{lemma:f0:vs:f1}, for every $p \in X$, 
\begin{displaymath}
  \int_{E_{p,n}} \tV_1(q,n) \, d\mu_0(q) \le  3\,\frac{\mu_0(E_{p,n})}{\PP_{n,p_0}({\widetilde N}_{\ooth, p_0}(p))} \int_{{\widetilde
  N}_{\ooth, p_0}(p)} \tV(\omega,n) \, d\PP_{n,p_0}(\omega).
\end{displaymath}
Note that by Lemma~\ref{prop:Doubling},
\begin{displaymath}
\mu_0(E_{p,n}) = e^{-wn\tau}\mu_n(B(p(n\tau),\ooh)) \le c_D^{-1}
e^{-wn\tau}  \mu_n(B(p(n\tau),\ooth)).
\end{displaymath}
Then, by Lemma~\ref{lemma:thick-tracking} we get
\begin{displaymath}
\frac{\mu_0(E_{p,n})}{\PP_{n,p_0}( {\widetilde N}_{\ooth, p_0}(p))} \le
c_D^{-(n+2)} \mu_0(B(p_0,\ooth)). 
\end{displaymath}
Thus, for every $p \in X$,
\begin{displaymath}
\frac{1}{\mu_0(B(p_0,\ooth))}  \int_{E_{p,n}} \tV_1(q,n) \, d\mu_0(q)
\le  3\,c_D^{-(n+2)} \int_{{\widetilde
  N}_{\ooth, p_0}(p)} \tV(\omega,n) \, d\PP_{n,p_0}(\omega).
\end{displaymath}
A walk tracking $p$ ends in $B(p(n\tau),\ooth)$, so the tracking
events for distinct $p\in X$ are disjoint.
Substituting into (\ref{eq:sum:over:Epn}), we get
\begin{displaymath}
\frac{1}{\mu_0(B(p_0,\ooth))} \int_{B(p_0,\ooth)} \tV_1(q,n) \,
 d\mu_0(q) \le 3\,c_D^{-(n+2)} \EE_{n,p_0}(\tV(\omega,n)). 
\end{displaymath}
Now, the statement follows from Proposition~\ref{prop:random-walk}.
(For example, one may take
$\alpha_2=\alpha_1+(\log 3+3\log(c_D^{-1}))/\tau$).  
\end{proof}

We record whether some interval $[j\tau,n\tau]$, with $0\leq j<n$
an integer, spends more than a proportion $\rho$ of its time in
$\{f_1=1\}$.
For $0<\rho<1$, $q\in\POT$, and $n\in\natls$, let
\begin{displaymath}
H_\rho(q,n)  = \max_{0 \le  j \le  n} \left\{\int_{j\tau}^{n\tau} f_1(q(t))
\, dt - \rho(n-j) \tau \right\}.
\end{displaymath}
Let
\begin{displaymath}
Y_\rho(q,n) = \begin{cases} 1 & \text{ if $H_\rho(q,n) > 0$}  \\
                                0 &  \text  {otherwise}. 
\end{cases}
\end{displaymath}

\begin{proposition}
\label{prop:short:range:good:statistics}
For every $0 < \rho < 1$ there  exists  $\tau_0 = \tau_0(\rho)$ and 
$\epsilon_\rho = \epsilon_\rho(\tau) \to 0$ as $\tau \to \infty$
such that the following holds.  For every  $\tau > \tau_0$, every
$p_0 \in \POT$ and every $n \ge 1$, 
\begin{displaymath}
\frac{1}{\mu_{p_0}(B(p_0,\ooth))} \int_{B(p_0,\ooth)} V(q(n\tau)) Y_\rho(q,n)  \,
d\mu_{p_0}(q)  \le e^{-(\rho - \epsilon_\rho(\tau))\tau} \cdot V(p_0). 
\end{displaymath}
\end{proposition}

\begin{proof}
{
Put $r=\ooth$ and $\gamma=1-2\alpha$, and write $\nu_t$ for the
conditional measure on the horosphere containing $p_0(t)$.  We first
claim that there is a uniform constant $C_\ast>0$ such that, for
$1\leq \ell\leq n$,
\begin{align}
&\frac{1}{\nu_0(B(p_0,r))}
\int_{B(p_0,r)} V(q(n\tau))
\exp\left(\int_{(n-\ell)\tau}^{n\tau}f_1(q(t))\,dt\right)
\,d\nu_0(q) \notag\\
&\qquad\leq
C_\ast C_\alpha e^{\alpha_2\ell\tau}
\left(e^{-\gamma(n-\ell)\tau}V(p_0)+1\right).
\label{eq:terminal:suffix:estimate}
\end{align}

Set $s=(n-\ell)\tau$ and $E=g_sB(p_0,r)$.  After changing variables
by $g_s$, the left hand side of
\eqref{eq:terminal:suffix:estimate} is
\begin{displaymath}
\frac{1}{\nu_s(E)}\int_E F_\ell(z)\,d\nu_s(z),
\qquad
F_\ell(z)=V(z(\ell\tau))
\exp\left(\int_0^{\ell\tau}f_1(z(t))\,dt\right).
\end{displaymath}
Choose a maximal $r$--separated set $\{x_i\}\subset E$.  The balls
$B(x_i,r)$ cover $E$, while the balls $B(x_i,r/2)$ are pairwise
disjoint.  Lemma~\ref{lemma:geodesic:expectation}, applied with base
point $x_i$, gives
\begin{displaymath}
\int_{B(x_i,r)}F_\ell(z)\,d\nu_s(z)
\leq e^{\alpha_2\ell\tau}V(x_i)\nu_s(B(x_i,r)).
\end{displaymath}
By Lemma~\ref{prop:Doubling} and the $1$--Lipschitz property of
$\log V$,
\begin{displaymath}
V(x_i)\nu_s(B(x_i,r))
\leq C_\ast\int_{B(x_i,r/2)}V(z)\,d\nu_s(z).
\end{displaymath}
Consequently,
\begin{displaymath}
\int_EF_\ell(z)\,d\nu_s(z)
\leq C_\ast e^{\alpha_2\ell\tau}
\int_{E^+}V(z)\,d\nu_s(z),
\qquad
E^+=\bigcup_iB(x_i,r/2).
\end{displaymath}
If $z\in E^+$, then the monotonicity of the AGY metric along a
horosphere gives
\begin{displaymath}
z(-s)\in B(p_0,3r/2)\subset B(p_0,2r)=B(p_0,\ooh).
\end{displaymath}
It follows from the horospherical volume-scaling formula and
Proposition~\ref{prop:AG} that
\begin{align*}
\int_{E^+}V(z)\,d\nu_s(z)
&\leq e^{ws}\int_{B(p_0,\ooh)}V(q(s))\,d\nu_0(q)\\
&\leq \tfrac12 C_\alpha e^{ws}\nu_0(B(p_0,\ooh))
\left(e^{-\gamma s}V(p_0)+1\right).
\end{align*}
Since
\begin{displaymath}
\nu_s(E)=e^{ws}\nu_0(B(p_0,r))
\end{displaymath}
and $\nu_0(B(p_0,\ooh))/\nu_0(B(p_0,r))$ is uniformly bounded by
Lemma~\ref{prop:Doubling}, this proves
\eqref{eq:terminal:suffix:estimate}.

If $Y_\rho(q,n)=1$, then for some $1\leq\ell\leq n$,
\begin{displaymath}
\int_{(n-\ell)\tau}^{n\tau}f_1(q(t))\,dt>\rho\ell\tau.
\end{displaymath}
Therefore,
\begin{displaymath}
Y_\rho(q,n)
\leq\sum_{\ell=1}^n e^{-\rho\ell\tau}
\exp\left(\int_{(n-\ell)\tau}^{n\tau}f_1(q(t))\,dt\right).
\end{displaymath}
Multiplying by $V(q(n\tau))$, integrating, and applying
\eqref{eq:terminal:suffix:estimate}, we obtain
\begin{align*}
&\frac{1}{\nu_0(B(p_0,r))}
\int_{B(p_0,r)}V(q(n\tau))Y_\rho(q,n)\,d\nu_0(q)\\
&\qquad\leq C_\ast C_\alpha
\sum_{\ell=1}^n e^{-(\rho-\alpha_2)\ell\tau}
\left(e^{-\gamma(n-\ell)\tau}V(p_0)+1\right)\\
&\qquad\leq
\frac{2C_\ast C_\alpha}{1-e^{-(\rho-\alpha_2)\tau}}
e^{-(\rho-\alpha_2)\tau}V(p_0),
\end{align*}
where in the last inequality we used $V(p_0)\geq1$ and took $\tau$
large enough that $\alpha_2(\tau)<\rho$.

Define
\begin{displaymath}
\epsilon_\rho(\tau)=\alpha_2(\tau)+\frac{1}{\tau}
\log\left(\frac{2C_\ast C_{\alpha(\tau)}}
{1-e^{-(\rho-\alpha_2(\tau))\tau}}\right).
\end{displaymath}
Then the preceding upper bound is
\begin{displaymath}
e^{-(\rho-\epsilon_\rho(\tau))\tau}V(p_0).
\end{displaymath}
Finally, $\alpha_2(\tau)\to0$,
$\log C_{\alpha(\tau)}/\tau\to0$, and
$e^{-(\rho-\alpha_2(\tau))\tau}\to0$.  Hence
$\epsilon_\rho(\tau)\to0$, as required.
}
\end{proof}

\subsection{An alternative random walk model}

The next walk has prescribed distributions on the horospheres, so
probability bounds will give horospherical measure bounds. Its
transition densities are Lipschitz, allowing us to apply mixing.

In this section, we let $\tau' = m \tau$, where $\tau$ is as in
\S\ref{sec:subsec:first:random:walk} and  $m \in \natls$ is a large
parameter.
{Let $T=T(K_0)$ be the constant in
Corollary~\ref{cor:AGY:endpoint:contraction}, and fix once and for all a
constant}
\[
{L>800T(w+1).}
\]

Let $U(0)$ denote $U_\alpha(0)$ where $\alpha \in
\cA$. Let $\Hor_0 = \Hor_\alpha$
be the horosphere containing $U(0)$, and let $\Hor_n =
\Hor_0(n \tau')$ be the horosphere containing $U(n\tau')$. For $q
\in \Hor_0$,  define
\begin{displaymath}
\lambda_0(q) = {e^{-L d_\AGY(q,U(0))}},
\end{displaymath}
and for $q \in \Hor_n$, let
\begin{displaymath}
\lambda_n(q) = \lambda_0(q(-n\tau')). 
\end{displaymath} 
Then, by \autoref{theorem:properties:AGY}, $\log
\lambda_n$ is {$L$-Lipschitz}, i.e.\ for any $p$, $q$ in $\Hor_n$, 
\begin{equation}
\label{eq:loglipschitz}
  {e^{-L d_\AGY(p,q)}} \le  \frac{\lambda_n(p)}{\lambda_n(q)} \le {e^{L d_\AGY(p,q)}}. 
\end{equation}
We use the abbreviation $\bar{\mu}_n$ to denote the
conditional measure $\mu_{p_n}$ where $p_n$ is any point of
$\Hor_n$.  (In particular, if $p_0 \in \Hor_0$, $\bar{\mu}_0 = \mu_0$). 
We then  define
\begin{displaymath}
\mu_n' =  \lambda_n \bar{\mu}_n = \lambda_n \mu_{p_n}. 
\end{displaymath}
Then, for $E \subset \Hor_0$, 
\begin{equation}
\label{eq:appendix:Volume-Growth}
\mu_n'(g_{n\tau'}(E))
=
e^{wn\tau'}\mu_0'(E).
\end{equation}
The weight $\lambda_0$ equals $1$ on $U(0)$ and decays away from it.
The next lemma proves that $\mu'_0$ is finite; then
\eqref{eq:appendix:Volume-Growth} gives finiteness of every $\mu'_n$.
These measures are not normalized.

\begin{lemma}
\label{lemma:V0:finite:exp}
We have
\begin{displaymath}
\int_{\Hor_0} V(q) \, d\mu_0'(q) <  \infty.
\end{displaymath}
\end{lemma}

\begin{proof}
{
Recall that $q_\alpha\in U(0)$ and that
$q_\alpha(t)\in K_0$ for every $t\in\reals$.  Choose $h>0$ sufficiently
small that $e^{2h}\leq2$, and put
\[
r_-:=\ooth,
\qquad
r_+:=e^{2h}r_-\leq\ooh.
\]
For $k\geq0$, define a shell in the horosphere through $q_\alpha(-kh)$ by
\[
S_k=\left\{z\in\Hor(q_\alpha(-kh)):
r_-\leq d_\AGY(z,q_\alpha(-kh))<r_+\right\},
\qquad
\mathcal E_k=g_{kh}S_k\subset\Hor_0.
\]

We first claim that
\[
\Hor_0\subset B(q_\alpha,r_-)\cup\bigcup_{k\geq0}\mathcal E_k.
\]
Indeed, for $q\in\Hor_0$, let
\[
D_q(t)=d_\AGY(q(-t),q_\alpha(-t)),
\qquad t\geq0.
\]
By \eqref{eq:Monotone}, $D_q$ is nonincreasing.  We use the proof of
Corollary~\ref{cor:AGY:endpoint:contraction} with initial path segments
of length $r_-=\ooth$ in place of $\ooh$.  More explicitly, if
$D_q(t)\geq r_-$, then $D_q(s)\geq r_-$ for $0\leq s\leq t$.
On each of the $\lfloor t/T\rfloor$ full blocks of length $T$,
Corollary~\ref{cor:AG43} halves the distance between the endpoints of
the initial segment of length $r_-$ of a close-to-minimizing path,
while \eqref{eq:Monotone} does not increase the length of the rest.
Here the base orbit stays in $K_0$ throughout every block.  Letting
the path-length error tend to zero, each block therefore decreases
the distance by at least $r_-/2=1/400$.  By monotonicity on the
remaining time interval, for $t\geq2T$ we obtain
\[
D_q(t)\leq D_q(0)-\frac{1}{400}\left\lfloor\frac{t}{T}\right\rfloor
\leq D_q(0)-\frac{t}{800T}
\qquad\text{if }D_q(t)\geq r_-.
\]
Consequently $D_q(t)$ must eventually be smaller than $r_-$;
otherwise the displayed upper bound would become negative.
If $D_q(0)\geq r_-$, let $k\geq1$ be the first integer
such that $D_q(kh)<r_-$.  Then $D_q((k-1)h)\geq r_-$, while
\eqref{eq:Monotone} gives
\[
D_q((k-1)h)\leq e^{2h}D_q(kh)<r_+.
\]
Thus $q\in\mathcal E_{k-1}$, proving the claim.

Lemma~\ref{lem:AGY-volume} and the choice $r_+\leq\ooh$ give
\[
\mu_{q_\alpha(-kh)}(S_k)\leq \sV_0.
\]
The horospherical volume-scaling formula therefore implies
\[
\bar\mu_0(\mathcal E_k)
=e^{wkh}\mu_{q_\alpha(-kh)}(S_k)
\leq \sV_0e^{wkh}.
\]
Since $q_\alpha(-kh)\in K_0$, the function $V$ is uniformly bounded along
the orbit of $q_\alpha$.  If $q=z(kh)\in\mathcal E_k$, the Lipschitz property of
$\log V$ and the exponential change of $V$ along the flow give
\[
V(q)\leq e^{kh}V(z)\leq C e^{kh}.
\]

For $q\in\mathcal E_k$ we have $D_q(kh)\geq r_-$.  Thus, when
$kh\geq2T$, the contraction estimate just proved yields
\[
d_\AGY(q,q_\alpha)
\geq d_\AGY(q(-kh),q_\alpha(-kh))+\frac{kh}{800T}
\geq\frac{kh}{800T}.
\]
Since $U(0)\subset W_\alpha(0)$ and
\eqref{eq:dagy:x:F:alpha:y:F:alpha} gives
$\sup_{u\in U(0)}d_\AGY(u,q_\alpha)\leq\ooth$, it follows that
\[
\lambda_0(q)
\leq e^{L\ooth}\exp\left(-\frac{Lkh}{800T}\right)
\qquad(q\in\mathcal E_k,\ kh\geq2T).
\]
Combining the preceding three estimates, we obtain
\[
\int_{\mathcal E_k}V(q)\,d\mu_0'(q)
=\int_{\mathcal E_k}V(q)\lambda_0(q)\,d\bar\mu_0(q)
\leq C\exp\left(\left(w+1-\frac{L}{800T}\right)kh\right).
\]
The finitely many indices for which $kh<2T$ contribute a finite amount by
the volume and $V$ estimates above.  The remaining series
over $k$ converges by the choice of $L$.  Finally, the integral
over $B(q_\alpha,r_-)$ is finite by Lemma~\ref{lem:AGY-volume}, the
$1$--Lipschitz property of $\log V$, and the inequality $0<\lambda_0\leq1$.
The covering claim now proves the lemma.
}
\end{proof}


{
Let $0<\delta\leq1/400$. We now construct a partition of unity
on each $\Hor_n$.
\par}
For each non-negative integer $n$, we
choose a maximal $\delta/8$-separated subset $\Delta_n$ of $\Hor_n$.  We
denote this set as $\{ x_\sigma\}_{\sigma \in \Delta_n}$. 
We now claim that
\begin{displaymath}
\bigcup_{\sigma \in \Delta_n} B(x_\sigma,\delta/8) = \Hor_n. 
\end{displaymath}
Indeed, if $y \in \Hor_n$, then there exists
$x_\sigma \in \Delta_n$ with $d_\AGY(x_\sigma,y) < \delta/8$  (or else
we would add $y$ to $\Delta_n$). 
Thus, $\{B(x_\sigma,\delta/4)\}_{\sigma \in \Delta_n}$ is an open cover
of $\Hor_n$. Let $M$ be the multiplicity of this cover, which is
bounded by a universal constant in view of Lemma~\ref{prop:Doubling}. 
(Indeed for all centers whose $\delta/4$--balls contain a fixed $p$, the
$\delta/16$--balls are disjoint and contained in $B(p,\delta/2)$.  Three
applications of Lemma~\ref{prop:Doubling} give
$\bar\mu_n(B(x_\sigma,\delta/16))\geq
c_D^3\bar\mu_n(B(p,\delta/2))$, so there are at most $c_D^{-3}$ such
centers.)
We now choose a partition of unity $\phi_\sigma$
subordinate to this cover as follows. For each $\sigma  \in \Delta_n$,
choose a smooth $\psi_\sigma: \Hor_n \to [0,1]$ with
\begin{displaymath}
  \psi_\sigma = 1 \text{ on } B(x_\sigma,\delta/8), \qquad \psi_\sigma
  = 0 \text{ outside $B(x_\sigma, \delta/4)$}, \qquad
  \cC^1(\psi_\sigma) = O(1/\delta).
\end{displaymath}
The bound on $\cC^1(\psi_\sigma)$ follows  from the proof
of \cite[Proposition~5.8]{AG}.
Now set
\begin{displaymath}
\phi_\sigma(q) = \frac{\psi_\sigma(q)}{\sum_{\rho \in \Delta_n}
  \psi_\rho(q)}. 
\end{displaymath}
On $B(x_\sigma,\delta/16)$, the numerator is $1$ and at most $M$
terms, each at most $1$, contribute to the denominator. Hence
\begin{displaymath}
\phi_\sigma(x) \ge \frac{1}{M} \qquad \text{ for
$x \in B(x_\sigma,\delta/16)$.}  
\end{displaymath}
Note that $\phi_\sigma = 0$ outside $B(x_\sigma, \delta/4)$.
From the construction, $\phi_\sigma$ is smooth with $\cC^1(\phi_\sigma) =
O(1/\delta)$. Here $\cC^1(\phi)$ denotes the $C^1$-norm of the
function $\phi$,  as defined in \cite{AG}. 

We now define a  random walk as follows. Define
\[
\Omega_n'= \Big\{ (q_0, q_1, \dots, q_n)  : q_j \in \Hor_j, 
\ 0\leq j \leq n, \qquad d_\AGY(q_j,q_{j+1}(-\tau')) \le\delta/2,\ 0
\le  j \le  n-1 \Big\}. 
\]
We equip $\Omega_{n}'$ with a measure $\PP_n'$ as follows. 
We define the law of $q_0 \in \Hor_0$ to be $\tfrac{1}{\mu_0'(\Hor_0)} \,
d\mu_0'(q_0)$.

For $p, q  \in \Hor_j$ define
\begin{displaymath}
\bar{P}_p(q) = \sum_{\sigma \in \Delta_j} \frac{\phi_\sigma(p)}{\int \phi_\sigma
\,  d\mu_j'}  \phi_\sigma(q), \qquad P_p(q) = \lambda_j(q)
\bar{P}_p(q).   
\end{displaymath}
Note that  for any  $p$, $q$, the sum has at most $M$ terms. 

For $0 \le j  \le n-1$, given $q_j \in \Hor_j$, 
we define the law of $q' \in \Hor_j$ to be
\begin{displaymath}
\bar{P}_{q_j}(q')  \, d\mu'_j(q') = P_{q_j}(q') \, d\bar{\mu}_j(q'),  
\end{displaymath}
and then define $q_{j+1} = q'(\tau') \in \Hor_{j+1}$. 

{
For a bounded Lipschitz function $f$ on a horosphere, write
\[
\|f\|_{\mathrm{Lip}}
=\|f\|_\infty+\operatorname{Lip}_{d_\AGY}(f).
\]
Since $d_\AGY$ is the Finsler length metric on the horosphere,
$\|f\|_{\mathrm{Lip}}\leq\cC^1(f)$ for $C^1$ functions $f$.
\par}
The immediate properties of this random walk are summarized in the
following:
\begin{lemma}
\label{lemma:random:walk:easy}
Under the measure $\PP_n'$, for any $0 \le j  \le n$, the following
hold:
\begin{itemize}
\item[{\rm (i)}] For $1 \le j \le  n$, we have $q_j(-\tau') \in B(q_{j-1},
  \delta/2)$.
\item[{\rm (ii)}] The law of $q_j$ is $\frac{1}{\mu_j'(\Hor_j)} \,
  d\mu_j'(q_j)$.
\item[{\rm (iii)}] For any $p \in \Hor_j$ and any
  $\sigma  \in \Delta_j$ with $\phi_\sigma(p) > 0$,
\begin{displaymath}
\int \phi_\sigma \, d\mu_j' \ge c_\delta \lambda_j(p) \bar{\mu}_j(B(p,\ooth)). 
\end{displaymath}
\item[{\rm (iv)}] For any $p \in \Hor_j$, 
\begin{displaymath}
\sup_{q} P_p(q) \le C_\delta \bar{\mu}_j(B(p,\ooth))^{-1}.  
\end{displaymath}
\item[{\rm (v)}] {For any $p\in\Hor_j$, the function $P_p$
is Lipschitz, and
\begin{displaymath}
\|P_p\|_{\mathrm{Lip}}\leq
C_\delta\bar\mu_j(B(p,\ooth))^{-1}.
\end{displaymath}
\par}
\end{itemize}
\end{lemma}

\begin{proof}
The statement (i) follows immediately from the definition of the
random walk law,  and the statement (ii) follows immediately
by induction on $j$. To see (iii) note that $\phi_\sigma  \ge 1/M$ on a
ball of radius {$\delta/16$}, and use Lemma~\ref{prop:Doubling} and (\ref{eq:loglipschitz}). Then (iv)
follows from (iii), (\ref{eq:loglipschitz}) and the definition of
$P_p(q)$.
{
For (v), put $A_p=\bar\mu_j(B(p,\ooth))^{-1}$. By (iii), the
bounded multiplicity of the partition, and its derivative bounds,
\[
\cC^1(\bar P_p)\leq C_\delta\lambda_j(p)^{-1}A_p.
\]
Every term in $\bar P_p$ is supported where $d_\AGY(p,q)\leq\delta/2$,
so \eqref{eq:loglipschitz} gives
$\lambda_j(q)\leq e^{L\delta/2}\lambda_j(p)$ there. Also,
$\|d\lambda_j\|_\AGY\leq L\lambda_j$ almost everywhere. The product
rule for locally Lipschitz functions and (iv) therefore give
\[
\|dP_p\|_\AGY
\leq L P_p+\lambda_j\|d\bar P_p\|_\AGY
\leq C_\delta A_p
\quad\text{almost everywhere}.
\]
The resulting local Lipschitz bound gives the same global bound for
the Finsler length metric. Together with (iv), this proves (v).
Here the fixed constant $L$ is absorbed into $C_\delta$.
\par}
\end{proof}

We also record some  easy estimates for later use. For $\tau'>0$, let
\begin{equation}
\label{eq:def:V0prime}
V_0'= {e^{3\tau'/4}}.
\end{equation}
{Then, Proposition~\ref{prop:AG}, applied with the fixed
choice $\alpha=1/8$, implies}
\begin{equation}
\label{eq:V_0:prime}
\frac{1}{\mu_q (B(q,\ooh))}
\int_{B(q,\ooh)} \hspace{-.5cm} V(y(\tau')) \, d\mu_q(y) \leq 
\begin{cases} 
{C_{1/8}e^{-3\tau'/4}} \cdot V(q) & \text{if} \quad V(q) \geq V_0'\\
{C_{1/8}}  & \text{if} \quad V(q) \leq V_0'.
\end{cases} 
\end{equation}

\begin{lemma}
\label{lemma:phi:sigma:V0}
There exists a constant  $C_\delta$ depending only  on  $\delta$ such  that  
for any $p \in \Hor_n$, 
\begin{equation}
\label{eq:phi:sigma:V_0}
\int P_p(q) V(q(\tau')) \, d\bar{\mu}_n(q) \leq 
\begin{cases} 
C_\delta {e^{-3\tau'/4}} \cdot V(p) & \text{if} \quad V(p) \geq V_0'\\
C_\delta & \text{if} \quad V(p) \leq V_0'.
\end{cases} 
\end{equation} 
\end{lemma}

\begin{proof}
Note that $P_p$ is supported in $B(p,\ooh)$. Then the lemma follows from 
Lemma~\ref{lemma:random:walk:easy}(iv), (\ref{eq:V_0:prime}) and
Lemma~\ref{prop:Doubling}. 
\end{proof}

\subsection{Tracking and good statistics}

For $\omega \in \Omega_{n}'$ and $0 \leq j \leq n$, define
\[
\tV_1(\omega,j) = V(q_j) \cdot \exp \left(  \sum_{i=0}^{j-1}  \int_{0}^{\tau'} f_1(q_{i+1}(-\tau'+t)) \, dt \right). 
\]
We note that the definition of $f_1$ involves $\tau$  (and  not
$\tau'$). 

Let $\cF_j'$ be the smallest $\sigma$-algebra for which $q_i = q_i(\omega)$ is measurable for all
$0 \leq i \leq j$.

\begin{proposition}
\label{prop:alt:random-walk}
For every  $\tau'  \in \tau \natls$ there  exists a constant $\alpha_3 =
\alpha_3(\tau,\delta,\tau') > 0$  with $\alpha_3 \to \alpha_2(\tau)$ as 
$\tau' \to \infty$ such that for any
$n \geq 1$, we have 
\[
\EE_{n}'(\tV_1(\omega, n) \mid \cF_0') \leq  e^{\alpha_3  n \tau'} \cdot V(q_0). 
\]
Here $\EE_{n}'$ denotes expectation with respect to the measure
$\PP_{n}'$. 
\end{proposition} 

\begin{proof}
We claim that there exists $\alpha_3 = \alpha_3(\tau,\delta,\tau')$
with $\alpha_3 \to \alpha_2(\tau)$ as $\tau' \to \infty$ such that
for $0 \leq j  <n$, 
\begin{equation}
\label{eq:alt:inductive-step}
\EE_{n}'(\tV_1(\omega, {j+1}) | \cF_j' )  \leq e^{\alpha_3 \tau'} \cdot \tV_1(\omega, j). 
\end{equation}
Note that 
\[
\frac{\EE_{n}'(\tV_1(\omega, {j+1}) \mid \cF_j' )}{\tV_1(\omega,j)} = 
\frac{1}{V(q_j)} \int P_{q_j}(q) V(q(\tau')) \exp\left( \int_0^{\tau'}
  f_1(q(t)) \,dt \right) \, d\bar{\mu}_j(q).
\]
Recall the definition $\tV_1(q,n)$ in (\ref{eq:def:tV1}),  and write  $\tau' =
m\tau$. Then, the right  hand side of  the above equation can
be rewritten as
\begin{displaymath}
\frac{1}{V(q_j)} \int P_{q_j}(q) \tV_1(q,m) \, d\bar{\mu}_j(q). 
\end{displaymath}
Note that $P_{q_j}(q) = 0$ if $q \not\in B(q_j,\ooth)$. 
Therefore,  in view of Lemma~\ref{lemma:random:walk:easy} (iv), 
\[
\frac{\EE_{n}'(\tV_1(\omega, {j+1}) \mid \cF_j' )}{\tV_1(\omega,j)} \le 
\frac{1}{V(q_j)} \frac{C_\delta}{\bar{\mu}_j(B(q_j,\ooth))} \int_{B(q_j,\ooth)}
\tV_1(q,m) \, d\bar{\mu}_j(q) \le C_\delta e^{\alpha_2 \tau'},
\]
where $\alpha_2 = \alpha_2(\tau)<1$ be as in
Lemma~\ref{lemma:geodesic:expectation}, 
and we have used Lemma~\ref{lemma:geodesic:expectation} with the base
point $p_0=q_j$ for the last estimate. 

For fixed $\tau$ and $\delta$, 
we now choose  $\alpha_3 = \alpha_3(\tau')$ to be the smallest
so that $C_\delta e^{\alpha_2 \tau'} \leq e^{\alpha_3 \tau'}$. Then,
$\alpha_3 \to \alpha_2(\tau)$ as $\tau' \to \infty$, and
(\ref{eq:alt:inductive-step}) holds. 
The rest of the proof is the same as the proof of
Proposition~\ref{prop:random-walk}. 
\end{proof}

\begin{lemma}
\label{lemma:longrange:good:statistics}
{
For every $\eta>0$ there exists $\tau_0=\tau_0(\eta)>0$ such that,
for every $0<\delta\leq1/400$ and every $\tau\geq\tau_0$, there exists
$\tau'_0=\tau'_0(\eta,\delta,\tau)>0$ such that if $\tau'>\tau'_0$
(and $\tau'\in\tau\natls$), the following holds:
\par}

{Let $n\geq0$ be an integer, and let} $G_+(\omega,\eta)$
denote the set of $k \in [0,n] \cap  \natls$ such that for all $k<j\le n$,
\begin{displaymath}
\sum_{i=k}^{j-1}\int_{0}^{\tau'} f_1(q_{i+1}(-\tau'+t)) \, dt \le \eta
(j-k) \tau'. 
\end{displaymath}
Then
\begin{displaymath}
\PP_n'\Big(|G_+(\omega,\eta)| <
(1-8\eta)n \Big) < {e^{-\tfrac{\eta^4}{4} n \tau'}}.
\end{displaymath}
\end{lemma}

\begin{proof}
{If $\eta\geq 1/8$, the event in the statement is empty, and
if $n=0$, the claim is immediate.  Thus, we may assume $0<\eta<1/8$ and
$n\geq1$.}
{
Let $\alpha_2=\alpha_2(\tau)$ be as in
Lemma~\ref{lemma:geodesic:expectation}. Since $\alpha_2(\tau)\to0$,
choose $\tau_0=\tau_0(\eta)$ so that
$\alpha_2(\tau)<\eta^4/16$ for every $\tau\geq\tau_0$.
Suppose $0<\delta\leq1/400$ and $\tau\geq\tau_0$.
\par}
In view of Proposition~\ref{prop:alt:random-walk} and
Lemma~\ref{lemma:V0:finite:exp}, for $\tau'\in \tau\natls$ large
enough {(depending on $\eta$, $\delta$, and $\tau$)}, we have
\begin{displaymath}
\EE_{n}'(\tV_1(\omega, n)) \leq
\mu_0'(\Hor_0)^{-1}\left(\int_{\Hor_0}V\,d\mu_0'\right) e^{{(\eta^4/16)} n \tau'} \le 
e^{{(\eta^4/8)}  n \tau'},
\end{displaymath}
where the last inequality is obtained by taking $\tau'$ sufficiently
large.  Then, by the Markov inequality,
\begin{displaymath}
\PP_n'\left(\sum_{i=0}^{n-1} \int_{0}^{\tau'} f_1(q_{i+1}(-\tau'+t)) \, dt
  \geq {\frac{\eta^4}{2}} n\tau' \right) \le
  {e^{-(\eta^4/4) n\tau'}}.
\end{displaymath}
In the rest of the proof, we assume that $\omega$ is such that
\begin{equation}
\label{eq:sum:small:enough}
  \sum_{i=0}^{n-1} \int_{0}^{\tau'} f_1(q_{i+1}(-\tau'+t)) \, dt
  < {\frac{\eta^4}{2}} n\tau'.
\end{equation}
Let $H_+(\omega,\eta) = ([0,n] \cap \natls)  \setminus G_+(\omega,\eta)$. If $k \in H_+(\omega,\eta)$
then there exists $k < j \le n$ such that
\begin{displaymath}
\sum_{i=k}^{j-1}\int_{0}^{\tau'} f_1(q_{i+1}(-\tau'+t)) \, dt >\eta (j-k)\tau'.
\end{displaymath}
We can rewrite this as
\begin{displaymath}
(j-k) <  \frac{1}{\eta \tau'} \sum_{i=k}^{j-1}\int_{0}^{\tau'}
f_1(q_{i+1}(-\tau'+t)) \, dt. 
\end{displaymath}
The sets $[k,j-1] \cap \natls$ cover $H_+(\omega,\eta)$ and we can choose a
subcover such that each point is covered at most twice. Then,
\begin{displaymath}
|H_+(\omega,\eta)| \le  \frac{2}{\eta \tau'}  \sum_{i=0}^{n-1} \int_{0}^{\tau'}
f_1(q_{i+1}(-\tau'+t)) \, dt \le {\eta^3 n},  
\end{displaymath}
where the last estimate follows from
(\ref{eq:sum:small:enough}). Therefore,
\begin{displaymath}
|G_+(\omega,\eta)| \ge {(1-\eta^3) n} \ge (1-8\eta)n. 
\end{displaymath}
\end{proof}

\begin{lemma}
\label{lemma:tracking}
For every $0 < \eta <1/4$, every {$0<\delta\leq1/400$} and  every $\tau > 0$  
there exists $\tau'_0 = \tau'_0(\eta,\delta,\tau)$ such
that if $\tau' > \tau'_0$ (and $\tau' \in \tau  \natls$) the following
holds:

Suppose $\omega = (q_0,\dots, q_n) \in \Omega_n'$ and $k \in
G_+(\omega,\eta)$. Then,
\begin{displaymath}
d_\AGY(q_k, q_n(-(n-k)\tau')) < \delta.
\end{displaymath}
\end{lemma}

Note that $q_k$ and $q_n(-(n-k)\tau')$ are both in the same
horosphere. 

\begin{proof}
Let
\begin{displaymath}
\hat{q}(t) =  q_{j+1}(-(j+1)\tau'+t)  \text{ for $j\tau' \le  t <  (j+1)\tau'$.}
\end{displaymath}
Let
\begin{displaymath}
F_j =  \{ t \in [k\tau',(k+j)\tau'] \st
f_1(\hat{q}(t)) = 0 \}.
\end{displaymath}
Since $k \in G_+(\omega,\eta)$, we have
\begin{equation}
\label{eq:Fj:large}
|F_j| \ge (1-\eta) j \tau'.   
\end{equation}
Consider the function
\begin{displaymath}
h_j(t) = d_\AGY(q_{k+j}(t-(k+j) \tau'),
\hat{q}(t)). 
\end{displaymath}
Note that $h_j$ may not be continuous at the points $i \tau'$, $i \in
\natls$. However,
the jump  $|h_j(i\tau'+) - h_j(i \tau'-)|$ is bounded by $\delta/2$.  Also,
$h_j$ is increasing on each interval $[i\tau',(i+1)\tau']$,  and
$h_j((k+j)\tau'-) = 0$.  

Let
\begin{displaymath}
E_j = \{  t \in [k\tau',(k+j)\tau'] \st h_j(t)
> \ooh \}.
\end{displaymath}
We claim that
\begin{equation}
\label{eq:Ej:cap:Fj}
|E_j \cap F_j| \le \eta j \tau'.
\end{equation}
Let $T = T(\supp{1-f_1})$, where $T(\cdot)$ is as in
Corollary~\ref{cor:AGY:endpoint:contraction}.  Let
\begin{displaymath}
\cJ  =  \{ i  \in [k,k+j-1]  \st |E_j \cap F_j \cap [i\tau', (i+1)\tau']| \ge
2T \},
\end{displaymath}
and for $i \in \cJ$ let
$t_i'$ denote the smallest $t \in [i \tau',(i+1)\tau']$ such  that
$h_j(t_i') \ge \ooh$. 

By Corollary~\ref{cor:AGY:endpoint:contraction} applied  to the  interval
$[t_i',(i+1)\tau']$, we get
\begin{displaymath}
|E_j \cap  F_j \cap [i \tau',(i+1)\tau']  \le 400 T  \left(
  h_j((i+1)\tau'-) - h_j(i\tau'+)\right).  
\end{displaymath}
For  $i \not\in \cJ$, we  have by definition,  
\begin{displaymath}
|E_j \cap  F_j \cap [i \tau',(i+1)\tau']  \le 2T.  
\end{displaymath}
Sum these bounds over all $j$ blocks. The increases of $h_j$
between jumps sum to at most $j\delta$, since $h_j\geq0$,
$h_j((k+j)\tau'-)=0$, and each jump has size at most $\delta/2$.
Thus
\begin{displaymath}
|E_j \cap F_j| \le  400 T j \delta + 2 T j.
\end{displaymath}
We choose $\tau_0'$ so  that for $\tau'  > \tau_0'$,  $400 T \delta
+2T  \le \eta
\tau'$. This  proves (\ref{eq:Ej:cap:Fj}).   
Then, in view of  (\ref{eq:Fj:large}), and since $\eta < 1/4$, 
\begin{equation}
\label{eq:Fj:setminus:Ej}
  |F_j\setminus E_j| \geq (1-2 \eta) j \tau' \ge \frac{1}{2} j\tau'. 
\end{equation}

Let $K$ be the $\Mod$-compact set 
\begin{displaymath}
K = \{ q \st \text{ there exists $p \in \Hor(q)$ with $f_1(p) = 0$ and
  $d_\AGY(p,q) \le \ooh$} \},
\end{displaymath}
and let $C_K$, $\beta_K$ be as in Corollary~\ref{cor:AGY:half}. Choose $\tau'
=\tau'(\eta,\delta,\tau)$ so that $C_K e^{-\beta_K \tau'} < 0.01$. 
We now claim that for $0 \le j \le n-k$,
\begin{equation}
\label{eq:dAGY:convergent:sum}
d_\AGY(q_k,q_{k+j}(-j\tau')) \le \frac{\delta}{2} \left( 1+ \sum_{i=1}^{j-1}
C_K e^{-\beta_K i \tau'}\right) \le \delta.   
\end{equation}
This will be done by induction on $j$.  For $j=0$ there is nothing to
prove, and for $j=1$ it follows from $d_\AGY(q_k,q_{k+1}(-\tau')) \leq
\delta/2$. Now suppose $j \ge 1$ and the statement holds for $j$. 

By the triangle inequality, it is enough to show that
\begin{displaymath}
d_\AGY(q_{k+j}(-j\tau'), q_{k+j+1}(-(j+1) \tau')) \le \frac{\delta}{2}
C_K e^{-\beta_K
   j \tau'}.  
\end{displaymath}
Note that for $t \in F_j \setminus E_j$,
\begin{displaymath}
d_\AGY(q_{k+j}(-(k+j)\tau' + t), \hat{q}(t)) \le \ooh, 
\end{displaymath}
which together with the definitions of $F_j$ and $K$
implies that $q_{k+j}(-(k+j)\tau' + t) \in K$. 
Also, note that
\begin{displaymath}
d_\AGY(q_{k+j+1}(-\tau'), q_{k+j}) \le \frac{\delta}{2}. 
\end{displaymath}
Then, by Corollary~\ref{cor:AGY:half} and
(\ref{eq:Fj:setminus:Ej}),  
\begin{displaymath}
d_\AGY(q_{k+j+1}(-(j+1)\tau'), q_{k+j}(-j\tau')) \le \frac{\delta}{2}
C_K e^{-\beta_K j   \tau'},  
\end{displaymath}
as required. 

The final estimate is just the case $j=n-k$ of
(\ref{eq:dAGY:convergent:sum}).
\end{proof}

\subsection{Good Statistics}

{
In this subsection, we take $\delta=1/400$. As in the previous
subsection, let $m$ be a large integer and put $\tau'=m\tau$.
\par}

For $\omega \in \Omega_{n}'$ and $0 \leq j \leq n$, define
\begin{displaymath}
V^\ast(\omega,j) = V(q_j) \exp\left( \sum_{i=0}^{j-1} Y_\rho(q_{i+1}(-\tau'),m)
\right). 
\end{displaymath}
We note that the definition of $Y_\rho$ involves $\tau$  (and  not
$\tau'$). 

Let $\cF_j'$ be the smallest $\sigma$-algebra for which $q_i = q_i(\omega)$ is measurable for all
$0 \leq i \leq j$.

\begin{proposition}
\label{prop:estimate:Vast}
For every $0<\rho<1$ there exist $\tau_0=\tau_0(\rho)>0$
a constant $C\geq 1$ and a function
$\beta_\rho(\tau)$ with $\beta_\rho(\tau)  \to 0$ as $\tau \to \infty$
such that the following holds.
For every $\tau>\tau_0$, every $m,n\in\natls$ with
$\tau'=m\tau$, we have, almost surely
\begin{displaymath}
\EE_{n}'(V^\ast(\omega, n) \mid \cF_0') \leq  C
e^{\beta_\rho(\tau)n} \cdot V(q_0), 
\end{displaymath}
Here $\EE_{n}'$ denotes expectation with respect to the measure
$\PP_{n}'$. 
\end{proposition}

\begin{proof}
{
Set
\[
\kappa_\rho(\tau)
=e^{-(\rho-\epsilon_\rho(\tau))\tau}.
\]
After increasing $\tau_0(\rho)$, we may assume that
$\epsilon_\rho(\tau)<\rho$ for $\tau>\tau_0$.

For $p\in\Hor_j$, define
\begin{align*}
T_jf(p)&=\int P_p(q)f(q(\tau'))\,d\bar\mu_j(q),\\
R_jf(p)&=\int P_p(q)Y_\rho(q,m)f(q(\tau'))\,d\bar\mu_j(q).
\end{align*}
Using Proposition~\ref{prop:AG} with the fixed choice $\alpha=1/8$,
there are absolute constants $b\geq1$ and $C_0\geq1$, and a number
$a_\tau=C_0e^{-3\tau/4}$, such that
\begin{equation}
\label{eq:Vast:ordinary-transition}
T_jV\leq a_\tau V+b.
\end{equation}
Here we used Lemma~\ref{lemma:random:walk:easy}(iv), the fact that
$P_p$ is supported in $B(p,\delta/2)\subset B(p,\ooh)$, and the
inequality $\tau'=m\tau\geq\tau$.

Increase $\tau_0$ once more so that $a_\tau\leq1/2$.  It follows from
\eqref{eq:Vast:ordinary-transition} and $T_j1=1$ that every string of
$\ell$ consecutive ordinary transitions satisfies
\begin{align}
T_jT_{j+1}\cdots T_{j+\ell-1}V
&\leq a_\tau^\ell V+b\sum_{s=0}^{\ell-1}a_\tau^s \notag\\
&\leq V+2b\leq C_1V,
\label{eq:Vast:ordinary-string}
\end{align}
where $C_1=1+2b$ and we used $V\geq1$.

On the other hand, Proposition~\ref{prop:short:range:good:statistics}
and Lemma~\ref{lemma:random:walk:easy}(iv) give
\begin{equation}
\label{eq:Vast:bad-transition}
R_jV\leq C_0\kappa_\rho(\tau)V.
\end{equation}
Indeed, since $P_p$ is supported in
$B(p,\delta/2)\subset B(p,\ooth)$,
\begin{align*}
R_jV(p)
&\leq
\frac{C_\delta}{\bar\mu_j(B(p,\ooh))}
\int_{B(p,\ooth)}V(q(m\tau))Y_\rho(q,m)\,d\bar\mu_j(q)\\
&\leq C_0\kappa_\rho(\tau)V(p).
\end{align*}

Since $Y_\rho$ takes values in $\{0,1\}$, we have
$e^{Y_\rho}=1+(e-1)Y_\rho$.  Thus, if
\[
K_jf(p)=\int P_p(q)e^{Y_\rho(q,m)}f(q(\tau'))\,d\bar\mu_j(q),
\]
then
\[
K_j=T_j+(e-1)R_j.
\]
By the Markov property and the definition of $V^\ast$,
\[
\EE_n'\big(V^\ast(\omega,n)\mid\cF_0'\big)
=(K_0K_1\cdots K_{n-1}V)(q_0).
\]

Expand the product $K_0K_1\cdots K_{n-1}$.  A term containing exactly
$k$ occurrences of the operators $R_j$ consists of those $k$ bad
transitions separated by at most $k+1$ strings of ordinary
transitions.  By \eqref{eq:Vast:ordinary-string} and
\eqref{eq:Vast:bad-transition}, each such term is bounded by
\[
C_1^{k+1}\big(C_0\kappa_\rho(\tau)\big)^kV(q_0).
\]
There are $\binom nk$ such terms.  Consequently,
\begin{align*}
\EE_n'\big(V^\ast(\omega,n)\mid\cF_0'\big)
&\leq C_1\sum_{k=0}^n\binom nk
\big((e-1)C_0C_1\kappa_\rho(\tau)\big)^kV(q_0)\\
&=C_1\big(1+(e-1)C_0C_1\kappa_\rho(\tau)\big)^nV(q_0)\\
&\leq C_1\exp\big((e-1)C_0C_1\kappa_\rho(\tau)n\big)V(q_0).
\end{align*}
The proposition follows by taking $C=C_1$ and
\[
\beta_\rho(\tau)
=(e-1)C_0C_1e^{-(\rho-\epsilon_\rho(\tau))\tau}.
\]
Clearly $\beta_\rho(\tau)\to0$ as $\tau\to\infty$.
}
\end{proof}


We need two backward counting conditions. The set $G_-$ controls
the frequency of blocks detected by $Y_\rho$ in every interval
ending at the chosen index. The set $G_\ast$ controls the frequency
of indices outside $G_+$, for which Lemma~\ref{lemma:tracking}
gives no tracking estimate.

\begin{lemma}
\label{lemma:Gminus:small}
Suppose $0<\eta<1/16$, and put
\begin{displaymath}
\rho=\frac{\eta^2}{100}.
\end{displaymath}
{
There exist $\tau_0=\tau_0(\eta)>0$ and $0<\kappa=\kappa(\eta)<1$
such that, for every $\tau\geq\tau_0$, there exists
$\tau'_0=\tau'_0(\eta,\tau)>0$ such that if $\tau'>\tau'_0$
(and $\tau'\in\tau\natls$), the following holds:

Let $n\geq0$ be an integer, put $m=\tau'/\tau$, and let
$G_-(\omega,\eta)$ denote the set of $k\in[1,n]\cap\natls$
such that, for every $0\leq j<k$,
\par}
\begin{displaymath}
 \sum_{i=j}^{k-1}
 Y_\rho(q_{i+1}(-\tau'),m)
 \leq {\eta^2(k-j)}.
\end{displaymath}
Let $G_\ast(\omega,\eta)$ be the set of
$k\in[1,n]\cap\natls$ such that, for every $0\leq j<k$,
\begin{displaymath}
 \big|
 \{\ell\in[j+1,k]\cap\natls:
       \ell\notin G_+(\omega,\eta)\}
 \big|
 \leq{\eta^2(k-j)}.
\end{displaymath}
{Then, for all sufficiently large $n$
(depending on $\eta$, $\tau$, and $\tau'$),}
{
\begin{displaymath}
 \PP_n'\Big(
 \big|G_+(\omega,\eta)\cap G_-(\omega,\eta)
       \cap G_\ast(\omega,\eta)\big|
 <(1-8\eta)n
 \Big)
 \leq\kappa^n.
\end{displaymath}
\par}
\end{lemma}

\begin{proof}
{
Choose $\tau_0=\tau_0(\eta)$ as in
Lemma~\ref{lemma:longrange:good:statistics}.
By Proposition~\ref{prop:estimate:Vast}, after increasing $\tau_0$
we may assume that
\[
\beta_\rho(\tau)<\frac{\eta^4}{4},
\qquad \rho=\frac{\eta^2}{100}.
\]
for every $\tau\geq\tau_0$. Suppose $\tau\geq\tau_0$, and choose
$\tau'_0=\tau'_0(\eta,\tau)\geq1$ so that
Lemma~\ref{lemma:longrange:good:statistics} applies with $\delta=1/400$.
\par}
By Proposition~\ref{prop:estimate:Vast} and
Lemma~\ref{lemma:V0:finite:exp}, there is a constant $C_\eta$ such that
\begin{displaymath}
 \EE_n'\big(V^\ast(\omega,n)\big)
 \leq C_\eta e^{\beta_\rho(\tau)n}.
\end{displaymath}
Since $V\geq1$, Markov's inequality gives
\begin{align*}
 &\PP_n'\left(
 \sum_{i=0}^{n-1}Y_\rho(q_{i+1}(-\tau'),m)
 \geq\frac{\eta^4}{2}n
 \right)                                                   \\
 &\qquad\leq
 C_\eta\exp\left(
 -\left(\frac{\eta^4}{2}-\beta_\rho(\tau)\right)n
 \right)
 \leq C_\eta e^{-\eta^4n/4}.
\end{align*}
Consequently, for some $\kappa_-=\kappa_-(\eta)<1$ and every
sufficiently large $n$, the probability above is at most
$\kappa_-^n$.

In the remainder of the proof, suppose that
\begin{equation}
\label{eq:Yrho:total:small}
 \sum_{i=0}^{n-1}Y_\rho(q_{i+1}(-\tau'),m)
 <\frac{\eta^4}{2}n.
\end{equation}
Let
\begin{displaymath}
 H_-(\omega,\eta)
 =([1,n]\cap\natls)\setminus G_-(\omega,\eta).
\end{displaymath}
For every $k\in H_-(\omega,\eta)$, there is $0\leq j<k$ such that
\begin{displaymath}
 \sum_{i=j}^{k-1}Y_\rho(q_{i+1}(-\tau'),m)>{\eta^2(k-j)}.
\end{displaymath}
The intervals $[j+1,k]\cap\natls$ cover $H_-(\omega,\eta)$, and we may
choose a subcover of multiplicity at most two.  Therefore,
{
\begin{displaymath}
 |H_-(\omega,\eta)|
 \leq\frac{2}{\eta^2}
 \sum_{i=0}^{n-1}Y_\rho(q_{i+1}(-\tau'),m)
 \leq\eta^2n.
\end{displaymath}
\par}

The proof of Lemma~\ref{lemma:longrange:good:statistics} shows that,
outside an event of probability at most
\begin{displaymath}
 e^{-\eta^4n\tau'/4},
\end{displaymath}
we have
\begin{displaymath}
 |H_+(\omega,\eta)|\leq\eta^3n,
 \qquad
 H_+(\omega,\eta)
 =([1,n]\cap\natls)\setminus G_+(\omega,\eta).
\end{displaymath}
Let
\begin{displaymath}
 H_\ast(\omega,\eta)
 =([1,n]\cap\natls)\setminus G_\ast(\omega,\eta).
\end{displaymath}
The same interval-covering argument gives
{
\begin{displaymath}
 |H_\ast(\omega,\eta)|
 \leq\frac{2}{\eta^2}|H_+(\omega,\eta)|
 \leq2\eta n.
\end{displaymath}
\par}
{
Combining these bounds, we obtain
\[
 |H_+\cup H_-\cup H_\ast|
 \leq(\eta^3+\eta^2+2\eta)n<8\eta n.
\]
\par}
{
Since $\tau'>1$, we have
$e^{-\eta^4n\tau'/4}\leq e^{-\eta^4n/4}$. Choose
\[
\max\{\kappa_-,e^{-\eta^4/4}\}<\kappa(\eta)<1.
\]
Then the probability in the statement is at most
$\kappa_-^n+e^{-\eta^4n/4}\leq\kappa(\eta)^n$ for all sufficiently
large $n$. The lemma follows.
\par}
\end{proof}

\begin{lemma}
\label{lemma:Gminus:implies:good:statistics}
For every $\tau>0$ there exists a $\Mod$--compact set
$K=K(\tau)$, containing $K_0$, with the following property.  For every
$0<\eta<1/16$, put $\rho=\eta^2/100$.  There exists
$\tau'_0=\tau'_0(\eta,\tau)>0$ such that, for every $m,n\in\natls$ with
$\tau'=m\tau\geq\tau'_0$, every $\omega=(q_0,\ldots,q_n)\in\Omega_n'$, and
every $k\in[1,n]\cap\natls$ satisfying
\begin{displaymath}
{k\in G_-(\omega,\eta)\cap G_\ast(\omega,\eta)},
\end{displaymath}
the following hold:
\begin{itemize}
\item[\rm (i)] The point $q_n(-(n-k)\tau')$ belongs to $K$, and for every
$0\leq a<k\tau'$,
\begin{displaymath}
 \big|
 \{v\in[a,k\tau']:\ q_n(-n\tau'+v)\in K\}
 \big|
 \geq(1-8\eta)(k\tau'-a).
\end{displaymath}
\item[\rm (ii)] For every $0  \le b < k \tau'$, we have
\begin{multline*}
 \big|\{u\in[b,k\tau']:\ q_n(-n\tau'+u)\in K\text{ and, for every }
 0\leq a<u,\\
 \big|\{v\in[a,u]:q_n(-n\tau'+v)\in K\}\big|
 \geq(1-8\eta)(u-a)\}\big|
 \geq(1-8\eta)(k\tau'-b).
\end{multline*}
\end{itemize}
\end{lemma}

{
\begin{proof}
Put $T=k\tau'$, write $p(u)=q_n(-n\tau'+u)$, and set
\[
K=K_0\cup\{q\in\POT:V(q)\leq9eV_0e^{4\tau}\}.
\]
This is a $\Mod$--compact set depending only on $\tau$.
Choose $\tau'_0$ sufficiently large that Lemma~\ref{lemma:tracking}
applies with $\delta=1/400$.

For $0\leq i<k$ and $0\leq r<m$, call
\[
I_{i,r}=[i\tau'+r\tau,\ i\tau'+(r+1)\tau]
\]
good if $f_1(q_{i+1}(-\tau'+t))=0$ for some
$t\in[r\tau,(r+1)\tau]$, and bad otherwise. If
$i+1\in G_+(\omega,\eta)$, tracking and backward monotonicity of
$d_\AGY$ give
\[
d_\AGY\big(q_{i+1}(-\tau'+t),p(i\tau'+t)\big)<\delta
\qquad(0\leq t\leq\tau').
\]
Since $f_1(q)=0$ implies $V(q)\leq9V_0e^{3\tau}$, the exponential
change of $V$ along the flow and the $1$--Lipschitz property of
$\log V$ imply that
\begin{equation}
\label{eq:good:tau:block:in:K}
p(u)\in K\qquad\text{for all }u\in I_{i,r}
\end{equation}
whenever $i+1\in G_+(\omega,\eta)$ and $I_{i,r}$ is good.

Call $i$ admissible if $Y_\rho(q_{i+1}(-\tau'),m)=0$ and
$i+1\in G_+(\omega,\eta)$. For an admissible $i$, the definition of
$Y_\rho$ and the fact that $f_1=1$ throughout every bad interval give
\begin{equation}
\label{eq:bad:tau:suffix}
\#\{r\in[a,m-1]\cap\zed:I_{i,r}\text{ is bad}\}
\leq\rho(m-a)\qquad(a\in[0,m]\cap\zed).
\end{equation}
Since $k\in G_-(\omega,\eta)\cap G_\ast(\omega,\eta)$, we also have
\begin{equation}
\label{eq:nonadmissible:pieces}
\#\{i\in[j,k-1]\cap\zed:i\text{ is not admissible}\}
\leq2\eta^2(k-j)\qquad(j\in[0,k-1]\cap\zed).
\end{equation}
In particular, $k-1$ is admissible since $2\eta^2<1$.

Let $B$ be the union of the nonadmissible $\tau'$--pieces and the
bad $\tau$--intervals in admissible pieces. By
\eqref{eq:good:tau:block:in:K}, $p(u)\in K$ outside $B$.
The function $a\mapsto|B\cap[a,T]|$ is affine between consecutive
grid points. Thus \eqref{eq:bad:tau:suffix} within each
admissible piece and \eqref{eq:nonadmissible:pieces} for the
nonadmissible pieces, followed by linear interpolation and addition,
give
\begin{equation}
\label{eq:A16:bad:suffix}
|B\cap[a,T]|\leq(2\eta^2+\rho)(T-a)
\qquad(0\leq a\leq T).
\end{equation}
The last small interval is good by \eqref{eq:bad:tau:suffix}, since
$k-1$ is admissible and $\rho<1$. Hence $p(T)\in K$, and
\eqref{eq:A16:bad:suffix} proves~(i).

For~(ii), define
\[
H=\{u\in[0,T]:|B\cap[a,u]|>8\eta(u-a)
\text{ for some }0\leq a<u\}.
\]
We claim that
\begin{equation}
\label{eq:A16:exceptional:suffix}
|H\cap[b,T]|\leq\frac{2(2\eta^2+\rho)}{8\eta}(T-b)
=\frac{201\eta}{400}(T-b)<8\eta(T-b)
\end{equation}
for every $0\leq b<T$. To see this, take a nonempty compact set
$A\subset H\cap[b,T]$. Slightly extending witnessing intervals at
their right endpoints, choose a finite cover of $A$ by intervals
$J\subset[0,T]$ of multiplicity at most two, each with right endpoint
at least $b$ and $|B\cap J|>8\eta|J|$. This is possible since
\eqref{eq:A16:bad:suffix} implies $T\notin H$.
Let $v$ be the minimum of $b$ and their left endpoints. Their union
contains $A\cup(v,b)$, so
\[
8\eta(|A|+b-v)\leq8\eta\sum_J|J|
\leq\sum_J|B\cap J|
\leq2(2\eta^2+\rho)(T-v).
\]
Since $2(2\eta^2+\rho)/(8\eta)<1$, rearranging and taking the
supremum over $A$ proves \eqref{eq:A16:exceptional:suffix}.
Every interior point of $B$ belongs to $H$ since $8\eta<1$.
Thus, apart from the finitely many subdivision endpoints, every
$u\notin H$ has $p(u)\in K$ and satisfies the condition in~(ii).
Equation~\eqref{eq:A16:exceptional:suffix} therefore proves~(ii).
\end{proof}
\par}

We introduce a variant of
the definition of good statistics. For a Mod-compact set $K$,
a constant $0<\eta<1/2$ and $t  > 0$
we say a quadratic differential $q\in g_t(\Hor_0)$ has  $(K,
\eta)$--good statistics relative to $\Hor_0$
if $q \in K$ and for all $s \in  [-t,0]$,
\begin{equation}
\label{eq:app:def:good:statistics}
|\{ u \in [s,0]  : q(u)\in K\}|  \geq (1-\eta)\cdot |s| .
\end{equation}

\begin{definition}
\label{def:K:eta:generic}
Let $K \subset \POT$ be $\Mod$--compact and let $0 < \eta < 1/2$. 
For $R>0$, we say $q \in \Hor_0$ is $(K, \eta)$-generic at $R$ if the following hold: 
\begin{itemize}
\item[{\rm (i)}] $q(R)$ has $(K, \eta)$--good statistics relative to
  $\Hor_0$, and 
\item[{\rm (ii)}] For $0 \le  s \le R$,
$$| \{ u \in (s, R)  : \text{$q(u)$ has $(K,
    \eta)$--good statistics relative to $\Hor_0$} \} | \geq (1-\eta) (R-s).$$
\end{itemize}
\end{definition} 

\begin{proposition}
\label{prop:good:statistics:final}
{
For every $0<\eta<1/2$ there exist $0<\kappa=\kappa(\eta)<1$
and $\tau_0=\tau_0(\eta)>0$ such that, for every $\tau\geq\tau_0$,
there exists $\tau'_0=\tau'_0(\eta,\tau)>0$ such that if
$\tau'>\tau'_0$ (and $\tau'\in\tau\natls$), the following holds:

Let $K=K(\tau)\supset K_0$ be the
$\Mod$--compact set supplied by
Lemma~\ref{lemma:Gminus:implies:good:statistics}, and set $R_j=j\tau'$.
\par}
Then, for
\[
\bar{k}_n(x) = \# \Big\{ j \in [1,n]  : \text{$x$ is $(K, \eta)$-generic at $R_j$} \Big\},
\]
we have, for all sufficiently large $n$
{(depending on $\eta$, $\tau$, and $\tau'$)}, that for any $\alpha \in
\cA$, 
\[
\mu_0(\{ x \in U_\alpha(0)  : \bar{k}_n(x) \le (1-\eta)n \}) \le \kappa^n
\mu_0(U_\alpha(0)).  
\]
\end{proposition}

\begin{proof}
Since $\cA$  is finite, it  is enough to prove the proposition for a
fixed $\alpha$.
Let
{
\begin{displaymath}
 H_n=G_+(\omega,\epsilon)\cap G_-(\omega,\epsilon)
     \cap G_\ast(\omega,\epsilon),\qquad\epsilon=\eta/16,
\end{displaymath}
\par}
and let $E_n$ be the event
\begin{displaymath}
|H_n| \ge  (1-8\epsilon)n \ge (1-\eta/2)n.
\end{displaymath}
{
Choose $\tau_0=\tau_0(\eta)$ and $0<\kappa_0=\kappa_0(\eta)<1$
as in Lemma~\ref{lemma:Gminus:small}, with $\eta/16$ in place of $\eta$.
Suppose $\tau\geq\tau_0$, and let $\tau'_0=\tau'_0(\eta,\tau)$ be as
in that lemma. Then, for $\tau'>\tau'_0$ (and $\tau'\in\tau\natls$),
$\PP_n'(E_n)\geq1-\kappa_0^n$ for all sufficiently large $n$.
\par}

Now  suppose $E_n$ holds. We can enlarge $\tau_0'$ so that
Lemma~\ref{lemma:Gminus:implies:good:statistics} applies. 
It follows from Lemma~\ref{lemma:Gminus:implies:good:statistics}
that for each $k \in H_n$, $q=q_n(-n\tau')$ is $(K,\eta)$--generic at $R_k
= k \tau'$. Therefore,
\begin{displaymath}
\bar{k}_n(q_n(-n\tau')) \ge (1-\eta/2)n >  (1-\eta) n. 
\end{displaymath}
Let $F_n$ be the event $q_n
\in U_\alpha(n \tau')$. Then, by
Lemma~\ref{lemma:random:walk:easy}(ii) and  (\ref{eq:appendix:Volume-Growth}), 
$\PP_n'(F_n) \ge c$ for some constant $c$  independent of $n$. Then,
$\PP_n'(E_n \mid F_n) \ge 1 - \kappa(\eta)^n$ for some $\kappa(\eta) <
1$. Now the proposition follows since $\lambda_0 = 1$ on $U_\alpha(0)$ and
the law of $q_n(-n\tau')$, conditioned on $F_n$, is
the normalized restriction of $\mu_0$ to $U_\alpha(0)$.

\end{proof}

\subsection{Mixing}
\label{sec:subsec:appendix:mixing}

{
In this subsection, let $0<\delta\leq1/400$, and use the random walk
constructed above with this value of $\delta$.
\par}

\begin{lemma}
\label{lemma:crude:lower:bound:on:volume}
{
There exist $c_1>0$ and $N_1>0$ such that for all $x\in\POT$,
\begin{displaymath}
\mu_x(B(x,\ooh))\geq c_1 V(x)^{-N_1}.
\end{displaymath}
\par}
\end{lemma}

\begin{proof}
{
Compare the return to a compact set in
\cite[Theorem~2.5 and the proof of Lemma~2.6]{mirzakhani:effective:counting}.
We record the adaptation to the present height and horospherical measure.
Fix $\alpha\in(0,1/4)$, put $b=1-2\alpha$, and take
$r=\ooh$ and $t=b^{-1}\log V(x)$. Proposition~\ref{prop:AG} and
Lemma~\ref{prop:Doubling} give
\[
\frac{1}{\mu_x(B(x,r/2))}
\int_{B(x,r/2)}V(y(t))\,d\mu_x(y)
\leq\frac{C_\alpha}{2c_D}\bigl(e^{-bt}V(x)+1\bigr)
=\frac{C_\alpha}{c_D}.
\]
Thus some $y\in B(x,r/2)$ has $z=y(t)$ in the fixed
$\Mod$--compact set $K=\{V\leq2\max(1,C_\alpha/c_D)\}$.
The compact-set lower bound in Lemma~\ref{lem:AGY-volume} and
doubling give $\mu_z(B(z,r/2))\geq c_1>0$ uniformly on $K$.
This bound follows from local period-coordinate comparison and
positivity of the horospherical volume form.
By backward monotonicity in \eqref{eq:Monotone},
\[
g_{-t}B(z,r/2)\subset B(y,r/2)\subset B(x,r).
\]
The conditional measures at $x$ and $y$ agree, so
\eqref{eq:Volume-Growth} yields
\[
\mu_x(B(x,r))
\geq\mu_x\bigl(g_{-t}B(z,r/2)\bigr)
=e^{-wt}\mu_z(B(z,r/2))
\geq c_1 V(x)^{-w/b}.
\]
Take $N_1=w/b$.
\par}
\end{proof}

Suppose $\Omega \subset \POT$ is a $\Mod$--compact set and
let $N_\delta(\Omega)$ be the $\delta$-neighborhood of $\Omega$ (in
the $D_\AGY$ metric). 

We recall a few lemmas from \cite{mirzakhani:effective:counting}.
{
As in \cite[Section~2.1]{mirzakhani:effective:counting}, we work with
embedded period boxes in a fixed finite cover without orbifold points.
More precisely, fix a torsion-free finite-index subgroup
$\Gamma_0<\Mod$, put $X_0=\POT/\Gamma_0$, and write
$\pi_0:\POT\to X_0$ for the projection. Let $\nu_0$ be the
Masur--Veech probability measure on $X_0$. The inverse image in $X_0$
of a measurable subset of $\POM$ has the same normalized measure as
that subset. Thus the estimates below for $\Mod$--invariant target sets
descend to $\POM$. The cover is fixed throughout, so constants may
depend on it.
\par}

\begin{lemma}[{\cite[Lemma~2.6]{mirzakhani:effective:counting}}]
\label{lemma:embedded:hor:ball}
{
There exist $c_3>0$ and $N_3\geq1$ such that, on setting
\[
r_3(x)=c_3V(x)^{-N_3},
\]
the restriction of $\pi_0$ to $B(x,r)$ is injective for every
$x\in\POT$ and $0<r\leq r_3(x)$. We may also require that the period
box of radius $10r_3(x)$ centered at $x$ embeds in $X_0$.
\par}
\end{lemma}

\begin{proof}
{
We spell out the argument of the cited lemma with the present height
$V$. Fix the parameters in \cite[Lemmas~6.2 and~6.4]{AG}. The height
comparison in Lemma~6.2 and the horocycle estimate in Lemma~6.4 give
constants $A\geq1$ and $b>0$ such that
\[
\int_0^1 V(g_tu_sx)\,ds\leq A e^{-bt}V(x)+A
\qquad(t\geq0),
\]
where $u_s=\left(\begin{smallmatrix}1&s\\0&1\end{smallmatrix}\right)$.
Taking $t=b^{-1}\log(2AV(x))$, we find $s\in[0,1]$ for which
$V(g_tu_sx)\leq A+1$. Consequently $\pi_0(g_tu_sx)$ belongs to a
fixed compact set $K\subset X_0$, and
$t\leq C_0+b^{-1}\log V(x)$.

Write $\mathcal B_a(x)$ for the period box of radius $a$ defined in
\cite[Section~2.9]{mirzakhani:effective:counting}. Since $X_0$ has no
orbifold points, compactness of $K$ gives $a_K>0$ such that the period
box $\mathcal B_{a_K}(z)$ embeds in $X_0$ whenever $\pi_0(z)\in K$.
The linear action in period coordinates and
\cite[Lemma~2.4]{mirzakhani:effective:counting} imply, with $h=g_tu_s$,
\[
h\mathcal B_a(x)\subset
\mathcal B_{C_1V(x)^{2/b}a}(hx)
\]
whenever the radius on the right is sufficiently small. Choose
$N_3\geq\max(1,2/b)$ and then choose $c_3>0$ so small that
$10C_1c_3<a_K$ and all these period boxes are defined. It follows that
$\mathcal B_{10r_3(x)}(x)$ embeds: any identification in that box would,
after applying $h$, give an identification in
$\mathcal B_{a_K}(hx)$, since the $SL(2,\reals)$ action commutes with
$\Gamma_0$.

Finally, \cite[Proposition~5.3]{AG} gives
$B(x,r)\subset\mathcal B_{2r}(x)$ for all sufficiently small $r$.
Decreasing $c_3$ once more, we may assume this for $r\leq r_3(x)$
and also that $r_3(x)\leq\ooh$. This proves the lemma. The fixed
factor $c_3$ is retained because $V$ can equal $1$.
\par}
\end{proof}

Let $\cC^1(\varphi)$ denote the $C^1$ norm of $\varphi$ as defined  in
\cite{AG}. The following is a rephrased version of
{\cite[Lemma~2.10]{mirzakhani:effective:counting}}. 

\begin{lemma}
\label{lemma:mirzakhani:partition}
There exists $N_5 > 0$ so that the following holds.
Let $x\in\POT$, and let
{$\epsilon=0.01e^{-N_3\ooh}r_3(x)$.}
There exists a finite collection $\{\psi_i\}$ of $C^\infty$ functions on $\Hor(x)$ with the following properties:
\begin{enumerate}
\item $0\leq\psi_i\leq 1$ for all $i$.
\item $\cC^1(\psi_i)\leq N_5 V(x)^{N_5}$.
\item For every $i$, $\psi_i$ is supported on $B(y_i,\epsilon)$ for some $y_i\in B(x,\ooh)$.
\item The covering $\{B(y_i,\epsilon)\}$ of $B(x,\ooh)$ has multiplicity at most $N_5$.
\item $\sum\psi_i\leq 1$, and the equality holds on a neighborhood of $B(x,\ooh)$. 
\end{enumerate}   
\end{lemma}
{
The extra fixed factor in $\epsilon$ ensures compatibility with
Lemma~\ref{lemma:embedded:hor:ball}: since $\log V$ is $1$--Lipschitz
on a horosphere and $y_i\in B(x,\ooh)$, we have
\[
r_3(y_i)\geq e^{-N_3\ooh}r_3(x),
\qquad \epsilon\leq0.01r_3(y_i).
\]
All fixed factors are absorbed in $N_5$ and the constants below.
\par}
Let $N(x)$ denote the number of functions $\psi_i$ in
Lemma~\ref{lemma:mirzakhani:partition}. Then, in view of
(iv), 
\begin{displaymath}
\sum_{i=1}^{N(x)} \mu_x(B(y_i, \epsilon)) \le N_5  \mu_x(B(x,\ooh+\epsilon))
\end{displaymath}
Thus, in view of Lemma~\ref{prop:Doubling},
$N(x) \le N_5 \epsilon^{-d}$, where $d$ depends only on the
dimension. WLOG, $N_5 >  d$, thus 
\begin{equation}
\label{eq:number:of:functions}
N(x) \le N_5 \epsilon^{-N_5}. 
\end{equation}

{
\begin{lemma}[Boundary layers of metric neighborhoods]
\label{lemma:neighborhood:boundary:layer}
Let $(X,d)$ be a proper geodesic metric space with a locally finite
Borel measure $\nu$. Suppose that, for some $D\geq1$ and $r_0>0$,
\[
0<\nu(B_d(z,2r))\leq D\nu(B_d(z,r))
\qquad(z\in X,\ 0<2r\leq r_0).
\]
For a nonempty set $F\subset X$, write
$F_a=\{z:d(z,F)<a\}$. Let $\delta>0$ and
$r_*=\min(\delta,r_0)$, and assume $\nu(F_{\delta+r_*})<\infty$.
Then, for $0<\epsilon<r_*$,
\begin{equation}
\label{eq:neighborhood:boundary:layer}
\nu(F_{\delta+\epsilon}\setminus F_{\delta-\epsilon})
\leq C\left(\frac{\epsilon}{r_*}\right)^\beta
\nu(F_{\delta+r_*}),
\end{equation}
where
\[
C=1+D^{-5},\qquad
\beta=\frac{\log(1+D^{-5})}{\log32}>0.
\]
In particular, if $\nu(X)=1$, then for every $0<\epsilon<\delta$,
\[
\nu(F_{\delta+\epsilon}\setminus F_{\delta-\epsilon})
\leq C_\delta\epsilon^\beta,
\qquad C_\delta=C\min(\delta,r_0)^{-\beta},
\]
uniformly in $F$.
\end{lemma}

\begin{proof}
Replace $F$ by its closure, and put
\[
\rho(z)=d(z,F),\qquad
A_s=\{z:\delta-s\leq\rho(z)<\delta+s\}.
\]
We claim that
\begin{equation}
\label{eq:boundary:layer:contraction}
\nu(A_s)\leq\lambda\nu(A_{32s}),\qquad
0<s\leq r_*/32,\qquad
\lambda=(1+D^{-5})^{-1}<1.
\end{equation}
Choose a maximal pairwise disjoint family
$B_i=B_d(z_i,8s)$ with $z_i\in A_s$. Properness implies separability,
so this family is countable. Maximality gives
\[
A_s\subset\bigcup_i B_d(z_i,16s).
\]
For each $i$, choose a nearest point of $F$ to $z_i$, and follow a
minimizing geodesic toward it for distance $4s$, reaching $w_i$.
This is possible since $\rho(z_i)\geq\delta-s\geq31s$.
Then $\rho(w_i)=\rho(z_i)-4s$. The $1$--Lipschitz property of
$\rho$ shows that the ball $H_i=B_d(w_i,s)$ satisfies
\[
H_i\subset\{\delta-6s<\rho<\delta-2s\}
\subset A_{32s}\setminus A_s.
\]
Also $H_i\subset B_i$, so the $H_i$ are pairwise disjoint. Since
\[
B_d(z_i,16s)\subset B_d(w_i,20s)\subset B_d(w_i,32s),
\]
five applications of doubling give
$\nu(B_d(z_i,16s))\leq D^5\nu(H_i)$; all radii are at most $r_0$.
Consequently,
\begin{align*}
\nu(A_{32s})
&\geq\nu(A_s)+\sum_i\nu(H_i)\\
&\geq\nu(A_s)+D^{-5}\sum_i\nu(B_d(z_i,16s))\\
&\geq(1+D^{-5})\nu(A_s).
\end{align*}
This proves \eqref{eq:boundary:layer:contraction}.

Set $n=\lfloor\log_{32}(r_*/\epsilon)\rfloor$. Iterating gives
\[
\nu(A_\epsilon)\leq\lambda^n\nu(A_{32^n\epsilon})
\leq\lambda^n\nu(F_{\delta+r_*}).
\]
Since $32^{-\beta}=\lambda$, we have
$\lambda^n\leq\lambda^{-1}(\epsilon/r_*)^\beta$, proving
\eqref{eq:neighborhood:boundary:layer}. If $\nu(X)=1$, this proves
the last assertion for $\epsilon<r_*$; for $\epsilon\geq r_*$ it
follows from the trivial bound by $1$.
\end{proof}
\par}

We need a mixing estimate for the indicator of $N_\delta(\Omega)$,
with constants independent of $\Omega$. The following variant of
\cite[Corollary~3.4]{mirzakhani:effective:counting} follows from
exponential mixing \cite{AGY} by approximating this indicator
with smooth functions.

Let $\chi$ be the  characteristic function of the
complement of $N_\delta(\Omega)$.
\begin{proposition}
\label{prop:hor-mixing}
For every $\delta > 0$ there exist   
$\kappa>0$ and $C_2>0$ (both depending on $\delta$)
such that for every $q \in \POT$ and every {Lipschitz function} $\varphi:
\Hor(q) \to \reals$ {with compact support in $B(q,0.1r_3(q))$},
for every $\Mod$-compact set $\Omega \subset \POT$ and every
$t  > 0$, we have
\begin{displaymath}
\Big| \int  \varphi(x) (1-\chi) (x(t)) \,  d\mu_q(x)-
\muMV(N_\delta(\Omega)/\Mod)\cdot \int \varphi \, d\mu_q
\Big| \leq {\|\varphi\|_{\mathrm{Lip}}} C_2 \cdot e^{-\kappa t}
\end{displaymath}
\end{proposition} 
\begin{proof}
{
We first bound the measure near the boundary of $N_\delta(\Omega)$,
uniformly in $\Omega$, using
Lemma~\ref{lemma:neighborhood:boundary:layer}. We then bracket
its indicator by smooth functions and apply mixing.

To apply the boundary estimate, we verify its local doubling hypothesis.
On the marked unit-area
stratum, the slow variation of the AGY norm
\cite[Proposition~5.5]{AG}, together with the period-box and
local-product coordinates of
\cite[Sections~2.9--2.10]{mirzakhani:effective:counting}, gives uniform
local comparisons with a normed vector space and Lebesgue measure.
In particular, writing $\widetilde\nu_0$ for the lift of $\nu_0$ and
$m=\dim_\reals\POT$, there are a fixed radius $r_1>0$ and numbers
$J(z)>0$ such that
\[
\widetilde\nu_0(B_\AGY(z,a))\asymp J(z)a^m
\qquad(0<a\leq r_1),
\]
with comparison constants independent of $z$ and $a$. Indeed, the
AGY norm is uniformly comparable to its value at the center in a
fixed-size period chart; in the unit-area local-product coordinates
the volume density is uniformly comparable to its value at the
center. Comparing the corresponding norm balls gives the displayed
estimate. Thus marked balls are uniformly locally doubling and have
uniform packing bounds at these scales.

These bounds also hold on $X_0$, even for balls larger than its local
injectivity radius. Here is the multiplicity argument. For a lift $z$
of $\bar z\in X_0$, let
\[
M_z(a)=\#\{\gamma\in\Gamma_0:
                D_\AGY(z,\gamma z)<a\}.
\]
For all sufficiently small $r>0$, the packing bound covers
$B_\AGY(z,2r)$ by at most $L$ balls of radius $r/2$, with $L$ uniform.
Each such ball contains at most $M_z(r)$ points of the orbit
$\Gamma_0z$: translate one orbit point in the ball back to $z$.
Hence $M_z(2r)\leq L M_z(r)$.
Counting lifts over quotient balls gives
\begin{align*}
\frac{\widetilde\nu_0(B_\AGY(z,r))}{M_z(2r)}
&\leq\nu_0(B_{X_0}(\bar z,r)),\\
\nu_0(B_{X_0}(\bar z,2r))
&\leq\frac{\widetilde\nu_0(B_\AGY(z,3r))}{M_z(r)}.
\end{align*}
For the first inequality, a point lifted to $B_\AGY(z,r)$ has at most
$M_z(2r)$ lifts there. For the second, a point lifted to
$B_\AGY(z,2r)$ has at least $M_z(r)$ lifts in $B_\AGY(z,3r)$.
Combining these inequalities with marked doubling proves that there
exist uniform $D$ and $r_0>0$ such that
\begin{equation}
\label{eq:ambient:quotient:doubling}
\nu_0(B_{X_0}(\bar z,2r))
\leq D\nu_0(B_{X_0}(\bar z,r))
\qquad(0<2r\leq r_0).
\end{equation}
The AGY length metric is complete and locally compact, so $X_0$ is
a proper geodesic space.

Let $F\subset X_0$ be the inverse image of $\Omega/\Mod$, and set
$E=N_\delta(F)$. We may replace $F$ by its closure. It is compact,
and $E$ is the inverse image of $N_\delta(\Omega)/\Mod$; in
particular,
\[
\nu_0(E)=\muMV(N_\delta(\Omega)/\Mod).
\]
By \eqref{eq:ambient:quotient:doubling} and
Lemma~\ref{lemma:neighborhood:boundary:layer}, there are $\beta>0$
and $C_\delta$ independent of $F$ such that
\begin{equation}
\label{eq:AGY:neighborhood:layer}
\nu_0(N_{\delta+\epsilon}(F)\setminus
              N_{\delta-\epsilon}(F))
\leq C_\delta\epsilon^\beta
\qquad(0<\epsilon<\delta).
\end{equation}

For $0<\epsilon<\min(\delta/10,r_0/10,1/10)$, choose smooth
compactly supported functions $f_\epsilon^-$ and $f_\epsilon^+$ on
$X_0$ such that
\begin{gather}
0\leq f_\epsilon^-\leq 1_E\leq f_\epsilon^+\leq1,
\qquad
\cC^1(f_\epsilon^\pm)\leq C\epsilon^{-1},
\label{eq:AGY:smooth:brackets}\\
f_\epsilon^-=1\text{ on }N_{\delta-3\epsilon}(F),
\qquad
\operatorname{supp}f_\epsilon^+\subset N_{\delta+3\epsilon}(F).
\notag
\end{gather}
For completeness, start with Lipschitz cutoffs of $d(\cdot,F)$:
the lower cutoff is $1$ for $d\leq\delta-2\epsilon$ and $0$ for
$d\geq\delta-\epsilon$, and the upper cutoff is $1$ for
$d\leq\delta+\epsilon$ and $0$ for $d\geq\delta+2\epsilon$.
Their Lipschitz constants are at most $\epsilon^{-1}$. Smooth in
period charts of sufficiently small diameter and combine with a
locally finite smooth partition of unity $\{\zeta_j\}$. If $h_j$
is a local smoothing of a cutoff $h$, choose its value error $a_j$
so small that
\[
\sum_j a_j\|d\zeta_j\|_\AGY\leq\epsilon^{-1},
\qquad a_j\leq1/10.
\]
Local norm comparison bounds $\|dh_j\|_\AGY$ by
$C\epsilon^{-1}$. The identity
\[
d\Big(\sum_j\zeta_j h_j\Big)
=\sum_j\zeta_j\,dh_j+\sum_j(h_j-h)\,d\zeta_j
\]
therefore gives the same uniform derivative bound after gluing.
Compose with a fixed smooth function equal to $0$ on
$(-\infty,1/5]$ and to $1$ on $[4/5,\infty)$, taking values in
$[0,1]$. Choosing the chart diameters and convolution radii small
enough preserves the indicated inner and outer margins. This proves
\eqref{eq:AGY:smooth:brackets}; the chart sizes may depend on $F$,
but the first-derivative bound does not. Properness gives compact
support. It follows from \eqref{eq:AGY:neighborhood:layer} that
\begin{equation}
\label{eq:AGY:bracket:mass}
0\leq\nu_0(f_\epsilon^+)-\nu_0(f_\epsilon^-)
\leq C_\delta\epsilon^\beta.
\end{equation}

By \cite[Proposition~3.2]{mirzakhani:effective:counting}, applied on
the fixed cover $X_0$, there exist $C>0$ and $\kappa_0>0$ such that,
for smooth compactly supported $f$ and smooth weights $\varphi$
with the stated support,
\begin{equation}
\label{eq:hor:mixing:smooth:target}
\left|\int_{\Hor(q)}\varphi(x)f(x(t))\,d\mu_q(x)
 -\nu_0(f)\int_{\Hor(q)}\varphi\,d\mu_q\right|
\leq C\cC^1(\varphi)\cC^1(f)e^{-\kappa_0t}.
\end{equation}
The support condition is admissible by
Lemma~\ref{lemma:embedded:hor:ball}. The fixed finite cover changes
only the mixing constants; compare also
\cite[Theorem~8.2 and its proof]{mirzakhani:effective:counting}.

The same estimate holds for compactly supported Lipschitz weights,
with $\|\varphi\|_{\mathrm{Lip}}$ in place of $\cC^1(\varphi)$.
Indeed, identify the horospherical period chart with a subset of a
real vector space and freeze the AGY norm at $q$. By
\cite[Proposition~5.3]{AG}, this norm and the horospherical metric
are uniformly comparable on the chart. Extend $\varphi$ by zero in
these coordinates and convolve with a nonnegative smooth approximate
identity. The resulting smooth functions $\varphi_a$ satisfy
\[
\varphi_a\longrightarrow\varphi\quad\text{uniformly},
\qquad
\cC^1(\varphi_a)\leq C\|\varphi\|_{\mathrm{Lip}},
\]
where $C$ is independent of $q$ and $a$. Since the original support
is a compact subset of $B(q,0.1r_3(q))$, sufficiently small convolution
radii keep the supports in this ball. Convolution also preserves
nonnegativity. For each fixed $t$, uniform convergence on this
common support, which has finite $\mu_q$-measure, lets us pass to the
limit in \eqref{eq:hor:mixing:smooth:target}. The first-derivative
bound is uniform as $a\to0$, so only a fixed factor is lost and the
decay exponent is unchanged. We use this Lipschitz version below.

First suppose $\varphi\geq0$. Apply
\eqref{eq:hor:mixing:smooth:target} with $f=f_\epsilon^-$ and
$f=f_\epsilon^+$ and use the pointwise brackets. Since
$\int\varphi\,d\mu_q\leq C\|\varphi\|_{\mathrm{Lip}}$ by
Lemma~\ref{lem:AGY-volume}, equations
\eqref{eq:AGY:smooth:brackets}--\eqref{eq:AGY:bracket:mass} give
\begin{equation}
\label{eq:hor:mixing:bracket:error}
\left|\int\varphi(x)1_E(x(t))\,d\mu_q(x)
       -\nu_0(E)\int\varphi\,d\mu_q\right|
\leq C_\delta\|\varphi\|_{\mathrm{Lip}}
       \big(\epsilon^{-1}e^{-\kappa_0t}+\epsilon^\beta\big).
\end{equation}
This also holds for signed Lipschitz $\varphi$. Apply the
nonnegative case to
$\varphi_+=\max(\varphi,0)$ and
$\varphi_-=\max(-\varphi,0)$ and subtract. Each has Lipschitz norm
at most $\|\varphi\|_{\mathrm{Lip}}$, so only a fixed factor is lost.

Fix $0<a_\delta<\min(\delta/10,r_0/10,1/10)$ and put
\[
\epsilon=a_\delta e^{-\kappa_0t/(1+\beta)},
\qquad \kappa=\frac{\beta\kappa_0}{1+\beta}>0.
\]
Both terms in \eqref{eq:hor:mixing:bracket:error} are bounded by a
constant depending only on $\delta$ times $e^{-\kappa t}$.
Finally, $1_E(x(t))=(1-\chi)(x(t))$ and
$\nu_0(E)=\muMV(N_\delta(\Omega)/\Mod)$, so the result follows.
\par}
\end{proof}

\begin{proposition}
\label{prop:all-generic}
Suppose $\theta > 0$, {$0<\delta\leq1/400$} and suppose $\Omega \subset \POT$ is a
Mod-compact set satisfying $\muMV(N_\delta(\Omega )/\Mod) \geq
(1-\theta)$.  Then there exists $t_\theta=t_\theta(\theta,\delta) >0$ and
$N_\ast  = N_\ast(\delta)$ 
such that 
for every $n\ge 0$, for every $q \in \Hor_n$,  if  
$t \geq t_\theta + N_\ast \log(V(q))$, then 
\begin{equation}
\label{eq:all-generic}
 \int P_q(x) \chi (x(t)) \,  d\bar{\mu}_n(x)
\leq   2 \theta.
\end{equation}
\end{proposition}

\begin{proof}
Let $\{\psi_i\}_{1 \le i \le N(q)}$ be the functions of
Lemma~\ref{lemma:mirzakhani:partition}. Let
$\varphi_i(x) = P_q(x) \psi_i(x)$.  By
Lemma~\ref{lemma:random:walk:easy}(v),
Lemma~\ref{lemma:crude:lower:bound:on:volume} and
Lemma~\ref{lemma:mirzakhani:partition}(2),  
\begin{displaymath}
{\|\varphi_i\|_{\mathrm{Lip}}}
\le C_\delta V(q)^{N_1+ N_3 N_5}.
\end{displaymath}
{
Here we use the product inequality
$\|P_q\psi_i\|_{\mathrm{Lip}}
\leq\|P_q\|_{\mathrm{Lip}}\|\psi_i\|_{\mathrm{Lip}}$
and $\|\psi_i\|_{\mathrm{Lip}}\leq\cC^1(\psi_i)$.
\par}
For each  $i$ we  apply
Proposition~\ref{prop:hor-mixing} with $\varphi$ replaced by
$\varphi_i$.  Summing over $i$ and noting that $\int P_q(x) \,
d\bar{\mu}_n(x) =1$, we get
\begin{multline*}
\Big| \int  P_q(x) (1-\chi) (x(t)) \,  d\bar{\mu}_n(x)-
\muMV(N_\delta(\Omega)/\Mod)
\Big| \leq \\ \leq N(q) (C_\delta V(q)^{N_1+ N_3 N_5}) C_2 \cdot
e^{-\kappa t}\leq  C_\delta' V(q)^{N_1 + 2 N_3 N_5} C_2
\cdot e^{-\kappa  t}. 
\end{multline*}
{
Let $N_\ast = 2 (N_1+N_3N_5)/\kappa$. Choose
$t_\theta=t_\theta(\theta,\delta)>0$ so that
$C_\delta'C_2e^{-\kappa t_\theta}\leq\theta$.
If $t\geq t_\theta+N_\ast\log V(q)$, then the right hand side of
the preceding estimate is at most $C_\delta'C_2e^{-\kappa t_\theta}$,
since $\kappa N_\ast\geq N_1+2N_3N_5$ and $V(q)\geq1$.
The proposition follows.
\par}
\end{proof}

Let $M$ be a large parameter. Let
\begin{displaymath}
V_M(q) = \max(V(q), M). 
\end{displaymath}
It follows from Lemma~\ref{lemma:phi:sigma:V0} that for any $n$ and
any $p \in \Hor_n$,
\begin{equation}
\label{eq:int:Vplus}
\frac{1}{V_M(p)}\int P_p(q) V_M(q(\tau')) \, d\bar{\mu}_n(q) \leq 
\begin{cases} 
C_\delta {e^{-3\tau'/4}} + \frac{M}{V(p)} & \text{if} \quad V(p) \geq V_0'\\
\frac{C_\delta+M}{V_M(p)} & \text{if} \quad V(p) \leq V_0'.
\end{cases} 
\end{equation}

For $\omega \in \Omega_n'$ and $0 \leq j \leq n$, define
\[
\widehat{V}(\omega,j) = V_M(q_j) \cdot \exp \left(  \sum_{i=1}^{j}  \chi(q_i) \right). 
\]

\begin{lemma}
\label{lemma:margulis:mixing}
Suppose $\theta > 0$, {$0<\delta\leq1/400$} and   
suppose $\Omega$ is a $\Mod$-compact subset of $\POT$ with
$\muMV(N_\delta(\Omega)/\Mod) \ge 1-\theta$.
{Then there exist $M=M(\theta,\delta)\geq1$ and
$\tau_0'=\tau_0'(\theta,\delta)>0$}
such that for all $\tau' >  \tau_0'$ with $\tau' \in \tau\natls$ and all
$n \geq 1$ we have 
\[
\EE_{n}'(\widehat{V}(\omega, n) \mid \cF_0') \leq e^{8 \theta  \, n}
V_M(q_0).  
\]
\end{lemma}

\begin{proof}
{Choose $M=M(\theta,\delta)\geq1$ so that}
\begin{equation}
\label{eq:app:choice:M}
  1+ \frac{e C_\delta}{M} + 2\theta(e-1) < e^{8\theta},  
\end{equation}
where $C_\delta$ is as in (\ref{eq:int:Vplus}).   
{Choose $\tau_0'=\tau_0'(\theta,\delta)>0$ large enough
that, for $\tau'>\tau_0'$, we have $V_0'>M$ and
$\tau'>t_\theta+N_\ast\log M$.}
  
We claim that
\begin{equation}
\label{eq:mixing:inductive-step}
\EE_{n}'(\widehat{V}(\omega, {j+1}) | \cF_j' )  \leq e^{8\theta} \cdot \widehat{V}(\omega, j). 
\end{equation}
Note that
\[
  \EE_{n}'(\widehat{V}(\omega, {j+1}) \mid \cF_j') = \exp \left(  \sum_{i=1}^{j} 
  \chi(q_i) \right)
  \int P_{q_j}(q)   e^{\chi(q(\tau'))}V_M(q(\tau'))\, d\bar{\mu}_j(q).
\]
Dividing both sides by $\widehat{V}(\omega,j)$ we get 
\begin{displaymath}
\frac{\EE_{n}'(\widehat{V}(\omega, {j+1}) \mid \cF_j'
  )}{\widehat{V}(\omega,j)} = \frac{1}{V_M(q_j)} \int P_{q_j}(q) e^{\chi(q(\tau'))}V_M(q(\tau'))\, d\bar{\mu}_j(q).
\end{displaymath} 
 
We consider three cases. First assume $V(q_j) \ge V_0'$. (We choose
$\tau'$ large enough  so that 
$V_0' > M$). Then, using (\ref{eq:int:Vplus}) and the trivial
bound $e^{\chi(q(\tau'))}\le e$, we get
\begin{align*}
\frac{1}{V_M(q_j)} \int P_{q_j}(q) e^{\chi(q(\tau'))}V_M(q(\tau'))\,
  d\bar{\mu}_j(q)  & \le e \cdot C_\delta {e^{-3\tau'/4}} +
                     \frac{e\cdot M}{V_0'}   \\
  & \le e \cdot (C_\delta +  M) {e^{-3\tau'/4}}. 
\end{align*}
Now assume  $V(q_j) <  V_0'$ but $t_\theta + N_\ast \log V(q_j) \ge
\tau'$. (In particular, with $\tau'$ chosen large enough, $V(q_j) \ge M$). 
Then, again using (\ref{eq:int:Vplus}) and the trivial bound
$e^{\chi(q(\tau'))}\le e$, we get, 
\begin{displaymath}
\frac{1}{V_M(q_j)} \int P_{q_j}(q) e^{\chi(q(\tau'))}V_M(q(\tau'))\,
d\bar{\mu}_j(q) \le 
  \frac{e \cdot (C_\delta+M)}{V(q_j)} \le e \cdot (C_\delta+M) e^{-(\tau' -t_\theta)/N_\ast}
\end{displaymath}
Now suppose $V(q_j) < V_0'$  but $t_\theta + N_\ast \log V(q_j) <
\tau'$. (This is the remaining case). 
We decompose
\begin{displaymath}
  \int P_{q_j}(q) e^{\chi(q(\tau'))}V_M(q(\tau'))\, d\bar{\mu}_j(q)
 = I_1 + I_2,
\end{displaymath}
where $I_1$ is the integral over $\{q \st V(q(\tau')) > M\}$ and $I_2$ is the integral over $\{ q \st V(q(\tau')) \le M\}$. 
We have, by Lemma~\ref{lemma:phi:sigma:V0}, 
\begin{displaymath}
\frac{I_1}{V_M(q_j)}\le \frac{e}{M} \int_{V(q(\tau')) \ge M}
P_{q_j}(q) V(q(\tau')) \, d\bar{\mu}_j(q) \le  \frac{e\cdot C_\delta}{M}. 
\end{displaymath}
Also,
\begin{align*}
\frac{I_2}{V_M(q_j)} & \le \int P_{q_j}(q) e^{\chi(q(\tau'))} \,
                       d\bar{\mu}_j(q) \\
                     &  = \int P_{q_j}(q) (1 + (e-1)\chi(q(\tau'))) \,  d\bar{\mu}_j(q) \\
  & \le  1 + 2\theta(e-1),
\end{align*}
where we have used Proposition~\ref{prop:all-generic} for the last
estimate.

Thus, combining the estimates in the three cases we get
\begin{multline*}
\frac{\EE_{n}'(\widehat{V}(\omega, {j+1}) | \cF_j' )}{\widehat{V}(\omega,j)} 
  \leq \max\Big( e \cdot  (C_\delta+M) {e^{-3\tau'/4}},  e\cdot
    (C_\delta+M) e^{-(\tau'-t_\theta)/N_\ast}, \\  \frac{ e\cdot C_\delta}{M} + 1
    + 2\theta(e-1)\Big). 
\end{multline*} 
{Thus, in view of (\ref{eq:app:choice:M}), after increasing
$\tau_0'$ if necessary, (\ref{eq:mixing:inductive-step}) holds for
every $\tau'>\tau_0'$.} The rest of the proof is the
same  as the proof of Proposition~\ref{prop:random-walk}.
\end{proof}

Let
\begin{displaymath}
  S(\omega) = \{ k  \in [0,n] \st q_k \in N_\delta(\Omega). \}.
\end{displaymath}

\begin{lemma}
\label{lemma:mixing:all:good}
{Let $\eta>0$ and $0<\delta\leq1/400$.}
Then,  for $\theta > 0$ sufficiently small depending on $\eta$, any
$\Mod$-compact $\Omega  \subset \POT$ with $\muMV(N_\delta(\Omega)/\Mod)
\ge 1-\theta$, any $\tau' \in \tau
\natls$ sufficiently large (depending on $\theta$, $\delta$) and any $n$
large enough,  
\begin{displaymath}
\PP_n'\Big( |S(\omega)| < (1-\eta) n \Big) < e^{-\tfrac{\eta}{2} n}.   
\end{displaymath}
\end{lemma}

\begin{proof}
{
Choose $0<\theta<\eta/32$. Let $M=M(\theta,\delta)$ and
$\tau_0'=\tau_0'(\theta,\delta)$ be as in
Lemma~\ref{lemma:margulis:mixing}, and suppose
$\tau'\in\tau\natls$ with $\tau'>\tau_0'$.
By Lemma~\ref{lemma:V0:finite:exp} and the finiteness of $\cA$,
choose $A\geq1$ so that $\EE_n'V(q_0)\leq A$ for every chart.
The initial law is independent of $\tau$, $\tau'$, $\theta$,
$\delta$, and $\Omega$, so $A$ can be chosen independently of
these parameters. Increasing $\tau_0'$ if necessary, we may assume
that $M<V_0'=e^{3\tau'/4}$.
Hence
\[
\EE_n'V_M(q_0)\leq A+M\leq A+e^{\tau'}.
\]
Then, by Lemma~\ref{lemma:margulis:mixing},
\begin{displaymath}
\EE_n'(\widehat{V}(\omega,n)) \le (A+e^{\tau'})e^{8\theta n}
\leq (A+e^{\tau'})
e^{\tfrac{\eta}{4} n}. 
\end{displaymath}
\par}
Note that if $|S(\omega)| < (1-\eta) n$, we have
$\widehat{V}(\omega,n) \ge e^{\eta n}$. 
Then, by the Markov inequality,
{
\begin{displaymath}
\PP_n'\Big( |S(\omega)| < (1-\eta) n \Big)
\leq (A+e^{\tau'})e^{-\tfrac{3\eta}{4}n}<e^{-\eta n/2}
\end{displaymath}
provided $n>\frac{4}{\eta}\log(A+e^{\tau'})$.
Thus the lower bound on $n$ can be chosen to depend only on
$\eta$ and $\tau'$. The lemma follows.
\par}
\end{proof}

\begin{lemma}
\label{lemma:all-generic}
Suppose  $0 < \eta < 1/2$ and {$0<\delta\leq1/400$} and $\tau > 0$
are arbitrary, and suppose $\Omega
\subset \POT$ is a $\Mod$-compact set. Then there exists $\tau_0' =
\tau_0'(\eta,\delta,\tau) > 0$ such that if $\tau' \in  \tau \natls$
and $\tau' >  \tau_0'$, and  if $k \in S(\omega) \cap
G_+(\omega,\eta/16)$, then 
\begin{displaymath}
q_n(-n \tau'+k\tau')  \in N_{2\delta}(\Omega). 
\end{displaymath}
\end{lemma}

\begin{proof}
Let $q = q_n(-n\tau')$. Since $k \in G_+(\omega,\eta/16)$, by
Lemma~\ref{lemma:tracking}, $D_\AGY(q_k,q(k\tau')) \le  d_\AGY(q_k, q(k
\tau')) \le \delta$. Since $k \in S(\omega)$, we have $q_k \in
N_\delta(\Omega)$. Therefore, $q(k \tau') \in N_{2\delta}(\Omega)$.
\end{proof}

{
\begin{corollary}[Geodesic returns]
\label{cor:geodesic:returns}
For every $0<\eta<1/2$ there exist $\tau_0=\tau_0(\eta)>0$,
$0<\kappa_{\rm v}=\kappa_{\rm v}(\eta)<1$, and
$\theta_0=\theta_0(\eta)>0$ such that, for every $\tau\geq\tau_0$,
every $0<\theta<\theta_0$,
every $0<\delta\leq1/400$, and every
$\Mod$--compact $\Omega\subset\POT$ satisfying
\[
\muMV(N_\delta(\Omega)/\Mod)\geq1-\theta,
\]
there exists $\tau'_0=\tau'_0(\eta,\tau,\theta,\Omega,\delta)>0$
such that if $\tau'>\tau'_0$ (and $\tau'\in\tau\natls$), the
following holds:

Let $n\geq0$ be an integer, and set
\[
v_n(x)=\#\{1\leq j\leq n:x(j\tau')\in N_{2\delta}(\Omega)\}.
\]
Then, for all sufficiently large $n$
(depending on $\eta$, $\tau$, $\theta$, $\Omega$, $\delta$, and $\tau'$)
and every $\alpha\in\cA$,
\[
\mu_0\bigl(\{x\in U_\alpha(0):v_n(x)\leq(1-\eta/2)n\}\bigr)
\leq\kappa_{\rm v}^{\,n}\mu_0(U_\alpha(0)).
\]
\end{corollary}

\begin{proof}
Let $\epsilon=\eta/64$, and choose $\tau_0=\tau_0(\eta)$ as in
Lemma~\ref{lemma:longrange:good:statistics}, with $\epsilon$ in place
of $\eta$. Choose $\theta_0=\theta_0(\eta)>0$ so that
Lemma~\ref{lemma:mixing:all:good} applies with $\eta/4$ in place of
$\eta$. Suppose $\tau\geq\tau_0$ and $0<\theta<\theta_0$.

Let
\[
E_n^+=\{|G_+(\omega,\epsilon)|\geq(1-\eta/8)n\},
\qquad
E_n^\Omega=\{|S(\omega)|\geq(1-\eta/4)n\}.
\]
Then, by Lemmas~\ref{lemma:longrange:good:statistics}
and~\ref{lemma:mixing:all:good}, for $\tau'\in\tau\natls$
sufficiently large,
\[
\PP_n'((E_n^+)^c)\leq e^{-\epsilon^4n\tau'/4},
\qquad
\PP_n'((E_n^\Omega)^c)\leq e^{-\eta n/8}
\]
for all sufficiently large $n$. We also assume that $\tau'>1$ and
that $\tau'$ is large enough for Lemma~\ref{lemma:all-generic}
to apply with $\eta/4$ in place of $\eta$.

Let $q=q_n(-n\tau')$. By Lemma~\ref{lemma:all-generic},
$q(j\tau')\in N_{2\delta}(\Omega)$ for every
$j\in G_+(\omega,\epsilon)\cap S(\omega)$. Hence, on
$E_n^+\cap E_n^\Omega$,
\[
v_n(q)\geq
\big|G_+(\omega,\epsilon)\cap S(\omega)\cap[1,n]\big|
\geq(1-3\eta/8)n-2>(1-\eta/2)n
\]
for sufficiently large $n$. Here we have omitted the index zero
from the two sets.

By Lemma~\ref{lemma:random:walk:easy}{\rm(ii)} and
\eqref{eq:appendix:Volume-Growth}, the law of $q$ is
\[
\frac{\mu_0'}{\mu_0'(\Hor_0)}.
\]
Fix $\alpha\in\cA$, and let $F_n=\{q\in U_\alpha(0)\}$.
Since $\lambda_0=1$ on
$U_\alpha(0)$, we have
\[
\PP_n'(F_n)=c_\alpha:=
\frac{\mu_0(U_\alpha(0))}{\mu_0'(\Hor_0)}>0,
\qquad
\PP_n'(q\in A\mid F_n)
=\frac{\mu_0(A)}{\mu_0(U_\alpha(0))}
\quad(A\subset U_\alpha(0)).
\]
Since $\mu_0'(\Hor_0)<\infty$ by Lemma~\ref{lemma:V0:finite:exp},
it follows that
\[
\frac{\mu_0(\{x\in U_\alpha(0):v_n(x)\leq(1-\eta/2)n\})}
     {\mu_0(U_\alpha(0))}
\leq c_\alpha^{-1}
\bigl(e^{-\epsilon^4n\tau'/4}+e^{-\eta n/8}\bigr).
\]
Choose $\kappa_{\rm v}=\kappa_{\rm v}(\eta)$ so that
\[
\max\{e^{-\epsilon^4/4},e^{-\eta/8}\}
<\kappa_{\rm v}<1.
\]
Since $\tau'>1$, the right hand side of the preceding estimate is
at most $\kappa_{\rm v}^{\,n}$ for all sufficiently large $n$.
Since $\cA$ is finite, this holds for all $\alpha\in\cA$ when $n$
is sufficiently large.

\end{proof}
\par}

In the next definition, good statistics has the full backward-time
meaning of Definition~\ref{def:good-statistics}.

\begin{definition}
\label{def:generic}
For $R>0$, we say $q \in \Hor_0$ is $(K, \eta,\delta,\Omega)$-generic at $R$ if the following hold: 
\begin{itemize}
\item[{\rm (i)}] $q(R) \in N_{2\delta}(\Omega)$ and has $(K, \eta)$--good statistics, and 
\item[{\rm (ii)}] For all  $0 \le s < R$,
\begin{displaymath}
| \{ u \in (s, R)  : \text{$q(u)$ has $(K, \eta)$--good statistics} \} | \geq (1-\eta) (R-s).
\end{displaymath}
\end{itemize}
\end{definition}

\begin{proposition}
\label{prop:appendix-final}
For any $1/2> \eta > 0$, there exists $\kappa = \kappa(\eta) < 1$ and
$\tau = \tau(\eta) > 0$ a $\Mod$--compact set $K=K(\tau) \supset K_0$
and $\theta_0 = \theta_0(\eta)>0$ such that
the following holds: 
For any $0  < \theta < \theta_0$ and any {$0<\delta\leq1/400$}, for any
$\Mod$--compact $\Omega \subset \POT$ with
\[
\muMV(N_\delta(\Omega)/\Mod) \geq (1 - \theta),
\]
there exists $\tau_0' =  \tau_0'(\eta,\theta,\Omega,\delta)$ such that if
$\tau' \in \tau \natls$ and $\tau' > \tau_0'$  the  following holds:
Let $R_j = j \tau'$.  
Then, for
\[
k_n''(x) = \# \Big\{ j \in [1,n]  : \text{$x$ is $(N_\delta(K), \eta, \delta,\Omega)$-generic at $R_j$} \Big\},
\]
we have, for all sufficiently large $n$ (depending on $\tau'$, $\eta$,
$\Omega$, $\delta$), and for  all $\alpha  \in \cA$, that
\[
\mu_0(\{ x \in U_\alpha(0)  : k_n''(x) \le (1-\eta)n \}) \le \kappa^n
\mu_0(U_\alpha(0)).  
\]
\end{proposition}

\begin{proof}
{
Choose $\tau=\tau(\eta)$ sufficiently large for
Proposition~\ref{prop:good:statistics:final}, with $\eta/4$ in place
of $\eta$, and Corollary~\ref{cor:geodesic:returns}.
Let $K=K(\tau)\supset K_0$ and $0<\kappa_{\rm g}=\kappa_{\rm g}(\eta)<1$
be as in the proposition, and let $\theta_0=\theta_0(\eta)$ and
$\kappa_{\rm v}=\kappa_{\rm v}(\eta)$ be as in the corollary.
Choose $\kappa=\kappa(\eta)$ so that
\[
\max\{\kappa_{\rm g},\kappa_{\rm v}\}<\kappa<1.
\]

Suppose $0<\theta<\theta_0$, $0<\delta\leq1/400$, and $\Omega$ satisfies
the measure hypothesis. Choose $\tau'_0$ large enough for both
Proposition~\ref{prop:good:statistics:final} and
Corollary~\ref{cor:geodesic:returns}, and suppose
$\tau'>\tau'_0$ (and $\tau'\in\tau\natls$). Put $R_j=j\tau'$.
Fix $\alpha\in\cA$. For $x\in U_\alpha(0)$, let
\[
\begin{aligned}
\mathcal I_n(x)
&=\{1\leq j\leq n:x\text{ is }(K,\eta/4)
       \text{--generic at }R_j\},\\
\mathcal J_n(x)
&=\{1\leq j\leq n:x(R_j)\in N_{2\delta}(\Omega)\}.
\end{aligned}
\]
Then, by Proposition~\ref{prop:good:statistics:final} and
Corollary~\ref{cor:geodesic:returns},
\[
\begin{aligned}
\mu_0(\{x\in U_\alpha(0):|\mathcal I_n(x)|\leq(1-\eta/4)n\})
&\leq\kappa_{\rm g}^{\,n}\mu_0(U_\alpha(0)),\\
\mu_0(\{x\in U_\alpha(0):|\mathcal J_n(x)|\leq(1-\eta/2)n\})
&\leq\kappa_{\rm v}^{\,n}\mu_0(U_\alpha(0))
\end{aligned}
\]
for all sufficiently large $n$.

Outside the union of these two exceptional sets,
\[
|\mathcal I_n(x)\cap\mathcal J_n(x)|
>(1-3\eta/4)n>(1-\eta)n.
\]
For every $j\in\mathcal I_n(x)\cap\mathcal J_n(x)$, we have
$x(R_j)\in N_{2\delta}(\Omega)$ and $x$ is $(K,\eta/4)$--generic
at $R_j$ in the sense of Definition~\ref{def:K:eta:generic}.
Since $x\in U_\alpha(0)$ and $K\supset K_0$, by
Lemma~\ref{lemma:good:statistics:beyond:zero}, the relative good
statistics in that definition imply good statistics on the full
backward geodesic, also at the intermediate times.
Since $\eta/4<\eta$ and $K\subset N_\delta(K)$, it follows that
$j$ is counted by $k_n''(x)$.

It follows that
\[
\begin{aligned}
&\mu_0(\{x\in U_\alpha(0):k_n''(x)\leq(1-\eta)n\})\\
&\qquad\leq
(\kappa_{\rm g}^{\,n}+\kappa_{\rm v}^{\,n})\mu_0(U_\alpha(0))
\leq\kappa^n\mu_0(U_\alpha(0))
\end{aligned}
\]
for all sufficiently large $n$. Since $\cA$ is finite, this holds
for all $\alpha\in\cA$ when $n$ is sufficiently large.
The proposition follows.
\par}
\end{proof}


\begin{thebibliography}{BDGGS24}

\bibitem[ABEM12]{athreya:bufetov:eskin:mirzakhani}
J.~Athreya, A.~Bufetov, A.~Eskin, and M.~Mirzakhani.
\newblock Lattice point asymptotics and volume growth on
  {T}eichm\"uller space.
\newblock {\em Duke Math. J.}, 161(6):1055--1111, 2012.

\bibitem[AG13]{AG}
A.~Avila and S.~Gou\"{e}zel.
\newblock Small eigenvalues of the {L}aplacian for algebraic measures in moduli
  space, and mixing properties of the {T}eichm\"{u}ller flow.
\newblock {\em Ann. of Math. (2)}, 178(2):385--442, 2013.

\bibitem[AGY06]{AGY}
A.~Avila, S.~Gou\"{e}zel, and J.-C.~Yoccoz.
\newblock Exponential mixing for the {T}eichm\"{u}ller flow.
\newblock {\em Publ. Math. Inst. Hautes \'{E}tudes Sci.}, 104:143--211, 2006.

\bibitem[AV07]{Avila-Viana}
A.~Avila and M.~Viana.
\newblock Simplicity of {L}yapunov spectra: proof of the
  {Z}orich-{K}ontsevich conjecture.
\newblock {\em Acta Math.}, 198(1):1--56, 2007.

\bibitem[BF25]{bader:furman:boundary}
U.~Bader and A.~Furman.
\newblock Lyapunov spectrum via boundary theory {I}.
\newblock Preprint, {\tt arXiv:2504.09731}, 2025.

\bibitem[BDGGS24]{bell:delecroix:gadre:gutierrezromo:schleimer:DF}
M.~Bell, V.~Delecroix, V.~Gadre, R.~Guti\'errez-Romo, and S.~Schleimer.
\newblock Diagonal flow detects topology of strata.
\newblock Preprint, {\tt arXiv:2101.12197}, 2024.

\bibitem[CM07a]{connell:muchnik:QC}
C.~Connell and R.~Muchnik.
\newblock Harmonicity of quasiconformal measures and {P}oisson boundaries of
  hyperbolic spaces.
\newblock {\em Geom. Funct. Anal.}, 17:707--769, 2007.

\bibitem[CM07b]{connell:muchnik:Gibbs}
C.~Connell and R.~Muchnik.
\newblock Harmonicity of {G}ibbs measures.
\newblock {\em Duke Math. J.}, 137(3):461--509, 2007.

\bibitem[CDR12]{choi:dumas:rafi:grafting}
Y.-E.~Choi, D.~Dumas, and K.~Rafi.
\newblock Grafting rays fellow travel {T}eichm\"uller geodesics.
\newblock {\em Int. Math. Res. Not. IMRN}, 2012(11):2445--2492, 2012.

\bibitem[CRS08]{rafi:LT}
Y.-E.~Choi, K.~Rafi, and C.~Series.
\newblock Lines of minima and {T}eichm\"uller geodesics.
\newblock {\em Geom. Funct. Anal.}, 18(3):698--754, 2008.

\bibitem[DDM14]{masur:SH}
S.~Dowdall, M.~Duchin, and H.~Masur.
\newblock Statistical hyperbolicity in {T}eichm\"uller space.
\newblock {\em Geom. Funct. Anal.}, 24(3):748--795, 2014.

\bibitem[DLR10]{rafi:LD}
M.~Duchin, C.~J. Leininger, and K.~Rafi.
\newblock Length spectra and degeneration of flat metrics.
\newblock {\em Invent. Math.}, 182(2):231--277, 2010.

\bibitem[Dum15]{dumas:skinning}
D.~Dumas.
\newblock Skinning maps are finite-to-one.
\newblock {\em Acta Math.}, 215(1):55--126, 2015.

\bibitem[ECPR]{eskin:cantrell:pan:rafi:HG}
A.~Eskin and S.~Cantrell, W.~Pan and K.~Rafi. 
\newblock Random walks on hyperbolic groups, 
\newblock in progress.

\bibitem[EM15]{eskin:matheus:LC}
A.~Eskin and C.~Matheus.
\newblock A coding-free simplicity criterion for the {L}yapunov exponents of
  {T}eichm\"uller curves.
\newblock {\em Geom. Dedicata}, 179:45--67, 2015.

\bibitem[EMM22]{mirzakhani:effective:counting}
A. Eskin, M. Mirzakhani and A. Mohammadi, Effective counting of simple
closed geodesics on hyperbolic surfaces, J. Eur. Math. Soc. (JEMS)
{\bf 24} (2022), no.~9, 3059--3108; MR4422208. 

\bibitem[EMR19]{eskin:mirzakhani:rafi:counting}
A.~Eskin, M.~Mirzakhani, and K.~Rafi.
\newblock Counting closed geodesics in strata.
\newblock {\em Invent. Math.}, 215(2):535--607, 2019.
  
\bibitem[FLP79]{flp:TTS}
A.~Fathi, F.~Laudenbach, and V.~Po\'enaru.
\newblock {\em Travaux de {T}hurston sur les surfaces}, volume 66--67 of
  {\em Ast\'erisque}.
\newblock Soci\'et\'e Math\'ematique de France, 1979.

\bibitem[FM12]{farb:MCG}
B.~Farb and D.~Margalit.
\newblock {\em A primer on mapping class groups}, volume~49 of
  {\em Princeton Mathematical Series}.
\newblock Princeton University Press, Princeton, NJ, 2012.

\bibitem[For02]{forni:deviation}
G.~Forni.
\newblock Deviation of ergodic averages for area-preserving flows on surfaces
  of higher genus.
\newblock {\em Ann. of Math. (2)}, 155(1):1--103, 2002.

\bibitem[Fra22]{frankel:meromorphic}
I.~Frankel.
\newblock Meromorphic {$L^2$} functions on flat surfaces.
\newblock {\em Geom. Funct. Anal.}, 32(4):832--860, 2022.

\bibitem[Fur63]{furstenberg:PF}
H.~Furstenberg.
\newblock A {P}oisson formula for semi-simple {L}ie groups.
\newblock {\em Ann. of Math. (2)}, 77:335--386, 1963.

\bibitem[Fur67]{furstenberg:PB}
H.~Furstenberg.
\newblock {P}oisson boundaries and envelopes of discrete groups.
\newblock {\em Bull. Amer. Math. Soc.}, 73(3):350--356, 1967.

\bibitem[Fur71]{furstenberg:RW}
H.~Furstenberg.
\newblock Random walks and discrete subgroups of {L}ie groups.
\newblock In {\em Advances in Probability and Related Topics}, Vol.~1,
  pages 1--63. Marcel Dekker, New York, 1971.

\bibitem[Gad14]{gadre:singular}
V.~Gadre.
\newblock Harmonic measures for distributions with finite support on the mapping
  class group are singular.
\newblock {\em Duke Math. J.}, 163(2):309--368, 2014.

\bibitem[GM91]{masur:EG}
F.~P. Gardiner and H.~A. Masur.
\newblock Extremal length geometry of {T}eichm\"uller space.
\newblock {\em Complex Variables Theory Appl.}, 16(2--3):209--237, 1991.

\bibitem[GMT17]{gadre:maher:tiozzo}
V.~Gadre, J.~Maher, and G.~Tiozzo.
\newblock Word length statistics for {T}eichm\"uller geodesics and singularity
  of harmonic measure.
\newblock {\em Comment. Math. Helv.}, 92(1):1--36, 2017.

\bibitem[Ham10]{hamenstadt:closing}
U.~Hamenst\"adt.
\newblock Dynamics of the {T}eichm\"uller flow on compact invariant sets.
\newblock {\em J. Mod. Dyn.}, 4(3):393--418, 2010.

\bibitem[HM79]{masur:QDF}
J.~Hubbard and H.~A. Masur.
\newblock Quadratic differentials and foliations.
\newblock {\em Acta Math.}, 142(3--4):221--274, 1979.

\bibitem[Hub06]{hubbard:TT}
J.~Hubbard.
\newblock {\em {T}eichm\"uller theory and applications to geometry, topology
  and dynamics. Volume 1: {T}eichm\"uller theory}.
\newblock Matrix Editions, Ithaca, NY, 2006.

\bibitem[IT92]{IT:Teich}
Y.~Imayoshi and M.~Taniguchi.
\newblock {\em An introduction to {T}eichm\"{u}ller spaces}.
\newblock Springer-Verlag, Tokyo, 1992.
\newblock Translated and revised from the Japanese by the authors.

\bibitem[Ker80]{kerckhoff:AG}
S.~P. Kerckhoff.
\newblock The asymptotic geometry of {T}eichm\"uller space.
\newblock {\em Topology}, 19(1):23--41, 1980.

\bibitem[KM94]{kaimanovich:masur:rennes}
V.~A. Kaimanovich and H.~Masur.
\newblock The {P}oisson boundary of the mapping class group and of
  {T}eichm\"uller space.
\newblock {\em Publications Math\'ematiques et Informatique de Rennes},
  1994, fascicule 2:1--69, 1994.

\bibitem[KM96]{kaimanovich:masur:PB}
V.~A. Kaimanovich and H.~Masur.
\newblock The {P}oisson boundary of the mapping class group.
\newblock {\em Invent. Math.}, 125(2):221--264, 1996.

\bibitem[KW04]{kleinbock:weiss:bounded}
D.~Kleinbock and B.~Weiss.
\newblock Bounded geodesics in moduli space.
\newblock {\em Int. Math. Res. Not.}, 2004(30):1551--1560, 2004.

\bibitem[LS04]{luo:stong:DT}
F.~Luo and R.~Stong.
\newblock Dehn--Thurston coordinates for curves on surfaces.
\newblock {\em Comm. Anal. Geom.}, 12:1--41, 2004.

\bibitem[Mas82]{masur:IE}
H.~Masur.
\newblock Interval exchange transformations and measured foliations.
\newblock {\em Ann. of Math. (2)}, 115(1):169--200, 1982.

\bibitem[PH92]{ph:TT}
R.~C. Penner and J.~L. Harer.
\newblock {\em Combinatorics of train tracks}, volume 125 of
  {\em Annals of Mathematics Studies}.
\newblock Princeton University Press, Princeton, NJ, 1992.

\bibitem[Raf07]{rafi:CM}
K.~Rafi.
\newblock A combinatorial model for the {T}eichm\"uller metric.
\newblock {\em Geom. Funct. Anal.}, 17(3):936--959, 2007.

\bibitem[Raf14]{rafi:hyperbolicity}
K.~Rafi.
\newblock Hyperbolicity in {T}eichm\"uller space.
\newblock {\em Geom. Topol.}, 18(5):3025--3053, 2014.

\bibitem[RS09]{rafi:schleimer:covers}
K.~Rafi and S.~Schleimer.
\newblock Covers and the curve complex.
\newblock {\em Geom. Topol.}, 13(4):2141--2162, 2009.

\bibitem[Sco83]{scott:geometries}
P.~Scott.
\newblock The geometries of {$3$}-manifolds.
\newblock {\em Bull. London Math. Soc.}, 15(5):401--487, 1983.

\bibitem[Shu26]{shuvaeva:AGY-volume}
E.~Shuvaeva.
\newblock {\em Geometry and dynamics in the Teichm\"uller and moduli spaces
of flat surfaces}.
\newblock Ph.D. dissertation, University of Chicago, 2026.
\newblock \url{https://doi.org/10.6082/uchicago.16955}.

\bibitem[SB01]{bonahon:WT}
Ya\c{s}ar S\"ozen and F.~Bonahon.
\newblock The {W}eil-{P}etersson and {T}hurston symplectic forms.
\newblock {\em Duke Math. J.}, 108(3):581--597, 2001.

\bibitem[Str84]{strebel:QD}
K.~Strebel.
\newblock {\em Quadratic differentials}, volume~5 of
  {\em Ergebnisse der Mathematik und ihrer Grenzgebiete}.
\newblock Springer-Verlag, Berlin, 1984.

\bibitem[Zor06]{zorich:FS}
A.~Zorich.
\newblock Flat surfaces.
\newblock In {\em Frontiers in number theory, physics, and geometry. {I}},
  pages 437--583. Springer, Berlin, 2006.

\bibitem[Zyk22]{zykoski:ID}
B.~Zykoski.
\newblock The $L^\infty$--isodelaunay decomposition of strata of abelian
  differentials.
\newblock Preprint, {\tt arXiv:2206.04143}, 2022.

\end{thebibliography}
\end{document}